\documentclass[sort&compress,times,preprint]{elsarticle}

\usepackage{url}
\usepackage{xcolor}

\definecolor{gold}{rgb}{0.85, 0.65, 0.13}
\global\long\def\red#1{\textcolor{black}{#1}}%
\global\long\def\violet#1{\textcolor{black}{#1}}%
\global\long\def\blue#1{\textcolor{black}{#1}}%

\usepackage[hidelinks, colorlinks]{hyperref}
\usepackage{bookmark}
\bookmarksetup{
	numbered,
	addtohook={%
		\ifnum\bookmarkget{level}>1 %
		\bookmarksetup{numbered=false}%
		\fi
	},
}

\usepackage{framed,multirow}
\usepackage[margin=25mm]{geometry}
\global\bibsep=0pt

\usepackage{amssymb}
\usepackage{latexsym}

\usepackage{amsfonts}
\usepackage{mathrsfs}
\usepackage{upgreek}

\usepackage{graphicx}
\usepackage{epstopdf}
\usepackage{algorithmic}
\ifpdf
\DeclareGraphicsExtensions{.eps,.pdf,.png,.jpg}
\else
\DeclareGraphicsExtensions{.eps}
\fi
\usepackage[export]{adjustbox}

\usepackage{amsmath, amsthm}
\usepackage{xspace}
\usepackage{bm}
\usepackage{physics}
\usepackage{stmaryrd}
\usepackage{cleveref}
\usepackage{enumerate}
\usepackage[labelfont=bf,textfont=normalfont]{caption}
\usepackage{subcaption}
\usepackage{mathtools}

\usepackage{booktabs}
\usepackage{multirow, multicol, makecell}
\usepackage{longtable}
\usepackage{threeparttable}

\newtheorem{theorem}{Theorem}
\newtheorem{proposition}{Proposition}
\newtheorem{lemma}{Lemma}
\newtheorem{assumption}{Assumption}
\newtheorem{remark}{Remark}
\newtheorem{corollary}{Corollary}
\newtheorem{definition}{Definition}

\crefname{assumption}{Assumption}{Assumptions}
\crefname{definition}{Definition}{Definitions}
\crefname{lemma}{Lemma}{Lemmas}
\crefname{remark}{Remark}{Remarks}
\crefname{theorem}{Theorem}{Theorems}
\crefname{proposition}{Proposition}{Propositions}
\crefname{section}{Section}{Sections}
\crefname{figure}{Fig.}{Figs.}
\crefname{equation}{}{}
\crefname{table}{Table}{Tables}
\crefname{appendix}{}{}

\newcommand{\fnc}[1]{\ensuremath{\mathcal{#1}}}
\newcommand{\vecfnc}[1]{\ensuremath{\boldsymbol{\mathcal{#1}}}} 

\newcommand{\uhk}[0]{\ensuremath{\bm{{u}}_{h,k}}}

\newcommand{\uk}[0]{\ensuremath{\bm{{u}}_{k}}}

\renewcommand{\H}[0]{\mathsf{H}}
\newcommand{\Dxk}[0]{\mathsf{D}_{xk}}
\newcommand{\Dyk}[0]{\mathsf{D}_{yk}}

\newcommand{\Sxk}[0]{\mathsf{S}_{xk}}
\newcommand{\Syk}[0]{\mathsf{S}_{yk}}

\newcommand{\Exk}[0]{\mathsf{E}_{xk}}
\newcommand{\Eyk}[0]{\mathsf{E}_{yk}}

\newcommand{\D}[0]{\mathsf{D}}
\newcommand{\Q}[0]{\mathsf{Q}}
\newcommand{\E}[0]{\mathsf{E}}
\newcommand{\M}[0]{\mathsf{M}}
\newcommand{\T}[0]{\mathsf{T}}
\newcommand{\F}[0]{\mathsf{F}}
\newcommand{\C}[0]{\mathsf{C}}
\newcommand{\A}[0]{\mathsf{A}}
\newcommand{\I}[0]{\mathsf{I}}
\newcommand{\G}[0]{\mathsf{G}}
\newcommand{\U}[0]{\mathsf{U}}
\newcommand{\K}[0]{\mathsf{K}}

\newcommand{\J}[0]{\mathsf{J}}
\newcommand{\X}[0]{\mathsf{X}}
\newcommand{\Y}[0]{\mathsf{Y}}
\renewcommand{\P}[0]{\mathsf{P}}

\newcommand{\Dxi}[0]{\hat{\mathsf{D}}_{\xi}}

\newcommand{\Sxi}[0]{\hat{\mathsf{S}}_{\xi}}

\newcommand{\Qxi}[0]{\hat{\mathsf{Q}}_{\xi}}
\newcommand{\Qeta}[0]{\hat{\mathsf{Q}}_{\eta}}

\newcommand{\Exi}[0]{\hat{\mathsf{E}}_{\xi}}

\newcommand{\Nxig}[0]{\hat{\mathsf{N}}_{\xi\gamma}}

\newcommand{\nxi}[0]{\hat{n}_{\xi}}
\newcommand{\neta}[0]{\hat{n}_{\eta}}

\newcommand{\B}[0]{\mathsf{B}}
\newcommand{\R}[0]{\mathsf{R}}
\newcommand{\N}[0]{\mathsf{N}}

\newcommand{\Bghat}[0]{\hat{\mathsf{B}}_{\gamma}}
\newcommand{\Nxgk}[0]{\mathsf{N}_{x\gamma k}}
\newcommand{\Nygk}[0]{\mathsf{N}_{y\gamma k}}

\newcommand{\Rg}[0]{\mathsf{R}_{\gamma}}
\newcommand{\Rgk}[0]{\mathsf{R}_{\gamma k}}
\newcommand{\Rgv}[0]{\mathsf{R}_{\gamma v}}

\newcommand{\Rbargk}[0]{\mathsf{\bar{\R}}_{\gamma k}}
\newcommand{\Dgk}[0]{\mathsf{D}_{\gamma k}}
\newcommand{\Dgv}[0]{\mathsf{D}_{\gamma v}}

\newcommand{\V}[0]{\mathsf{V}}

\DeclareMathOperator{\mydiag}{diag}

\newcommand{\etal}[0]{{et~al.\@}\xspace}
\newcommand{\eg}[0]{{e.g.\@}\xspace}
\newcommand{\ie}[0]{{i.e.\@}\xspace}

\newcommand{\ignore}[1]{} 

\newcommand{\polyref}[1]{\ensuremath{\mathbb{P}^{#1}}(\hat{\Omega})}

\newcommand{\contref}[1]{\ensuremath{\mathcal{C}^{#1}}(\hat{\Omega})}

\newcommand{\poly}[1]{\ensuremath{\mathbb{P}^{#1}}({\Omega}_k)}
\newcommand{\vecpoly}[2]{[\ensuremath{\mathbb{P}^{#1}}({\Omega}_k)]^{#2}}
\newcommand{\cont}[1]{\ensuremath{\mathcal{C}^{#1}}({\Omega}_k)}

\newcommand{\vecLtwo}[1]{[\ensuremath{L^2}({\Omega}_k)]^{#1}}

\newcommand{\IR}[1]{\mathbb{R}^{#1}}%
\newcommand{\IRtwo}[2]{\mathbb{R}^{{#1}\times{#2}}}%

\newcommand{\jump}[1]{\left\llbracket {#1}\right\rrbracket }%
\newcommand{\avg} [1]{\left\{  {#1}\right\}  }%

\newcommand{\sumfhat}[0] {\sum_{\gamma \subset \hat{\Gamma}}}

\newcommand{\sumfk}[0] {\sum_{\gamma \subset \Gamma_k}}

\newcommand{\liftVglob}[1]{{{\mathcal{L}}}\left({#1} \right)}
\newcommand{\liftVloc}[1]{{{\mathcal{L}}}^\gamma\left({#1} \right)}

\newcommand{\liftSglob}[1]{{\mathcal{S}}\left({#1} \right)}
\newcommand{\liftSloc}[1]{{\mathcal{S}}^\gamma\left({#1} \right)}

\newcommand{\innerprod}[2]{\left\langle {#1}, {#2}\right\rangle }

\newcommand{\myabs}[1]{|{#1}|}

\begin{document}


\begin{frontmatter}

\title{Simultaneous approximation terms and functional accuracy for diffusion problems discretized with multidimensional summation-by-parts operators} 

\author{Zelalem Arega Worku \corref{cor1}}
\cortext[cor1]{Corresponding author: } 
\ead{zelalem.worku@mail.utoronto.ca}

\author{David W. Zingg}
\ead{dwz@oddjob.utias.utoronto.ca}

\address{Institute for Aerospace Studies, University of Toronto, Toronto, Ontario, M3H 5T6, Canada}


\addcontentsline{toc}{section}{Abstract}
\begin{abstract}
Several types of simultaneous approximation term (SAT) for diffusion problems discretized with \blue{diagonal-norm} multidimensional summation-by-parts (SBP) operators are analyzed based on a common framework. Conditions under which the SBP-SAT discretizations are consistent, conservative, adjoint consistent, and energy stable are presented. For SATs leading to primal and adjoint consistent discretizations, the error in output functionals is shown to be of order $h^{2p}$ when a degree $p$ multidimensional SBP operator is used to discretize the spatial derivatives. SAT penalty coefficients corresponding to various discontinuous Galerkin fluxes developed for elliptic partial differential equations are identified. We demonstrate that the original method of Bassi and Rebay, the modified method of Bassi and Rebay, and the symmetric interior penalty method are equivalent when implemented with SBP diagonal-E operators that have diagonal norm matrix, \eg, the Legendre-Gauss-Lobatto SBP operator in one space dimension. Similarly, the local discontinuous Galerkin and the compact discontinuous Galerkin schemes are equivalent for this family of operators. The analysis remains valid on curvilinear grids if a degree $\le p+1$ bijective polynomial mapping from the reference to physical elements is used. Numerical experiments with the two-dimensional Poisson problem support the theoretical results.
\end{abstract}


\end{frontmatter}
\section{Introduction}
High-order methods can provide superior solution accuracy for a given computational cost. Furthermore, when used with unstructured and discontinuous elements they enable efficient $ hp $-adaptation and high code parallelization while still being consistent, locally conservative, and stable for wide range of fluid flow problems. Many of these powerful features can be attributed to the solution discontinuity between adjacent elements. The manner in which elements are coupled affects most essential properties of discretizations, such as accuracy, consistency, conservation, stability, adjoint consistency, functional convergence, conditioning, stiffness, sparsity, symmetry, and so on. Therefore, the coupling procedure at interfaces between adjacent elements is a critical aspect of discontinuous high-order methods. In this paper, we analyze the numerical properties of discretizations arising from the use of one such coupling procedure, simultaneous approximation terms (SATs) \cite{carpenter1994time}, for diffusion problems.

Discontinuous high-order methods developed in the past few decades include summation-by-parts (SBP) methods coupled with SATs, discontinuous Galerkin (DG), and flux reconstruction (FR) methods. In the DG and FR methods, element coupling and boundary conditions are enforced via numerical fluxes. A unified framework of DG fluxes for elliptic problems is analyzed by Arnold \etal \cite{arnold2002unified}, \blue{a review of the SBP-SAT method is presented in \cite{fernandez2014review,svard2014review}, and the connections between SBP-SAT, DG, and FR methods can be found, for example, in \cite{gassner2013skew,ranocha2016summation,chan2018discretely,montoya2021unifying}}. Motivated by developments in the DG method, Carpenter, Nordstr\"{o}m, and Gottlieb \cite{carpenter1999stable} introduced the Carpenter-Nordstr\"{o}m-Gottlieb (CNG) SAT to solve the multi-domain problem for high-order finite difference methods \cite{carpenter2007revisiting}. In later works \cite{carpenter2010revisiting, carpenter2007revisiting}, they showed that SATs are closely related to DG fluxes and introduced the Baumann-Oden (BO)\cite{baumann1999discontinuous} and local discontinuous Galerkin (LDG) \cite{cockburn1998local} type SATs for one-dimensional classical SBP operators. Although these SATs are consistent, conservative, and energy stable, not all of them possess other desired properties such as symmetry and adjoint consistency. Hicken and Zingg \cite{hicken2011superconvergent} presented conditions that SATs must satisfy for SBP-SAT discretizations to be adjoint consistent. Furthermore, they showed that, under mild assumptions, linear functionals superconverge for adjoint consistent discretizations. Adjoint consistency and functional superconvergence properties are further studied in \cite{berg2012superconvergent,hicken2014dual,hicken2012output,eriksson2018dual,eriksson2018finite,nikkar2019dual,nordstrom2017relation,ghasemi2020conservation}. Recently, Craig Penner and Zingg \cite{penner2020superconvergent} showed that functional superconvergence is retained in curvilinear coordinates for adjoint consistent discretizations of hyperbolic PDEs with generalized SBP operators \cite{fernandez2014generalized}. 

Multidimensional SBP operators were introduced by Hicken, Del Rey Fern{\'a}ndez, and Zingg \cite{hicken2016multidimensional}. The SBP operators constructed in \cite{hicken2016multidimensional} are classified as SBP-$ \Gamma $ operators -- a family of operators that have $ {p+d-1 \choose d-1} $ volume nodes on each facet, where $ d $ is the spatial dimension of the problem. Later, the SBP-$ \Omega $ \cite{fernandez2018simultaneous} and SBP diagonal-E\footnote{Abbreviated as SBP-E in all figures and tables.} \cite{chen2017entropy} operator families were introduced. SBP-$ \Omega $ operators have volume nodes strictly in the interior domain of the element, while the SBP diagonal-E operators are characterized by two features: facet nodes that are collocated with volume nodes and diagonal surface integral matrices. A broader classification of the operators that is based on the dimensions spanned by the volume to facet node extrapolation matrices\footnote{Also known as the interpolation/extrapolation matrix.}, $ \R $, categorizes the SBP-$ \Omega $, SBP-$ \Gamma $, and SBP diagonal-E operators under the $ \R^{d} $, $ \R^{d-1} $, and $ \R^{0} $ operator families, respectively, where the superscript on $ \R $ indicates the dimensions spanned by the extrapolation matrices \cite{marchildon2020optimization}. For a degree $ p $ multidimensional SBP operator that has a diagonal norm matrix\footnote{The norm matrix is known as the mass matrix in the DG literature.}, the diagonal entries of the norm matrix and the corresponding volume nodes define a degree $ 2p-1 $ quadrature rule, and this connection simplifies the construction of multidimensional SBP operators as quadrature rules are readily available in the literature \cite{fernandez2018simultaneous}. The analysis presented in this paper is restricted to multidimensional SBP operators that have a diagonal norm matrix.

SATs for hyperbolic problems discretized with SBP-$ \Gamma $ and SBP-$ \Omega $ operators were studied in \cite{hicken2016multidimensional,fernandez2018simultaneous}. A framework to implement SATs with multidimensional SBP operators for second-order partial differential equations (PDEs) was subsequently proposed by Yan, Crean, and Hicken \cite{yan2018interior}. The framework presented in \cite{yan2018interior} is flexible enough to construct compact\footnote{If SATs couple only immediate neighbor elements, they are referred to as compact stencil SATs; otherwise, they are referred to as wide or extended stencil SATs.} stencil SATs that lead to consistent, conservative, adjoint consistent, and energy stable SBP-SAT discretizations. Furthermore, it was shown in \cite{yan2018interior,yan2020immersed} that the modified method of Bassi and Rebay (BR2) \cite{bassi1997highbr2}, the symmetric interior penalty (SIPG) \cite{douglas1976interior}, and the compact discontinuous Galerkin (CDG) \cite{peraire2008compact} methods fall under this framework. Numerical properties of discretizations of the two-dimensional heat equation with SBP-$ \Gamma $ and SBP-$ \Omega $ operators coupled with the BR2 and SIPG SATs were also investigated in \cite{yan2018interior}. \red{For tensor-product SBP discretizations in multiple dimensions, wide stencil DG fluxes, such as the LDG method, are widely used for coupling of viscous terms \cite{carpenter2014entropy,parsani2015entropy,gassner2018br1}. However, many numerical properties of discretizations resulting from the use of wide stencil SATs and multidimensional SBP operators have not been analyzed so far. In light of this, we study properties of compact and wide stencil SATs under a general SAT framework for multidimensional SBP operators.}

The three main objectives of the present work are: (1) to extend the framework in \cite{yan2018interior} to allow construction of wide stencil SATs and study their numerical properties, (2) to demonstrate that when diffusion problems are discretized with degree $ p $ multidimensional SBP operators in a primal and adjoint consistent manner, the error in output functionals is of order $ h^{2p} $, and (3) to show the equivalence of different types of DG-based SATs when implemented with SBP diagonal-E operators that have a diagonal norm matrix. We also specify SAT coefficients that correspond to the consistent DG fluxes in \cite{arnold2002unified,peraire2008compact} and provide stability analysis for the SATs that are not studied in \cite{yan2018interior}. All results are presented in two space dimensions; however, generalization to three space dimensions is straightforward.

The paper is organized as follows: In \cref{sec:notation}, we introduce our notation and present important definitions. After describing the model problem in \cref{sec:model problem}, the SBP-SAT discretization and the generic SAT framework are provided in \cref{sec:SBP-SAT discretization}. We analyze properties of SBP-SAT discretizations in \cref{sec:Properties of SBP-SAT discretization} and present SATs corresponding to popular DG methods in \cref{sec:Existing and DG SATs}. In \cref{sec:Practial issues}, we demonstrate the equivalence of various types of SATs when implemented with the diagonal-norm $ \R^{0} $ SBP family and study the sparsity of system matrices arising from SBP-SAT discretizations. In \cref{sec:Numerical Results}, we investigate numerical properties of various SBP-SAT discretizations of the steady version of the model problem. Finally, we present concluding remarks in \cref{sec:conclusions}. 

\section{Notation and definitions}
\label{sec:notation}

We closely follow the notation in \cite{hicken2016multidimensional, yan2018interior, crean2018entropy}. A $d$-dimensional compact domain and its boundary are denoted by $\Omega \subset \mathbb{R}^d$ and $\partial\Omega$, respectively. The domain is tessellated into $n_e$ non-overlapping elements, $ {\mathcal T}_h \equiv \{\{ \Omega_k\}_{k=1}^{n_e}: \Omega=\cup_{k=1}^{n_e} {\Omega}_k\}$. The boundaries of each element will be referred to as facets or interfaces, and we denote their union by $ \Gamma_k \equiv \partial\Omega_k $. A reference element,  $\hat{\Omega}$, and its boundary, $\hat{\Gamma}$, are used to construct SBP operators which are then mapped to each physical element. The reference triangle is a right angle triangle with vertices $ \hat{v}_1 =(-1,-1) $, $ \hat{v}_2 = (1, -1) $, and  $ \hat{v}_3 = (-1,1) $ and facets $ f_1=\overrightarrow{\hat{v}_2 \hat{v}_3} $, $ f_2=\overrightarrow{\hat{v}_3 \hat{v}_1} $, and $ f_3=\overrightarrow{\hat{v}_1 \hat{v}_2} $. The boundaries $\partial\Omega$, $\Gamma_k$, and $\hat{\Gamma}$ are assumed to be piecewise smooth. The set of all interior interfaces is denoted by $\Gamma ^I \equiv \{\Gamma_k \cap \Gamma_v : k,v=1,\dots,n_e, k\neq v \}$. The set of facets of $ \Omega_k $ that are also interior facets is denoted by $ \Gamma^I_k \equiv \Gamma^I\cap\Gamma_k $, and $ \Gamma ^B \equiv \{\partial\Omega \cap \Gamma_k : k = 1,\dots, n_e\}$ delineates the set of all boundary facets. Finally, the set containing all facets is denoted by $ \Gamma \equiv \Gamma^I \cup \Gamma^B$. The set of $n_p$ volume nodes in element $ \Omega_k $ is represented by $ S_k=\{(x_i,y_i)\}_{i=1}^{n_p} $, while the set of $n_f$ nodes on facet $ \gamma \in {\Gamma_k}$ is denoted by $S_\gamma=\{(x_i,y_i)\}_{i=1}^{n_f} $. Similarly, we represent the set of volume nodes in the reference element, $ \hat{\Omega} $, and facet nodes on $ \gamma \in \hat{\Gamma} $ by $ \hat{S}=\{(\xi_i,\eta_i)\}_{i=1}^{n_p} $, and  $\hat{S}_\gamma=\{(\xi_i,\eta_i)\}_{i=1}^{n_f} $, respectively. Operators associated with the reference element bear a hat, $ (\hat\cdot) $.

Scalar functions are written in uppercase script type, \eg, $\fnc{U}_k \in \cont{\infty}$, and vector-valued functions of dimension $ d $ are represented by boldface uppercase script letters, \eg, $\vecfnc{W}_k \in \vecLtwo{d}$. The space of polynomials of total degree $ p $ is denoted by $\polyref{p} $, and $ n_p^* = {p+d \choose d} $ is the cardinality of the polynomial space. The restrictions of functions to grid points are denoted by bold letters, \eg, $ \uk \in \IR{n_p}$ is the evaluation of $ \fnc{U}_k $ at grid points $ S_k $, while vectors resulting from numerical approximations have subscript $ h $, \eg, $ \uhk \in \IR{n_p}$. When dealing with error estimates, we define $h\equiv \max_{a, b \in S_k} \norm{a - b}_2$ as the diameter of an element. Matrices are denoted by sans-serif uppercase letters, \eg, $\V \in \IRtwo{n_p}{n_p^*}$; $ \bm{1} $ denotes a vector consisting of all ones, $ \bm {0} $ denotes a vector or matrix consisting of all zeros. The sizes of  $ \bm{1} $ and $ \bm {0} $ should be clear from context. Finally, $ \I_{n} $ represents the identity matrix of size $ n \times n $ unless specified otherwise. 

The definition of multidimensional SBP operators first appeared in \cite{hicken2016multidimensional}, and is presented below on the reference element. 
\begin{definition}[Two-dimensional SBP operator \cite{hicken2016multidimensional}] \label{def:sbp}
	The matrix $\Dxi \in \IRtwo{n_p}{n_p}$ is a degree $ p $ SBP operator approximating the first derivative $ \pdv{\xi} $ on the set of nodes $ \hat{S}=\{(\xi_i,\eta_i)\}_{i=1}^{n_p} $ if
	\begin{enumerate}
		\item $ \Dxi \bm{p} = \pdv{\fnc{P}}{\xi}$ for all $\fnc{P} \in \polyref{p} $
		\item $ \Dxi=\hat{\H}^{-1} \Qxi $ where $ \hat{\H} $ is a symmetric positive definite (SPD) matrix, and 
		\item $ \Qxi = \Sxi + \frac{1}{2} \Exi $ where $ \Sxi = - \Sxi^T $, $ \Exi = \Exi^T $ and $ \Exi $ satisfies
		\[ \bm{p}^T \Exi  \bm {q} = \int_{\hat{\Gamma}} \fnc{P}\fnc{Q} \nxi \dd{\Gamma}\]
		for all $ \fnc{P},\fnc{Q} \in \polyref{r} $, where $ r \ge p $, and $ \nxi $ is the $ \xi $-component of the outward pointing unit normal vector, $ \hat{\bm{n}} = [\nxi, \neta]^T $, on $ \hat{\Gamma} $.
	\end{enumerate}  
\end{definition} 
An analogous definition applies for operators in the $ \eta $ direction. Properties $ 2 $  and $ 3 $ in \cref{def:sbp} give 
\begin{equation} \label{eq:SBP property}
	\begin{aligned}
		\Qxi + \Qxi^T &= \Exi,
	\end{aligned}
\end{equation}
which will be referred to as the SBP property. Throughout this paper, the matrix $ \hat{\H} $ is assumed to be diagonal. The set of nodes $ \hat{S} $ and the diagonal entries of $ \hat{\H} $ define a quadrature rule of at least degree $ 2p-1 $; thus, the inner product of two functions $ \fnc{P} $ and $ \fnc{Q} $ is approximated by \cite{hicken2016multidimensional, fernandez2018simultaneous}
\[ \bm{p}^T\hat{\H} \bm{q} = \int_{\hat{\Omega}}\fnc{P}\fnc{Q}\dd{\Omega} +\order{h^{2p}}.\]
Together with the fact that $ \hat{\H} $ is SPD, the above approximation can be used to define the norm
\[ \bm{u}^T\hat{\H} \bm{u} = \norm{\bm {u}}_{\hat{\H}}^2= \int_{\hat{\Omega}} \fnc{U}^2\dd{\Omega} + \order{h^{2p}}, \]
which is a degree $ 2p-1 $ approximation of the $ L^2 $ norm.

Under the assumption that a quadrature rule exists on $ \gamma \in \hat{\Gamma} $ with nodes $ \hat{S}_\gamma $ the surface integral matrix $ \Exi $ can be decomposed as \cite{fernandez2018simultaneous}
\begin{equation} \label{eq: Exi}
\Exi = \sumfhat \Rg^T \Bghat \Nxig \Rg,
\end{equation}
where, $ \Bghat \in \IRtwo{n_f}{n_f} $ is a diagonal matrix containing a minimum of degree $ 2p $ positive quadrature weights on $ \gamma $ along its diagonal, $ \Nxig \in \IRtwo{n_f}{n_f} $ contains the $ \xi $ component of $ \hat{\bm{n}}_\gamma $ along its diagonal, and $ \Rg \in \IRtwo{n_f}{n_p} $ is a matrix that extrapolates the solution from the volume nodes to the facet nodes. The quadrature accuracy requirement on $ \Bghat $ is a sufficient but not necessary condition to construct SBP operators \cite{shadpey2020entropy}. In this paper, we consider SBP operators with facet quadrature based on the Legendre-Gauss (LG) rule which offers a degree $ 2p+1 $ accuracy. The extrapolation matrix, $ \Rg $, is exact for polynomials of degree $ p $ on the reference element. \violet{For SBP-$ \Omega $ operators}, it is constructed as \cite{fernandez2018simultaneous} 
\begin{equation} \label{eq:extrapolation matrix}
\Rg=\hat{\V}_\gamma \hat{\V}_\Omega^+ = \hat{\V}_\gamma \left(\hat{\V}_\Omega^T \hat{\V}_\Omega\right)^{-1}\hat{\V}_{\Omega}^T,
\end{equation} 
where $ \hat{\V}_\gamma \in \IRtwo{n_f}{n_p^*} $ and $ \hat{\V}_\Omega \in \IRtwo{n_p}{n_p^*} $ are Vandermonde matrices constructed using \blue{the orthonormalized canonical basis discussed in \cref{sec:Curvilinear Transformation}} and the set of nodes $ \hat{S}_\gamma $ and $ \hat{S} $, respectively, and $ (\cdot)^+$ represents the Moore-Penrose pseudoinverse. \violet{For SBP-$ \Gamma $ operators, $ \R_{\gamma} $ is obtained using a $\hat{\V}_\Omega $ matrix constructed using the basis evaluated at the $ p+1 $ volume nodes that lie on facet $ \gamma $ \cite{fernandez2018simultaneous}. Finally, for SBP diagonal-E operators, $ \Rg $ contains unity at each entry $ (i,j) $ if $ i+(n - 1)n_f = j $, where $ n=\{1, 2, 3\} $ is the facet number; all other entries are zero \cite{shadpey2020entropy}.} 

Some definitions that are used in DG formulations of diffusion problems will prove useful for later discussions. Following \cite{arnold2002unified}, we introduce the broken finite element spaces associated with the tessellation $ \fnc{T}_h$ of $ \Omega $. The spaces of scalar and vector functions, $ V_h $ and $ \Sigma_h $ respectively, whose restrictions to each element, $ \Omega_k $, belong to the space of polynomials are defined by 
\begin{equation} \label{FEM space}
\begin{aligned}
V_h & \equiv \{\fnc{P}\in L^2(\Omega): \fnc{P}|_{\Omega_k} \in \poly{p}  \;\;  \forall\Omega_k \in \fnc{T}_h\},\\
\Sigma_h & \equiv \{\vecfnc{V}\in [L^2(\Omega)]^2: \vecfnc{V}|_{\Omega_k} \in \vecpoly{p}{2}  \;\; \forall\Omega_k \in \fnc{T}_h\},\\
\end{aligned}
\end{equation}
and the set in which traces\footnote{Traces define the restriction of functions along the boundaries of each element; thus, functions in $ T(\Gamma) $ are double-valued on $ \Gamma^I $ and single valued on $ \Gamma^B$ \cite{arnold2002unified}. See \cite{riviere2008discontinuous} for trace theorems which affect the function spaces in which the solution and test functions are sought.} of the functions in $\fnc{T}_h$ lie is defined by
\begin{equation} \label{FEM space trace}
	T(\Gamma) \equiv \Pi_{\Omega_k\in\fnc{T}_h} L^2(\Gamma_k). 
\end{equation}
The jump, $ \jump{\cdot}$, and average, $ \avg{\cdot} $, operators for scalar and vector-valued functions are defined as
\begin{equation} \label{eq:Jump and Average}
\begin{aligned}
\jump{\fnc{P}} &=\fnc{P}_k\bm{n}_k + \fnc{P}_v \bm{n}_v,  & \avg{\fnc{P}} &= \frac{1}{2} (\fnc{P}_k + \fnc{P}_v), & &  \forall \fnc{P}\in T(\Gamma),\\
\jump{\vecfnc{V}} &=\vecfnc{V}_k\cdot\bm{n}_k + \vecfnc{V}_v \cdot\bm{n}_v,  & \avg{\vecfnc{V}} &= \frac{1}{2} (\vecfnc{V}_k + \vecfnc{V}_v), & & \forall \vecfnc{V}\in [T(\Gamma)]^2.\\
\end{aligned}
\end{equation}
At the boundaries, $ \jump{\fnc{P}}=\fnc{P}_k\bm{n}_k $ and $ \avg{\vecfnc{V}}=\vecfnc{V}_k $, and the $ \avg{\fnc{P}} $ and $ \jump{\vecfnc{V}} $ are left undefined \cite{arnold2002unified}. Surface integral terms that appear in the DG flux formulation\footnote{The flux formulation is obtained by transforming second-order PDEs into a system of first-order PDEs.} are converted to volume integrals via lifting operators. For vector-valued functions, the global lifting operator for interior facets, $ \fnc{L} : [L^2(\Gamma^I)]^2 \rightarrow \Sigma_h$, and the local lifting operator for interior facets, $ \fnc{L}^\gamma : [L^2(\gamma)]^2 \rightarrow \Sigma_h$ are defined by \cite{peraire2008compact}
\begin{align}
\int_{\Omega}\liftVglob{\vecfnc{V}}\cdot \vecfnc{Z} \dd{\Omega} 
&= -\int_{\Gamma^I}\vecfnc{V}\cdot \avg{\vecfnc{Z}} \dd{\Gamma} 
&& \forall \vecfnc{Z}\in \Sigma_h,  \label{eq: lift global vector} \\
\int_{\Omega_\gamma}\liftVloc{\vecfnc{V}}\cdot \vecfnc{Z} \dd{\Omega} 
&= -\int_{\gamma}\vecfnc{V}\cdot \avg{\vecfnc{Z}} \dd{\Gamma}
&& \forall \vecfnc{Z}\in \Sigma_h,\; \gamma\in\Gamma^I, \label{eq: lift local vector} 
\end{align}
where $ \Omega_{\gamma} = \Gamma_k \cup \Gamma_v$. Similarly, for scalar functions, the global lifting operator, $ \fnc{S}: L^2(\Gamma^I) \rightarrow \Sigma_h$, the local lifting operator, $ \fnc{S}^\gamma: L^2(\gamma) \rightarrow \Sigma_h $,  and the lifting operators at Dirichlet boundary facets, $ \fnc{S}^D : L^2(\gamma) \rightarrow \Sigma_h$, are defined by 
\begin{align}
\int_{\Omega} \liftSglob{\fnc{P}} \cdot \vecfnc{Z} \dd{\Omega} 
&= -\int_{\Gamma^I}\fnc{P}\jump{\vecfnc{Z}} \dd{\Gamma} && \forall \vecfnc{Z}\in \Sigma_h, \label{eq: lift global scalar} 
\\
\int_{\Omega_\gamma} \liftSloc{\fnc{P}} \cdot \vecfnc{Z} \dd{\Omega} 
&= -\int_{\gamma}\fnc{P}\jump{\vecfnc{Z}} \dd{\Gamma} && \forall \vecfnc{Z}\in \Sigma_h,\; \gamma\in \Gamma^I, \label{eq: lift local scalar}
\\
\int_{\Omega_\gamma}\fnc{S}^D({\fnc{P}})\cdot \vecfnc{Z} \dd{\Omega} 
&= -\int_{\gamma}\fnc{P}\vecfnc{Z}\cdot\bm{n} \dd{\Gamma} 
&& \forall \vecfnc{Z}\in \Sigma_h,\; \gamma\in \Gamma^D \label{eq: lift Dirichlet},
\end{align}
respectively. \violet{Note that the surface integrals on the right-hand side (RHS) of \cref{eq: lift global vector} and \cref{eq: lift global scalar} do not include boundary facets; hence, these global lifting operators differ from similar definitions, \eg, in \cite{arnold2002unified, bassi2005discontinuous}. The consequence of such definitions of the global lifting operators is that the boundary conditions are enforced using compact SATs only, \ie, extended stencil SATs are applied exclusively on interior facets. This is important for adjoint consistency of discretizations of problems with non-homogeneous Dirichlet boundary conditions, as explained in \cref{sec:Adjoint Consistency}.} Moreover, the lifting operator at Dirichlet boundaries is defined locally; however, a global lifting operator definition would give the same final SBP-SAT discretization of the PDEs we are interested in.

\section{Model problem} \label{sec:model problem}
We consider the two-dimensional diffusion equation, 
\begin{equation} \label{eq:diffusion problem}
\begin{aligned}\pdv{\fnc{U}}{t}+L(\fnc{U})& =\fnc{F}\;\text{in}\;\Omega, &  & \fnc{U}=\fnc{U}_{0}\;\text{at}\;t=0, & & 
B_D(\fnc{U}) =\fnc{U}_{D}\;\text{on}\;\Gamma^{D}, &  & B_N(\fnc{U})=\fnc{U}_{N}\;\text{on}\;\Gamma^{N},
\end{aligned}
\end{equation}
where the linear differential operators in $ \Omega $,  on $ \Gamma^D $, and on $ \Gamma^N $ are given, respectively, by $ L(\fnc{U})=-\nabla\cdot\left(\lambda\nabla\fnc{U}\right) $, $ B_D(\fnc{U}) = \fnc{U} $, and $ B_N(\fnc{\fnc{U}})= \left(\lambda\nabla\fnc{U}\right)\cdot\bm{n} $, $ \fnc{F}\in L^2(\Omega) $ is the source term, $ \lambda= \bigl[ \begin{smallmatrix} \lambda_{xx} & \lambda_{xy} \\
\lambda_{yx} & \lambda_{yy} \end{smallmatrix} \bigr]$ is an SPD tensor with diffusivity coefficients in each combination of directions, and we assume that $ \Gamma \cap \Gamma^D \neq \emptyset$. \blue{The energy stability analysis presented in this work applies to SBP-SAT discretizations of the unsteady model problem given in \cref{eq:diffusion problem}.}

\blue{In order to study adjoint consistency and superconvergence of functionals, we will consider the steady version of \cref{eq:diffusion problem}, the Poisson problem}. We also consider a linear functional of the form 
\begin{equation} \label{eq:linear functional}
\fnc{I}(\fnc{U})=\innerprod{\fnc{G}_\Omega}{\fnc{U}}_{\Omega} 
+\innerprod{\fnc{G}_{\Gamma^{D}}}{C_D(\fnc{U})}_{\Gamma^{D}}
+\innerprod{\fnc{G}_{\Gamma^{N}}}{C_N(\fnc{U})}_{\Gamma^{N}},
\end{equation}
where $ \fnc{G}_\Omega \in L^2(\Omega)$, $ \fnc{G}_{\Gamma^D} \in L^2({\Gamma^D})$, $ \fnc{G}_{\Gamma^N} \in L^2({\Gamma^N})$, \blue{$ C_D $ and $ C_N $ are linear differential operators at the Dirichlet and Neumann boundaries, respectively}, and $\innerprod{\cdot}{\cdot}_{\Omega} $, $ \innerprod{\cdot}{\cdot}_{\Gamma^D} $, and $ \innerprod{\cdot}{\cdot}_{\Gamma^D} $ represent the $ L^2(\Omega) $, $ L^2(\Gamma^D) $, and $ L^2(\Gamma^N) $ inner products, respectively. Such a functional is said to be compatible with the steady version of \cref{eq:diffusion problem} if \cite{hartmann2007adjoint}
\begin{equation}\label{eq:compatibility condition}
	\innerprod{L(\fnc{U})}{\psi}_{\Omega} + \innerprod{B_D(\fnc{U})}{C_D^*(\psi)}_{\Gamma^{D}} + 
	\innerprod{B_N(\fnc{U})}{C_N^*(\psi)}_{\Gamma^{N}}
	= 
	\innerprod{\fnc{U}}{L^*(\psi)}_{\Omega} + \innerprod{C_D (\fnc{U})}{B_D^*(\psi)}_{\Gamma^{D}} + 
	\innerprod{C_N (\fnc{U})}{B_N^*(\psi)}_{\Gamma^{N}},	
\end{equation}
where $ L^* $, $ B^*_D $, $ C^*_D $, $ B^*_N $, and $ C^*_N $ are the adjoint operators to the linear differential operators $ L $, $ B_D $, $ C_D $, $ B_N $, and $ C_N $, respectively, and \blue{$ \psi $ is a unique adjoint variable in an appropriate function space, \eg, we assume $ \fnc{U},\fnc{\psi} \in H^2 $}. A compatible linear functional satisfies the relations \cite{giles1997adjoint,hartmann2007adjoint}
\begin{equation}\label{eq:Adjoint relation}
\begin{aligned}
\fnc{I}(\fnc{U})
&=\innerprod{\fnc{U}}{\fnc{G}_\Omega}_{\Omega}
+\innerprod{C_{D}(\fnc{U})}{\fnc{G}_{\Gamma^{D}}}_{\Gamma^{D}}
+\innerprod{C_{N}(\fnc{U})}{\fnc{G}_{\Gamma^{N}}}_{\Gamma^{N}}
=\innerprod{\fnc{U}}{L^*(\psi)}_\Omega
+\innerprod{C_D(\fnc{U})}{B_D^*(\psi)}_{\Gamma^{D}}
+\innerprod{C_N(\fnc{U})}{B_N^*(\psi)}_{\Gamma^{N}}\\
&=\innerprod{L(\fnc{U})}{\psi}_\Omega
+\innerprod{B_D(\fnc{U})}{C_D^*(\psi)}_{\Gamma^D}
+\innerprod{B_N(\fnc{U})}{C_N^*(\psi)}_{\Gamma^N}
=\innerprod{\fnc{F}}{\psi}_\Omega
+\innerprod{B_D(\fnc{U})}{C_D^*(\psi)}_{\Gamma^D}
+\innerprod{B_N(\fnc{U})}{C_N^*(\psi)}_{\Gamma^N}.
\end{aligned}
\end{equation}
In the subsequent analysis, we will consider the compatible linear functional 
\begin{equation} \label{eq:Functional}
\fnc{I} (\fnc{U})= \int_{\Omega}\fnc{G}\fnc{U}\dd{\Omega} 
- \int_{\Gamma^D}\psi_{D} (\lambda\nabla\fnc{U})\cdot\bm{n}\dd{\Gamma}
+ \int_{\Gamma^N}\psi_{N}\fnc{U}\dd{\Gamma},
\end{equation}
where $ \fnc{G}\in L^2{(\Omega)}$, $\psi_{N}=[\bm{n}\cdot(\lambda\nabla\psi)] \in L^2{(\Gamma^N)}$, and $\psi_{D}\in L^2{(\Gamma^D)}$. The functional given in \cref{eq:Functional} is obtained by substituting $ \fnc{G}_\Omega = \fnc{G} $, $ \fnc{G}_{\Gamma^{D}} = \psi_D $, $ \fnc{G}_{\Gamma^{N}} = \psi_N $, $ C_D(\fnc{U})=-(\lambda\nabla\fnc{U})\cdot\bm{n} $, $ C_N(\fnc{U})=\fnc{U} $, $ B_D(\fnc{U}) = \fnc{U}_D $, and $ B_N(\fnc{U}) = (\lambda\nabla\fnc{U})\cdot\bm{n} = \fnc{U}_N$ in \cref{eq:Adjoint relation}, \ie,
\begin{align}
\fnc{I}(\fnc{U})&=\innerprod{\fnc{G}}{\fnc{U}}_{\Omega} 
+\innerprod{\psi_{D}}{-(\lambda\nabla\fnc{U})\cdot\bm{n}}_{\Gamma^D}
+\innerprod{\psi_{N}}{\fnc{U}}_{\Gamma^N} 
\label{eq:Adjoint relation 2-1}
\\
&= \innerprod{L(\fnc{U})}{\psi}_{\Omega} + \innerprod{\fnc{U}_D}{C^*_D(\psi)}_{\Gamma^D} 
+\innerprod{\fnc{U}_N}{C^*_N(\psi)}_{\Gamma^N} 
\label{eq:Adjoint relation 2-2}
\\
&= \innerprod{\fnc{U}}{L^*(\psi)}_{\Omega} 
+\innerprod{-(\lambda\nabla\fnc{U})\cdot\bm{n}}{B^*_D(\psi)}_{\Gamma^D}
+\innerprod{\fnc{U}}{{B^*_N(\psi)}}_{\Gamma^N}.
\label{eq:Adjoint relation 2-3}
\end{align}
Following \cite{hartmann2007adjoint}, we apply integration by parts to $ \innerprod{L(\fnc{U})}{\psi}_{\Omega} $ twice and rearrange terms to find
\begin{equation} \label{eq:IBP on L(U)}
\begin{aligned}
\innerprod{L(\fnc{U})}{\psi}_{\Omega}
&-\innerprod{\fnc{U}}{(\lambda\nabla\psi)\cdot\bm{n}}_{\Gamma^D}
+\innerprod{(\lambda\nabla\fnc{U})\cdot \bm{n}}{\psi}_{\Gamma^N}
\\ 
&= -\innerprod{\fnc{U}}{\nabla\cdot(\lambda \nabla\psi)}_{\Omega}
-\innerprod{(\lambda\nabla\fnc{U})\cdot \bm{n}}{\psi}_{\Gamma^D}
+\innerprod{\fnc{U}}{(\lambda\nabla\psi)\cdot\bm{n}}_{\Gamma^N},
\end{aligned}
\end{equation}
where symmetry of the inner product is used assuming the problem is real-valued. Equations \cref{eq:Adjoint relation 2-2,eq:Adjoint relation 2-3,eq:IBP on L(U)} imply
\begin{equation}\label{eq:Adjoint operators}
\begin{aligned}
L^*(\psi) &= -\nabla\cdot\left(\lambda\nabla{\fnc{\psi}}\right), 
&& B^*_D(\psi) =  \psi,
&& B^*_N(\psi) = (\lambda\nabla\psi)\cdot\bm{n},
&& C^*_D(\psi) = -(\lambda\nabla\psi)\cdot\bm{n},
&& C^*_N(\psi) = \psi.	
\end{aligned}
\end{equation}
Using the result in \cref{eq:Adjoint operators} and subtracting \cref{eq:Adjoint relation 2-3} from \cref{eq:Adjoint relation 2-1} we have 
\begin{equation}
\innerprod{\fnc{G}+\nabla\cdot(\lambda\nabla\psi)}{\fnc{U}}_{\Omega}  
+\innerprod{\psi-\psi_{D}}{(\lambda\nabla\fnc{U})\cdot\bm{n}}_{\Gamma^D}
+\innerprod{\psi_{N} - (\lambda\nabla\fnc{U})\cdot\bm{n}}{\fnc{U}}_{\Gamma^N} = 0.
\end{equation}
Thus, the adjoint for the model problem satisfies
\begin{equation}\label{eq:Adjoint problem}
\begin{aligned} 
L^*(\psi)- \fnc{G} &= 0 \; \text{in} \; \Omega, &&
\psi = \psi_{D} \; \text{on} \; \Gamma^D, && (\lambda\nabla\psi)\cdot\bm{n} = \psi_{N} \; \text{on} \; \Gamma^N.
\end{aligned}
\end{equation}

\section{SBP-SAT discretization} \label{sec:SBP-SAT discretization}
In this section, the discretization of the model equation \cref{eq:diffusion problem} with multidimensional SBP operators is presented. Notation and definition of operators follow \cite{yan2018interior}. The following assumption is used in the construction of SBP operators on curved elements which is presented in \cref{sec:Curvilinear Transformation}. 
\begin{assumption} \label{assu: mapping}
	We assume that there exists a bijective and time-invariant polynomial mapping, $ \fnc{M}_k:{ \hat{\Omega}}\rightarrow{ {\Omega}_k} $, of degree $ p_{\rm map} \le p+1 $ for all $ \Omega_k \in \fnc{T}_h $. Furthermore, volume and facet quadrature rules with the set of nodes in the reference element exist such that \blue{diagonal-norm} SBP operators satisfying \cref{def:sbp} can be constructed on the reference element.
\end{assumption}

The extrapolation matrix is exact for constant functions in the physical element,$ \Omega_k $, particularly $ \Rgk \bm{1} = \bm{1} $. Polynomials in $ \Omega_k $ are not necessarily polynomials in the reference element, $ \hat{\Omega} $; thus, SBP operators in the physical domain are not exact for polynomials in $ \Omega_k $. However, under \cref{assu: mapping} the accuracy of the derivative operators in the physical elements is not compromised \cite{crean2018entropy}. We state, without proof, Theorem 9 in \cite{crean2018entropy} that establishes the accuracy of SBP derivative operators on physical elements. 
\begin{theorem}\label{thm:Accuracy of Dx}
	Let \cref{assu: mapping} hold and the metric terms be computed exactly using \cref{eq:grid metrics volume} and \cref{eq:grid metrics facet}. Then for $ \bm{u}_k\in\IR{n_p} $ holding the values of $ \fnc{U}\in\contref{p+1} $ at the nodes in $ \hat{S} $, the derivative operators given by \cref{eq:Dx and Dy} are order $ p $ accurate, \ie,
	\begin{equation*}
	[\Dxk\bm{u}_k]_i = {\pdv{\fnc{U}}{x}} (\xi_i,\eta_i) +\order{h^{p}}, \quad \text{and} \quad [\Dyk\bm{u}_k]_i = {\pdv{\fnc{U}}{y}} (\xi_i,\eta_i) +\order{h^{p}}.
	\end{equation*}
\end{theorem}
Furthermore, if the SBP operators on physical elements are constructed as described in \cref{sec:Curvilinear Transformation} and \cref{assu: mapping} holds, then $ \Dxk \bm{1} = \bm{0} $ and $ \Dyk \bm{1} = \bm{0} $, which are the conditions required to satisfy the discrete metric identity/freestream preservation condition \cite{crean2018entropy, shadpey2020entropy}. 

\begin{figure}
	\begin{centering}
		\includegraphics[width=0.27\columnwidth]{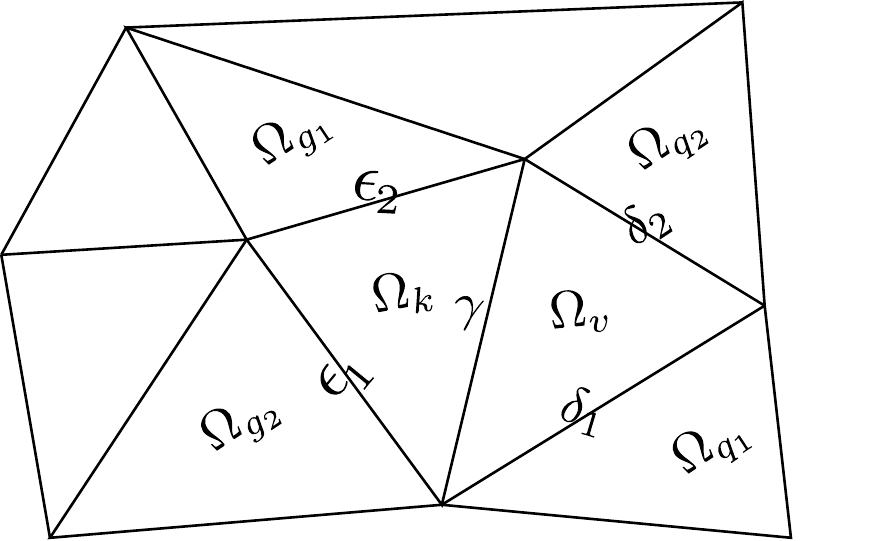}
	\caption{\label{fig:element and facet label}Element and facet labeling. For convenience common labeling is often used, \eg, $\sum_{\delta \subset \Gamma_v}(\cdot) $ implies $\sum_{\delta_{1}}(\cdot) + \sum_{\delta_{2}}(\cdot)$.}
	\end{centering}
\end{figure}

The diffusivity coefficients are evaluated at the volume nodes and stored in an SPD block matrix,
\begin{equation}\label{eq:Lambda}
\Lambda = \begin{bmatrix}
\Lambda_{xx} & \Lambda_{xy}\\
\Lambda_{yx} & \Lambda_{yy}
\end{bmatrix},
\end{equation} 
where each block is diagonal, \eg, $\Lambda_{xx} = \mydiag(\lambda_{xx}(x_1,y_1), \dots, \lambda_{xx}(x_{n_p}, y_{n_p}))$.
The second derivative is approximated by applying the first derivative twice, 
\begin{equation} \label{eq:D2 1st form}
[\nabla \cdot (\lambda \nabla)]_k \approx \D^{(2)}_k = \left[\begin{array}{cc}
\Dxk & \Dyk\end{array}\right] \Lambda_k \left[\begin{array}{c}
\Dxk\\
\Dyk
\end{array}\right],
\end{equation} 
and the normal derivative at facet $ \gamma $ is given by 
\begin{equation} \label{eq:D_gamma k}
[\bm{n}\cdot(\lambda\nabla)]_k \approx \Dgk = \N_{\gamma k}^T \Rbargk \Lambda_k\left[\begin{array}{c}
\Dxk\\
\Dyk
\end{array}\right],
\end{equation}
where 
\begin{equation*}
\begin{aligned}
\N_{\gamma k} &= \left[\begin{array}{c}
\N_{x\gamma k} \\
\N_{y\gamma k}
\end{array}\right], \quad \text{and}
& &
\Rbargk = \left[\begin{array}{cc}
\Rgk\\
& \Rgk
\end{array}\right].
\end{aligned}
\end{equation*}
Furthermore, a discrete analogue of application of integration by parts to the term $ \int_{\Omega_k}\fnc{V}\nabla\cdot(\lambda\nabla\fnc{U}) \dd{\Omega} $ yields the relation (see \cite[Proposition 1]{yan2018interior}),
\begin{equation} \label{eq:D2 identity}
\D^{(2)}_k=-\H_{k}^{-1}\M_{k}+ \H_{k}^{-1} \sum_{\gamma \subset \Gamma_k} \R_{\gamma k}^{T}\B_{\gamma}\Dgk,
\end{equation}
where $ \M_k $ is a symmetric positive semidefinite matrix given by
\begin{equation} \label{eq:M_k matrix}
\M_{k}=\left[\begin{array}{cc}
\D_{xk}^{T} & \D_{yk}^{T}\end{array}\right]\bar{\H}_k\Lambda_k\left[\begin{array}{c}
\Dxk\\
\Dyk
\end{array}\right],
\end{equation}
and 
\begin{equation*} 
\bar{\H}_{k}=\left[\begin{array}{cc}
\H_{k}\\
& \H_{k}
\end{array}\right].
\end{equation*}
A further decomposition of the $ \D^{(2)}_k $ matrix can be obtained by applying the the SBP property twice:
\begin{proposition}
	If SBP operators on physical elements are constructed as discussed in \cref{sec:Curvilinear Transformation}, then the second derivative operator in \cref{eq:D2 1st form}, which is constructed by applying the first derivative twice, has the decomposition 
	\begin{equation} \label{eq:D2 identity 2}
	\D_{k}^{(2)}=\H_{k}^{-1}\left(\D_{k}^{(2)}\right)^{T}\H_{k}-\H_{k}^{-1}\sum_{\gamma\subset\Gamma_{k}}\D_{\gamma k}^{T}\B_{\gamma}\R_{\gamma k}+\H_{k}^{-1}\sum_{\gamma\subset\Gamma_{k}}\R_{\gamma k}^{T}\B_{\gamma}\D_{\gamma k}.
	\end{equation}
\end{proposition}
\begin{proof}
	Using the SBP property in \cref{eq:SBP property}, we have
	\begin{equation*}
	\begin{aligned}
	\D_{xk}=\H_{k}^{-1}\H_{k}\D_{xk}\H_{k}^{-1}\H_{k}&=\H_{k}^{-1}\Q_{xk}\H_{k}^{-1}\H_{k} = \H_{k}^{-1}\left(-\Q_{xk}^{T}+\E_{xk}\right)\H_{k}^{-1}\H_{k},
	\end{aligned}
	\end{equation*}
	therefore, 
	\begin{equation} \label{eq:Dxk decompostion}
	\D_{xk}=-\H_{k}^{-1}\D_{xk}^{T}\H_{k}+\H_{k}^{-1}\E_{xk}, \quad \text{and} \quad \D_{xk}^{T}=-\H_{k}\D_{xk}\H_{k}^{-1}+\E_{xk}\H_{k}^{-1}.
	\end{equation}
	Furthermore, 
	\begin{equation} \label{eq:Exk Dgk decomposition}
	\begin{bmatrix}
	\E_{xk} & \E_{yk}\end{bmatrix}\Lambda_{k}\begin{bmatrix}
	\D_{xk}\\
	\D_{yk}
	\end{bmatrix} 
	= \sum_{\gamma \subset \Gamma_{k}} \Rgk^T\B_{\gamma}\N_{\gamma k}^T \bar{\R}_{\gamma k} \Lambda_k \begin{bmatrix}
	\D_{xk}\\
	\D_{yk}
	\end{bmatrix} 
	= \sum_{\gamma \subset \Gamma_{k}} \Rgk^T\B_{\gamma} \Dgk.
	\end{equation} 
	Using \cref{eq:Dxk decompostion} and the fact that $ \bar{\H}_k\Lambda_k\bar{\H}^{-1}_k = \Lambda_k$ (since $ \H_k $ and components of $ \Lambda_k $ are diagonal), we arrive at
	\begin{align}
		\D_k^{(2)} &= \begin{bmatrix}
		\D_{xk} & \D_{yk}\end{bmatrix}\Lambda_{k}\begin{bmatrix}
		\D_{xk}\\
		\D_{yk}
		\end{bmatrix} 
		= -\H_{k}^{-1}\begin{bmatrix}
		\D_{xk}^{T} & \D_{yk}^{T}\end{bmatrix}\Lambda_{k}\begin{bmatrix}
		\D_{xk}\\
		\D_{yk}
		\end{bmatrix}\H_{k}
		+\H_{k}^{-1}\begin{bmatrix}
		\E_{xk} & \E_{yk}\end{bmatrix}\Lambda_{k}\begin{bmatrix}
		\D_{xk}\\
		\D_{yk}
		\end{bmatrix} \nonumber\\
		&
		=\H_{k}^{-1}\begin{bmatrix} 
		\D_{xk}^{T} & \D_{yk}^{T}\end{bmatrix}\Lambda_{k}\begin{bmatrix}
		\D_{xk}^{T}\\
		\D_{yk}^{T}
		\end{bmatrix}\H_{k}
		-\H_{k}^{-1}\begin{bmatrix}
		\D_{xk}^{T} & \D_{yk}^{T}\end{bmatrix}\Lambda_{k}\begin{bmatrix}
		\E_{xk}\\
		\E_{yk}
		\end{bmatrix}
		+\H_{k}^{-1}\begin{bmatrix}
		\E_{xk} & \E_{yk}\end{bmatrix}\Lambda_{k}\begin{bmatrix}
		\D_{xk}\\
		\D_{yk}
		\end{bmatrix}.\label{eq:D2 last line}
	\end{align}
	Noting that $ \Lambda_{k} = \Lambda_{k}^T $, $ \E_{xk}=\E_{xk}^T $, $ \E_{yk}=\E_{yk}^T $, and substituting \cref{eq:Exk Dgk decomposition} in \cref{eq:D2 last line} yields the desired result.
\end{proof}
\begin{remark}
	Identity \cref{eq:D2 identity 2} mimics application of integration by parts twice on $ \int_{\Omega}\fnc{V}\nabla\cdot\left(\lambda\nabla\fnc{U}\right) \dd{\Omega} $.
\end{remark}
The SBP-SAT semi-discretization of \cref{eq:diffusion problem} for element $ \Omega_k $ can now be written as 
\begin{equation} \label{eq:SBP-SAT discretization}
\dv{\uhk}{t}=\D^{(2)}_{k}\uhk+\bm{f}_{k}-\H_{k}^{-1}\bm{s}_{k}^{I}(\uhk) -\H_{k}^{-1}\bm{s}_{k}^{B}(\uhk,\bm{u}_{\gamma k}, \bm{w}_{\gamma k}),
\end{equation}
where the interior facet SATs and the boundary SATs are given, respectively, by 
\begin{equation} \label{eq:Interface SATs}
\begin{aligned}
\bm{s}_{k}^{I}(\uhk)&=\sum_{\gamma\subset\Gamma_{k}^{I}}\left[\begin{array}{cc}
\R_{\gamma k}^{T} & \D_{\gamma k}^{T}\end{array}\right]\left[\begin{array}{cc}
\T_{\gamma  k}^{(1)} & \T_{\gamma k}^{(3)}
\\
\T_{\gamma k}^{(2)} & \T_{\gamma k}^{(4)}
\end{array}\right]\left[\begin{array}{c}
\R_{\gamma k}\bm{u}_{h,k}-\R_{\gamma v}\bm{u}_{h,v}
\\
\D_{\gamma k}\bm{u}_{h,k}+\D_{\gamma v}\bm{u}_{h,v}
\end{array}\right]
\\
& \quad +\sum_{\gamma\subset\Gamma_{k}^{I}}
\left\{ \sum_{\epsilon\subset\Gamma_{k}^{I}}\R_{\gamma k}^{T}\T_{\gamma\epsilon k}^{(5)}\left[\R_{\epsilon k}\bm{u}_{h,k}-\R_{\epsilon g}\bm{u}_{h,g}\right]
+\sum_{\delta\subset\Gamma_{v}^{I}}\R_{\gamma k}^{T}\T_{\gamma\delta v}^{(6)}\left[\R_{\delta v}\bm{u}_{h,v}-\R_{\delta q}\bm{u}_{h,q}\right]\right\},
\end{aligned}
\end{equation}
and 
\begin{equation} \label{eq:Boundary SATs}
\begin{aligned}
\bm{s}_{k}^{B}(\uhk,\bm{u}_{\gamma k}, \bm{w}_{\gamma k}) &=\sum_{\gamma\subset\Gamma^{D}}\left[\begin{array}{cc}
\R_{\gamma k}^{T} & \D_{\gamma k}^{T}\end{array}\right]\left[\begin{array}{c}
\T_{\gamma}^{(D)}\\
-\B_{\gamma}
\end{array}\right](\Rgk\uhk  -\bm{u}_{\gamma k})
+\sum_{\gamma\subset\Gamma^{N}}\R_{\gamma k}^{T}\B_{\gamma}\left(\D_{\gamma k}\bm{u}_{h,k}-\bm{w}_{\gamma k}\right),
\end{aligned}
\end{equation}
here, $ \bm{u}_{\gamma k} $ is the restriction of $ \fnc{U} $ on $ S_{\gamma}$, $ \bm{w}_{\gamma k} $ is the restriction of $ \bm{n}\cdot(\lambda\nabla\fnc{U}) $ on $ S_{\gamma}$, $ \epsilon \in \{\epsilon_{1},\epsilon_{2}\} $, $ \delta \in \{\delta_{1},\delta_{2}\} $, and the matrices $\T_{\gamma k}^{(1)}, \T_{\gamma k}^{(2)}, \T_{\gamma k}^{(3)}, \T_{\gamma k}^{(4)}, \T_{\gamma \epsilon k}^{(5)}, \T_{\gamma \delta k}^{(6)} \in \IRtwo{n_f}{n_f}$ are SAT penalty/coefficient matrices. Elements and facets are labeled as shown in \cref{fig:element and facet label}. \blue{To avoid calculating the gradient of the solution, $ [\begin{smallmatrix}	\Dxk\bm{u}_{h,k} \\ \Dyk\bm{u}_{h,k}\end{smallmatrix}] $, multiple times to find terms such as $ \D_{\gamma k}\bm{u}_{h,k} $ in \cref{eq:Interface SATs}, one can compute and store the gradient of the solution in a vector.}

The structure of the \violet{interior facet} SATs given by \cref{eq:Interface SATs} differs from the form considered in \cite{yan2018interior} by the inclusion of the last two terms, which enable the study of wide stencil SATs that couple a target element with second neighbors, \eg, BR1 and LDG type SATs. \blue{We point out that, unlike wide stencil DG fluxes, the boundary SATs do not include extended stencil terms. This facilitates the design of adjoint consistent schemes for problems with non-homogeneous Dirichlet boundary conditions; however, it bears an adverse effect on energy stability of some the DG-based SATs. The connection between the SATs and DG fluxes as well as the stability issues due to the form of the boundary SATs are discussed in \cref{sec:Existing and DG SATs}.} All of the SATs considered in this work have $ \T_{\gamma k}^{(4)}=\bm{0}$. While it is possible to construct SATs with a nonzero $ \T_{\gamma k}^{(4)}$ coefficient, this can decrease the global accuracy of the numerical solution and increase the condition number and stiffness of the arising system matrix \cite{carpenter2010revisiting,eriksson2018finite}. \blue{Indeed, for most of the SBP-SAT discretizations studied in \cref{sec:Numerical Results}, we observed that setting $ \T_{\gamma k}^{(4)}= 1/2 \B_\gamma$ decreases the accuracy of the solution and increases the condition number of the arising system matrix by two to three orders of magnitude}. For the analyses that follow, we do not make the assumption that $ \T_{\gamma k}^{(4)}$ is zero. The assumptions we make regarding the SAT coefficients are stated below.
\begin{assumption}\label{assu:Coefficient matrices}
	For any element $ \Omega_k $ and facets $ a,b\in \{\gamma,\epsilon,\delta\} $, we assume that the coefficient matrices $ \T^{(1)}_{ak} $, $\T^{(2)}_{ak}, \T^{(3)}_{ak}, \T^{(4)}_{ak}$, and $ \T^{(D)}_{a} $ are symmetric, $ \T^{(5)}_{abk} = (\T^{(5)}_{bak}) ^T  $, and $ \T^{(6)}_{abk} = (\T^{(6)}_{bak}) ^T  $.
\end{assumption}
Premultiplying \cref{eq:SBP-SAT discretization} by $ \bm{v}^T_k \H_k $ and employing identity \cref{eq:D2 identity}, the weak form of the SBP-SAT discretization reads
\begin{equation} \label{eq:Weak form 1 elem}
\begin{aligned}
\bm{v}_{k}^{T}\H_{k}\dv{\uhk}{t}&=-\bm{v}_{k}^{T}\M_k\uhk 
+\sum_{\gamma \subset \Gamma_k} \bm{v}_{k}^{T}\R_{\gamma k}^{T}\B_{\gamma}\Dgk \uhk
+ \bm{v}_{k}^{T}\H_{k}\bm{f}_{k} 
-\bm{v}_{k}^{T}\bm{s}_{k}^{I} (\uhk)
-\bm{v}_{k}^{T}\bm{s}_{k}^{B}(\uhk,\bm{u}_{\gamma k}, \bm{w}_{\gamma k}),
\end{aligned}
\end{equation}
for all $ \bm{v}_k \in \IR{n_p}$. Summing \cref{eq:Weak form 1 elem} over all elements yields
\begin{equation} \label{eq:Discretization summed over all elements}
\sum_{\Omega_k\subset {\mathcal T}_h}\bm{v}_{k}^{T}\H_{k}\dv{\uhk}{t} = R_{h}(\bm{u}_{h},\bm{v}), \quad \forall\bm{v}\in\mathbb{R}^{\Sigma n_{e}},
\end{equation}
where 
\begin{align}
		R_{h}(\bm{u}_{h},\bm{v})&=-\sum_{\Omega_{k}\subset\mathcal{T}_{h}}\bm{v}_{k}^{T}\M_{k}\bm{u}_{h,k}
		+\sum_{\Omega_{k}\subset\mathcal{T}_{h}}\bm{v}_{k}^{T}\H_{k}\bm{f}_{k} -\sum_{\gamma\subset\Gamma^{I}}\left[\begin{array}{c}
		\R_{\gamma k}\bm{v}_{k}\\
		\R_{\gamma v}\bm{v}_{v}\\
		\D_{\gamma k}\bm{v}_{k}\\
		\D_{\gamma v}\bm{v}_{v}
		\end{array}\right]^{T}\left[\begin{array}{cccc}
		\T_{\gamma   k}^{(1)} & -\T_{\gamma   k}^{(1)} & \T_{\gamma k}^{(3)}-\B_{\gamma} & \T_{\gamma k}^{(3)}\\
		-\T_{\gamma   v}^{(1)} & \T_{\gamma   v}^{(1)} & \T_{\gamma v}^{(3)} & \T_{\gamma v}^{(3)}-\B_{\gamma}\\
		\T_{\gamma k}^{(2)} & -\T_{\gamma k}^{(2)} & \T_{\gamma k}^{(4)} & \T_{\gamma k}^{(4)}\\
		-\T_{\gamma v}^{(2)} & \T_{\gamma v}^{(2)} & \T_{\gamma v}^{(4)} & \T_{\gamma v}^{(4)}
		\end{array}\right]\left[\begin{array}{c}
		\R_{\gamma k}\bm{u}_{h,k}\\
		\R_{\gamma v}\bm{u}_{h,v}\\
		\D_{\gamma k}\bm{u}_{h,k}\\
		\D_{\gamma v}\bm{u}_{h,v}
		\end{array}\right] 
		\nonumber
		\\& \quad-\sum_{\gamma\subset\Gamma^{I}}\left\{\sum_{\epsilon\subset\Gamma_{k}^{I}} \left[\begin{array}{c}
		\R_{\gamma k}\bm{v}_{k}\\
		\R_{\gamma v}\bm{v}_{v}
		\end{array}\right]^{T}\left[\begin{array}{cc}
		\T_{\gamma\epsilon k}^{(5)} & -\T_{\gamma\epsilon k}^{(5)}\\
		\T_{\gamma\epsilon k}^{(6)} & -\T_{\gamma\epsilon k}^{(6)}
		\end{array}\right]\left[\begin{array}{c}
		\R_{\epsilon k}\bm{u}_{h,k}\\
		\R_{\epsilon g}\bm{u}_{h,g}
		\end{array}\right]
		-\sum_{\delta\subset\Gamma_{v}^{I}}\left[\begin{array}{c}
		\R_{\gamma k}\bm{v}_{k}\\
		\R_{\gamma v}\bm{v}_{v}
		\end{array}\right]^{T}\left[\begin{array}{cc}
		\T_{\gamma\delta v}^{(6)} & -\T_{\gamma\delta v}^{(6)}\\
		\T_{\gamma\delta v}^{(5)} & -\T_{\gamma\delta v}^{(5)}
		\end{array}\right]\left[\begin{array}{c}
		\R_{\delta v}\bm{u}_{h,v}\\
		\R_{\delta q}\bm{u}_{h,q}
		\end{array}\right] \right\}
		\nonumber
		\\& 
		\quad
		-\sum_{\gamma\subset\Gamma^{D}}\left[\begin{array}{c}
		\R_{\gamma k}\bm{v}_{k}\\
		\D_{\gamma k}\bm{v}_{k}
		\end{array}\right]^{T}\left[\begin{array}{cc}
		\T_{\gamma}^{(D)} & -\B_{\gamma}\\
		-\B_{\gamma} & \bm{0}
		\end{array}\right]\left[\begin{array}{c}
		\R_{\gamma k}\bm{u}_{h,k}-\bm{u}_{\gamma k}\\
		\D_{\gamma k}\bm{u}_{h,k}
		\end{array}\right]+\sum_{\gamma\subset\Gamma^{N}}\bm{v}_{k}^{T}\R_{\gamma k}^{T}\B_{\gamma}\bm{w}_{\gamma k}.\label{eq:Residual 1st form}
\end{align}
Instead of adding the terms responsible for extending the stencil (terms containing coefficients $ \T^{(5)} $ and $ \T^{(6)} $) facet by facet, we can add them element by element. That is, we regroup facet terms associated with  $ \T^{(5)} $ and $ \T^{(6)} $ by element and rewrite the residual as 
\begin{equation} \label{eq:Residual 2nd form}
\begin{aligned}
	R_{h}(\bm{u}_{h},\bm{v})&=-\sum_{\Omega_{k}\subset\mathcal{T}_{h}}\bm{v}_{k}^{T}\M_{k}\bm{u}_{h,k}+\sum_{\Omega_{k}\subset\mathcal{T}_{h}}\bm{v}_{k}^{T}\H_{k}\bm{f}_{k} -\sum_{\gamma \subset \Gamma^I} (\bm{v}^{\star})^T \T^{\star} \bm{u}_h^{\star} 
	-\sum_{\Omega_{k}\subset\mathcal{T}_{h}}\sum_{\gamma,\epsilon\subset \Gamma_k^I} (\bm{v}^{\diamond})^T \T^{\diamond} \bm{u}_h^{\diamond}\\ &
	\quad -\sum_{\gamma\subset\Gamma^{D}}\left[\begin{array}{c}
	\R_{\gamma k}\bm{v}_{k}\\
	\D_{\gamma k}\bm{v}_{k}
	\end{array}\right]^{T}\left[\begin{array}{cc}
	\T_{\gamma}^{(D)} & -\B_{\gamma}\\
	-\B_{\gamma} & \bm{0}
	\end{array}\right]\left[\begin{array}{c}
	\R_{\gamma k}\bm{u}_{h,k}-\bm{u}_{\gamma k}\\
	\D_{\gamma k}\bm{u}_{h,k}
	\end{array}\right] +\sum_{\gamma\subset\Gamma^{N}}\bm{v}_{k}^{T}\R_{\gamma k}^{T}\B_{\gamma}\bm{w}_{\gamma k}, 
	\end{aligned}
	\end{equation}
	where $ \sum_{\gamma \subset \Gamma^I} (\bm{v}^{\star})^T \T^{\star} \bm{u}_h^{\star}  $ is equal to the third term on the RHS of \cref{eq:Residual 1st form}, and 
	\begin{equation} 
	\begin{aligned}
	\sum_{\gamma,\epsilon\subset \Gamma_k^I} (\bm{v}^{\diamond})^T \T^{\diamond} \bm{u}_h^{\diamond} &= 
	\left[\begin{array}{c}
	\R_{\gamma k}\bm{v}_{k}\\
	\R_{\gamma v}\bm{v}_{v}\\
	\R_{\epsilon_{1}k}\bm{v}_{k}\\
	\R_{\epsilon_{1}g_{1}}\bm{v}_{g_{1}}\\
	\R_{\epsilon_{2}k}\bm{v}_{k}\\
	\R_{\epsilon_{2}g_{2}}\bm{v}_{g_{2}}
	\end{array}\right]^T
	\left[\begin{array}{cccccc}
	\bm{0} & \bm{0} & \T_{\gamma\epsilon_{1}k}^{(5)} & -\T_{\gamma\epsilon_{1}k}^{(5)} & \T_{\gamma\epsilon_{2}k}^{(5)} & -\T_{\gamma\epsilon_{2}k}^{(5)}\\
	\bm{0} & \bm{0} & \T_{\gamma\epsilon_{1}k}^{(6)} & -\T_{\gamma\epsilon_{1}k}^{(6)} & \T_{\gamma\epsilon_{2}k}^{(6)} & -\T_{\gamma\epsilon_{2}k}^{(6)}\\
	\T_{\epsilon_{1}\gamma k}^{(5)} & -\T_{\epsilon_{1}\gamma k}^{(5)} & \bm{0} & \bm{0} & \T_{\epsilon_{1}\epsilon_{2}k}^{(5)} & -\T_{\epsilon_{1}\epsilon_{2}k}^{(5)}\\
	\T_{\epsilon_{1}\gamma k}^{(6)} & -\T_{\epsilon_{1}\gamma k}^{(6)} & \bm{0} & \bm{0} & \T_{\epsilon_{1}\epsilon_{2}k}^{(6)} & -\T_{\epsilon_{1}\epsilon_{2}k}^{(6)}\\
	\T_{\epsilon_{2}\gamma k}^{(5)} & -\T_{\epsilon_{2}\gamma k}^{(5)} & \T_{\epsilon_{2}\epsilon_{1}k}^{(5)} & -\T_{\epsilon_{2}\epsilon_{1}k}^{(5)} & \bm{0} & \bm{0}\\
	\T_{\epsilon_{2}\gamma k}^{(6)} & -\T_{\epsilon_{2}\gamma k}^{(6)} & \T_{\epsilon_{2}\epsilon_{1}k}^{(6)} & -\T_{\epsilon_{2}\epsilon_{1}k}^{(6)} & \bm{0} & \bm{0}
	\end{array}\right]
	\left[\begin{array}{c}
	\R_{\gamma k}\bm{u}_{h,k}\\
	\R_{\gamma v}\bm{u}_{h,v}\\
	\R_{\epsilon_{1}k}\bm{u}_{h,k}\\
	\R_{\epsilon_{1}g_{1}}\bm{u}_{h,g_{1}}\\
	\R_{\epsilon_{2}k}\bm{u}_{h,k}\\
	\R_{\epsilon_{2}g_{2}}\bm{u}_{h,g_{2}}
	\end{array}\right].
\end{aligned}
\end{equation}

Yet another form of the residual, and one that is important for energy analysis, is obtained by employing the ``borrowing trick" of \cite{carpenter1999stable} which is generalized for multidimensional SBP operators in \cite{yan2018interior}. The approach allows to find conditions for energy stability by writing the volume term on the RHS of \cref{eq:Residual 1st form}, $ -\sum_{\Omega_{k}\subset\mathcal{T}_{h}}\bm{v}_{k}^{T}\M_{k}\bm{u}_{h,k}, $ as a facet contribution. The following lemma is as an extension of Lemma 1 in \cite{yan2018interior}.
\begin{lemma}
	Given a facet based weight $ \alpha_{\gamma k} > 0$ satisfying $ \sum_{\gamma \subset \Gamma_k} \alpha_{\gamma k} = 1$ for each facet $ \gamma $, the residual of the SBP-SAT discretization for the homogeneous version  of \cref{eq:diffusion problem}, \ie, $ \fnc{F}=0$, $\fnc{U}_D = 0 $ , and $ \fnc{U}_N =0$, can be written as
	\begin{equation} \label{eq:Residual 3rd form}
	\begin{aligned}
	R_{h}(\bm{u}_{h},\bm{u}_{h})&=-\sum_{\gamma\subset\Gamma^{I}} \begin{bmatrix}
	\D_{\gamma k}\bm{u}_{h,k}\\
	\D_{\gamma v}\bm{u}_{h,v}
	\end{bmatrix}^{T}\begin{bmatrix}
	\T_{\gamma k}^{(4)} & \T_{\gamma k}^{(4)}\\
	\T_{\gamma v}^{(4)} & \T_{\gamma v}^{(4)}
	\end{bmatrix}
	\begin{bmatrix}
	\D_{\gamma k}\bm{u}_{h,k}\\
	\D_{\gamma v}\bm{u}_{h,v}
	\end{bmatrix} 
	-\sum_{\gamma\subset\Gamma^{I}} X_1
	-\sum_{\Omega_{k}\in\mathcal{T}_{h}}\sum_{\gamma,\epsilon\in\Gamma^I_k} X_2
	\\
	&\quad
	-\sum_{\gamma\subset\Gamma^{D}}\begin{bmatrix}
	\R_{\gamma}\bm{v}_{k}\\
	\F_{k}\bm{v}_{k}
	\end{bmatrix}^{T}
	\begin{bmatrix}
	\T_{\gamma}^{(D)} & -\B_{\gamma}\C_{\gamma k}\\
	-\C_{\gamma k}^{T}\B_{\gamma} & \alpha_{\gamma k}\Lambda_{k}^{*}
	\end{bmatrix}
	\begin{bmatrix}
	\R_{\gamma k}\bm{u}_{h,k}\\
	\F_{k}\bm{u}_{h,k}
	\end{bmatrix},
	\end{aligned}
	\end{equation}	
	where 
	\begin{equation*}
	\begin{aligned}
	X_1 &=\begin{bmatrix}
	\R_{\gamma k}\bm{u}_{h,k}\\
	\R_{\gamma v}\bm{u}_{h,v}\\
	\F_{k}\bm{u}_{h,k}\\
	\F_{v}\bm{u}_{h,v}
	\end{bmatrix}^{T}
	\begin{bmatrix}
	\frac{1}{2}\T_{\gamma k}^{(1)} & -\frac{1}{2}\T_{\gamma  k}^{(1)} & (\T_{\gamma k}^{(3)}-\B_{\gamma})\C_{\gamma k} & \T_{\gamma k}^{(3)}\C_{\gamma v}\\
	-\frac{1}{2}\T_{\gamma  v}^{(1)} & \frac{1}{2}\T_{\gamma  v}^{(1)} & \T_{\gamma v}^{(3)}\C_{\gamma k} & (\T_{\gamma v}^{(3)}-\B_{\gamma})\C_{\gamma v}\\
	\C_{\gamma k}^{T}\T_{\gamma k}^{(2)} & -\C_{\gamma k}^{T}\T_{\gamma k}^{(2)} & \alpha_{\gamma k}\Lambda_{k}^{*} & \bm{0}\\
	-\C_{\gamma v}^{T}\T_{\gamma v}^{(2)} & \C_{\gamma v}^{T}\T_{\gamma v}^{(2)} & \bm{0} & \alpha_{\gamma v}\Lambda_{v}^{*}
	\end{bmatrix}
	\begin{bmatrix}
	\R_{\gamma k}\bm{u}_{h,k}\\
	\R_{\gamma v}\bm{u}_{h,v}\\
	\F_{k}\bm{u}_{h,k}\\
	\F_{v}\bm{u}_{h,v}
	\end{bmatrix}
	\end{aligned}
	\end{equation*}
	\begin{equation*}
	\begin{aligned}
	X_2 &=
	\begin{bmatrix}
	\R_{\gamma k}\bm{u}_{h,k}\\
	\R_{\gamma v}\bm{u}_{h,v}\\
	\R_{\epsilon_{1}k}\bm{u}_{h,k}\\
	\R_{\epsilon_{1}g_{1}}\bm{u}_{h,g_{1}}
	\end{bmatrix}^{T}
	\begin{bmatrix}
	\frac{1}{8}\T_{\gamma k}^{(1)} & -\frac{1}{8}\T_{\gamma k}^{(1)} & \T_{\gamma\epsilon_{1}k}^{(5)} & -\T_{\gamma\epsilon_{1}k}^{(5)}\\
	-\frac{1}{8}\T_{\gamma v}^{(1)} & \frac{1}{8}\T_{\gamma v}^{(1)} & \T_{\gamma\epsilon_{1}k}^{(6)} & -\T_{\gamma\epsilon_{1}k}^{(6)}\\
	\T_{\epsilon_{1}\gamma k}^{(5)} & -\T_{\epsilon_{1}\gamma k}^{(5)} & \frac{1}{8}\T_{\epsilon_{1}k}^{(1)} & -\frac{1}{8}\T_{\epsilon_{1}k}^{(1)}\\
	\T_{\epsilon_{1}\gamma k}^{(6)} & -\T_{\epsilon_{1}\gamma k}^{(6)} & -\frac{1}{8}\T_{\epsilon_{1}g_1}^{(1)} & \frac{1}{8}\T_{\epsilon_{1}g_1}^{(1)}
	\end{bmatrix}
	\begin{bmatrix}
	\R_{\gamma k}\bm{u}_{h,k}\\
	\R_{\gamma v}\bm{u}_{h,v}\\
	\R_{\epsilon_{1}k}\bm{u}_{h,k}\\
	\R_{\epsilon_{1}g_{1}}\bm{u}_{h,g_{1}}
	\end{bmatrix}
	\\
	&\quad
	+\begin{bmatrix}
	\R_{\gamma k}\bm{u}_{h,k}\\
	\R_{\gamma v}\bm{u}_{h,v}\\
	\R_{\epsilon_{2}k}\bm{u}_{h,k}\\
	\R_{\epsilon_{2}g_{2}}\bm{u}_{h,g_{2}}
	\end{bmatrix}^{T}
	\begin{bmatrix}
	\frac{1}{8}\T_{\gamma k}^{(1)} & -\frac{1}{8}\T_{\gamma k}^{(1)} & \T_{\gamma\epsilon_{2}k}^{(5)} & -\T_{\gamma\epsilon_{2}k}^{(5)}\\
	-\frac{1}{8}\T_{\gamma v}^{(1)} & \frac{1}{8}\T_{\gamma v}^{(1)} & \T_{\gamma\epsilon_{2}k}^{(6)} & -\T_{\gamma\epsilon_{2}k}^{(6)}\\
	\T_{\epsilon_{2}\gamma k}^{(5)} & -\T_{\epsilon_{2}\gamma k}^{(5)} & \frac{1}{8}\T_{\epsilon_{2} k}^{(1)} & -\frac{1}{8}\T_{\epsilon_{2} k}^{(1)}\\
	\T_{\epsilon_{2}\gamma k}^{(6)} & -\T_{\epsilon_{2}\gamma k}^{(6)} & -\frac{1}{8}\T_{\epsilon_{2}g_2}^{(1)} & \frac{1}{8}\T_{\epsilon_{2}g_2}^{(1)}
	\end{bmatrix}
	\begin{bmatrix}
	\R_{\gamma k}\bm{u}_{h,k}\\
	\R_{\gamma v}\bm{u}_{h,v}\\
	\R_{\epsilon_{2}k}\bm{u}_{h,k}\\
	\R_{\epsilon_{2}g_{2}}\bm{u}_{h,g_{2}}
	\end{bmatrix}\\
	&\quad
	+\begin{bmatrix}
	\R_{\epsilon_{1}k}\bm{u}_{h,k}\\
	\R_{\epsilon_{1}g_{1}}\bm{u}_{h,g_{1}}\\
	\R_{\epsilon_{2}k}\bm{u}_{h,k}\\
	\R_{\epsilon_{2}g_{2}}\bm{u}_{h,g_{2}}
	\end{bmatrix}^{T}
	\begin{bmatrix}
	\frac{1}{8}\T_{\epsilon_{1}k}^{(1)} & -\frac{1}{8}\T_{\epsilon_{1}k}^{(1)} & \T_{\epsilon_{1}\epsilon_{2}k}^{(5)} & -\T_{\epsilon_{1}\epsilon_{2}k}^{(5)}\\
	-\frac{1}{8}\T_{\epsilon_{1}g_1}^{(1)} & \frac{1}{8}\T_{\epsilon_{1}g_1}^{(1)} & \T_{\epsilon_{1}\epsilon_{2}k}^{(6)} & -\T_{\epsilon_{1}\epsilon_{2}k}^{(6)}\\
	\T_{\epsilon_{2}\epsilon_{1}k}^{(5)} & -\T_{\epsilon_{2}\epsilon_{1}k}^{(5)} & \frac{1}{8}\T_{\epsilon_{2} k}^{(1)} & -\frac{1}{8}\T_{\epsilon_{2} k}^{(1)}\\
	\T_{\epsilon_{2}\epsilon_{1}k}^{(6)} & -\T_{\epsilon_{2}\epsilon_{1}k}^{(6)} & -\frac{1}{8}\T_{\epsilon_{2}g_2}^{(1)} & \frac{1}{8}\T_{\epsilon_{2}g_2}^{(1)}
	\end{bmatrix}
	\begin{bmatrix}
	\R_{\epsilon_{1}k}\bm{u}_{h,k}\\
	\R_{\epsilon_{1}g_{1}}\bm{u}_{h,g_{1}}\\
	\R_{\epsilon_{2}k}\bm{u}_{h,k}\\
	\R_{\epsilon_{2}g_{2}}\bm{u}_{h,g_{2}}
	\end{bmatrix}
	\end{aligned}
	\end{equation*}
	and
\begin{equation*}
	\begin{aligned}
	\F_{k}&=\left[\begin{array}{cc}
	\Lambda_{xx} & \Lambda_{xy}\\
	\Lambda_{yx} & \Lambda_{yy}
	\end{array}\right]_{k}\left[\begin{array}{c}
	\D_{xk}\\
	\D_{yk}
	\end{array}\right],&\Lambda_{k}^{*}&=\left[\begin{array}{cc}
	\Lambda_{xx} & \Lambda_{xy}\\
	\Lambda_{yx} & \Lambda_{yy}
	\end{array}\right]_{k}^{-1}\left[\begin{array}{cc}
	\H_{k}\\
	& \H_{k}
	\end{array}\right],&\C_{\gamma k}&=\left[\begin{array}{cc}
	\N_{x \gamma k}\R_{\gamma k} & \N_{y \gamma k}\R_{\gamma k}\end{array}\right],\\\F_{v}&=\left[\begin{array}{cc}
	\Lambda_{xx} & \Lambda_{xy}\\
	\Lambda_{yx} & \Lambda_{yy}
	\end{array}\right]_{v}\left[\begin{array}{c}
	\D_{xv}\\
	\D_{yv}
	\end{array}\right],&\Lambda_{v}^{*}&=\left[\begin{array}{cc}
	\Lambda_{xx} & \Lambda_{xy}\\
	\Lambda_{yx} & \Lambda_{yy}
	\end{array}\right]_{v}^{-1}\left[\begin{array}{cc}
	\H_{v}\\
	& \H_{v}
	\end{array}\right],&\C_{\gamma v}&=\left[\begin{array}{cc}
	\N_{x\gamma v}\R_{\gamma v} & \N_{y\gamma v}\R_{\gamma v}\end{array}\right],
	\end{aligned}
\end{equation*}
\end{lemma}
\begin{proof}
	Since the result follows from simple algebraic manipulations, the complete proof is omitted but we state some of the steps. We made use of the decomposition provided in \cite{yan2018interior},
	\begin{equation}
	\M_{k} = \sum_{\gamma \subset \Gamma_{k}} \alpha_{\gamma k}\F_{k}^T\Lambda_k^* \F_{k},
	\end{equation}
	applied the relations $ \C_{\gamma k}\F_{k}=\D_{\gamma k} $, $ \C_{\gamma v}\F_{v} = \D_{\gamma v} $, subtracted $ \frac{1}{2}\T^{(1)} $ terms from the first $ 2\times 2 $ block of $ \T^{\star}$ and added $ \frac{1}{4}\T^{(1)} $ terms in the $ 2\times 2 $ diagonal blocks of $ \T^{(\diamond)} $, which is then decomposed into the three terms in $ X_2 $. Note that $ \T^{(\diamond)} $ is summed element by element; therefore, we recover  $ \frac{1}{2}\T^{(1)} $ terms at each interior facet due to contributions from two abutting elements.    
\end{proof}

\section{Properties of the SBP-SAT discretization} \label{sec:Properties of SBP-SAT discretization}
In this section, some numerical properties of the SBP-SAT discretization given in \cref{eq:SBP-SAT discretization,eq:Interface SATs,eq:Boundary SATs} are analyzed. We will make use of the three equivalent forms of the residual in \cref{eq:Residual 1st form}, \cref{eq:Residual 2nd form}, and \cref{eq:Residual 3rd form} depending on the property under consideration. 

\subsection{Consistency} Consistency is a requirement that the SBP-SAT discretization of the primal problem represents the continuous PDE accurately. Consider the steady version of the model problem \cref{eq:diffusion problem}; then consistency requires that the SBP-SAT discretization be at least approximately satisfied by the exact solution \cite{arnold2015stability}. More precisely, we require that 
\begin{equation} 
\lim_{h\rightarrow0}\sum_{\Omega_{k}\in\mathcal{T}_{h}}\norm{L_{h,k}(\bm{u}_k)-\bm{f}_k}_{\H_k}=0,
\end{equation}
where $ \bm{u}_k \in \IR{n_p} $ is a grid function representing the exact solution, $ \fnc{U}_k\in\cont{p+1}$, 
\begin{equation}\label{eq:Lhk}
L_{h,k}( \bm{u}_k)= - \D_{k}^{(2)}\bm{u}_k 
+  \H^{-1}_k\bm{s}_k^I(\bm{u}_k) 
+  \H^{-1}_k\bm{s}_k^B(\bm{u}_k,\bm{u}_{\gamma k}, \bm{w}_{\gamma k}) 
\end{equation}
is the discrete counterpart of the linear operator $ L_k $ applied on $ \fnc{U}_k $, and $ \norm{\cdot}_{\H_k} $ is the norm defined by $ \H_k $ matrix on $ \IR{n_p} $.
\begin{theorem} \label{thm:Consistency}
	Let \cref{assu: mapping} hold, the SBP operators be constructed as described in \cref{sec:Curvilinear Transformation}, the solution to the PDE \cref{eq:diffusion problem} be $ \fnc{U}\in\fnc{C}^{p+1}(\Omega) $, and $ \bm{u}_k\in\IR{n_p} $ be a grid function representing $ \fnc{U}_k $. Then, the discretization \cref{eq:SBP-SAT discretization,eq:Interface SATs,eq:Boundary SATs} is a consistent approximation of the diffusion problem given by \cref{eq:diffusion problem}, and $\norm{L_{h,k} \bm{u}_k - \bm{f}_k}_\infty = \fnc{O}(h^{p-1})$.
\end{theorem}
\begin{proof}
	The result is a simple consequence of \cref{thm:Accuracy of Dx}, and the accuracy of the extrapolation matrix. It follows by substituting $ \bm{u} $ in place of $ \bm{u}_h $ in \cref{eq:SBP-SAT discretization,eq:Interface SATs,eq:Boundary SATs}, and noting that $ \D_{k}^{(2)} \bm{u}_k = [\nabla\cdot(\lambda\nabla\fnc{U})]_{S_{k}} + \order{h^{p-1}}$, $ \Dgk\bm{u}_k = -\Dgv\bm{u}_v = [\bm{n}_k\cdot(\lambda\nabla\fnc{U})]_{S_\gamma} +\order{h^{p}} $, and the extrapolation matrices are order $ p+1 $ accurate, \eg, $ \Rgk\bm{u}_k = \fnc{U}|_{S_\gamma} + \order*{h^{p+1}} $. 
\end{proof}
\subsection{Conservation} A conservative discretization needs to satisfy the divergence theorem discretely. Multiplying \cref{eq:diffusion problem} by a test function, $ \fnc{V}\in H^{2}(\Omega) $, and integrating by parts we find 
\begin{equation} \label{eq:Conservation IBP}
\int_{\Omega}\fnc{V}\pdv{\fnc{U}}{t} \dd\Omega=-\int_{\Omega}\nabla \fnc{V}\cdot(\lambda\nabla \fnc{U})\dd\Omega+\int_{\Gamma^B}\fnc{V}(\lambda\nabla \fnc{U})\cdot \bm{n}\dd\Gamma+\int_{\Omega} \fnc{V}\fnc{F}\dd\Omega,
\end{equation}
which is equivalent to applying the divergence theorem if we set $ \fnc{V}=1 $. Thus, for a problem with no source term and $ \bm{v} = \bm{1} $, we require all except the boundary terms on the RHS of \cref{eq:Residual 2nd form} to vanish. 
\begin{theorem}\label{thm:Conservation}
	Let \cref{assu: mapping} hold and the metric terms be evaluated exactly, then the SBP-SAT discretization given in \cref{eq:SBP-SAT discretization,eq:Interface SATs,eq:Boundary SATs} is conservative if the penalty matrices satisfy
	\begin{equation} \label{eq:Conditions for conservation}
	\begin{aligned}
	\T_{\gamma   k}^{(1)}&=\T_{\gamma   v}^{(1)},
	&&
	\T_{\gamma k}^{(3)}+\T_{\gamma v}^{(3)}=\B_{\gamma},
	&&
	\T_{abk}^{(5)}=-\T_{abk}^{(6)},
	\end{aligned}
	\end{equation}
	where $ a,b\in\left\{ \gamma,\epsilon_{1},\epsilon_{2} \right\}$. 
\end{theorem}
\begin{proof}
	In \cref{eq:Discretization summed over all elements} and \cref{eq:Residual 2nd form}, we set $ \bm{v} = \bm {1} $, $ \bm{f}_k =\bm{0} $. Applying $ \T_{\gamma   v}^{(1)} =\T_{\gamma   k}^{(1)}$, and $ \T_{\gamma k}^{(3)} - \B_{\gamma} = -\T_{\gamma v}^{(3)} $ in $ \T^{(\star)} $, $ \T_{abk}^{(5)}=-\T_{abk}^{(6)} $ for $ a,b\in\left\{ \gamma,\epsilon_{1},\epsilon_{2}\right\}$ in $ \T^{(\diamond)} $, and using the exactness of the extrapolation matrix for constants along with the operations $ \Dxk \bm{1} = \Dyk\bm{1}=\bm{0}$ and $\Dgk \bm{1} = \Dgv \bm{1} = \bm{0} $, we obtain
	$ \sum_{\Omega_{k}\subset\mathcal{T}_{h}}\bm{v}_{k}^{T}\M_{k}\bm{u}_{k} = \bm{0}$, $\sum_{\gamma \subset \Gamma^I} (\bm{v}^{\star})^T \T^{\star} \bm{u}^{\star} = \bm{0}$, and $\sum_{\gamma,\epsilon\in \Gamma_k^I} (\bm{v}^{\diamond})^T \T^{\diamond} \bm{u}^{\diamond} = \bm{0}$.
	Therefore, \cref{eq:Discretization summed over all elements} becomes
	\begin{align}
	\sum_{\Omega_k\subset {\mathcal T}_h}\bm{1}^{T}\H_{k}\dv{\uhk}{t} =
	\sum_{\gamma\subset\Gamma^{N}}\bm{1}^{T}\Rgk^T\B_{\gamma}\bm{w}_{\gamma k}  
	+\sum_{\gamma\subset\Gamma^{D}}\qty[\bm{1}^{T}\Rgk^T\B_{\gamma}\D_{\gamma k}\bm{u}_{h,k}
	-\bm{1}^{T}\Rgk^T\T_{\gamma}^{(D)}\left(\R_{\gamma k}\bm{u}_{h,k}-\bm{u}_{\gamma k}\right)],
	\end{align}
	\ie, all interface terms in the residual vanish and $ \sum_{\Omega_k\subset {\mathcal T}_h}\bm{1}^{T}\H_{k}\dv{\uhk}{t} $ depends only on boundary terms, as desired.
\end{proof}
\subsection{Adjoint Consistency} \label{sec:Adjoint Consistency}
Adjoint consistency refers to an accurate discrete representation of the continuous adjoint problem, \ie, the exact solution to the continuous adjoint problem \cref{eq:Adjoint problem} needs to satisfy
\begin{equation} \label{eq:Adjoint consitstency definition}
\lim_{h\rightarrow0}\sum_{\Omega_{k}\subset\mathcal{T}_{h}}\norm{L_{h,k}^* (\bm{\psi})-\bm{g}_k}_{\H_k}=0,
\end{equation}
where $ L_{h,k}^* $ is the discrete adjoint operator (see \cite{hicken2011superconvergent,hicken2014dual} for similar definitions). The discretization of the linear functional \cref{eq:Functional} needs to be modified in a consistent manner to attain an adjoint consistent discretization \cite{hartmann2007adjoint}. One possible modification is
\begin{equation}\label{eq:Functional discrete}
I_h(\bm{u}_h) = \sum_{\Omega_{k}\subset\mathcal{T}_{h}} \bm{g}^T_k \H_k \uhk - \sum_{\gamma \subset \Gamma^D} \bm{\psi}_{\gamma k}^T\B_{\gamma}\Dgk \uhk 
+ \sum_{\gamma \subset \Gamma^N} \bm{z}_{\gamma k}^T\B_{\gamma}\Rgk\uhk + \sum_{\gamma \subset \Gamma^D} \bm{\psi}_{\gamma k}^T\T_{\gamma}^{(D)}(\Rgk\uhk - \bm{u}_{\gamma k}),
\end{equation} 
\violet{where $ \bm{\psi}_{\gamma k} $ and $ \bm{z}_{\gamma k} $ are restrictions of $ \psi $ and $ \bm{n}\cdot(\lambda\nabla\psi) $ on $ S_{\gamma}$, respectively}. The last term in \cref{eq:Functional discrete} is added for adjoint consistency \cite{hartmann2007adjoint,hartmann2013higher,yan2018interior,hicken2011superconvergent,hicken2012output}. Similarly, we modify the discretization of another form of the functional that is given by the last equality in \cref{eq:Adjoint relation} as
\begin{equation}\label{eq:Functional discrete 2}
I_h(\bm{\psi}_h)
= \sum_{\Omega_{k}\subset\mathcal{T}_{h}} \bm{f}_k^{T} \H_k \bm{\psi}_{h,k}
- \sum_{\gamma \subset \Gamma^D} \bm{u}_{\gamma k}^T\B_{\gamma}\Dgk \bm{\psi}_{h,k} 
+ \sum_{\gamma \subset \Gamma^N} \bm{w}_{\gamma k}^T\B_{\gamma}\Rgk\bm{\psi}_{h,k}
+ \sum_{\gamma \subset \Gamma^D} \bm{u}_{\gamma k}^T\T_{\gamma}^{(D)}(\violet{\Rgk\bm{\psi}_{h,k}} - \bm{\psi}_{\gamma k}).
\end{equation}
\blue{A general procedure to modify the functional for adjoint consistency of discretizations of problems with non-homogeneous Dirichlet boundary conditions is discussed in \cite{hartmann2013higher}. If the boundary SATs contain extended stencil terms, it is not clear whether a similar modification is applicable to retain adjoint consistency.}

\violet{To find the discrete adjoint operator, we require that $ I_h(\bm{u}_h) - I_h(\bm{\psi}_h) =0$, which is a discrete analogue of $ \fnc{I}(\fnc{U})-\fnc{I}(\fnc{\psi})=0 $}. Subtracting $ \sum_{\Omega_{k}\subset\mathcal{T}_{h}}\bm{\psi}^T_{h,k}\H_k(L_{h,k}(\uhk)-\bm{f}_k) = 0 $ from \cref{eq:Functional discrete} gives
\begin{equation}
\begin{aligned}
I_h(\bm{u}_h)
&= \sum_{\Omega_{k}\subset\mathcal{T}_{h}} \bm{g}^T_k \H_k \uhk 
- \sum_{\gamma \subset \Gamma^D} \bm{\psi}_{\gamma k}^T\B_{\gamma}\Dgk \uhk 
+ \sum_{\gamma \subset \Gamma^N} \bm{z}_{\gamma k}^T\B_{\gamma}\Rgk\uhk\\
&\quad+ \sum_{\gamma \subset \Gamma^D} \bm{\psi}_{\gamma k}^T\T_{\gamma}^{(D)}(\Rgk\uhk - \bm{u}_{\gamma k}) - \sum_{\Omega_{k}\subset\mathcal{T}_{h}}\bm{\psi}^T_{h,k}\H_k( L_{h,k}(\uhk) - \bm{f}_k).\\
\end{aligned}
\end{equation} 
Rearranging, adding, and subtracting terms we have 
\begin{equation}\label{eq:Functional discrete 3}
\begin{aligned}
I_h(\bm{u}_h)
&=\sum_{\Omega_{k}\subset\mathcal{T}_{h}}\bm{\psi}^T_{h,k}\H_k\bm{f}_k
-\sum_{\gamma \subset \Gamma^D}\bm{u}_{\gamma k}^T\B_{\gamma}\Dgk\bm{\psi}_{h,k}
+\sum_{\gamma \subset \Gamma^N}\bm{w}_{\gamma k}^T\B_{\gamma}\Rgk \bm{\psi}_{h,k}
+ \sum_{\gamma \subset \Gamma^D} \bm{u}_{\gamma k}^T\T_{\gamma}^{(D)}(\Rgk\psi_{\gamma k} 
- \bm{\psi}_{\gamma k})
\\
&
\quad
-\sum_{\Omega_{k}\subset\mathcal{T}_{h}}\bm{\psi}^T_{h,k}\H_k L_{h,k}(\uhk)
+\sum_{\Omega_{k}\subset\mathcal{T}_{h}}\uhk^T\H_k\bm{g}_k
+\sum_{\gamma \subset \Gamma^D}\bm{u}_{\gamma k}^T\B_{\gamma}\Dgk\bm{\psi}_{h,k}
-\sum_{\gamma \subset \Gamma^D} \bm{u}_{\gamma k}^T\T_{\gamma}^{(D)}(\Rgk\bm{\psi}_{h,k} - \bm{\psi}_{\gamma k})
\\
&
\quad
-\sum_{\gamma \subset \Gamma^N}\bm{w}_{\gamma k}^T\B_{\gamma}\Rgk \bm{\psi}_{h,k}	
- \sum_{\gamma \subset \Gamma^D} \bm{\psi}_{\gamma k}^T\B_{\gamma}\Dgk \uhk 
+ \sum_{\gamma \subset \Gamma^N} \bm{z}_{\gamma k}^T\B_{\gamma}\Rgk\uhk
+ \sum_{\gamma \subset \Gamma^D} \bm{\psi}_{\gamma k}^T\T_{\gamma}^{(D)}(\Rgk\uhk - \bm{u}_{\gamma k}),
\end{aligned}
\end{equation}
where the sum of the first four terms on the RHS is equal to \violet{$ I_h(\bm{\psi}_h) $} due to \cref{eq:Functional discrete 2}. \violet{Rearranging terms, we find}
\begin{equation} \label{eq:Vanishing term for adjoint consistency}
 -\sum_{\Omega_{k}\subset\mathcal{T}_{h}}\bm{\psi}^T_{h,k}\H_k L_{h,k}(\uhk) 
+ \sum_{\Omega_{k}\subset\mathcal{T}_{h}}\uhk^T\H_k\bm{g}_k + B_{\rm terms} + {I_h(\bm{\psi}_h) - I_h(\bm{u}_h)} = 0,
\end{equation}
where $ B_{\rm terms} $ is the sum of the boundary terms in the last two lines of \cref{eq:Functional discrete 3}. Using \cref{eq:Lhk} and \cref{eq:D2 identity 2} we write
\begin{align}
	-\bm{\psi}^T_{h,k}\H_k L_{h,k}( \bm{u}_{h,k}) &=  \bm{\psi}^T_{h,k} \H_k \D_{k}^{(2)}\uhk
	-  \bm{\psi}^T_{h,k}\bm{s}_k^I(\bm{u}_{h,k}) 
	-  \bm{\psi}^T_{h,k}\bm{s}_k^B(\bm{u}_{h,k},\bm{u}_{\gamma k}, \bm{w}_{\gamma k}) 
	\nonumber\\
	\quad
	&=\bm{\psi}^T_{h,k}\left(\D_{k}^{(2)}\right)^{T}\H_{k}\uhk
	-\sum_{\gamma\subset\Gamma_{k}}\bm{\psi}^T_{h,k}\D_{\gamma k}^{T}\B_{\gamma}\R_{\gamma k}\uhk
	+\sum_{\gamma\subset\Gamma_{k}}\bm{\psi}^T_{h,k}\R_{\gamma k}^{T}\B_{\gamma}\D_{\gamma k}\uhk 
	\nonumber \\
	&\quad
	-  \bm{\psi}^T_{h,k}\bm{s}_k^I(\bm{u}_{h,k}) 
	-  \bm{\psi}^T_{h,k}\bm{s}_k^B(\bm{u}_{h,k},\bm{u}_{\gamma k}, \bm{w}_{\gamma k}).		
\end{align} 
Summing over all elements and transposing the result, we find
\begin{equation} \label{eq:Sum Lhk}
\begin{aligned}
	-\sum_{\Omega_{k}\subset\mathcal{T}_{h}}\left(\bm{\psi}^T_{h,k}\H_k L_{h,k}( \bm{u}_{h,k})\right)^T &=
	\sum_{\Omega_{k}\subset\mathcal{T}_{h}}\uhk^T\H_k\D_{k}^{(2)}\bm{\psi}_{h,k} 
	-\sum_{\gamma \subset \Gamma^I} (\bm{u}_h^{\star})^T \tilde{\T}^{\star} \psi_h^{\star} 
	-\sum_{\Omega_{k}\subset\mathcal{T}_{h}}\sum_{\gamma,\epsilon\in \Gamma_k^I} (\bm{u}_h^{\diamond})^T \tilde{\T}^{\diamond} \bm{\psi}_h^{\diamond}
	\\
	& \quad
	-\sum_{\gamma \subset \Gamma^D}(\Rgk\uhk - \bm{u}_{\gamma k})^T \T_{\gamma}^{(D)}\Rgk\bm{\psi}_{h,k} 
	-\sum_{\gamma \subset \Gamma^D} \bm{u}_{\gamma k}^T \B_\gamma \Dgk \bm{\psi}_{h,k}
	+\sum_{\gamma \subset \Gamma^D} \uhk^T\Dgk^T\B_{\gamma}\Rgk\bm{\psi}_{h,k}
	\\
	&\quad
	+\sum_{\gamma \subset \Gamma^N} \bm{w}_{\gamma k}^T \B_\gamma \Rgk \bm{\psi}_{h,k}  
	-\sum_{\gamma \subset \Gamma^N} \uhk^T\Rgk^T\B_{\gamma}\Dgk \bm{\psi}_{h,k},
\end{aligned}
\end{equation}
where $ \tilde{\T}^{\diamond} = (\T^{\diamond})^T$, and 
\begin{equation}
	\tilde{\T}^{\star}=\begin{bmatrix}
	\T_{\gamma   k}^{(1)} & -\T_{\gamma   v}^{(1)} & \T_{\gamma k}^{(2)}+\B_{\gamma} & -\T_{\gamma v}^{(2)}\\
	-\T_{\gamma   k}^{(1)} & \T_{\gamma   v}^{(1)} & -\T_{\gamma k}^{(2)} & \T_{\gamma v}^{(2)}+\B_{\gamma}\\
	\T_{\gamma k}^{(3)}-\B_{\gamma} & \T_{\gamma v}^{(3)} & \T_{\gamma k}^{(4)} & \T_{\gamma v}^{(4)}\\
	\T_{\gamma k}^{(3)} & \T_{\gamma v}^{(3)}-\B_{\gamma} & \T_{\gamma k}^{(4)} & \T_{\gamma v}^{(4)}
	\end{bmatrix}.
\end{equation}
Substituting \cref{eq:Sum Lhk} \violet{into \cref{eq:Vanishing term for adjoint consistency}, enforcing $ I_h(\bm{u}_h) - I_h(\bm{\psi}_h) =0$, and simplifying, we obtain
\begin{equation} \label{eq:L^*hk + g}
\begin{aligned}
	&\sum_{\Omega_{k}\subset{\mathcal T}_{h}}\bm{u}_{h,k}^T\H_{k}\left(\D_{k}^{(2)}\bm{\psi}_{h,k}+\bm{g}_{k}\right)
	-\sum_{\gamma\subset\Gamma^{I}}(\bm{u}_{h}^{\star})^{T}\tilde{\T}^{\star}\bm{\psi}_{h}^{\star}
	-\sum_{\Omega_{k}\subset\mathcal{T}_{h}}\sum_{\gamma,\epsilon\in\Gamma_{k}^{I}}(\bm{u}_{h}^{\diamond})^{T}\tilde{\T}^{\diamond}\bm{\psi}_{h}^{\diamond}
	\\&
	\quad
	-\sum_{\gamma\subset\Gamma^{D}}\left[\begin{array}{c}
	\R_{\gamma k}\bm{u}_{h,k}\\
	\D_{\gamma k}\bm{u}_{h,k}
	\end{array}\right]^{T}\left[\begin{array}{c}
	\T_{\gamma}^{(D)}\\
	-\B_{\gamma}
	\end{array}\right]\left(\R_{\gamma k}\bm{\psi}_{h,k}-\bm{\psi}_{\gamma k}\right)
	-\sum_{\gamma\subset\Gamma^{N}}\bm{u}_{h,k}^{T}\R_{\gamma k}^{T}\B_{\gamma}\left(\D_{\gamma k}\bm{\psi}_{h,k}-\bm{z}_{\gamma k}\right) = 0.
\end{aligned}
\end{equation} 
Rewriting \cref{eq:L^*hk + g} as 
$
	\sum_{\Omega_{k}\subset{\mathcal T}_{h}} \bm{u}_{h,k}^T\H_k \left(L_{h,k}^*\left(\bm{\psi}_h\right) - \bm{g}_k \right) = 0
$, 
we identify the discrete adjoint operator, $ L_{h,k}^* $, satisfying}
\begin{equation} \label{eq:Discrete adjoint problem}
	L_{h,k}^*(\bm{\psi}_h) - \bm{g}_k = \bm{0},
\end{equation}
\violet{to be}
\begin{equation} \label{eq:Discrete Adjoint operator}
L_{h,k}^*(\bm{\psi}_{h}) = -\D_k^{(2)}\bm{\psi}_{h,k}
+  \H^{-1}_k(\bm{s}_k^I)^*(\bm{\psi}_{h,k}) 
+  \H^{-1}_k(\bm{s}_k^B)^*(\bm{\psi}_{h,k},\bm{\psi}_{\gamma k}, \bm{z}_{\gamma k}),
\end{equation}
\violet{where} the adjoint interior facet SATs and boundary SATs are given, respectively, by
\begin{equation} \label{eq:sI*}
\begin{aligned}
	\left(\bm{s}_{k}^{I}\right)^{*}  &=\sum_{\gamma\subset\Gamma_{k}^{I}}\begin{bmatrix}
	\R_{\gamma k}^{T} & \D_{\gamma k}^{T}\end{bmatrix}
	\begin{bmatrix}
	\T_{\gamma   k}^{(1)} & -\T_{\gamma   v}^{(1)} & \T_{\gamma k}^{(2)}+\B_{\gamma} & -\T_{\gamma v}^{(2)}\\
	\T_{\gamma k}^{(3)}-\B_{\gamma} & \T_{\gamma v}^{(3)} & \T_{\gamma k}^{(4)} & \T_{\gamma v}^{(4)}
	\end{bmatrix}
	\begin{bmatrix}
	\R_{\gamma k}\bm{\psi}_{h,k}\\
	\R_{\gamma v}\bm{\psi}_{h,v}\\
	\D_{\gamma k}\bm{\psi}_{h,k}\\
	\D_{\gamma v}\bm{\psi}_{h,v}
	\end{bmatrix}\\
	&\quad
	+\sum_{\gamma\subset\Gamma_{k}^{I}}\left\{\sum_{\epsilon\subset\Gamma_{k}^{I}}\R_{\gamma k}^{T}\left[\T_{\gamma\epsilon k}^{(5)}\R_{\epsilon k}\bm{\psi}_{h,k}+\T_{\gamma\epsilon k}^{(6)}\R_{\epsilon g}\bm{\psi}_{h,g}\right]
	-\sum_{\delta\subset\Gamma_{v}^{I}}\R_{\gamma k}^{T}\left[\T_{\gamma\delta v}^{(5)}\R_{\delta v}\bm{\psi}_{h,v}+\T_{\gamma\delta v}^{(6)}\R_{\delta q}\bm{\psi}_{h,q}\right] \right\},
\end{aligned}
\end{equation}
\begin{equation}
\begin{aligned}
\left(\bm{s}_{k}^{B}\right)^{*}&=\sum_{\gamma\subset\Gamma^{D}}\left[\begin{array}{cc}
\R_{\gamma k}^{T} & \D_{\gamma k}^{T}\end{array}\right]\left[\begin{array}{c}
\T_{\gamma}^{(D)}\\
-\B_{\gamma}
\end{array}\right]\left[\begin{array}{cc}
\R_{\gamma k}\bm{\psi}_{h,k}-\bm{\psi}_{\gamma k}\end{array}\right]
+\sum_{\gamma\subset\Gamma^{N}}\R_{\gamma k}^{T}\B_{\gamma}\left(\D_{\gamma k}\bm{\psi}_{h,k}-\bm{z}_{\gamma k}\right).
\end{aligned}
\end{equation}
\violet{Note that the last term in $(\bm{s}_{k}^{I})^{*}$ is obtained by rewriting the extended interior facet SATs by facet contribution, regrouping by element, \ie, 
\begin{equation}
	\begin{aligned}
		\sum_{\Omega_{k}\subset\mathcal{T}_{h}}\sum_{\gamma,\epsilon\subset\Gamma_{k}^{I}}(\bm{\psi}_{h}^{\diamond})^{T}\T^{\diamond}\bm{u}_{h}^{\diamond}&=\sum_{\gamma\subset\Gamma^{I}}\bigg\{\sum_{\delta\subset\Gamma_{v}^{I}}\bm{\psi}_{h,v}^{T}\R_{\delta v}^{T}\T_{\delta\gamma v}^{(5)}[\R_{\gamma v}\bm{u}_{h,v}-\R_{\gamma k}\bm{u}_{h,k}]+\sum_{\epsilon\subset\Gamma_{k}^{I}}\bm{\psi}_{h,k}^{T}\R_{\epsilon k}^{T}\T_{\epsilon\gamma k}^{(5)}[\R_{\gamma k}\bm{u}_{h,k}-\R_{\gamma v}\bm{u}_{h,v}]\\&\quad+\sum_{\delta\subset\Gamma_{v}^{I}}\bm{\psi}_{h,q}^{T}\R_{\delta q}^{T}\T_{\delta\gamma v}^{(6)}[\R_{\gamma v}\bm{u}_{h,v}-\R_{\gamma k}\bm{u}_{h,k}]+\sum_{\epsilon\subset\Gamma_{k}^{I}}\bm{\psi}_{h,g}^{T}\R_{\epsilon g}^{T}\T_{\epsilon\gamma k}^{(6)}[\R_{\gamma k}\bm{u}_{h,k}-\R_{\gamma v}\bm{u}_{h,v}]\bigg\}\\&=\sum_{\Omega_{k}\subset\mathcal{T}_{h}}\sum_{\gamma\subset\Gamma_{k}^{I}}\bigg\{\sum_{\epsilon\subset\Gamma_{k}^{I}}\left[\bm{\psi}_{h,k}^{T}\R_{\epsilon k}^{T}\T_{\epsilon\gamma k}^{(5)}+\bm{\psi}_{h,g}^{T}\R_{\epsilon g}^{T}\T_{\epsilon\gamma k}^{(6)}\right]\R_{\gamma k}\bm{u}_{h,k} -\sum_{\delta\subset\Gamma_{v}^{I}}\left[\bm{\psi}_{h,v}^{T}\R_{\delta v}^{T}\T_{\delta\gamma v}^{(5)}+\bm{\psi}_{h,q}^{T}\R_{\delta q}^{T}\T_{\delta\gamma v}^{(6)}\right]\R_{\gamma k}\bm{u}_{h,k}\bigg\},
	\end{aligned}
\end{equation}
and transposing to find 
\begin{equation}
	\begin{aligned}
		\sum_{\Omega_{k}\subset\mathcal{T}_{h}}\sum_{\gamma,\epsilon\in\Gamma_{k}^{I}}(\bm{u}_{h}^{\diamond})^{T}\tilde{\T}^{\diamond}\bm{\psi}_{h}^{\diamond}&=\sum_{\Omega_{k}\subset\mathcal{T}_{h}}\sum_{\gamma\subset\Gamma_{k}^{I}}\bigg\{\sum_{\epsilon\subset\Gamma_{k}^{I}}\bm{u}_{h,k}^{T}\R_{\gamma k}^{T}\left[\T_{\gamma\epsilon k}^{(5)}\R_{\epsilon k}\bm{\psi}_{h,k}+\T_{\gamma\epsilon k}^{(6)}\R_{\epsilon g}\bm{\psi}_{h,g}\right]-\sum_{\delta\subset\Gamma_{v}^{I}}\bm{u}_{h,k}^{T}\R_{\gamma k}^{T}\left[\T_{\gamma\delta v}^{(5)}\R_{\delta v}\bm{\psi}_{h,v}+\T_{\gamma\delta v}^{(6)}\R_{\delta q}\bm{\psi}_{h,q}\right]\bigg\}.
	\end{aligned}
\end{equation}
It is also possible to regroup the extended interior facet SATs as}
\begin{equation} \label{eq:Adjoint Interface SAT extended terms}
\sum_{\gamma\subset\Gamma_{k}^{I}}\sum_{\epsilon\subset\Gamma_{k}^{I}}\left[\R_{\epsilon k}^{T}\left(\T_{\epsilon\gamma k}^{(5)}\R_{\gamma k}\bm{\psi}_{h,k}+\T_{\epsilon\gamma k}^{(6)}\R_{\gamma v}\bm{\psi}_{h,v}\right)-\R_{\epsilon g}^{T}\left(\T_{\epsilon\gamma k}^{(5)}\R_{\gamma k}\bm{\psi}_{h,k}+\T_{\epsilon\gamma k}^{(6)}\R_{\gamma v}\bm{\psi}_{h,v}\right)\right], 
\end{equation}
\violet{which can replace the last term in \cref{eq:sI*}}. We now state a theorem which is an extension of Theorem 1 in \cite{yan2018interior}.
\begin{theorem}\label{thm:Adjoint consistency}
	Let \cref{assu: mapping,assu:Coefficient matrices} hold, and the metric terms be evaluated exactly. Then the SBP-SAT discretization given in \cref{eq:SBP-SAT discretization,eq:Interface SATs,eq:Boundary SATs} and the discrete functional \cref{eq:Functional discrete} are adjoint consistent with respect to the steady version of the continuous PDE \cref{eq:diffusion problem} and functional \cref{eq:Functional}, \ie, \cref{eq:Adjoint consitstency definition} holds, and $\norm{ L_{h,k}^* \bm{\psi}_k - \bm{g}_k}_\infty ={\cal{O}}(h^{p-1})$ if $ \psi_k \in \cont{p+1} $ and the coefficient matrices satisfy the following relations
	\begin{equation} \label{eq:Coefficients for adjoint consistency}
	\begin{aligned}
	\T_{\gamma   k}^{(1)}&=\T_{\gamma   v}^{(1)}, 
	&& \T_{\gamma k}^{(2)}+\T_{\gamma v}^{(2)}=-\B_{\gamma}, 
	&&\T_{\gamma k}^{(3)}+\T_{\gamma v}^{(3)}=\B_{\gamma},
	&&\T_{\gamma k}^{(4)}=\T_{\gamma v}^{(4)},
	&& \T_{abk}^{(5)}=-\T_{abk}^{(6)},
	\end{aligned}
	\end{equation}
	where $ a,b\in\left\{ \gamma,\epsilon_{1},\epsilon_{2}\right\}$. 
\end{theorem}
\begin{proof}
	The result is a consequence of the accuracies of the derivative and extrapolation operators, \eg, for the exact adjoint solution $ \D_{k}^{(2)}\bm{\psi}_{k}+\bm{g}_{k} = [\nabla\cdot(\lambda\nabla\psi_k)]_{S_k}+\fnc{G}_k|_{S_{k}} + {\mathcal{O}}(h^{p-1}) = {\mathcal{O}}(h^{p-1})  $. Similarly, $ \Dgk\bm{\psi}_{k} = -\Dgk\bm{\psi}_{v}=[\bm{n}\cdot(\lambda\nabla\psi_k)]_{S_\gamma} +\order{h^{p}}$ and $ \Rgk\bm{\psi}_{k}=\psi_k|_{S_\gamma} + {\cal O}(h^{p+1}) $ which also holds for the extrapolations at the other facets. The desired result is obtained by substituting these approximations in \cref{eq:Discrete Adjoint operator} and using the coefficients in \cref{eq:Coefficients for adjoint consistency}.
\end{proof}

\violet{All but the second and fourth} conditions in \cref{thm:Adjoint consistency} are required for conservation (see \cref{thm:Conservation}). Similar analysis in \cite{hicken2011superconvergent} shows that an additional condition is required for a conservative discretization to be adjoint consistent, and such requirement can also be inferred from \cite{arnold2002unified,hartmann2013higher}. The following corollary follows from a comparison of the conditions presented in \cref{thm:Adjoint consistency,thm:Conservation}. 
\begin{corollary}
	Adjoint consistency is a sufficient but not necessary condition for conservation.
\end{corollary}
\begin{remark}
	In the DG literature, conservative properties of numerical fluxes associated with the solution and with the gradient of the solution are used to define conservation and adjoint consistency. A numerical discretization is conservative if the numerical flux associated with the gradient of the solution is conservative, and the discretization is adjoint consistent if both numerical fluxes associated with the solution and the gradient of the solution are conservative \cite{arnold2002unified,hartmann2013higher}. Conservation of the numerical flux associated with the gradient amounts to applying integration by parts once, as in \cref{eq:Conservation IBP}, and requiring that all SATs except those corresponding to $ \int_{\Gamma^B}(\lambda\nabla\fnc{U})\cdot\bm{n} \dd{\Gamma}$ vanish for $ \bm{v=1} $ in the discretized equation. On the other hand, conservation of both numerical fluxes amounts to applying integration by parts twice on the diffusive term and requiring that all SATs in the discretization except those corresponding to $ \int_{\Gamma^B}\fnc{V}(\lambda\nabla\fnc{U})\cdot\bm{n} - \fnc{U}(\lambda\nabla\fnc{V})\cdot\bm{n}\dd{\Gamma} $ vanish for smooth test function, $ \fnc{V}$. This gives $ (\tilde{\T}^{\star})^T $ instead of $ \T^{\star} $ in \cref{eq:Residual 2nd form} which in turn yields the same conditions for conservation as those stated for adjoint consistency in \cref{eq:Coefficients for adjoint consistency}. The ambiguity on the definition of conservation, \ie, whether to require the SATs to vanish after a single or double applications of integration by parts, is discussed in \cite{carpenter2010revisiting}. Furthermore, in \cite{ghasemi2020conservation,nordstrom2017relation} the latter definition is used to conclude that adjoint consistency is a necessary and sufficient condition for conservation.
\end{remark}
\begin{remark}
	The discrete adjoint operator, \cref{eq:Discrete Adjoint operator}, is derived for compact \violet{stencil} SATs in \cite{yan2018interior} using the variational relation $ J_{h}^{\prime}[\bm{u}_h](\delta \bm{u})+R_{h}^{\prime}[\bm{u}_h]\left(\delta \bm{u},\psi\right) = 0$, where $ J_{h}^{\prime}[\bm{u}_h] $ and $ R_{h}^{\prime}[\bm{u}_h] $ denote the Fr\'{e}chet derivatives of $ J_h $ and $ R_h $ with respect to $ \bm{u}_h $, respectively, and $ \delta\bm{u}_h $ is the variation of $ \bm{u}_h $. This suggests that the modification of the functional in \cref{eq:Functional discrete 2} is necessary for adjoint consistency, since our approach to obtain \cref{eq:Discrete Adjoint operator} relies on this modification.
\end{remark}
\begin{remark}
	\cref{thm:Adjoint consistency} implies adjoint consistent schemes need not have a symmetric stiffness matrix resulting from the residual \cref{eq:Residual 3rd form}. However, as in \cite{yan2018interior}, we consider adjoint consistent SATs that yield a symmetric stiffness matrix by requiring $ \T_{\gamma k}^{(3)}-\T_{\gamma k}^{(2)}=\B_{\gamma} $.
\end{remark}
\subsection{Functional accuracy} \label{sec:Functional accuracy}
The accuracy of the target functional depends on the primal and adjoint consistency of the SBP-SAT discretization of the underlying PDE. We establish functional error estimates for primal and adjoint consistent SBP-SAT discretizations of the Poisson problem on unstructured curvilinear grids. The following assumption will be necessary to proceed with the analysis.
\begin{assumption} \label{assu:Solution error}
	Unique numerical solutions to the primal and adjoint equations, $ \bm{u}_{h,k} $ and $ \bm{\psi}_{h,k} $ respectively, exist, and as $ h\rightarrow 0 $ they approximate the exact primal and adjoint solutions, $ \bm{u}_{k} $ and $ \bm{\psi}_{k} $ respectively, to order $ h^{p+1} $ in the infinity norm, \ie, 
	\begin{equation*}
	\begin{aligned}
	\norm{{\bm{u}_{h,k}-\bm{u}_{k}}}_\infty&=\order{h^{p+1}},
	\;\text{\;and}	
	&&
	\norm{\bm{\psi}_{h,k}-\bm{\psi}_{k}}_\infty=\order{h^{p+1}}.
	\end{aligned}
	\end{equation*}
\end{assumption}
The reason to invoke \cref{assu:Solution error} is related to difficulties in showing that the SBP-SAT discretization is pointwise stable (see, \eg, \cite{hicken2012output,hicken2011superconvergent,penner2020superconvergent} for similar assumptions, and \cite{gustafsson1981convergence,svard2006order,svard2019convergence} for analysis of convergence rates in one dimensional framework). Despite the fact that the truncation error for the primal and adjoint discretizations is $ {\mathcal{O}}(h^{p-1})$, numerical experiments show that the estimates in \cref{assu:Solution error} are usually attained.
In fact, the functional error estimate in \cref{thm:Functional accuracy} below holds even for order $ h^p $ primal and adjoint solution error values. To simplify the analysis, we first consider the case where the discretization has only one element and later show that the error estimate holds for more general cases.
\subsubsection{Functional error estimate on a single element}
The residual for the discrete Poisson problem, RHS of \cref{eq:SBP-SAT discretization}, on a single element premultiplied by $ \bm{\psi}_{k}^T\H_k $ reads
\begin{equation} \label{eq:Residual 1elem 1st form}
\begin{aligned}
\bm{\psi}_{k}^{T}\H_{k}R_{h,u}
&=-\bm{\psi}_{k}^{T}\H_{k}\D_{k}^{(2)}\left(\bm{u}_{h,k}-\bm{u}_{k}\right)
-\bm{\psi}_{k}^{T}\H_{k}\D_{k}^{(2)}\bm{u}_{k}
-\bm{\psi}_{k}^{T}\H_{k}\bm{f}_{k}
+\sum_{\gamma\subset\Gamma^{D}}\bm{\psi}_{k}^{T}\R_{\gamma k}^{T}\T_{\gamma}^{(D)}\left(\R_{\gamma k}\bm{u}_{h,k}-\bm{u}_{\gamma k}\right)
\\
&\quad
-\sum_{\gamma\subset\Gamma^{D}}\bm{\psi}_{k}^{T}\D_{\gamma k}^{T}\B_{\gamma}\left(\R_{\gamma k}\bm{u}_{h,k}-\bm{u}_{\gamma k}\right)
+\sum_{\gamma\subset\Gamma^{N}}\bm{\psi}_{k}^{T}\R_{\gamma k}^{T}\B_{\gamma}\left(\D_{\gamma k}\bm{u}_{h,k}-\bm{w}_{\gamma k}\right) = 0,
\end{aligned}
\end{equation}
where we have added and subtracted the term $ \bm{\psi}_{k}^{T}\H_{k}\D_{k}^{(2)}\bm{u}_{k} $. Similarly, the residual of the discrete adjoint problem, \cref{eq:Discrete adjoint problem}, premultiplied with $ \bm{u}^T_{k}\H_k $ can be written as
\begin{equation}\label{eq:Residual 1elem 2nd form}
\begin{aligned}
\bm{u}_{k}^{T}\H_{k}R_{h,\psi}&=-\bm{u}_{k}^{T}\H_{k}\D_{k}^{(2)}\left(\bm{\psi}_{h,k}-\bm{\psi}_{k}\right)-\bm{u}_{k}^{T}\H_{k}\D_{k}^{(2)}\bm{\psi}_{k}-\bm{u}_{k}^{T}\H_{k}\bm{g}_{k}+\sum_{\gamma\subset\Gamma^{D}}\bm{u}_{k}^{T}\R_{\gamma k}^{T}\T_{\gamma}^{(D)}\left(\R_{\gamma k}\bm{\psi}_{h,k}-\bm{\psi}_{\gamma k}\right)\\&\quad-\sum_{\gamma\subset\Gamma^{D}}\bm{u}_{k}^{T}\D_{\gamma k}^{T}\B_{\gamma}\left(\R_{\gamma k}\bm{\psi}_{h,k}-\bm{\psi}_{\gamma k}\right)+\sum_{\gamma\subset\Gamma^{N}}\bm{u}_{k}^{T}\R_{\gamma k}^{T}\B_{\gamma}\left(\D_{\gamma k}\bm{\psi}_{h,k}-\bm{z}_{\gamma k}\right)=0.
\end{aligned}
\end{equation}
We will use the above forms of the residuals, \cref{eq:Residual 1elem 1st form} and \cref{eq:Residual 1elem 2nd form}, to prove the following \violet{theorem}.
\begin{theorem} \label{thm:Functional accuracy}
	If \cref{assu: mapping,assu:Coefficient matrices,assu:Solution error} hold, \violet{$\fnc{U}, \psi\in \fnc{C}^{2p+2}(\Omega) $, $\lambda \in \fnc{C}^{2p+1}(\Omega) $}, and $ \bm{u}_{h,k} \in \IR{n_p}$ is a solution to a consistent and adjoint consistent SBP-SAT discretization of the form \cref{eq:SBP-SAT discretization}, then the discrete functional \cref{eq:Functional discrete} is an order $ h^{2p} $ approximation to the compatible linear functional \cref{eq:Functional}, \ie,
	\begin{equation} \label{eq:Functional error estimate}
	\fnc{I}(\fnc{U}) -  I_h(\bm{u}_{h}) = \order{h^{2p}}.
	\end{equation}
\end{theorem}
\begin{proof}
	The \blue{diagonal} norm matrices $ \H_k $ and  $ \B_{\gamma} $ contain quadrature weights of at least order $h^{2p}$ and $ h^{2p+1} $ accuracy, respectively. We discretize \cref{eq:Functional} using these quadratures as 
	\begin{align} \label{eq:Functional 1}
		{\mathcal I}\left({\mathcal U}\right)&=\int_{\Omega}{\mathcal G}{\mathcal U}\dd{\Omega}-\int_{\Gamma^{D}}\psi_{D}\left(\lambda\nabla{\mathcal U}\right)\cdot\bm{n}\dd{\Gamma}+\int_{\Gamma^{N}}\psi_{N}{\mathcal {\mathcal U}}\dd{\Gamma}
		\nonumber
		\\&
		=\bm{g}_k^{T}\H_{k}\bm{u}_{k}-\sum_{\gamma\subset\Gamma^{D}}\bm{\psi}_{\gamma k}^{T}\B_{\gamma}\bm{w}_{\gamma k}+\sum_{\gamma \subset \Gamma^{N}}\bm{z}_{\gamma k}^{T}\B_{\gamma}\bm{u}_{\gamma k}+{\mathcal O}\left(h^{2p}\right).
	\end{align}
	From \cref{eq:Adjoint relation}, the compatibility of the linear functional implies
	\begin{align}
			{\cal I}\left({\cal U}\right)=\fnc{I}\left(\psi\right)&=\int_{\Omega}{\cal F}\psi\dd{\Omega}-\int_{\Gamma^{D}}\mathcal{U}_{D}\left(\lambda\nabla{\cal \psi}\right)\cdot\bm{n}\dd{\Gamma}+\int_{\Gamma^{N}}{\cal {\cal U}}_{N}\psi\dd{\Gamma} \nonumber
			\\&=\bm{f}_{k}^{T}\H_{k}\bm{\psi}_{k}-\sum_{\gamma\subset\Gamma^{D}}\bm{u}_{\gamma k}^{T}\B_{\gamma}\bm{z}_{\gamma k}+\sum_{\gamma\subset\Gamma^{N}}\bm{w}_{\gamma k}^{T}\B_{\gamma}\bm{\psi}_{\gamma k}+{\cal O}\left(h^{2p}\right),\label{eq:Compatiblility in proof}
	\end{align}
	Subtracting \cref{eq:Functional discrete} from \cref{eq:Functional 1} and rearranging, we obtain
	\begin{equation} \label{eq:Functional in proof}
	\begin{aligned}
	{\cal I}\left({\cal U}\right)&={\cal I}_{h}\left(\bm{u}_{h}\right)-\bm{g}_{k}^{T}\H_{k}\left(\bm{u}_{h,k}-\bm{u}_{k}\right)+\sum_{\gamma\subset\Gamma^{D}}\bm{\psi}_{\gamma k}^{T}\B_{\gamma}\left(\D_{\gamma k}\bm{u}_{h,k}-\bm{w}_{\gamma k}\right)
	\\&\quad -\sum_{\gamma\subset\Gamma^{N}}\bm{z}_{\gamma k}^{T}\B_{\gamma}\left(\R_{\gamma k}\bm{u}_{h,k}-\bm{u}_{\gamma k}\right)
	-\sum_{\gamma\subset\Gamma^{D}}\bm{\psi}_{\gamma k}^{T}\T_{\gamma}^{(D)}\left(\R_{\gamma k}\bm{u}_{h,k}-\bm{u}_{\gamma k}\right)+{\cal O}\left(h^{2p}\right).
	\end{aligned}
	\end{equation}
	\violet{Adding and subtracting terms, we can rewrite \cref{eq:Functional in proof} as 
	\begin{equation} \label{eq:Functional in proof 1}
		\begin{aligned}
			{\cal I}\left({\cal U}\right)&={\cal I}_{h}\left(\bm{u}_{h}\right)-\bm{g}_{k}^{T}\H_{k}\left(\bm{u}_{h,k}-\bm{u}_{k}\right)+\sum_{\gamma\subset\Gamma^{D}}\bm{\psi}_{\gamma k}^{T}\B_{\gamma}\left(\D_{\gamma k}\bm{u}_{h,k}-\D_{\gamma k}\bm{u}_{k}\right)+\sum_{\gamma\subset\Gamma^{D}}\bm{\psi}_{\gamma k}^{T}\B_{\gamma}\left(\D_{\gamma k}\bm{u}_{k}-\bm{w}_{\gamma k}\right)\\&\quad-\sum_{\gamma\subset\Gamma^{N}}\bm{z}_{\gamma k}^{T}\B_{\gamma}\left(\R_{\gamma k}\bm{u}_{h,k}-\R_{\gamma k}\bm{u}_{k}\right)-\sum_{\gamma\subset\Gamma^{N}}\bm{z}_{\gamma k}^{T}\B_{\gamma}\left(\R_{\gamma k}\bm{u}_{k}-\bm{u}_{\gamma k}\right)\\&\quad-\sum_{\gamma\subset\Gamma^{D}}\bm{\psi}_{\gamma k}^{T}\T_{\gamma}^{(D)}\left(\R_{\gamma k}\bm{u}_{h,k}-\R_{\gamma k}\bm{u}_{k}\right)-\sum_{\gamma\subset\Gamma^{D}}\bm{\psi}_{\gamma k}^{T}\T_{\gamma}^{(D)}\left(\R_{\gamma k}\bm{u}_{k}-\bm{u}_{\gamma k}\right)+{\cal O}\left(h^{2p}\right).
		\end{aligned}
	\end{equation}
	Adding \cref{eq:Residual 1elem 1st form} to \cref{eq:Functional in proof 1} and simplifying, we have
	\begin{equation} \label{eq:Functional in proof 2}
	\begin{aligned} 
		{\cal I}\left({\cal U}\right)&={\cal I}_{h}\left(\bm{u}_{h}\right)+\bigg\{\sum_{\gamma\subset\Gamma^{D}}\bm{\psi}_{\gamma k}^{T}\B_{\gamma}\D_{\gamma k}\H_{k}^{-1}-\bm{g}_{k}^{T}-\sum_{\gamma\subset\Gamma^{N}}\bm{z}_{\gamma k}^{T}\B_{\gamma}\R_{\gamma k}\H_{k}^{-1}-\sum_{\gamma\subset\Gamma^{D}}\bm{\psi}_{k}^{T}\D_{\gamma k}^{T}\B_{\gamma}\R_{\gamma k}\H_{k}^{-1}
		\\&\quad-\bm{\psi}_{k}^{T}\H_{k}\D_{k}^{(2)}\H_{k}^{-1}+\sum_{\gamma\subset\Gamma^{D}}\left(\bm{\psi}_{k}^{T}\R_{\gamma k}^{T}-\bm{\psi}_{\gamma k}^{T}\right)\T_{\gamma}^{(D)}\R_{\gamma k}\H_{k}^{-1}+\sum_{\gamma\subset\Gamma^{N}}\bm{\psi}_{k}^{T}\R_{\gamma k}^{T}\B_{\gamma}\D_{\gamma k}\H_{k}^{-1}\bigg\}\H_{k}\left(\bm{u}_{h,k}-\bm{u}_{k}\right)
		\\&\quad+\sum_{\gamma\subset\Gamma^{D}}\bm{\psi}_{\gamma k}^{T}\B_{\gamma}\left(\D_{\gamma k}\bm{u}_{k}-\bm{w}_{\gamma k}\right)-\sum_{\gamma\subset\Gamma^{N}}\bm{z}_{\gamma k}^{T}\B_{\gamma}\left(\R_{\gamma k}\bm{u}_{k}-\bm{u}_{\gamma k}\right)-\sum_{\gamma\subset\Gamma^{D}}\bm{\psi}_{k}^{T}\D_{\gamma k}^{T}\B_{\gamma}\left(\R_{\gamma k}\bm{u}_{k}-\bm{u}_{\gamma k}\right)
		\\&\quad+\sum_{\gamma\subset\Gamma^{N}}\bm{\psi}_{k}^{T}\R_{\gamma k}^{T}\B_{\gamma}\left(\D_{\gamma k}\bm{u}_{k}-\bm{w}_{\gamma k}\right)-\bm{\psi}_{k}^{T}\H_{k}\D_{k}^{(2)}\bm{u}_{k}-\bm{\psi}_{k}^{T}\H_{k}\bm{f}_{k}+{\cal O}\left(h^{2p}\right).
	\end{aligned}
	\end{equation}
	Applying identity \cref{eq:D2 identity 2} in \cref{eq:Functional in proof 2} and simplifying gives 
	\begin{align} \label{eq:Functional in proof 3}
		{\cal I}\left({\cal U}\right)&={\cal I}_{h}\left(\bm{u}_{h}\right)+\bigg\{-\bm{g}_{k}^{T}-\bm{\psi}_{k}^{T}\left(\D_{k}^{(2)}\right)^{T}-\sum_{\gamma\subset\Gamma^{D}}\left(\bm{\psi}_{k}^{T}\R_{\gamma k}^{T}-\bm{\psi}_{\gamma k}^{T}\right)\B_{\gamma}\D_{\gamma k}\H_{k}^{-1}+\sum_{\gamma\subset\Gamma^{N}}\left(\bm{\psi}_{k}^{T}\D_{\gamma k}^{T}-\bm{z}_{\gamma k}^{T}\right)\B_{\gamma}\R_{\gamma k}\H_{k}^{-1}
		\nonumber
		\\&\quad+\sum_{\gamma\subset\Gamma^{D}}\left(\bm{\psi}_{k}^{T}\R_{\gamma k}^{T}-\bm{\psi}_{\gamma k}^{T}\right)\T_{\gamma}^{(D)}\R_{\gamma k}\H_{k}^{-1}\bigg\}\H_{k}\left(\bm{u}_{h,k}-\bm{u}_{k}\right)+\sum_{\gamma\subset\Gamma^{D}}\bm{\psi}_{\gamma k}^{T}\B_{\gamma}\left(\D_{\gamma k}\bm{u}_{k}-\bm{w}_{\gamma k}\right)-\sum_{\gamma\subset\Gamma^{N}}\bm{z}_{\gamma k}^{T}\B_{\gamma}\left(\R_{\gamma k}\bm{u}_{k}-\bm{u}_{\gamma k}\right)
		\nonumber
		\\&\quad-\sum_{\gamma\subset\Gamma^{D}}\bm{\psi}_{k}^{T}\D_{\gamma k}^{T}\B_{\gamma}\left(\R_{\gamma k}\bm{u}_{k}-\bm{u}_{\gamma k}\right)+\sum_{\gamma\subset\Gamma^{N}}\bm{\psi}_{k}^{T}\R_{\gamma k}^{T}\B_{\gamma}\left(\D_{\gamma k}\bm{u}_{k}-\bm{w}_{\gamma k}\right)-\bm{\psi}_{k}^{T}\H_{k}\D_{k}^{(2)}\bm{u}_{k}-\bm{\psi}_{k}^{T}\H_{k}\bm{f}_{k}+{\cal O}\left(h^{2p}\right).
	\end{align}
	The sum of the terms in the curly braces is the transpose of $ L_{h,k}^* (\bm{\psi}_k) - \bm{g}_k$, which is order $ h^{p-1} $ due to \cref{thm:Adjoint consistency}, $ \H_{k}(\bm{u}_{h,k}-\bm{u}_{k}) = \fnc{O} (h^{p+3})$ since  $ \H_k = \J_k\hat{\H} = \fnc{O}(h^2)$, and $ \norm*{\bm{u}_{h,k}-\bm{u}_{k}}_\infty = \fnc{O} (h^{p+1}) $ due to \cref{assu:Solution error}. Hence, we have 
	\begin{equation} \label{eq:Functional in proof 4}
		\begin{aligned}
			{\cal I}\left({\cal U}\right)&={\cal I}_{h}\left(\bm{u}_{h}\right)+\sum_{\gamma\subset\Gamma^{D}}\bm{\psi}_k^T \R_{\gamma k}^T\B_{\gamma}\left(\D_{\gamma k}\bm{u}_{k}-\bm{w}_{\gamma k}\right)-\sum_{\gamma\subset\Gamma^{N}}\bm{\psi}_k^T \D_{\gamma k}^T\B_{\gamma}\left(\R_{\gamma k}\bm{u}_{k}-\bm{u}_{\gamma k}\right)-\sum_{\gamma\subset\Gamma^{D}}\bm{\psi}_{k}^{T}\D_{\gamma k}^{T}\B_{\gamma}\left(\R_{\gamma k}\bm{u}_{k}-\bm{u}_{\gamma k}\right)\\&\quad+\sum_{\gamma\subset\Gamma^{N}}\bm{\psi}_{k}^{T}\R_{\gamma k}^{T}\B_{\gamma}\left(\D_{\gamma k}\bm{u}_{k}-\bm{w}_{\gamma k}\right)-\bm{\psi}_{k}^{T}\H_{k}\D_{k}^{(2)}\bm{u}_{k}-\bm{\psi}_{k}^{T}\H_{k}\bm{f}_{k}+{\cal O}\left(h^{2p}\right),
		\end{aligned}
	\end{equation}
	where the relations $ (\bm{\psi}_{\gamma k}^{T}- \bm{\psi}_k^T \R_{\gamma k}^T) \B_{\gamma}(\D_{\gamma k}\bm{u}_{k}-\bm{w}_{\gamma k}) = \fnc{O}(h^{p+2})$ and $ (\bm{z}_{\gamma k}^{T}- \bm{\psi}_k^T \D_{\gamma k}^T) \B_{\gamma}(\R_{\gamma k}\bm{u}_{k}-\bm{u}_{\gamma k}) = \fnc{O}(h^{p+2})$ are used to obtain the second and third terms on the RHS, respectively. We can further simplify  \cref{eq:Functional in proof 4} as
	\begin{equation} \label{eq:Functional in proof 5}
		{\cal I}\left({\cal U}\right)={\cal I}_{h}\left(\bm{u}_{h}\right)-\bigg\{\bm{\psi}_{k}^{T}\H_{k}\bm{f}_{k}+\bm{\psi}_{k}^{T}\H_{k}\D_{k}^{(2)}\bm{u}_{k}-\sum_{\gamma\subset\Gamma^{B}}\bm{\psi}_{k}^{T}\R_{\gamma k}^T\B_{\gamma}\left(\D_{\gamma k}\bm{u}_{k}-\bm{w}_{\gamma k}\right)+\sum_{\gamma\subset\Gamma^{B}}\bm{\psi}_{k}^{T}\D_{\gamma k}^{T}\B_{\gamma}\left(\R_{\gamma k}\bm{u}_{k}-\bm{u}_{\gamma k}\right)\bigg\}+{\cal O}\left(h^{2p}\right).
	\end{equation}
	Moreover, using \cref{eq:Functional 1}, \cref{eq:Compatiblility in proof}, and \cref{eq:D2 identity 2}, a straightforward algebraic manipulation of \cref{eq:Functional in proof 4} yields 
	\begin{equation} \label{eq:Functional in proof 6}
		{\cal I}\left({\cal U}\right)={\cal I}_{h}\left(\bm{u}_{h}\right)-\bigg\{\bm{u}_{k}^{T}\H_{k}\bm{g}_{k}+\bm{u}_{k}^{T}\H_{k}\D_{k}^{(2)}\bm{\psi}_{k}-\sum_{\gamma\subset\Gamma^{B}}\bm{u}_{k}^{T}\R_{\gamma k}^{T}\B_{\gamma}\left(\D_{\gamma k}\bm{\psi}_{k}-\bm{z}_{\gamma k}\right)+\sum_{\gamma\subset\Gamma^{B}}\bm{u}_{k}^{T}\D_{\gamma k}^{T}\B_{\gamma}\left(\R_{\gamma k}\bm{\psi}_{k}-\bm{\psi}_{\gamma k}\right)\bigg\}+{\cal O}\left(h^{2p}\right).
	\end{equation}
	Subtracting \cref{eq:Functional in proof 6} from \cref{eq:Functional in proof 5}, we find 
	\begin{equation} \label{eq:Functional difference}
		\tau_{u}-\tau_{\psi}={\cal O}(h^{2p}), 
	\end{equation}
	where 
	\begin{align}
		\tau_u &\coloneqq  \bm{\psi}_{k}^{T}\H_{k}\bm{f}_{k}+\bm{\psi}_{k}^{T}\H_{k}\D_{k}^{(2)}\bm{u}_{k}-\sum_{\gamma\subset\Gamma^{B}}\bm{\psi}_{k}^{T}\R_{\gamma k}^T\B_{\gamma}\left(\D_{\gamma k}\bm{u}_{k}-\bm{w}_{\gamma k}\right)+\sum_{\gamma\subset\Gamma^{B}}\bm{\psi}_{k}^{T}\D_{\gamma k}^{T}\B_{\gamma}\left(\R_{\gamma k}\bm{u}_{k}-\bm{u}_{\gamma k}\right), 
		\label{eq:tau_u}
		\\
		\tau_{\psi} &\coloneqq \bm{u}_{k}^{T}\H_{k}\bm{g}_{k}+\bm{u}_{k}^{T}\H_{k}\D_{k}^{(2)}\bm{\psi}_{k}-\sum_{\gamma\subset\Gamma^{B}}\bm{u}_{k}^{T}\R_{\gamma k}^{T}\B_{\gamma}\left(\D_{\gamma k}\bm{\psi}_{k}-\bm{z}_{\gamma k}\right)+\sum_{\gamma\subset\Gamma^{B}}\bm{u}_{k}^{T}\D_{\gamma k}^{T}\B_{\gamma}\left(\R_{\gamma k}\bm{\psi}_{k}-\bm{\psi}_{\gamma k}\right).  
		\label{eq:tau_psi}
	\end{align}
}
	
	\violet{For an affine mapping, the derivative and extrapolation operators in \cref{eq:tau_u} and \cref{eq:tau_psi} are exact for polynomials of degree $ p $, and since all the terms in $ \tau_{u} $ and $ \tau_{\psi} $ are integrals approximated by quadrature rules of degree $ \ge 2p-1 $, it is sufficient if we can show that \cref{eq:Functional error estimate} holds for polynomial integrands of total degree $ 2p-1 $. A similar technique is used to show quadrature accuracy in \cite{hicken2013summation}. If we set $ \fnc{U}\in \poly{p} $ such that $ (\psi \fnc{U}) \in \poly{2p+1} $, then we have $ \tau_{u}=\fnc{O}(h^{2p}) $ due to the accuracy of the SBP operators and the primal PDE. Similarly, if we set $ \psi \in \poly{p} $ such that $ (\psi \fnc{U}) \in \poly{2p+1} $, then we obtain $ \tau_{\psi}=\fnc{O}(h^{2p}) $ due to the accuracy of the SBP operators and the adjoint PDE. Hence, we obtain $ \fnc{I}(\fnc{U}) -  I_h(\bm{u}_{h}) = \fnc{O}(h^{2p}) $. On curved elements, the SBP operators are not exact for polynomials of degree greater than zero. However, since \cref{eq:Functional difference} must hold for all combinations of $ \cal{U} $ and $ \psi $, we conclude that each of the error terms, $\tau_{u}$ and $ \tau_{\psi} $, must be $ \fnc{O}(h^{2p}) $, which gives $ \fnc{I}(\fnc{U}) -  I_h(\bm{u}_{h}) = \fnc{O}(h^{2p}) $, as desired.}
\end{proof}

\subsubsection{Functional error estimate with interior facet SATs} 
In this subsection, we will show that the functional estimate established in \cref{thm:Functional accuracy} holds for primal and adjoint consistent interior facet SATs. We consider two elements $ \Omega_k $ and $ \Omega_v $ sharing interface $ \gamma $ and introduce the following vectors and block matrices 
\begin{equation}
\begin{aligned}
\bm{u}_{h}&=\begin{bmatrix}
\bm{u}_{h,k}^{T} & \bm{u}_{h,v}^{T} & \bm{u}_{h,g_{1}}^{T} & \bm{u}_{h,g_{2}}^{T} & \bm{u}_{h,q_{1}}^{T} & \bm{u}_{h,q_{2}}^{T}\end{bmatrix}^{T},&\mathbb{H}_{11}&=\H_{k},\\\bm{u}&=\left[\begin{array}{cccccc}
\bm{u}_{k}^{T} & \bm{u}_{v}^{T} & \bm{u}_{g_{1}}^{T} & \bm{u}_{g_{2}}^{T} & \bm{u}_{q_{1}}^{T} & \bm{u}_{q_{2}}^{T}\end{array}\right]^{T},&\mathbb{H}_{22}&=\H_{v},
\\\bm{\psi}_{h}&=\begin{bmatrix}
\bm{\psi}_{h,k}^{T} & \bm{\psi}_{h,v}^{T} & \bm{\psi}_{h,g_{1}}^{T} & \bm{\psi}_{h,g_{2}}^{T} & \bm{\psi}_{h,q_{1}}^{T} & \bm{\psi}_{h,q_{2}}^{T}\end{bmatrix}^{T},&\mathbb{D}_{11}&=\D_{k},
\\\bm{\psi}&=\begin{bmatrix}
\bm{\psi}_{k}^{T} & \bm{\psi}_{v}^{T} & \bm{\psi}_{g_{1}}^{T} & \bm{\psi}_{g_{2}}^{T} & \bm{\psi}_{q_{1}}^{T} & \bm{\psi}_{q_{2}}^{T}\end{bmatrix}^{T},&\mathbb{D}_{22}&=\D_{v},
\\\bm{f}&=\begin{bmatrix}
\bm{f}_{k}^{T} & \bm{f}_{v}^{T} & \bm{0} & \bm{0} & \bm{0} & \bm{0}\end{bmatrix}^{T},&\left(\mathbb{D}_{\Lambda}\right)_{11}&=\D_{k}\Lambda_{k},
\\\bm{g}&=\begin{bmatrix}
\bm{g}_{k}^{T} & \bm{g}_{v}^{T} & \bm{0} & \bm{0} & \bm{0} & \bm{0}\end{bmatrix}^{T},&\left(\mathbb{D}_{\Lambda}\right)_{22}&=\D_{v}\Lambda_{v},
\end{aligned}
\end{equation}
where the vectors in the left column are of size $ \IR{6n_p} $, the matrices in the right column have $ 6\times 6 $ blocks which are of size $ \IRtwo{n_p}{n_p} $ each. Except for the blocks specified, the entries of the matrices in the right column are all zeros. 

We can write the primal residual for the two elements premultiplied by $ \bm{\psi}^{T}\mathbb{H} $ as 
\begin{align} \label{eq:Residual 1st form 2elem}
\bm{\psi}^{T}\mathbb{H} R_{h,u}&=-\bm{\psi}^{T}\mathbb{H}\mathbb{D}\mathbb{D}_{\Lambda}\left(\bm{u}_{h}-\bm{u}\right)-\bm{\psi}^{T}\mathbb{H}\mathbb{D}\mathbb{D}_{\Lambda}\bm{u}-\bm{\psi}^{T}\mathbb{H}\bm{f}
+\bm{\psi}^{T}\mathbb{A}\left(\bm{u}_{h}-\bm{u}\right)
+\violet{\bm{\psi}^{T}\mathbb{A}\bm{u}}
+\bm{\psi}^{T}\mathbb{B}\left(\bm{u}_{h}-\bm{u}\right) +\violet{\bm{\psi}^{T}\mathbb{B}\bm{u}} = 0,
\end{align} 
where we have drop all SATs except those at the facet $ \gamma $, and the nonzero entries of block matrices $ \mathbb{A},\mathbb{B}\in \IRtwo{6n_p}{6n_p} $ are given by
\begin{equation}
\begin{aligned}
\mathbb{A}_{11}&=\R_{\gamma k}^{T}\T_{\gamma k}^{(1)}\R_{\gamma k}+\R_{\gamma k}^{T}\T_{\gamma\xi_{1}k}^{(5)}\R_{\xi_{1}k}+\R_{\gamma k}^{T}\T_{\gamma\xi_{2}k}^{(5)}\R_{\xi_{2}k},&\mathbb{A}_{13}&=-\R_{\gamma k}^{T}\T_{\gamma\xi_{1}g_{1}}^{(5)}\R_{\xi_{1}g_{1}},\\\mathbb{A}_{12}&=-\R_{\gamma k}^{T}\T_{\gamma k}^{(1)}\R_{\gamma v}+\R_{\gamma k}^{T}\T_{\gamma\delta_{1}v}^{(6)}\R_{\delta_{1}v}+\R_{\gamma k}^{T}\T_{\gamma\delta_{2}v}^{(6)}\R_{\delta_{2}v},&\mathbb{A}_{14}&=-\R_{\gamma k}^{T}\T_{\gamma\xi_{2}g_{2}}^{(5)}\R_{\xi_{2}g_{2}},\\\mathbb{A}_{21}&=-\R_{\gamma v}^{T}\T_{\gamma v}^{(1)}\R_{\gamma k}+\R_{\gamma v}^{T}\T_{\gamma\xi_{1}k}^{(6)}\R_{\xi_{1}k}+\R_{\gamma v}^{T}\T_{\gamma\xi_{2}k}^{(6)}\R_{\xi_{2}k},&\mathbb{A}_{15}&=-\R_{\gamma k}^{T}\T_{\gamma\delta_{1}v}^{(6)}\R_{\delta_{1}q_{1}},\\\mathbb{A}_{22}&=\R_{\gamma v}^{T}\T_{\gamma v}^{(1)}\R_{\gamma v}+\R_{\gamma v}^{T}\T_{\gamma\delta_{1}v}^{(5)}\R_{\delta_{1}v}+\R_{\gamma v}^{T}\T_{\gamma\delta_{2}v}^{(5)}\R_{\delta_{2}v},&\mathbb{A}_{16}&=-\R_{\gamma k}^{T}\T_{\gamma\delta_{2}v}^{(6)}\R_{\delta_{2}q_{2}},\\\mathbb{B}_{11}&=\D_{\gamma k}^{T}\T_{\gamma k}^{(2)}\R_{\gamma k}+\R_{\gamma k}^{T}\T_{\gamma k}^{(3)}\D_{\gamma k}+\D_{\gamma k}^{T}\T_{\gamma k}^{(4)}\D_{\gamma k},&\mathbb{A}_{23}&=-\R_{\gamma v}^{T}\T_{\gamma\xi_{1}k}^{(6)}\R_{\xi_{1}g_{1}},\\\mathbb{B}_{12}&=-\D_{\gamma k}^{T}\T_{\gamma k}^{(2)}\R_{\gamma v}+\R_{\gamma k}^{T}\T_{\gamma k}^{(3)}\D_{\gamma v}+\D_{\gamma k}^{T}\T_{\gamma k}^{(4)}\D_{\gamma v},&\mathbb{A}_{24}&=-\R_{\gamma v}^{T}\T_{\gamma\xi_{2}k}^{(6)}\R_{\xi_{2}g_{2}},\\\mathbb{B}_{21}&=-\D_{\gamma v}^{T}\T_{\gamma v}^{(2)}\R_{\gamma k}+\R_{\gamma v}^{T}\T_{\gamma v}^{(3)}\D_{\gamma k}+\D_{\gamma k}^{T}\T_{\gamma v}^{(4)}\D_{\gamma v},&\mathbb{A}_{25}&=-\R_{\gamma v}^{T}\T_{\gamma\delta_{1}v}^{(5)}\R_{\delta_{1}q_{1}},\\\mathbb{B}_{22}&=\D_{\gamma v}^{T}\T_{\gamma v}^{(2)}\R_{\gamma v}+\R_{\gamma v}^{T}\T_{\gamma v}^{(3)}\D_{\gamma v}+\D_{\gamma v}^{T}\T_{\gamma v}^{(4)}\D_{\gamma v},&\mathbb{A}_{26}&=-\R_{\gamma v}^{T}\T_{\gamma\delta_{2}v}^{(5)}\R_{\delta_{2}q_{2}}.
\end{aligned}
\end{equation}
For the discrete adjoint problem, \cref{eq:Discrete adjoint problem}, with extended interior facet SATs grouped as in \cref{eq:Adjoint Interface SAT extended terms}, the truncation error reads
\begin{equation} \label{eq:Residual adjoint 2elem}
e_{\psi}\equiv\mathbb{D}\mathbb{D}_{\Lambda}\bm{\psi}+\bm{g}-\mathbb{H}^{-1}\mathbb{K}\bm{\psi}-\mathbb{H}^{-1}\mathbb{L}\bm{\psi},
\end{equation}
where 
\begin{equation}
\begin{aligned}
	\mathbb{K}_{11}&=\R_{\gamma k}^{T}\T_{\gamma k}^{(1)}\R_{\gamma k}+\R_{\xi_{1}k}^{T}\T_{\xi_{1}\gamma k}^{(5)}\R_{\gamma k}+\R_{\xi_{2}k}^{T}\T_{\xi_{2}\gamma k}^{(5)}\R_{\gamma k},&\mathbb{K}_{31}&=-\R_{\xi_{1}g_{1}}^{T}\T_{\xi_{1}\gamma k}^{(5)}\R_{\gamma k},\\\mathbb{K}_{12}&=-\R_{\gamma k}^{T}\T_{\gamma v}^{(1)}\R_{\gamma v}+\R_{\xi_{1}k}^{T}\T_{\xi_{1}\gamma k}^{(6)}\R_{\gamma v}+\R_{\xi_{2}k}^{T}\T_{\xi_{2}\gamma k}^{(6)}\R_{\gamma v},&\mathbb{K}_{32}&=-\R_{\xi_{1}g_{1}}^{T}\T_{\xi_{1}\gamma k}^{(6)}\R_{\gamma v},\\\mathbb{K}_{21}&=-\R_{\gamma v}^{T}\T_{\gamma k}^{(1)}\R_{\gamma k}+\R_{\delta_{1}v}^{T}\T_{\delta_{1}\gamma v}^{(6)}\R_{\gamma k}+\R_{\delta_{2}v}^{T}\T_{\delta_{2}\gamma v}^{(6)}\R_{\gamma k},&\mathbb{K}_{41}&=-\R_{\xi_{2}g_{2}}^{T}\T_{\xi_{2}\gamma k}^{(5)}\R_{\gamma k},\\\mathbb{K}_{22}&=\R_{\gamma v}^{T}\T_{\gamma v}^{(1)}\R_{\gamma v}+\R_{\delta_{1}v}^{T}\T_{\delta_{1}\gamma v}^{(5)}\R_{\gamma v}+\R_{\delta_{2}v}^{T}\T_{\delta_{2}\gamma v}^{(5)}\R_{\gamma v},&\mathbb{K}_{42}&=-\R_{\xi_{2}g_{2}}^{T}\T_{\xi_{2}\gamma k}^{(6)}\R_{\gamma v},\\\mathbb{L}_{11}&=\R_{\gamma k}^{T}\left(\T_{\gamma k}^{(2)}+\B_{\gamma}\right)\D_{\gamma k}+\D_{\gamma k}^{T}\left(\T_{\gamma k}^{(3)}-\B_{\gamma}\right)\R_{\gamma k}+\D_{\gamma k}^{T}\T_{\gamma k}^{(4)}\D_{\gamma k},&\mathbb{K}_{51}&=-\R_{\delta_{1}q_{1}}^{T}\T_{\delta_{1}\gamma v}^{(6)}\R_{\gamma k},\\\mathbb{L}_{12}&=-\R_{\gamma k}^{T}\T_{\gamma v}^{(2)}\D_{\gamma v}+\D_{\gamma k}^{T}\T_{\gamma v}^{(3)}\R_{\gamma v}+\D_{\gamma k}^{T}\T_{\gamma v}^{(4)}\D_{\gamma v},&\mathbb{K}_{52}&=-\R_{\delta_{1}q_{1}}^{T}\T_{\delta_{1}\gamma v}^{(5)}\R_{\gamma v},\\\mathbb{L}_{21}&=-\R_{\gamma v}^{T}\T_{\gamma k}^{(2)}\D_{\gamma k}+\D_{\gamma v}^{T}\T_{\gamma k}^{(3)}\R_{\gamma k}+\D_{\gamma v}^{T}\T_{\gamma k}^{(4)}\D_{\gamma k},&\mathbb{K}_{61}&=-\R_{\delta_{2}q_{2}}^{T}\T_{\delta_{2}\gamma v}^{(6)}\R_{\gamma k},\\\mathbb{L}_{22}&=\R_{\gamma v}^{T}\left(\T_{\gamma v}^{(2)}+\B_{\gamma}\right)\D_{\gamma v}+\D_{\gamma v}^{T}\left(\T_{\gamma v}^{(3)}-\B_{\gamma}\right)\R_{\gamma v}+\D_{\gamma v}^{T}\T_{\gamma v}^{(4)}\D_{\gamma v},&\mathbb{K}_{62}&=-\R_{\delta_{2}q_{2}}^{T}\T_{\delta_{2}\gamma v}^{(5)}\R_{\gamma v}.
\end{aligned}
\end{equation}
\cref{thm:Adjoint consistency} ensures that $ e_{\psi} = {\cal O}(h^{p-1}) $ for adjoint consistent SATs .
Neglecting the boundary terms, we write the functional as
\begin{equation} \label{eq:Functional interface SATs}
{\mathcal I}\left({\mathcal U}\right)=\int_{\Omega}{\mathcal G}{\mathcal U}\dd{\Omega}=\int_{\Omega}{\mathcal \psi}{\mathcal F}\dd{\Omega}=\bm{g}^{T}\mathbb{H}\bm{u}+ \order{h^{2p}} =\bm{f}^{T}\mathbb{H}\bm{\psi}+ \order{h^{2p}} .
\end{equation}
Adding \cref{eq:Residual 1st form 2elem} to $ \fnc{I}(\fnc{U})=\bm{g}^T\mathbb{H}\bm{u} + \fnc{O}(h^{2p}) = \bm{g}^T\mathbb{H}\bm{u}_h - \bm{g}^T\mathbb{H}(\bm{u}_h - \bm{u}) + \fnc{O}(h^{2p})$ and simplifying, we obtain 
\begin{equation}\label{eq:Superconv with SATs}
\begin{aligned}
{\mathcal I}\left({\mathcal U}\right)&= I_{h}\left(\bm{u}_{h}\right)
+\left\{\bm{\psi}^{T}\mathbb{B}\mathbb{H}^{-1} -\bm{\psi}^{T}\mathbb{H}\mathbb{D}\mathbb{D}_{\Lambda}\mathbb{H}^{-1}-\bm{g}^{T}+\bm{\psi}^{T}\mathbb{A}\mathbb{H}^{-1}\right\} \mathbb{H}\left(\bm{u}_{h}-\bm{u}\right)
\\
&\quad
-\bm{\psi}^{T}\mathbb{H}\mathbb{D}\mathbb{D}_{\Lambda}\bm{u}-\bm{\psi}^{T}\mathbb{H}\bm{f} + \bm{\psi}^T\mathbb{A}\bm{u} + \bm{\psi}^T\mathbb{B}\bm{u} +\order{h^{2p}}.
\end{aligned}
\end{equation}
Identity \cref{eq:D2 identity 2} gives 
\begin{equation} \label{eq:D2 idnetity 2 on SAT}
\mathbb{H}\mathbb{D}\mathbb{D}_{\Lambda}\mathbb{H}^{-1}=\mathbb{D}_{\Lambda}^{T}\mathbb{D}^{T}-\mathbb{C}\mathbb{H}^{-1},
\end{equation}
where the nonzero entries of $ \mathbb{C} $ are 
\begin{equation}\label{eq: C matrix}
\begin{aligned}
\mathbb{C}_{11} & =\D_{\gamma k}^{T}\B_{\gamma}\R_{\gamma k}-\R_{\gamma k}^{T}\B_{\gamma}\D_{\gamma k}+\sum_{\epsilon\subset\Gamma_{k}}\left(\D_{\epsilon k}^{T}\B_{\epsilon}\R_{\epsilon k}-\R_{\epsilon k}^{T}\B_{\epsilon}\D_{\epsilon k}\right),\\
\mathbb{C}_{22} & =\D_{\gamma v}^{T}\B_{\gamma}\R_{\gamma v}-\R_{\gamma v}^{T}\B_{\gamma}\D_{\gamma v}+\sum_{\delta\subset\Gamma_{v}}\left(\D_{\delta v}^{T}\B_{\delta}\R_{\delta v}-\R_{\delta v}^{T}\B_{\delta}\D_{\delta v}\right).
\end{aligned}
\end{equation}
Substituting \cref{eq:D2 idnetity 2 on SAT} into \cref{eq:Superconv with SATs}, we find
\begin{equation} \label{eq:Functional interior SAT}
\begin{aligned}
{\mathcal I}\left({\mathcal U}\right)&={I}_{h}\left(\bm{u}_{h}\right)+
\left\{\bm{\psi}^{T}\left(\mathbb{B}+\mathbb{C}\right)\mathbb{H}^{-1}-\bm{\psi}^{T}\mathbb{D}_{\Lambda}^{T}\mathbb{D}^{T}-\bm{g}^{T}+\bm{\psi}^{T}\mathbb{A}\mathbb{H}^{-1}\right\} \mathbb{H}\left(\bm{u}_{h}-\bm{u}\right)
\\&\quad
-\bm{\psi}^{T}\mathbb{H}\mathbb{D}\mathbb{D}_{\Lambda}\bm{u}-\bm{\psi}^{T}\mathbb{H}\bm{f}
+\bm{\psi}^T\mathbb{A}\bm{u} + \bm{\psi}^T\mathbb{B}\bm{u}
+\order{h^{2p}}.
\end{aligned}
\end{equation}
We define $ \mathbb{C}_\gamma $ as the matrix $ \mathbb{C} $ without the terms on facets other than $ \gamma $, \violet{and replace $ \mathbb{C} $ in \cref{eq:Functional interior SAT} by $ \mathbb{C}_\gamma $}. We then note that $ \mathbb{A}=\mathbb{K}^{T} $ and $ \mathbb{B}+\mathbb{C}_{\gamma}=\mathbb{L}^{T} $, and the sum of the terms in the curly braces is equal to $ e_{\psi} ^T$, which is $\fnc{O}(h^{p-1}) $ for adjoint consistent interior facet SATs. \violet{Hence, the second term on the RHS of \cref{eq:Functional interior SAT} is order $ h^{2p+2} $. Therefore, the terms remaining due to the inclusion of interior SATs are $ \bm{\psi}^T\mathbb{A}\bm{u}$ and $ \bm{\psi}^T\mathbb{B}\bm{u} $. These terms are added to $ \tau_u $ given in \cref{eq:tau_u}. Similarly, using \cref{eq:Functional interface SATs} and \cref{eq:D2 idnetity 2 on SAT}, it is possible to show that $ \tau_{\psi} $ must include the adjoint interior SATs, $ \bm{u}^T\mathbb{K}\bm{\psi}$ and $ \bm{u}^T\mathbb{L}\bm{\psi} $. Then, applying the same argument used to establish the order of $\tau_u $ in the proof of \cref{thm:Functional accuracy}, we arrive at $ \fnc{I}(\fnc{U}) -  I_h(\bm{u}_{h}) = \fnc{O}(h^{2p}) $ for discretizations with interior facet SATs.}

\begin{remark}
	Dropping terms in the matrix $ \mathbb{C} $ that are associated with facets other than $ \gamma $  is not necessary if one considers all facets for the analysis, but this would require working with bigger matrices. 
\end{remark}

\subsection{Energy stability analysis}\label{sec:Energy stability} In general, energy stability of SBP-SAT discretizations implies  
\begin{equation}
\dv{}{t}\left(\norm{\bm{u}_h}^2_\H\right) =\bm{u}_h^T \H\dv{\bm{u}_h}{t}+\dv{\bm{u}_h^T}{t}\H\bm{u}_h = 2 R_{h}\left(\bm{u}_{h},\bm{u}_{h}\right) \le 0.
\end{equation} 
We analyze the time stability of the SBP-SAT discretization \cref{eq:SBP-SAT discretization,eq:Interface SATs,eq:Boundary SATs} of the homogeneous diffusion problem. A class of adjoint inconsistent SATs is considered first, and later we present conditions for stability of a class of adjoint consistent SATs. The following theorem, \violet{whose proof can be found in \cite{albert1969conditions, gallier2010schur}, is} useful for the energy analysis.
\begin{theorem} \label{thm:Positive semi-definiteness}
	For a symmetric matrix of the form 
	$\Y =\bigl[\begin{smallmatrix}
	\Y_{11} & \Y_{12} \\ \Y_{12}^T & \Y_{22}
	\end{smallmatrix}\bigr]$,
	\begin{enumerate}[i)]
		\item $ \Y\succeq 0 $ if and only if $ \Y_{11}\succeq 0 $, $ (\I -\Y_{11}\Y_{11}^{+})\Y_{12} = \bm{0} $, and $ \Y_{22} - \Y_{12}^T \Y_{11}^{+}\Y_{12}\succeq 0 $,
		\item $ \Y\succeq 0 $ if and only if $ \Y_{22}\succeq 0 $, $ (\I-\Y_{22}\Y_{22}^{+})\Y_{12}^T  = \bm{0}$, and $ \Y_{11} - \Y_{12}\Y_{22}^{+}\Y_{12}^T \succeq 0$,
	\end{enumerate}
	where $ \Y\succeq 0 $ indicates that $ \Y $ is positive semidefinite.
\end{theorem}

\subsection{Energy analysis for adjoint inconsistent SATs} All of the adjoint inconsistent SATs presented in this work do not couple second neighbor elements. Therefore, we focus on compact adjoint inconsistent SATs and prove the following statement. 
\begin{theorem}\label{thm:Stability Adjoint Inconsistent}
	A conservative but adjoint inconsistent SBP-SAT discretization, \cref{eq:SBP-SAT discretization,eq:Interface SATs,eq:Boundary SATs}, of the homogeneous diffusion problem \cref{eq:diffusion problem}, \ie, $ \fnc{F}=0$, $\fnc{U}_D = 0 $, and $ \fnc{U}_N =0$, is energy stable with respect to the \blue{diagonal norm matrix}, $ \H $, if 
	\begin{equation} \label{eq:Conditions for stabilty of adjoint inconsistent SATs}
	\begin{aligned}
		\T_{\gamma k}^{(3)}+\T_{\gamma k}^{(2)}-\B_{\gamma}&=\bm{0} ,&\T_{\gamma v}^{(3)}-\T_{\gamma k}^{(2)}&=\bm{0},&\T_{\gamma k}^{(1)}&\succeq0,&\T_{\gamma}^{(D)}-\frac{1}{\alpha_{\gamma k}}\B_{\gamma}\Upsilon_{\gamma \gamma k}\B_{\gamma}&\succeq0,\\\T_{\gamma v}^{(3)}+\T_{\gamma v}^{(2)}-\B_{\gamma}&=\bm{0},&\T_{\gamma k}^{(3)}-\T_{\gamma v}^{(2)}&=\bm{0},&\T_{\gamma k}^{(4)}&=\T_{\gamma v}^{(4)}\succeq 0,&\T_{abk}^{(5)}=-\T_{abk}^{(6)}&=\bm{0},
	\end{aligned}
	\end{equation}
	where for element $ \Omega_k $, and facets $ a,b\in\{\gamma,\epsilon_{1},\epsilon_{2},\delta_{1},\delta_{2}\} $
	\begin{equation} \label{eq:Upsilon definition}
	\Upsilon_{abk}\equiv\C_{ak}\left(\Lambda_{k}^{*}\right)^{-1}\C_{bk}^{T} 
	=\N_{ak}^{T}\bar{\R}_{ak}\bar{\H}_{k}^{-1}\Lambda_{k}\bar{\R}_{bk}^{T}\N_{bk}.
	\end{equation}
\end{theorem}
\begin{proof}
	We have to show that for the conditions given in \cref{eq:Conditions for stabilty of adjoint inconsistent SATs} the residual satisfies $ 2R_{h}\left(\bm{u}_{h},\bm{u}_{h}\right) =  R_{h}\left(\bm{u}_{h},\bm{u}_{h}\right)+R_{h}^{T}\left(\bm{u}_{h},\bm{u}_{h}\right)\le0 $, which is the case if all the $ 4\times 4 $ block matrices on the RHS of \cref{eq:Residual 3rd form} are positive semidefinite. The positive semidefiniteness of the $ 4\times 4 $ block matrix in the first and last terms of \cref{eq:Residual 3rd form} are analyzed in \cite{yan2018interior}, and it is shown that these block matrices are positive semidefinite if $ \T_{\gamma k}^{(4)} \succeq 0$, $ \T_{\gamma}^{(D)} \succeq 0 $, and the Schur complement $ \T_{\gamma}^{(D)} - (1/\alpha_{\gamma k})\B_{\gamma}\Upsilon_{\gamma \gamma k}\B_{\gamma} \succeq 0$. Substituting $ \T_{abk}^{(5)}=-\T_{abk}^{(6)}=\bm{0} $ and $ \T_{\gamma k}^{(1)} =\T_{\gamma v}^{(1)}$ (due to conservation) in \cref{eq:Residual 3rd form}, and regrouping $ \T_{\gamma k}^{(1)} $ terms, the last $ 4\times 4 $ block matrix that we need to show is positive semidefinite is $ \A + \A^T $, where 
	\begin{equation} \label{eq:Matrix A for compact SATs}
	\A = \begin{bmatrix}
	\T_{\gamma k}^{(1)}&-\T_{\gamma k}^{(1)}&\left(\T_{\gamma k}^{(3)}-\B_{\gamma}\right)\C_{\gamma k}&\T_{\gamma k}^{(3)}\C_{\gamma v}\\-\T_{\gamma k}^{(1)}&\T_{\gamma k}^{(1)}&\T_{\gamma v}^{(3)}\C_{\gamma k}&\left(\T_{\gamma\gamma v}^{(3)}-\B_{\gamma}\right)\C_{\gamma v}\\\C_{\gamma k}^{T}\T_{\gamma k}^{(2)}&-\C_{\gamma k}^{T}\T_{\gamma k}^{(2)}&\alpha_{\gamma k}\Lambda_{k}^{*}& \bm{0} \\-\C_{\gamma v}^{T}\T_{\gamma v}^{(2)}&\C_{\gamma v}^{T}\T_{\gamma v}^{(2)}& \bm{0} &\alpha_{\gamma v}\Lambda_{v}^{*}
	\end{bmatrix}.
	\end{equation}
	The off-diagonal block matrices of $ \A + \A^T $ vanish for the conditions given in \cref{eq:Conditions for stabilty of adjoint inconsistent SATs}. Moreover, $ \bigl[\begin{smallmatrix}
	\T_{\gamma k}^{(1)}&-\T_{\gamma k}^{(1)}\\-\T_{\gamma k}^{(1)}&\T_{\gamma k}^{(1)}
	\end{smallmatrix}\bigr] = \bigl[ \begin{smallmatrix}
	1&-1\\-1&1 
	\end{smallmatrix}\bigr] \otimes \T_{\gamma k}^{(1)}$ is positive semidefinite if $ \T_{\gamma k}^{(1)} \succeq 0 $ due to properties of Kronecker products and because $ \bigl[\begin{smallmatrix}
	1&-1\\-1&1
	\end{smallmatrix}\bigr] \succeq 0 $. Finally, $ \bigl[\begin{smallmatrix}
	\alpha_{\gamma k}\Lambda_{k}^{*}& \bm{0}\\\bm{0}&\alpha_{\gamma v}\Lambda_{v}^{*}
	\end{smallmatrix}\bigr] \succeq 0$ since $ \alpha_{\gamma k}>0 $, $ \alpha_{\gamma v}>0 $, $ \Lambda^*_k \succeq 0 $, and $ \Lambda^*_v \succeq 0 $. Therefore, $ \A + \A^T \succeq 0 $ which completes the proof.
\end{proof}

\subsection{Energy analysis for adjoint consistent SATs} We consider a class of adjoint consistent SATs for which $ \T_{\gamma k}^{(3)}-\T_{\gamma k}^{(2)}=\B_{\gamma} $ and $ \T_{a k}^{(1)} $ is SPD, where $ a\in\{\gamma, \epsilon_{1},\epsilon_{2}, \delta_{1},\delta_{2}\} $. This class covers all types of adjoint consistent SAT studied in this work, \red{and \cref{thm:Stability for Adjoint Consistent SATs} provides sufficient conditions for energy stability of discretizations involving these SATs.} A class of adjoint consistent SATs for which $ \T_{abk}^{(5)}=\T_{abk}^{(6)}=\bm{0} $ in addition to the above two conditions is studied in \cite{yan2018interior}. 
\begin{theorem} \label{thm:Stability for Adjoint Consistent SATs}
	An adjoint consistent SBP-SAT discretization, \cref{eq:SBP-SAT discretization,eq:Interface SATs,eq:Boundary SATs}, of the homogeneous diffusion problem \cref{eq:diffusion problem}, \ie, $ \fnc{F}=0$, $\fnc{U}_D = 0 $, and $ \fnc{U}_N =0$, for which \violet{\cref{assu:Coefficient matrices} holds}, $ \T_{\gamma k}^{(3)}-\T_{\gamma k}^{(2)}=\B_{\gamma} $, and $ \T_{a k}^{(1)} \succ 0 $  is energy stable with respect to the \blue{diagonal norm} matrix, $ \H $, if the SAT coefficient matrices satisfy
	\begin{align}
	\T_{\gamma k}^{(1)}-\frac{2}{\zeta}\left(\frac{1}{\alpha_{\gamma k}}\T_{\gamma k}^{(2)}\Upsilon_{\gamma\gamma k}\T_{\gamma k}^{(2)}+\frac{1}{\alpha_{\gamma v}}\T_{\gamma v}^{(2)}\Upsilon_{\gamma\gamma v}\T_{\gamma v}^{(2)}\right) &\succeq 0, 
	\label{eq:T1 stablility 1}
	\\
	\T_{a k}^{(1)}-64\T_{abk}^{(5)}\left(\T_{bk}^{(1)}\right)^{-1}\T_{ba k}^{(5)}&\succeq 0, 
	\label{eq:T1 stablility 2}
	\\
	\T_{\gamma}^{(D)}-\frac{1}{\alpha_{\gamma k}}\B_{\gamma}\Upsilon_{\gamma k}\B_{\gamma} &\succeq 0,
	\label{eq:TD stablility}
	\\
	\T_{\gamma k}^{(4)} &\succeq 0,
	\label{eq:T4 stablility}	 	
	\end{align}
	where $ \Upsilon_{abk} $ is defined in \cref{eq:Upsilon definition}, $ a,b\in\{\gamma,\epsilon_{1},\epsilon_{2}\} $, $ \zeta = 2 $ \violet{for compact stencil SATs, \ie, SATs with $ \T_{abk}^{(5)}=\T_{abk}^{(6)}=\bm{0} $,} otherwise $ \zeta=1 $, and $ \T\succ 0 $ indicates $ \T $ is positive definite.
\end{theorem}
\begin{proof}
	We apply the conditions for adjoint consistency given in \cref{eq:Coefficients for adjoint consistency} on the residual \cref{eq:Residual 3rd form}. Because the residual is symmetric under these conditions, it is sufficient to show that $ R_h({\bm{u}_h, \bm{u}_h}) \le 0$. The conditions in \cref{eq:TD stablility,eq:T4 stablility} are established in \cref{thm:Stability Adjoint Inconsistent}. We rewrite the $ 4\times 4 $ block matrix in $ X_1 $, the second term in \cref{eq:Residual 3rd form}, as $
	\A= \bigl[\begin{smallmatrix}
	\A_{11} & \A_{12}\\
	\A_{12}^{T} & \A_{22}
	\end{smallmatrix}\bigr],
	$
	where 
	\begin{equation}
	\A_{11} = \frac{1}{2}\begin{bmatrix}
	\T_{\gamma k}^{(1)}&-\T_{\gamma k}^{(1)}\\-\T_{\gamma k}^{(1)}&\T_{\gamma k}^{(1)}
	\end{bmatrix},
	\;
	\A_{12} =\begin{bmatrix}
	\T_{\gamma k}^{(2)}\C_{\gamma 		k}&-\T_{\gamma v}^{(2)}\C_{\gamma v}\\-\T_{\gamma k}^{(2)}\C_{\gamma k}&\T_{\gamma v}^{(2)}\C_{\gamma v}
	\end{bmatrix},
	\;
	\A_{22} =\begin{bmatrix}
	\alpha_{\gamma k}\Lambda_{k}^{*}&\\&\alpha_{\gamma v}\Lambda_{v}^{*}
	\end{bmatrix}.
	\end{equation}
	Since $ \A_{22} $ is positive definite, \cref{thm:Positive semi-definiteness} ensures $ \A \succeq 0 $ if and only if $ \A_{11}-\A_{12}\A_{22}^{-1}\A_{12}^{T} \succeq 0 $, \ie,
	\begin{equation}
	\begin{bmatrix}
	1 & -1\\
	-1 & 1
	\end{bmatrix}\otimes\left[\frac{1}{2}\T_{\gamma k}^{(1)}-\left(\frac{1}{\alpha_{\gamma k}}\T_{\gamma k}^{(2)}\Upsilon_{\gamma\gamma k}\T_{\gamma k}^{(2)}+\frac{1}{\alpha_{\gamma v}}\T_{\gamma v}^{(2)}\Upsilon_{\gamma\gamma v}\T_{\gamma v}^{(2)}\right)\right]\succeq 0,
	\end{equation}
	which gives the condition for stability in \cref{eq:T1 stablility 1} with $ \zeta=1 $.	Setting $ \T^{(5)}_{abk} = \bm{0}$ and regrouping $ \T^{(1)}_{ak} $ terms as in \cref{eq:Matrix A for compact SATs}, the terms with extended SATs in \cref{eq:Residual 3rd form} vanish. Imposing $ \T_{\gamma k}^{(3)}-\T_{\gamma k}^{(2)}=\B_{\gamma} $ and the adjoint consistency conditions, we obtain $ \A_{11} = \bigl[\begin{smallmatrix}
	\T_{\gamma k}^{(1)}&-\T_{\gamma k}^{(1)}\\-\T_{\gamma k}^{(1)}&\T_{\gamma k}^{(1)}
	\end{smallmatrix}\bigr] $ while $ \A_{12} $ and $ \A_{22} $ remain unchanged. This yields the stability condition \cite{yan2018interior}
	\begin{equation}\label{eq:Stability condition on Tgk1 compact SATs}
	\T_{\gamma k}^{(1)}-\left(\frac{1}{\alpha_{\gamma k}}\T_{\gamma k}^{(2)}\Upsilon_{\gamma\gamma k}\T_{\gamma k}^{(2)}+\frac{1}{\alpha_{\gamma v}}\T_{\gamma v}^{(2)}\Upsilon_{\gamma\gamma v}\T_{\gamma v}^{(2)}\right) \succeq 0.
	\end{equation}
	After applying the conditions for adjoint consistency, the first $ 4 \times 4 $ block matrix in $ X_2 $, the third term in \cref{eq:Residual 3rd form}, reads
	\begin{equation}
	\begin{aligned}
	\G = \begin{bmatrix}
	\G_{11}&\G_{12}\\\G_{12}^{T}&\G_{22} 
	\end{bmatrix}
	\equiv\begin{bmatrix}
	\frac{1}{8}\T_{\gamma k}^{(1)}&-\frac{1}{8}\T_{\gamma k}^{(1)}&\T_{\gamma\epsilon_{1}k}^{(5)}&-\T_{\gamma\epsilon_{1}k}^{(5)}\\-\frac{1}{8}\T_{\gamma k}^{(1)}&\frac{1}{8}\T_{\gamma k}^{(1)}&-\T_{\gamma\epsilon_{1}k}^{(5)}&\T_{\gamma\epsilon_{1}k}^{(5)}\\\T_{\epsilon_{1}\gamma k}^{(5)}&-\T_{\epsilon_{1}\gamma k}^{(5)}&\frac{1}{8}\T_{\epsilon_{1}k}^{(1)}&-\frac{1}{8}\T_{\epsilon_{1}k}^{(1)}\\-\T_{\epsilon_{1}\gamma k}^{(5)}&\T_{\epsilon_{1}\gamma k}^{(5)}&-\frac{1}{8}\T_{\epsilon_{1}k}^{(1)}&\frac{1}{8}\T_{\epsilon_{1}k}^{(1)}
	\end{bmatrix},
	\end{aligned}
	\end{equation}
	where $ \G_{11} $, $ \G_{12} $, and $ \G_{22} $ are $ 2\times 2 $ block matrices. Energy stability requires that $ \G\succeq 0 $ which, using \cref{thm:Positive semi-definiteness}, implies we need to find conditions such that $ \left({\I}-\G_{22}\G_{22}^{+}\right)\G_{12}^{T}=\bm{0} $ and $ \G_{11}-\G_{12}\G_{22}^{+}\G_{12}^{T} \succeq 0 $ since $ \G_{22} \succeq 0 $. 
	But 
	\begin{equation}
	\G_{22}^{+}=\left(\begin{bmatrix}
	1 & -1\\
	-1 & 1
	\end{bmatrix}\otimes\frac{1}{8}\T_{\epsilon_{1}k}^{(1)}\right)^{+}
	=\begin{bmatrix}
	1 & -1\\
	-1 & 1
	\end{bmatrix}^{+}\otimes\left(\frac{1}{8}\T_{\epsilon_{1}k}^{(1)}\right)^{+}
	=\begin{bmatrix}
	1 & -1\\
	-1 & 1
	\end{bmatrix}\otimes\left(\frac{1}{2}\T_{\epsilon_{1}k}^{(1)}\right)^{-1},
	\end{equation}
	where we have used $ \bigl[\begin{smallmatrix}
	1&-1\\-1&1
	\end{smallmatrix}\bigr]^{+} = \frac{1}{4}\bigl[\begin{smallmatrix}
	1&-1\\-1&1
	\end{smallmatrix}\bigr]$ 
	and the fact that $ \T_{\epsilon_{1}k}^{(1)} $ is invertible for the class of SATs under consideration. Therefore,
	\begin{equation}
	\left(\I-\G_{22}\G_{22}^{+}\right)\G_{12}^{T}
	=\left\{\begin{bmatrix}
	\I\\
	& \I
	\end{bmatrix}-\left(\begin{bmatrix}
	1 & -1\\
	-1 & 1
	\end{bmatrix}\otimes\frac{1}{8}\T_{\epsilon_{1}k}^{(1)}\right)\left(\begin{bmatrix}
	1 & -1\\
	-1 & 1
	\end{bmatrix}\otimes\left(\frac{1}{2}\T_{\epsilon_{1}k}^{(1)}\right)^{-1}\right)\right\}\begin{bmatrix}
	\T_{\epsilon_{1}\gamma k}^{(5)}&-\T_{\epsilon_{1}\gamma k}^{(5)}\\-\T_{\epsilon_{1}\gamma k}^{(5)}&\T_{\epsilon_{1}\gamma k}^{(5)}
	\end{bmatrix}=\bm{0}.
	\end{equation}
	Furthermore, using properties of Kronecker products it can be shown that   
	\begin{equation}
	\begin{aligned}
	\G_{11}-\G_{12}\G_{22}^{+}\G_{12}^{T}&=\begin{bmatrix}
	1 & -1\\
	-1 & 1
	\end{bmatrix}\otimes\left(\frac{1}{8}\T_{\gamma k}^{(1)}\right)
	-\begin{bmatrix}
	1 & -1\\
	-1 & 1
	\end{bmatrix}^{3}\otimes\left(\T_{\gamma\epsilon_{1}k}^{(5)}\left(\frac{1}{2}\T_{\epsilon_{1}k}^{(1)}\right)^{-1}\T_{\epsilon_{1}\gamma k}^{(5)}\right)
	\\&=\begin{bmatrix}
	1 & -1\\
	-1 & 1
	\end{bmatrix}\otimes\left(\frac{1}{8}\T_{\gamma k}^{(1)}-8\T_{\gamma\epsilon_{1}k}^{(5)}\left(\T_{\epsilon_{1}k}^{(1)}\right)^{-1}\T_{\epsilon_{1}\gamma k}^{(5)}\right),
	\end{aligned}
	\end{equation}
	which implies $ \G_{11}-\G_{12}\G_{22}^{+}\G_{12}^{T}\succeq 0 $ if 
	\begin{equation}
	\T_{\gamma k}^{(1)}-64\T_{\gamma\epsilon_{1}k}^{(5)}\left(\T_{\epsilon_{1}k}^{(1)}\right)^{-1}\T_{\epsilon_{1}\gamma k}^{(5)}\succeq 0.
	\end{equation}
	Similar energy analyses for the rest of the $ 4\times 4 $ block matrices in $ X_2 $ (or simple geometric arguments) reveal that all the $ 4\times 4 $ block matrices in $ X_2 $ are positive semidefinite if 
	\begin{equation} \label{eq:Stability condition on Tak1 extended SAT}
	\T_{a k}^{(1)}-64\T_{abk}^{(5)}\left(\T_{bk}^{(1)}\right)^{-1}\T_{ba k}^{(5)}\succeq 0,
	\end{equation}
	for $ a,b\in\{\gamma,\epsilon_{1},\epsilon_{2}\} $. We have shown that the conditions stated in \cref{thm:Stability for Adjoint Consistent SATs} are sufficient for all the $ 4\times 4 $ block matrices in the residual \cref{eq:Residual 3rd form} to be positive semidefinite; therefore, $R_h(\bm{u}_h, \bm{u}_h) \le 0$ as desired.
\end{proof}

\section{Existing and DG based SATs} \label{sec:Existing and DG SATs}
The SAT coefficients associated with different types of DG fluxes are obtained by discretizing the residual of the DG primal formulation of the Poisson problem, which has the general form \cite{arnold2002unified,peraire2008compact} 
\begin{equation} \label{eq:RHS of DG primal formulation}
	\begin{aligned}
		\fnc{R}(\fnc{U}_{h},\fnc{V})
		&=-\int_{\Omega}\lambda\nabla\fnc{U}_{h}\cdot\nabla\fnc{V}\dd{\Omega}
		+\int_{\Omega}\fnc{V}\fnc{F}\dd{\Omega}
		-\int_{\Gamma^{I}}\jump{\widehat{\fnc{U}}-\fnc{U}_{h}}\cdot\avg{\lambda\nabla\fnc{V}}+\avg{\widehat{\fnc{U}}-\fnc{U}_{h}}\jump{\lambda\nabla\fnc{V}}\dd{\Gamma}
		\\&\quad
		+\int_{\Gamma^{I}}\jump{\fnc{V}}\cdot\avg{\widehat{\vecfnc{W}}}+\avg{\fnc{V}}\jump{\widehat{\vecfnc{W}}}\dd{\Gamma}
		+\int_{\Gamma^{D}}(\fnc{U}_{h}-\fnc{U}_{D})\lambda\nabla\fnc{V}\cdot\bm{n}+\fnc{V}\widehat{\vecfnc{W}}\cdot\bm{n}\dd{\Gamma}
		+\int_{\Gamma^{N}}\fnc{V}\fnc{U}_{N}\dd{\Gamma},
	\end{aligned}
\end{equation}
where $ \widehat{\fnc{U}} $ and $ \widehat{\vecfnc{W}} $ are numerical fluxes of the solution, $ \fnc{U}_h $, and the auxiliary variable $ \vecfnc{W}_h$, respectively. Equation \cref{eq:RHS of DG primal formulation} is obtained after setting the numerical fluxes of the solution as $ \widehat{\fnc{U}}= \fnc{U}_D $ on $ \Gamma^D $ and $ \widehat{\fnc{U}}= \fnc{U}_h  $ on $ \Gamma^N $, and the numerical fluxes of the auxiliary variable on $ \Gamma^N $ as $ \widehat{\fnc{W}} = \fnc{U}_N $. For schemes with global lifting operators, the auxiliary variable, the solution, and the flux of the solution are related by 
\begin{equation} \label{eq:Auxiliary variable W}
	\vecfnc{W}_{h}=\lambda\nabla\fnc{U}_{h}-\lambda\fnc{L}\left(\jump{\widehat{\fnc{U}}-\fnc{U}_{h}}\right)-\lambda\fnc{S}\left(\avg{\widehat{\fnc{U}}-\fnc{U}_{h}}\right)-\lambda\fnc{S}^{D}(\widehat{\fnc{U}}-\fnc{U}_{h}).
\end{equation}
For compact SATs, the global lifting operators in \cref{eq:Auxiliary variable W} are replaced by local lifting operators, \ie,
\begin{equation} \label{eq:Auxiliary variable W local}
\vecfnc{W}_{h}^{\gamma}=\lambda\nabla\fnc{U}_{h}-\lambda\fnc{L}^\gamma\left(\jump{\widehat{\fnc{U}}-\fnc{U}_{h}}\right)-\lambda\fnc{S}^\gamma \left(\avg{\widehat{\fnc{U}}-\fnc{U}_{h}}\right)-\lambda\fnc{S}^{D}(\widehat{\fnc{U}}-\fnc{U}_{h}).
\end{equation}

\blue{The forms of the interior facet SATs in \cref{eq:Interface SATs} and boundary SATs in \cref{eq:Boundary SATs} are closely related to the integral terms on the interior and boundary facets in the DG primal formulation. For example, to see how the boundary SATs and boundary integral terms in \cref{eq:RHS of DG primal formulation} are related, integrate by parts the first term on the RHS of \cref{eq:RHS of DG primal formulation} and substitute $\widehat{\vecfnc{W}}= \vecfnc{W}^{\gamma}_h $ on $ \Gamma^D $, to obtain
\begin{equation} \label{eq:RHS of DG primal formulation 2}
\begin{aligned}
\fnc{R}(\fnc{U}_{h},\fnc{V})
&=\int_{\Omega}\fnc{V}\nabla\cdot(\lambda\nabla\fnc{U}_{h}) \dd\Omega
+\int_{\Omega}\fnc{V}\fnc{F}\dd\Omega - \int_{\Gamma^{I}}\fnc{V}(\lambda\nabla\fnc{U}_{h})\cdot \bm{n}\dd \Gamma
\\&\quad -\int_{\Gamma^{I}}\jump{\widehat{\fnc{U}}-\fnc{U}_{h}}\cdot\avg{\lambda\nabla\fnc{V}}+\avg{\widehat{\fnc{U}}-\fnc{U}_{h}}\jump{\lambda\nabla\fnc{V}}\dd \Gamma
+\int_{\Gamma^{I}}\jump{\fnc{V}}\cdot\avg{\widehat{\vecfnc{W}}}+\avg{\fnc{V}}\jump{\widehat{\vecfnc{W}}}\dd\Gamma
\\&\quad
+\int_{\Gamma^{D}}(\fnc{U}_{h}-\fnc{U}_{D})\lambda\nabla\fnc{V}\cdot\bm{n}- \fnc{V}\lambda\fnc{S}^{D}(\widehat{\fnc{U}}-\fnc{U}_{h})\cdot \bm{n}\dd\Gamma
+\int_{\Gamma^{N}}\fnc{V}(\fnc{U}_{N} -(\lambda\nabla\fnc{U}_{h})\cdot \bm{n}) \dd\Gamma.
\end{aligned}
\end{equation}
The discrete analogue of the boundary integral terms in the last line of \cref{eq:RHS of DG primal formulation 2} is of the same form as $ \bm{v}_k^T\bm{s}_k^B $. The structure of the boundary SATs remains unchanged for DG fluxes based on global lifting operators due to the definition of the global lifting operators, \cref{eq: lift global vector,eq: lift global scalar}, which involve only interior facet integrals.} \red{The connection between the interface coupling terms in the DG formulation and the interior facet SATs can be shown by discretizing the interior surface integrals in \cref{eq:RHS of DG primal formulation 2}, which requires discretization of the lifting operators that appear in the numerical fluxes.} To find the discrete analogues of the global lifting operator for vector functions, we first write \cref{eq: lift global vector} for $ \Omega_k $ as 
\begin{equation} \label{eq:Global lifting on Omegak}
	\int_{\Omega_{k}}\lambda_{k}{\mathcal L}_{k}\left(\jump{\fnc{U}_h}\right)\cdot{\vecfnc{Z}}_{k}\dd{\Omega}=-\frac{1}{2}\int_{\Gamma_{k}^I}\lambda_{k}\jump{\fnc{U}_h}\cdot{\vecfnc{Z}}_{k}\dd{\Gamma}=-\frac{1}{2}\int_{\Gamma_{k}^I}\jump{\fnc{U}_h}\cdot\lambda_{k}^{T}{\vecfnc{Z}}_{k}\dd{\Gamma}.
\end{equation} 
Note that the sum of the lifting operators defined by \cref{eq:Global lifting on Omegak} at an interface shared by two elements is
\begin{equation}
	-\frac{1}{2}\int_{\gamma}\jump{\fnc{U}_h}\cdot\lambda_{k}^{T}\vecfnc{Z}{k}\dd{\Gamma}-\frac{1}{2}\int_{\gamma}\jump{\fnc{U}_h}\cdot\lambda_{v}^{T}\vecfnc{Z}_{v}\dd{\Gamma}=-\int_{\gamma}\jump{\fnc{U}_h}\cdot\frac{1}{2}\left(\lambda_{k}^{T}\vecfnc{Z}_{k}+\lambda_{v}^{T}\vecfnc{Z}_{v}\right)\dd{\Gamma}=-\int_{\gamma}\jump{\fnc{U}_h}\cdot\avg{\lambda^{T}\vecfnc{Z}}\dd{\Gamma},
\end{equation}
which enables \cref{eq: lift global vector} to be recovered upon summing over all interfaces. Neglecting truncation error, the discretization of \cref{eq:Global lifting on Omegak} follows as 
\begin{equation}
	\begin{bmatrix}
	{\bm{z}}_{x,k}\\
	{\bm{z}}_{y,k}
	\end{bmatrix}^{T}
	\begin{bmatrix}
	\H_{k}\\
	& \H_{k}
	\end{bmatrix}
	\begin{bmatrix}
	\mathscr{L}_{x,k}\\
	\mathscr{L}_{y,k}
	\end{bmatrix}=-\frac{1}{2}\sum_{\gamma\subset\Gamma_{k}^I}
	\begin{bmatrix}
	{\bm{z}}_{x,k}\\
	{\bm{z}}_{y,k}
	\end{bmatrix}^{T}
	\begin{bmatrix}
	\Lambda_{xx} & \Lambda_{xy}\\
	\Lambda_{yx} & \Lambda_{yy}
	\end{bmatrix}_{k}
	\begin{bmatrix}
	\R_{\gamma k}^{T}\\
	& \R_{\gamma k}^{T}
	\end{bmatrix}
	\begin{bmatrix}
	\N_{x,\gamma}\\
	\N_{y,\gamma}
	\end{bmatrix}\B_{\gamma}\left(\R_{\gamma k}\bm{u}_{k}-\R_{\gamma v}\bm{u}_{v}\right),
\end{equation}
and thus, the $ x $-coordinate discrete global lifting operator for vector functions, $ \mathscr{L}_{x,k} $, is given by 
\begin{equation} \label{eq:Lift global vector}
	\begin{aligned}
		\mathscr{L}_{x,k}&=-\frac{1}{2}\sum_{\gamma\subset\Gamma_{k}^I}\H_{k}^{-1}\left(\Lambda_{xx}\R_{\gamma k}^{T}\N_{x,\gamma}+\Lambda_{xy}\R_{\gamma k}^{T}\N_{y,\gamma}\right)\B_{\gamma}\left(\R_{\gamma k}\bm{u}_{k}-\R_{\gamma v}\bm{u}_{v}\right).
	\end{aligned}
\end{equation}
The $ y $-coordinate discrete global lifting operator, $ \mathscr{L}_{y,k} $, has analogues expression. For the same reason, we will state only the $ x $-coordinate discrete lifting operators for the other types of lifting operators presented below. The local lifting operator for vector functions on element $ \Omega_k $ and facet $ \gamma\in \Gamma_k^I $ is defined as
\begin{equation}  \label{eq:Local lifting on Omegak}
	\int_{\Omega_{k}} \lambda_{k}{\cal L}_{k}^{\gamma}\left(\jump{{\fnc{U}_h}}\right)\cdot{\vecfnc{Z}}_{k}\dd{\Omega}=-\frac{1}{2}\int_{\gamma}\jump{{\fnc{U}_h}}\cdot\lambda_{k}^{T}{\vecfnc{Z}}_{k}\dd{\Gamma},
\end{equation}
which upon discretization gives the $ x $-coordinate local lifting operator \cite{yan2018interior}
\begin{equation} \label{eq:Lift local vector}
		\mathscr{L}_{x,k}^{\gamma}=-\frac{1}{2}\H_{k}^{-1}\left(\Lambda_{xx}\R_{\gamma k}^{T}\N_{x,\gamma}+\Lambda_{xy}\R_{\gamma k}^{T}\N_{y,\gamma}\right)\B_{\gamma}\left(\R_{\gamma k}\bm{u}_{k}-\R_{\gamma v}\bm{u}_{v}\right).
\end{equation}
Applying a similar approach, we write the global lifting operator for scalar valued functions, \cref{eq: lift global scalar}, on a single element as  
\begin{align}
	\int_{\Omega_{k}}\lambda_k{\cal S}_{k}\left(\jump{\fnc{U}_h}\cdot\bm{n}_{k}\right)\cdot\vecfnc{Z}_{k}\dd{\Omega}=-\int_{\Gamma_{k}^{I}}\left(\jump{\fnc{U}_h}\cdot\bm{n}_{k}\right)\lambda_k^T\vecfnc{Z}_{k}\cdot\bm{n}_{k}\dd{\Gamma}=-\sum_{\gamma\subset\Gamma_{k}^{I}}\int_{\gamma}\left({\fnc{U}}_{h,k}-{\fnc{U}}_{h,v}\right)\lambda_k^T\vecfnc{Z}_{k}\cdot\bm{n}_{k}\dd{\Gamma}, 
\end{align}
which gives the $ x $-coordinate discrete analogue of the global lifting operator for scalar functions as,
\begin{equation} \label{eq:Lift global scalar}
	\mathscr{S}_{x,k}=-\sum_{\gamma\subset\Gamma_{k}^I}\H_{k}^{-1}\left(\Lambda_{xx}\R_{\gamma k}^{T}\N_{x,\gamma}+\Lambda_{xy}\R_{\gamma k}^{T}\N_{y,\gamma}\right)\B_{\gamma}\left(\R_{\gamma k}\bm{u}_{k}-\R_{\gamma v}\bm{u}_{v}\right).
\end{equation}
Moreover, the discretization of the local lifting operator for scalar functions at interior facet $ \gamma\in\Gamma_k^I $ gives
\begin{equation} \label{eq:Lift local scalar}
	\mathscr{S}_{x,k}^{\gamma}=-\H_{k}^{-1}\left(\Lambda_{xx}\R_{\gamma k}^{T}\N_{x,\gamma}+\Lambda_{xy}\R_{\gamma k}^{T}\N_{y,\gamma}\right)\B_{\gamma}\left(\R_{\gamma k}\bm{u}_{k}-\R_{\gamma v}\bm{u}_{v}\right).
\end{equation}
Finally, at Dirichlet boundary facets, the $ x $-coordinate discrete lifting operators is given by
\begin{equation} \label{eq:Lift Dirichlet}
		\mathscr{S}^D_{x,k}=-\H_{k}^{-1}\left(\Lambda_{xx}\R_{\gamma k}^{T}\N_{x,\gamma}+\Lambda_{xy}\R_{\gamma k}^{T}\N_{y,\gamma}\right)\B_{\gamma}\left(\R_{\gamma k}\bm{u}_{k}-\bm{u}_{\gamma k}\right).		
\end{equation}

Before we proceed with identification of SATs pertaining to known DG methods, we state the following two lemmas which will be useful to analyze energy stability of some of the schemes studied in the following subsections. 
\begin{lemma}\label{lem:Inverse of sum of SPD matrices}
	Let $ \X \in \IRtwo{n}{n} $ and $ \Y \in \IRtwo{n}{n} $ be two SPD matrices, then the inverse of the sum of the matrices satisfies
	\begin{equation}\label{eq:Inverse of sum of SPD matrics}
		\begin{aligned}
			\X^{-1}+\Y^{-1}-(\X+\Y)^{-1} &\succ 0, 
			&&
			\X^{-1}-(\X+\Y)^{-1} \succ 0,
			&& \text{and} \quad
			\Y^{-1}-(\X+\Y)^{-1} \succ 0,
		\end{aligned}
	\end{equation}	
	where $ \Y \succ 0 $ indicates that $ \Y $ is positive definite.
\end{lemma} 
\begin{proof}
	We start from the following result in \cite{henderson1981deriving}, 
	\begin{equation} 
		\begin{aligned}
			\left(\X+\Y\right)^{-1}&=\X^{-1}-\X^{-1}\Y\left(\X+\Y\right)^{-1}, && \text{or}\quad
			\left(\X+\Y\right)^{-1}=\Y^{-1}-\Y^{-1}\X\left(\X+\Y\right)^{-1},
		\end{aligned}
	\end{equation}
	and write 
	\begin{equation}\label{eq:Inverse of sum of SPD matrics 2}
		\X^{-1}+\Y^{-1}-\left(\X+\Y\right)^{-1}=\X^{-1}+\Y^{-1}-\X^{-1}+\X^{-1}\Y\left(\X+\Y\right)^{-1}=\Y^{-1}+\X^{-1}\Y\left(\X+\Y\right)^{-1}.
	\end{equation}
	Note that $ \X^{-1}\Y\left(\X+\Y\right)^{-1} $ is symmetric because $ \X^{-1}-\left(\X+\Y\right)^{-1} $ is symmetric (since the sum of two symmetric matrices is symmetric, and the inverse of a symmetric matrix is symmetric). Furthermore, $ \X^{-1}\Y\left(\X+\Y\right)^{-1} $ is positive definite because $ \X^{-1}\succ0 $, $ \Y^{-1}\succ0 $, $ \left(\X+\Y\right)^{-1}\succ0 $, and their product $ \X^{-1}\Y\left(\X+\Y\right)^{-1} $ is symmetric, \ie, we used the fact that the product of two SPD matrices is positive definite if their product is symmetric as well. Therefore, we obtain $ \Y^{-1}+\X^{-1}\Y\left(\X+\Y\right)^{-1}\succ0 $ which implies that \cref{eq:Inverse of sum of SPD matrics 2} yields the first inequality in \cref{eq:Inverse of sum of SPD matrics}. By a similar argument we can write
	\begin{align}
		\X^{-1}-\left(\X+\Y\right)^{-1}&=\X^{-1}-\X^{-1}+\X^{-1}\Y\left(\X+\Y\right)^{-1}=\X^{-1}\Y\left(\X+\Y\right)^{-1}\succ0, 
		\\
		\Y^{-1}-\left(\X+\Y\right)^{-1}&=\Y^{-1}-\Y^{-1}+\Y^{-1}\X\left(\X+\Y\right)^{-1}=\Y^{-1}\X\left(\X+\Y\right)^{-1}\succ0.
	\end{align}
\end{proof}
\begin{lemma} \label{lem:I-YXY}
	Given an SPD matrix $ \X \in \IRtwo{n}{n} $ and a rectangular matrix $ \Y \in \IRtwo{n}{m} $ such that $ \X=\Y\Y^{T} $, we have
	\begin{equation}
		\I_{m}-\Y^{T}\X^{-1}\Y\succeq0,
	\end{equation} 
	where $ \I_{m} $ is an $ m\times m $ identity matrix.
\end{lemma}
\begin{proof}
	Consider the singular value decomposition $ \Y=\U\Sigma_{r}\V^{T} $, then
	\begin{align}
		\Y^{T}\X^{-1}\Y&=\Y^{T}\left(\Y\Y^{T}\right)^{-1}\Y=\Y^{T}\left(\Y^{+T}\Y^{+}\right)\Y=\left(\Y^{T}\Y^{+T}\right)\left(\Y^{+}\Y\right) \label{eq:SPD proof 2}\\
		&=\left(\V\Sigma_{r}\U^{T}\U\Sigma_{r}^{+}\V^{T}\right)\left(\V\Sigma_{r}^{+}\U^{T}\U\Sigma_{r}\V^{T}\right)=\left(\V\I_{r}\V^{T}\right)\left(\V\I_{r}\V^{T}\right)=\V\I_{r}\V^{T},
	\end{align}
	where $ \I_r $ is a diagonal matrix containing unity in its diagonal only up to the $ n $-th row and column indices, \ie, up to the rank of $ \Y $. In the second equality in \cref{eq:SPD proof 2} we made use of the property $ \left(\Y\Y^T\right)^{-1}=\left(\Y\Y^T\right)^{+}=\Y^{+T}\Y^{+} $ since $ \X=\Y\Y^T $ is invertible. Noting that the identity matrix can be written as $ \I_{m}=\V\I_{m} \V^{T} $, we have
	\begin{equation}
		\I_{m}-\Y^{T}\X^{-1}\Y=\V\I_{m}\V^{T}-\V\I_{r}\V^{T}=\V\left(\I_{m}-\I_{r}\right)\V^{T}\succeq0,
	\end{equation}
	which is the desired result.
\end{proof}
\subsection{BR1 SAT: The first method of Bassi and Rebay}
The numerical fluxes for the BR1 method \cite{bassi1997highnavier} are $ \widehat{\fnc{U}}=\avg{\fnc{U}_h} $ and $ \widehat{\fnc{W}}=\avg{\fnc{W}_h} $. Substituting these fluxes in \cref{eq:RHS of DG primal formulation} and simplifying, the residual for the BR1 method becomes \cite{arnold2002unified}
\begin{equation} \label{eq:BR1 residual}
	\begin{aligned}
		{\cal R}\left({\fnc{U}_h},{\cal V}\right)&=-\int_{\Omega}\lambda\nabla{\fnc{U}_h}\cdot\nabla{\cal V}\dd{\Omega}
		+\int_{\Omega}{\cal V}{\cal F}\dd{\Omega}
		+\int_{\Gamma^I}\jump{{\fnc{U}_h}}\cdot\avg{\lambda\nabla{\cal V}}
		+\avg{\lambda\nabla{\fnc{U}_h}}\cdot\jump{{\cal V}}\dd{\Gamma}
		\\
		& \quad
		-\int_{\Omega} \lambda{\cal L}(\jump{{\fnc{U}_h}}){\cal L}(\jump{{\cal V}})\dd{\Omega} +\int_{\Gamma^{D}}(\fnc{U}_{h}-\fnc{U}_{D})\lambda\nabla\fnc{V}\cdot\bm{n} \dd{\Gamma}
		+\int_{\Gamma^{D}}\fnc{V}(\lambda_k\nabla\fnc{U}_k)\cdot\bm{n}\dd{\Gamma}
		\\
		& \quad
		+\int_{\Gamma^{D}}\fnc{V}\lambda\fnc{S}^D(\fnc{U}_{h}-\fnc{U}_D)\cdot\bm{n}\dd{\Gamma}
		+\int_{\Gamma^{N}}\fnc{V}\fnc{U}_{N}\dd{\Gamma}.
	\end{aligned}
\end{equation}

\blue{If the surface integrals on the RHS of the global lifting operators in \cref{eq: lift global vector,eq: lift global scalar} include all facets, then the discretization of the BR1 primal formulation gives the boundary SATs:
	\begin{equation}
	\begin{aligned}
	\bm{s}_{k}^{B}(\uhk,\bm{u}_{\gamma k}, \bm{w}_{\gamma k}) &=\sum_{\gamma\subset\Gamma^{D}}\left[\begin{array}{cc}
	\R_{\gamma k}^{T} & \D_{\gamma k}^{T}\end{array}\right]\left[\begin{array}{c}
	\T_{\gamma}^{(D)}\\
	-\B_{\gamma}
	\end{array}\right](\Rgk\uhk  -\bm{u}_{\gamma k})
	+\sum_{\gamma\subset\Gamma^{N}}\R_{\gamma k}^{T}\B_{\gamma}\left(\D_{\gamma k}\bm{u}_{h,k}-\bm{w}_{\gamma k}\right)
	\\
	& \quad
	+\frac{1}{2}\sum_{\gamma\subset\Gamma_{k}^{I}}\sum_{\epsilon\subset\Gamma_{k}^{D}}\R_{\gamma k}^{T}\B_{\gamma}\Upsilon_{\gamma\epsilon k}\B_{\epsilon}\left(\R_{\epsilon k}\bm{u}_{h,k}-\bm{u}_{\epsilon k}\right)
	-\frac{1}{2}\sum_{\gamma\subset\Gamma_{v}^{I}}\sum_{\epsilon\subset\Gamma_{k}^{D}}\R_{\gamma v}^{T}\B_{\gamma}\Upsilon_{\gamma\epsilon k}\B_{\epsilon}\left(\R_{\epsilon k}\bm{u}_{h,k}-\bm{u}_{\epsilon k}\right)
	\\
	&\quad +\sum_{\gamma\subset\Gamma_{k}^{D}}\sum_{\epsilon\subset\Gamma_{k}^{D}}\R_{\gamma k}^{T}\B_{\gamma}\Upsilon_{\gamma\epsilon k}\B_{\epsilon}\left(\R_{\epsilon k}\bm{u}_{h,k}-\bm{u}_{\epsilon k}\right).
	\end{aligned}
	\end{equation} 
	Furthermore, the following term must be added to the interior facet SATs given by \cref{eq:Interface SATs}, 
	\begin{equation*}
	\frac{1}{2}\sum_{\gamma\subset\Gamma_{k}^{D}}\sum_{\epsilon\subset\Gamma_{k}^{I}}\bm{v}_{k}^{T}\R_{\gamma k}^{T}\B_{\gamma}\Upsilon_{\gamma\epsilon k}\B_{\epsilon}\left(\R_{\epsilon k}\bm{u}_{h,k}-\R_{\epsilon g}\bm{u}_{g}\right).
	\end{equation*}
	With these terms added to the interior and boundary facet SATs, it is possible to show that the SBP-SAT discretizations based on the primal and flux formulations of the BR1 method are identical. The BR1 SATs based on the flux formulation can be found, \eg, in Theorem 6.2  of \cite{chen2020review} by setting $ \beta=\alpha=0 $ therein. The extended boundary SATs affect adjoint consistency (and functional superconvergence) adversely, as discussed in \cref{sec:Adjoint Consistency}; however, not including them compromises the energy stability of the scheme. We now propose a modified BR1 type SAT that is stable but does not have extended boundary terms.}
\begin{proposition} \label{prop:BR1 SAT}
	A stabilized version of the BR1 type SAT is recovered if the coefficient matrices in \cref{eq:Residual 1st form} are set as
	\begin{equation*} \label{eq:BR1 coefficients}
		\begin{aligned}
			\T_{\gamma k}^{(1)}&=\T_{\gamma v}^{(1)}=\frac{1}{2}\B_{\gamma}\left[\frac{1}{\alpha_{\gamma k}}\Upsilon_{\gamma\gamma k}+\frac{1}{\alpha_{\gamma v}}\Upsilon_{\gamma\gamma v}\right]\B_{\gamma},&\T_{\gamma k}^{(3)}&=\T_{\gamma v}^{(3)}=\frac{1}{2}\B_{\gamma},&\T_{\gamma k}^{(2)}&=\T_{\gamma v}^{(2)}=-\frac{1}{2}\B_{\gamma},&\T_{\gamma k}^{(4)}&=\T_{\gamma v}^{(4)}=\bm{0},\\\T_{\gamma\epsilon k}^{(5)}&=-\T_{\gamma\epsilon k}^{(6)}=\frac{1}{16}\B_{\gamma}\Upsilon_{\gamma\epsilon k}\B_{\epsilon},&\T_{\gamma\delta v}^{(5)}&=-\T_{\gamma\delta v}^{(6)}=\frac{1}{16}\B_{\gamma}\Upsilon_{\gamma\delta v}\B_{\delta},&\T_{\gamma}^{(D)}&=\frac{1}{\alpha_{\gamma k}}\B_{\gamma}\Upsilon_{\gamma\gamma k}\B_{\gamma}.&&
		\end{aligned}
	\end{equation*}
	Moreover, the  BR1 SAT produces a consistent, conservative and adjoint consistent discretization.
\end{proposition}
\begin{proof}
	Discretizing \cref{eq:BR1 residual} using SBP operators and the discrete lifting operators \cref{eq:Lift global vector} and \cref{eq:Lift Dirichlet}, and comparing the result with \cref{eq:Residual 1st form} yields all the coefficients in \cref{prop:BR1 SAT} except $ \T_{\gamma k}^{(1)} $, $ \T_{\gamma \epsilon k}^{(5)} $, $ \T_{\gamma \delta v}^{(5)} $, and $ \T_{\gamma}^{(D)} $ which are modified for stability reasons. Before modification these coefficients read $ \T_{\gamma k}^{(1)}= (1/4)\B_{\gamma}[\Upsilon_{\gamma\gamma k}+\Upsilon_{\gamma\gamma v}]\B_{\gamma} $, $ \T_{\gamma\epsilon k}^{(5)}=(1/4)\B_{\gamma}\Upsilon_{\gamma\epsilon k}\B_{\epsilon} $, $ \T_{\gamma\delta v}^{(5)}= (1/4)\B_{\gamma}\Upsilon_{\gamma\delta v}\B_{\delta} $, and $ \T_{\gamma}^{(D)}=\B_{\gamma}\Upsilon_{\gamma \gamma k}\B_{\gamma}, $ which do not lead to stable discretization according to \cref{thm:Stability for Adjoint Consistent SATs}. In order to prove that the coefficients presented in \cref{prop:BR1 SAT} lead to stable discretization, we have to show that all the conditions in \cref{thm:Stability for Adjoint Consistent SATs} are met. From \cref{eq:TD stablility} and \cref{eq:T4 stablility}, we immediately see that the conditions on $ \T^{(D)}_{\gamma} $ and $ \T_{\gamma k}^{(4)} $ are satisfied. Substituting $ \T_{\gamma k}^{(2)} $, $ \T_{\gamma v}^{(2)} $, and the modified $ \T_{\gamma k}^{(1)} $ in \cref{eq:T1 stablility 1} we have 
	\begin{equation}
		\frac{1}{2}\B_{\gamma}\left[\frac{1}{\alpha_{\gamma k}}\Upsilon_{\gamma\gamma k}+\frac{1}{\alpha_{\gamma v}}\Upsilon_{\gamma\gamma v}\right]\B_{\gamma}-\frac{2}{4}\left(\frac{1}{\alpha_{\gamma k}}\B_{\gamma}\Upsilon_{\gamma\gamma k}\B_{\gamma}+\frac{1}{\alpha_{\gamma v}}\B_{\gamma}\Upsilon_{\gamma\gamma v}\B_{\gamma}\right)=\bm{0}.
	\end{equation}
	It remains to show that \cref{eq:T1 stablility 2} is satisfied. In \cref{thm:Stability for Adjoint Consistent SATs} we assumed $ \T_{ak}^{(1)} \succ 0 $ for $ a\in \{\gamma,\epsilon_{1},\epsilon_{2}\} $, which implies that $ \T_{ak}^{(1)} $ is invertible. This is achieved by the proposed $ \T_{\gamma k}^{(1)} $ coefficient since $ \Upsilon_{\gamma\gamma k}\succ 0 $. Note that in  \cref{eq:Upsilon definition} we have $ \bar{\H}_k^{-1}\Lambda_k \succ 0 $, and the normals in both $ x $ and $ y $ directions cannot be zero simultaneously. Using the proposed coefficients, we have
	\begin{align}
		\T_{\gamma\epsilon_{1}k}^{(5)}\left(\T_{\epsilon_{1}k}^{(1)}\right)^{-1}\T_{\epsilon_{1}\gamma k}^{(5)}&=\frac{1}{16}\B_{\gamma}\Upsilon_{\gamma\epsilon_{1}k}\B_{\epsilon_{1}}\left(\frac{1}{2}\B_{\epsilon_{1}}\left[\frac{1}{\alpha_{\epsilon_{1}k}}\Upsilon_{\epsilon_{1}\epsilon_{1}k}+\frac{1}{\alpha_{\epsilon_{1}g_{1}}}\Upsilon_{\epsilon_{1}\epsilon_{1}g_{1}}\right]\B_{\epsilon_{1}}\right)^{-1}\frac{1}{16}\B_{\epsilon_{1}}\Upsilon_{\epsilon_{1}\gamma k}\B_{\gamma} 
		\nonumber
		\\&=\frac{1}{128}\B_{\gamma}\Upsilon_{\gamma\epsilon_{1}k}\left[\frac{1}{\alpha_{\epsilon_{1}k}}\Upsilon_{\epsilon_{1}\epsilon_{1}k}+\frac{1}{\alpha_{\epsilon_{1}g_{1}}}\Upsilon_{\epsilon_{1}\epsilon_{1}g_{1}}\right]^{-1}\Upsilon_{\epsilon_{1}\gamma k}\B_{\gamma}
		\preceq\frac{1}{128}\B_{\gamma}\Upsilon_{\gamma\epsilon_{1}k}\left[\frac{1}{\alpha_{\epsilon_{1}k}}\Upsilon_{\epsilon_{1}\epsilon_{1}k}\right]^{-1}\Upsilon_{\epsilon_{1}\gamma k}\B_{\gamma}, \label{eq:T5 T1 T5}
	\end{align} 
	where we applied \cref{lem:Inverse of sum of SPD matrices} in the last step. But we can write 
	\begin{align}
		\Upsilon_{\gamma\epsilon_{1}k}\left[\frac{1}{\alpha_{\epsilon_{1}k}}\Upsilon_{\epsilon_{1}\epsilon_{1}k}\right]^{-1}\Upsilon_{\epsilon_{1}\gamma k}&=\N_{\gamma k}^{T}\bar{\R}_{\gamma k}\bar{\H}_{k}^{-1}\Lambda_{k}\bar{\R}_{\epsilon_{1}k}^{T}\N_{\epsilon_{1}k}\left(\frac{1}{\alpha_{\epsilon_{1}k}}\N_{\epsilon_{1}k}^{T}\bar{\R}_{\epsilon_{1}k}\bar{\H}_{k}^{-1}\Lambda_{k}\bar{\R}_{\epsilon_{1}k}^{T}\N_{\epsilon_{1}k}\right)^{-1}\N_{\epsilon_{1}k}^{T}\bar{\R}_{\epsilon_{1}k}\bar{\H}_{k}^{-1}\Lambda_{k}\bar{\R}_{\gamma k}^{T}\N_{\gamma k}
		\nonumber
		\\&
		=\alpha_{\epsilon_{1}k}\P^{T}\Y^{T}\left[\Y\Y^{T}\right]^{-1}\Y\P,
	\end{align}
	where $ \P=\left[\bar{\H}_{k}^{-1}\Lambda_{k}\right]^{\frac{1}{2}}\bar{\R}_{\gamma k}^{T}\N_{\gamma k} $ and $ \Y=\N_{\epsilon_{1}k}^{T}\bar{\R}_{\epsilon_{1}k}\left[\bar{\H}_{k}^{-1}\Lambda_{k}\right]^{\frac{1}{2}} $. \cref{lem:I-YXY} implies
	$
		(\I-\Y^{T}[\Y\Y^{T}]^{-1}\Y)\succeq0,
	$
	which gives
	\begin{align} \label{eq:Lemma 3 applied}
		\alpha_{\epsilon_{1}k}\P^{T}\I\P-\Upsilon_{\gamma\epsilon_{1}k}\left[\frac{1}{\alpha_{\epsilon_{1}k}}\Upsilon_{\epsilon_{1}\epsilon_{1}k}\right]^{-1}\Upsilon_{\epsilon_{1}\gamma k}=\alpha_{\epsilon_{1}k}\Upsilon_{\gamma\gamma k}-\Upsilon_{\gamma\epsilon_{1}k}\left[\frac{1}{\alpha_{\epsilon_{1}k}}\Upsilon_{\epsilon_{1}\epsilon_{1}k}\right]^{-1}\Upsilon_{\epsilon_{1}\gamma k}\succeq0.
	\end{align}
	Since \violet{$ 0 < \alpha_{\epsilon_{1}k} < 1 $},  $ \Upsilon_{\gamma\gamma k} \succ 0$, and $ \Upsilon_{\gamma\gamma v} \succ 0 $, we write
	\begin{equation}
		\left[\frac{1}{\alpha_{\gamma k}}\Upsilon_{\gamma\gamma k}+\frac{1}{\alpha_{\gamma v}}\Upsilon_{\gamma\gamma v}\right]-\Upsilon_{\gamma\epsilon_{1}k}\left[\frac{1}{\alpha_{\epsilon_{1}k}}\Upsilon_{\epsilon_{1}\epsilon_{1}k}\right]^{-1}\Upsilon_{\epsilon_{1}\gamma k}\succeq0,
	\end{equation}
	which, together with \cref{eq:T5 T1 T5}, yields the inequality
	\begin{equation}
		\T_{\gamma\epsilon_{1}k}^{(5)}\left(\T_{\epsilon_{1}k}^{(1)}\right)^{-1}\T_{\epsilon_{1}\gamma k}^{(5)}\preceq\frac{1}{128}\B_{\gamma}\Upsilon_{\gamma\epsilon_{1}k}\left[\frac{1}{\alpha_{\epsilon_{1}k}}\Upsilon_{\epsilon_{1}\epsilon_{1}k}\right]^{-1}\Upsilon_{\epsilon_{1}\gamma k}\B_{\gamma}\preceq\frac{1}{128}\B_{\gamma}\left[\frac{1}{\alpha_{\gamma k}}\Upsilon_{\gamma\gamma k}+\frac{1}{\alpha_{\gamma v}}\Upsilon_{\gamma\gamma v}\right]\B_{\gamma}.
	\end{equation}
	 Therefore, 
	 \begin{equation}
	 	\T_{\gamma k}^{(1)}-64\T_{\gamma\epsilon_{1}k}^{(5)}\left(\T_{\epsilon_{1}k}^{(1)}\right)^{-1}\T_{\epsilon_{1}\gamma k}^{(5)}
	 	\succeq
	 	\T_{\gamma k}^{(1)}-\frac{1}{2}\B_{\gamma}\left[\frac{1}{\alpha_{\gamma k}}\Upsilon_{\gamma\gamma k}+\frac{1}{\alpha_{\gamma v}}\Upsilon_{\gamma\gamma v}\right]\B_{\gamma}=\bm{0},
	 \end{equation}
	 which is the result required for \cref{eq:T1 stablility 2} to hold. Note that the same analysis can be done for any combination of facets  $ a,b \in \{\gamma,\epsilon_{1},\epsilon_{2}\} $, in \cref{eq:T1 stablility 2}. Finally, from \cref{thm:Consistency,thm:Conservation,thm:Adjoint consistency} it can easily be seen that the BR1 SAT satisfies all the conditions required for consistency, conservation, and adjoint consistency.	 
\end{proof}
\begin{remark}
	The proposed \violet{interior facet} BR1 SAT is equivalent to the consistent method of Brezzi \etal \cite{brezzi1999discontinuous}, the modified BR1 method in \cite{alhawwary2018accuracy}, the stabilized central flux in \cite{hesthaven2007nodal}, and the penalty approach in \cite{kannan2009study} in the sense that all of these methods can be reproduced by considering $ \sigma_1 \T_{\gamma k}^{(1)} $, $ \sigma_5 \T_{\gamma \epsilon k}^{(5)} $, and $ \sigma_D \T_{\gamma}^{(D)} $ in \cref{prop:BR1 SAT} for $ \sigma_1,\sigma_5,\sigma_D > 0 $.  
\end{remark}

\begin{remark} \label{rem:Sigma reduced}
	Assuming the source term is zero, it can be shown that the continuous energy estimate satisfies \cite{gassner2018br1}
	\begin{equation} \label{eq:Energy BR1 continuous}
		\frac{1}{2}\frac{{\rm d}}{{\rm d}t}\norm{\fnc{U}_{h}}{}^{2}=\fnc{R}(\fnc{U}_{h},\fnc{U}_{h})\le\sum_{\gamma\subset\Gamma}\int_{\gamma}\widehat{\fnc{U}}\jump{\vecfnc{W}_{h}}+\widehat{\vecfnc{W}}\cdot\jump{\fnc{U}_{h}}-\jump{\fnc{U}_{h}\vecfnc{W}_{h}}\dd{\Gamma}.
	\end{equation}
	Substituting the BR1 fluxes into \cref{eq:Energy BR1 continuous} and using the identity 
	\begin{equation}\label{eq:Jump-avg identity}
		\avg{\fnc{U}_{h}}\jump{\vecfnc{W}_{h}}+\avg{\vecfnc{W}_{h}}\cdot\jump{\fnc{U}_{h}}-\jump{\fnc{U}_{h}\vecfnc{W}_{h}}=0
	\end{equation}
	gives $ {\rm d}/{{\rm d}t}(\norm{\fnc{U}_h}^2) \le 0 $, which establishes the energy stability of the BR1 method for diffusion problems. \red{A discrete energy stability analysis of the SBP-SAT discretization based on the flux formulation leads to a similar conclusion. Such a proof follows the same technique used to show the entropy stability of the LDG method in \cite{chen2020review} (see, proof of Theorem 6.2 therein). If the BR1 SAT is applied only on the interior facets, however, the identity \cref{eq:Jump-avg identity} cannot be applied, and the energy stability of the discretization is compromised unless the SAT coefficients are modified, \eg, as in \cref{prop:BR1 SAT}. Exceptions apply when the BR1 and LDG SATs are implemented with the SBP diagonal-E operators for which all the extended SATs vanish (see, \cref{sec:Equivalence of SATs}). In this case, the discrete form of \cref{eq:Jump-avg identity} can be used to show the energy stability of the SBP-SAT discretization based on the flux formulation, which yields the same discretization as the primal formulation when implemented with the SBP diagonal-E operators.}
\end{remark}

\subsection{LDG SAT: The local discontinuous Galerkin method}
The LDG scheme \cite{cockburn1998local} is obtained by choosing the DG numerical fluxes as $ \fnc{\widehat{U}} = \avg{\fnc{U}_h} - \bm{\beta}\cdot\jump{\fnc{U}_h}$ and $ \widehat{\vecfnc{W}} = \avg{\vecfnc{W}_h} + \bm{\beta}\jump{\vecfnc{W}_h} - \mu h^{-1}\jump{\fnc{U}_h} $ on interior facets, and $ \widehat{\vecfnc{W}} = \vecfnc{W}_h - \mu h^{-1}(\fnc{U}_h - \fnc{U}_D)\bm{n} $ on $ \Gamma^D $ \cite{arnold2002unified,peraire2008compact}. The switch function, $ \bm{\beta} $, is defined on each interface as 
\begin{equation}\label{eq:LDG switch}
	\bm{\beta} = \frac{1}{2}(\beta_{\gamma k} \bm{n}_k + \beta_{\gamma v} \bm{n}_v),
\end{equation}
where $ \beta_k,\beta_v\in\{0, 1\} $ are switches defined for $ \Omega_k $ and $ \Omega_v $ at their shared interface. Furthermore, the switches satisfy  
\begin{equation}
\beta_{\gamma k} + \beta_{\gamma v} = 1.
\end{equation}
The values of the switches \violet{are set to zero at boundary facets, \ie, $ \beta_{\gamma k}=\beta_{\gamma v}=0 $ for $ \gamma\subset \Gamma^B $. For interior facets, the values} are determined based on the sign of the dot product $ \bm{n}\cdot \bm{g} $, where $ \bm{g} $ is an arbitrary global vector \cite{sherwin20062d}, \ie,
\begin{equation} \label{eq:LDG switch with g}
	\begin{aligned}
		\beta_{\gamma k}=\begin{cases}
		1 & {\rm if}\;\bm{n}_{k}\cdot\bm{g}\ge0,\\
		0 & {\rm if}\;\bm{n}_{k}\cdot\bm{g}<0.
		\end{cases}
	\end{aligned}
\end{equation}
Although it is possible to use other vectors as switch functions, the form in \cref{eq:LDG switch} is necessary to avoid wider stencil width \cite{sherwin20062d,peraire2008compact}. For instance, if we set $ \bm{\beta} = \bm{0} $ and $ \mu=0 $, we recover the BR1 fluxes. \violet{On curved elements, the normal vector varies along the facets; hence, $ \beta_{\gamma k} $ is not constant in general. This leads to cases where \cref{assu:Coefficient matrices} does not hold; particularly, $ \T_{\gamma k}^{(1)} \neq (\T_{\gamma k}^{(1)})^T$, $ \T_{abk}^{(5)} \neq (\T_{bak}^{(5)})^T $, $ \T_{abk}^{(6)} \neq (\T_{bak}^{(6)})^T $. Additionally, it increases the number of elements that are coupled, resulting in a denser system matrix. To remedy this, we calculate $ \beta_{\gamma k} $ using straight facets in 2D (or flat facets in 3D) regardless of whether or not the physical elements are curved.}

Substituting the numerical fluxes in \cref{eq:RHS of DG primal formulation} and simplifying, the residual of the LDG method reads
\begin{equation} \label{eq:LDG residual}
	\begin{aligned}
		\fnc{R}\left(\fnc{U}_{h},\fnc{V}\right)&=-\int_{\Omega}\lambda\nabla\fnc{U}_{h}\cdot\nabla\fnc{V}\dd{\Omega}+\int_{\Omega}\fnc{V}\fnc{F}\dd{\Omega}+\int_{\Gamma^{I}}\jump{\fnc{U}_{h}}\cdot\left\{ \lambda\nabla\fnc{V}\right\} +\jump{\fnc{V}}\cdot\avg{\lambda\nabla\fnc{U}_{h}}\dd{\Gamma}
		\\&\quad
		+\int_{\Gamma^{I}}\bm{\beta}\cdot\jump{\fnc{U}_{h}}\jump{\lambda\nabla\fnc{V}}+\jump{\lambda\nabla\fnc{U}_{h}}\bm{\beta}\cdot\jump{\fnc{V}}\dd{\Gamma}-\mu h^{-1}\int_{\Gamma^{I}}\jump{\fnc{V}}\cdot\jump{\fnc{U}_{h}}\dd{\Gamma}
		\\&\quad-\int_{\Omega}\left[\fnc{L}\left(\jump{\fnc{V}}\right)+\fnc{S}\left(\bm{\beta}\cdot\jump{\fnc{V}}\right)\right]\cdot\left[\lambda\fnc{L}\left(\jump{\fnc{U}_{h}}\right)+\lambda\fnc{S}\left(\bm{\beta}\cdot\jump{\fnc{U}_{h}}\right)\right]\dd{\Omega}-\mu h^{-1}\int_{\Gamma^{D}}\fnc{V}\left(\fnc{U}_{h}-\fnc{U}_D\right)\dd{\Gamma}
		\\&\quad +\int_{\Gamma^{D}}\left(\fnc{U}_{h}-\fnc{U}_{D}\right)\lambda\nabla\fnc{V}\cdot\bm{n}\dd{\Gamma}+\int_{\Gamma^{D}}\fnc{V}(\lambda\nabla\fnc{U}_{h})\cdot \bm{n}\dd{\Gamma}+\int_{\Gamma^{D}}\fnc{V}\lambda\fnc{S}^D(\fnc{U}_{h}-\fnc{U}_D)\cdot\bm{n}\dd{\Gamma}+\int_{\Gamma^{N}}\fnc{V}\fnc{U}_{N}\dd{\Gamma}.
	\end{aligned}
\end{equation}
\red{The boundary terms resulting from the discretization of \cref{eq:LDG residual} are different from the LDG boundary coupling terms obtained using global lifting operators defined on all interfaces. We have also used $ -\mu h^{-1}\int_{\Gamma^D}\fnc{V}(\fnc{U}_h - \fnc{U}_D) \dd{\Gamma}$ as the boundary stabilizing term instead of $ -\mu h^{-1}\int_{\Gamma^D}\fnc{V}\fnc{U}_h \dd{\Gamma}$. If these changes are not applied, the LDG boundary SATs would include extended stencil terms, \ie, 
\begin{equation}
\begin{aligned}
	\bm{s}_{k}^{B}(\uhk,\bm{u}_{\gamma k}, \bm{w}_{\gamma k}) &=\sum_{\gamma\subset\Gamma^{D}}\left[\begin{array}{cc}
	\R_{\gamma k}^{T} & \D_{\gamma k}^{T}\end{array}\right]\left[\begin{array}{c}
	\T_{\gamma}^{(D)}\\
	-\B_{\gamma}
	\end{array}\right](\Rgk\uhk  -\bm{u}_{\gamma k})
	+\sum_{\gamma\subset\Gamma^{N}}\R_{\gamma k}^{T}\B_{\gamma}\left(\D_{\gamma k}\bm{u}_{h,k}-\bm{w}_{\gamma k}\right)
	\\
	& \quad
	+\frac{1+\beta_{\gamma k}-\beta_{\gamma v}}{2}\sum_{\gamma\subset\Gamma_{k}^{I}}\sum_{\epsilon\subset\Gamma_{k}^{D}}\R_{\gamma k}^{T}\B_{\gamma}\Upsilon_{\gamma\epsilon k}\B_{\epsilon}\left(\R_{\epsilon k}\bm{u}_{h,k}-\bm{u}_{\epsilon k}\right)
	+\sum_{\gamma\subset\Gamma_{k}^{D}}\sum_{\epsilon\subset\Gamma_{k}^{D}}\R_{\gamma k}^{T}\B_{\gamma}\Upsilon_{\gamma\epsilon k}\B_{\epsilon}\left(\R_{\epsilon k}\bm{u}_{h,k}-\bm{u}_{\epsilon k}\right)
	\\
	&\quad
	-\frac{1+\beta_{\gamma k}-\beta_{\gamma v}}{2}\sum_{\gamma\subset\Gamma_{v}^{I}}\sum_{\epsilon\subset\Gamma_{k}^{D}}\R_{\gamma v}^{T}\B_{\gamma}\Upsilon_{\gamma\epsilon k}\B_{\epsilon}\left(\R_{\epsilon k}\bm{u}_{h,k}-\bm{u}_{\epsilon k}\right) + \mu h^{-1}\R_{\gamma k}^T\B_{\gamma}\bm{u}_{h,k},
\end{aligned}
\end{equation} 
where $ \T_{\gamma k}^{(D)}=\B_{\gamma}\Upsilon_{\gamma\gamma k}\B_{\gamma} $. Moreover, the term 
\begin{equation*}
\frac{1+\beta_{\epsilon k}-\beta_{\epsilon v}}{2}\sum_{\gamma\subset\Gamma_{k}^{D}}\sum_{\epsilon\subset\Gamma_{k}^{I}}\bm{v}_{k}^{T}\R_{\gamma k}^{T}\B_{\gamma}\Upsilon_{\gamma\epsilon k}\B_{\epsilon}\left(\R_{\epsilon k}\bm{u}_{h,k}-\R_{\epsilon g}\bm{u}_{g}\right).
\end{equation*}
must be added to the interior facet SATs given by \cref{eq:Interface SATs}. Similar to the unmodified BR1 SAT, it can be shown that the unmodified LDG SAT based on the primal and flux formulations are the same. As explained in \cref{rem:Sigma reduced}, the LDG SAT written in flux formulation is energy stable (even when $ \mu=0 $); hence, it follows that the unmodified LDG SAT based on the primal formulation is energy stable.} The penalty coefficients corresponding to the \violet{interior LDG SATs} can be obtained by discretizing \cref{eq:LDG residual} and comparing the result with \cref{eq:Residual 1st form}. The coefficients $ \T_{\gamma k}^{(2)} $, $ \T_{\gamma v}^{(2)} $, $ \T_{\gamma k}^{(3)} $, $ \T_{\gamma v}^{(3)} $, $ \T_{\gamma k}^{(4)} $ and $ \T_{\gamma v}^{(4)} $ are the same as those presented in \cref{prop:LDG SAT} below. The rest of the coefficients are 
\begin{equation} \label{eq:LDG original coefficients}
	\begin{aligned}
		\T_{\gamma k}^{(1)}&=\T_{\gamma v}^{(1)}=\violet{\B_{\gamma}\bigg[\frac{1+\left(\beta_{k\gamma}-\beta_{v\gamma}\right)^{2}+2(\beta_{k\gamma}-\beta_{v\gamma})}{4}\Upsilon_{\gamma\gamma k}+\frac{1+\left(\beta_{k\gamma}-\beta_{v\gamma}\right)^{2}-2(\beta_{k\gamma}-\beta_{v\gamma})}{4}\Upsilon_{\gamma\gamma v}\bigg]\B_{\gamma}}+\mu h^{-1}\B_{\gamma},
		\\
		\T_{\gamma\epsilon k}^{(5)}&=-\T_{\gamma\epsilon k}^{(6)}=\frac{\left[1+\beta_{\gamma k}-\beta_{\gamma v}\right]\left[1+\beta_{\epsilon k}-\beta_{\epsilon g}\right]}{4}\B_{\gamma}\Upsilon_{\gamma\epsilon k}\B_{\epsilon},\\\T_{\gamma\delta v}^{(6)}&=-\T_{\gamma \delta v}^{(5)}=\frac{\left[\beta_{\gamma k}-\beta_{\gamma v}-1\right]\left[1+\beta_{\delta v}-\beta_{\delta q}\right]}{4}\B_{\gamma}\Upsilon_{\gamma\delta v}\B_{\delta},\\\T_{\gamma k}^{(D)}&=\B_{\gamma}\Upsilon_{\gamma\gamma k}\B_{\gamma}+\mu h^{-1}\B_{\gamma}.
	\end{aligned}
\end{equation}
As expected, the coefficients in  \cref{eq:LDG original coefficients} are identical to the unmodified BR1 SAT when $ \bm{\beta}$ and $ \mu $ are set to zero. For mesh independent SAT penalties, \ie, $ \mu=0$, the coefficients in \cref{eq:LDG original coefficients} do not guarantee energy stability. To see this, consider the case where $ \beta_{\gamma k} = 1 $ and $ \beta_{\gamma v} = 0 $; the stability requirement in \cref{eq:T1 stablility 1} demands positive semidefiniteness of
\violet{ 
\begin{equation*}
	\T_{\gamma k}^{(1)}-2\left(\frac{1}{\alpha_{\gamma k}}\T_{\gamma k}^{(2)}\Upsilon_{\gamma\gamma k}\T_{\gamma k}^{(2)}+\frac{1}{\alpha_{\gamma v}}\T_{\gamma v}^{(2)}\Upsilon_{\gamma\gamma v}\T_{\gamma v}^{(2)}\right)=\B_{\gamma}\Upsilon_{\gamma\gamma k}\B_{\gamma}-\frac{2}{\alpha_{\gamma k}}\B_{\gamma}\Upsilon_{\gamma\gamma k}\B_{\gamma}.
\end{equation*}
However, this cannot be achieved since $ 0< \alpha_{\gamma k}<1 $}. It is also clear from \cref{eq:TD stablility} that $ \T_{\gamma}^{(D)} $ is not large enough to ensure energy stability. To remedy this, we propose a stabilized form of the LDG scheme that does not have a mesh dependent stability parameter.
\begin{proposition}\label{prop:LDG SAT}
	A consistent, conservative, adjoint consistent, and stable LDG type SAT with no mesh dependent stabilization parameter, \ie, $ \mu=0 $, is obtained if the penalty coefficients in \cref{eq:Residual 1st form} are chosen such that 
	\begin{equation*}
		\begin{aligned}
			\T_{\gamma k}^{(1)}&=\T_{\gamma v}^{(1)}=\B_{\gamma}\bigg[\frac{\beta_{\gamma k}-\beta_{\gamma v}+1}{\alpha_{\gamma k}}\Upsilon_{\gamma\gamma k}+\frac{\beta_{\gamma v}-\beta_{\gamma k}+1}{\alpha_{\gamma v}}\Upsilon_{\gamma\gamma v}\bigg]\B_{\gamma},&\T_{\gamma}^{(D)}&=\frac{1}{\alpha_{\gamma k}}\B_{\gamma}\Upsilon_{\gamma\gamma k}\B_{\gamma},&\T_{\gamma k}^{(4)}&=\T_{\gamma v}^{(4)}=\bm{0},\\\T_{\gamma\epsilon k}^{(5)}&=-\T_{\gamma\epsilon k}^{(6)}=\frac{\left[1+\beta_{\gamma k}-\beta_{\gamma v}\right]\left[1+\beta_{\epsilon k}-\beta_{\epsilon g}\right]}{16}\B_{\gamma}\Upsilon_{\gamma\epsilon k}\B_{\epsilon},&\T_{\gamma k}^{(2)}&=\frac{\beta_{\gamma v}-\beta_{\gamma k}-1}{2}\B_{\gamma},&\T_{\gamma k}^{(3)}&=\frac{\beta_{\gamma v}-\beta_{\gamma k}+1}{2}\B_{\gamma},\\\T_{\gamma\delta v}^{(6)}&=-\T_{\gamma\delta v}^{(5)}=\frac{\left[\beta_{\gamma k}-\beta_{\gamma v}-1\right]\left[1+\beta_{\delta v}-\beta_{\delta q}\right]}{16}\B_{\gamma}\Upsilon_{\gamma\delta v}\B_{\delta},&\T_{\gamma v}^{(2)}&=\frac{\beta_{\gamma k}-\beta_{\gamma v}-1}{2}\B_{\gamma},&\T_{\gamma v}^{(3)}&=\frac{\beta_{\gamma k}-\beta_{\gamma v}+1}{2}\B_{\gamma}.
		\end{aligned}
	\end{equation*}
\end{proposition}
\begin{proof}
	The proofs for consistency, conservation, and adjoint consistency are straightforward. Moreover, the stiffness matrix arising from the discretization is symmetric since $ \T_{\gamma k}^{(3)}-\T_{\gamma k}^{(2)} = \B_{\gamma} $. We see that the energy stability conditions in \cref{eq:TD stablility,eq:T4 stablility} are met. It remains to show that the coefficients satisfy the energy stability requirements in \cref{eq:T1 stablility 1,eq:T1 stablility 2}. Note that if either or both $ \beta_{\gamma v}=1 $ and $ \beta_{\epsilon g}=1 $ then $ \T_{\gamma\epsilon k}^{(5)}=\T_{\epsilon\gamma k}^{(5)}=\bm{0} $, and the scheme is stable since \cref{eq:T1 stablility 1} is satisfied, \ie,
	\begin{equation} \label{eq:LDG stability with betak=0}
		\violet{\T_{\gamma k}^{(1)}-2\left(\frac{1}{\alpha_{\gamma k}}\T_{\gamma k}^{(2)}\Upsilon_{\gamma\gamma k}\T_{\gamma k}^{(2)}+\frac{1}{\alpha_{\gamma v}}\T_{\gamma v}^{(2)}\Upsilon_{\gamma\gamma v}\T_{\gamma v}^{(2)}\right)
		=\frac{2}{\alpha_{\gamma v}}\B_{\gamma}\Upsilon_{\gamma\gamma v}\B_{\gamma}-\frac{2}{\alpha_{\gamma v}}\B_{\gamma}\Upsilon_{\gamma\gamma v}\B_{\gamma}= \bm{0}}.
	\end{equation}
	 Thus, we only need to consider the case where $ \beta_{\gamma k}=\beta_{\epsilon k}=1 $, which gives $ \T_{\gamma\epsilon k}^{(5)}=(1/4)\B_{\gamma}\Upsilon_{\gamma\epsilon k}\B_{\epsilon} $, $ \T_{\epsilon\gamma k}^{(5)}=(1/4)\B_{\epsilon}\Upsilon_{\epsilon\gamma k}\B_{\gamma} $, and $ \T_{\epsilon_{1}k}^{(1)}=(2/\alpha_{\epsilon_{1}k})\B_{\epsilon_{1}}\Upsilon_{\epsilon_{1}\epsilon_{1}k}\B_{\epsilon_{1}} $. Hence, we have
	 \begin{equation}
	 	\T_{\gamma k}^{(1)}-2\left(\frac{1}{\alpha_{\gamma k}}\T_{\gamma k}^{(2)}\Upsilon_{\gamma\gamma k}\T_{\gamma k}^{(2)}+\frac{1}{\alpha_{\gamma v}}\T_{\gamma v}^{(2)}\Upsilon_{\gamma\gamma v}\T_{\gamma v}^{(2)}\right)
	 	=\frac{2}{\alpha_{\gamma k}}\B_{\gamma}\Upsilon_{\gamma\gamma k}\B_{\gamma}-\frac{2}{\alpha_{\gamma k}}\B_{\gamma}\Upsilon_{\gamma\gamma k}\B_{\gamma}=\bm{0}.
	 \end{equation}
	 Furthermore, we find 
	 \begin{equation}
	 	\T_{\gamma\epsilon_{1}k}^{(5)}\left(\T_{\epsilon_{1}k}^{(1)}\right)^{-1}\T_{\epsilon_{1}\gamma k}^{(5)}=\frac{1}{16}\B_{\gamma}\Upsilon_{\gamma\epsilon_{1}k}\B_{\epsilon_{1}}\left[\frac{2}{\alpha_{\epsilon_{1}k}}\B_{\epsilon_{1}}\Upsilon_{\epsilon_{1}\epsilon_{1}k}\B_{\epsilon_{1}}\right]^{-1}\B_{\epsilon_{1}}\Upsilon_{\epsilon_{1}\gamma k}\B_{\gamma}=\frac{1}{32}\B_{\gamma}\Upsilon_{\gamma\epsilon_{1}k}\left[\frac{1}{\alpha_{\epsilon_{1}k}}\Upsilon_{\epsilon_{1}\epsilon_{1}k}\right]^{-1}\Upsilon_{\epsilon_{1}\gamma k}\B_{\gamma},
	 \end{equation}
	 but, as in \cref{eq:Lemma 3 applied}, application of \cref{lem:I-YXY} yields $ \frac{1}{\alpha_{\gamma k}}\Upsilon_{\gamma\gamma k}-\Upsilon_{\gamma\epsilon_{1}k}[\frac{1}{\alpha_{\epsilon_{1}k}}\Upsilon_{\epsilon_{1}\epsilon_{1}k}]^{-1}\Upsilon_{\epsilon_{1}\gamma k}\succeq0 $, which implies that 
	 \begin{equation}
	 	\T_{\gamma\epsilon_{1}k}^{(5)}\left(\T_{\epsilon_{1}k}^{(1)}\right)^{-1}\T_{\epsilon_{1}\gamma k}^{(5)} \preceq \frac{1}{32}\B_{\gamma}\left[\frac{1}{\alpha_{\gamma k}}\Upsilon_{\gamma\gamma k}\right]\B_{\gamma}.
	 \end{equation}
	 The last condition, \cref{eq:T1 stablility 2}, that we need to show for energy stability follows as
	 \begin{equation}
	 	\T_{\gamma k}^{(1)}-64\T_{\gamma\epsilon_{1}k}^{(5)}\left(\T_{\epsilon_{1}k}^{(1)}\right)^{-1}\T_{\epsilon_{1}\gamma k}^{(5)}\succeq\T_{\gamma k}^{(1)}-2\B_{\gamma}\left[\frac{1}{\alpha_{\gamma k}}\Upsilon_{\gamma\gamma k}\right]\B_{\gamma}=\bm{0}.
	 \end{equation}
	 Therefore, the LDG SAT coefficients in \cref{prop:LDG SAT} lead to an energy stable SBP-SAT discretization.
\end{proof}

\subsection{CDG SAT: The compact discontinuous Galerkin method}
The CDG method \cite{peraire2008compact} has similar numerical fluxes as the LDG scheme, but uses local instead of global lifting operators. More precisely, the numerical fluxes for the CDG method are $ \fnc{\widehat{U}} = \avg{\fnc{U}_h} - \bm{\beta}\cdot\jump{\fnc{U}_h}$ and $ \widehat{\vecfnc{W}} = \avg{\vecfnc{W}_h^{\gamma}} + \bm{\beta}\jump{\vecfnc{W}_h^{\gamma}} - \mu h^{-1}\jump{\fnc{U}_h} $ on interior facets, and $ \widehat{\vecfnc{W}} = \vecfnc{W}_h^{\gamma} - \mu h^{-1}(\fnc{U}_h - \fnc{U}_D)\bm{n} $ on $ \Gamma^D $ \cite{peraire2008compact}. For this reason, the residual of the CDG method can be obtained from \cref{eq:LDG residual} by replacing the global lifting operators by the corresponding local lifting operators. The implication of this on the SBP-SAT discretization is the nullification of SAT coefficients that lead to extended stencils. 
\begin{proposition}\label{prop:CDG SAT}
	A consistent, conservative, adjoint consistent, and energy stable version of the CDG scheme with no mesh dependent stabilization parameter, \ie, $ \mu=0 $, has the SAT coefficients 
	\begin{equation*}
		\begin{aligned}
			\T_{\gamma k}^{(1)}&=\T_{\gamma v}^{(1)}=\frac{1}{2}\B_{\gamma}\bigg[\frac{\beta_{\gamma k}-\beta_{\gamma v}+1}{\alpha_{\gamma k}}\Upsilon_{\gamma\gamma k}+\frac{\beta_{\gamma v}-\beta_{\gamma k}+1}{\alpha_{\gamma v}}\Upsilon_{\gamma\gamma v}\bigg]\B_{\gamma},&\T_{\gamma k}^{(3)}&=\frac{\beta_{\gamma v}-\beta_{\gamma k}+1}{2}\B_{\gamma},&\T_{\gamma k}^{(4)}&=\T_{\gamma v}^{(4)}=\bm{0},\\\T_{\gamma k}^{(2)}&=\frac{\beta_{\gamma v}-\beta_{\gamma k}-1}{2}\B_{\gamma},&\T_{\gamma v}^{(3)}&=\frac{\beta_{\gamma k}-\beta_{\gamma v}+1}{2}\B_{\gamma},&\T_{\gamma\epsilon k}^{(5)}&=\T_{\gamma\epsilon k}^{(6)}=\bm{0},\\\T_{\gamma v}^{(2)}&=\frac{\beta_{\gamma k}-\beta_{\gamma v}-1}{2}\B_{\gamma},&\T_{\gamma}^{(D)}&=\frac{1}{\alpha_{\gamma k}}\B_{\gamma}\Upsilon_{\gamma\gamma k}\B_{\gamma},&\T_{\gamma\delta v}^{(6)}&=\T_{\gamma\delta v}^{(5)}=\bm{0}.
		\end{aligned}
	\end{equation*}
\end{proposition}
\begin{proof}
	The proofs for consistency, conservation, and adjoint consistency are straightforward. Energy stability follows if we can show that \cref{eq:T1 stablility 1} holds (note that \cref{eq:TD stablility} is satisfied). If $ \beta_{\gamma k} = 1 $ then $ \beta_{\gamma v}=0 $, and
	\begin{equation}
		\T_{\gamma k}^{(1)}-\left(\frac{1}{\alpha_{\gamma k}}\T_{\gamma k}^{(2)}\Upsilon_{\gamma\gamma k}\T_{\gamma k}^{(2)}+\frac{1}{\alpha_{\gamma v}}\T_{\gamma v}^{(2)}\Upsilon_{\gamma\gamma v}\T_{\gamma v}^{(2)}\right)=\frac{1}{\alpha_{\gamma k}}\B_{\gamma}\Upsilon_{\gamma\gamma k}\B_{\gamma}-\frac{1}{\alpha_{\gamma k}}\B_{\gamma}\Upsilon_{\gamma\gamma k}\B_{\gamma} = \bm{0}.
	\end{equation}
	Similarly, if $ \beta_{\gamma v} = 1$ then $ \beta_{\gamma k} = 0 $, and we have $ \T_{\gamma k}^{(1)}-((1/{\alpha_{\gamma k}})\T_{\gamma k}^{(2)}\Upsilon_{\gamma\gamma k}\T_{\gamma k}^{(2)}+(1/{\alpha_{\gamma v}})\T_{\gamma v}^{(2)}\Upsilon_{\gamma\gamma v}\T_{\gamma v}^{(2)}) = \bm{0} $. Hence, the stability condition in \cref{eq:T1 stablility 1} is satisfied.
\end{proof}
The SAT coefficients in \cref{prop:CDG SAT}, except for $ \T_{\gamma k}^{(1)} $ and $ \T_{\gamma}^{(D)} $, are found by discretizing the residual resulting from the CDG method and comparing the result with \cref{eq:Residual 1st form}. The coefficients $ \T_{\gamma k}^{(1)} $ and $ \T_{\gamma}^{(D)} $ for the original CDG method are the same as those stated in \cref{eq:LDG original coefficients}. Similar SAT coefficients for the CDG method are proposed in \cite{yan2020immersed}, where the stability issue with the original CDG method is discussed. In \cite{brdar2012compact}, numerical studies revealed that the original CDG method can be unstable for variable coefficient diffusion problems and for discretizations on quadrilateral grids. 
\begin{remark}
	As noted in \cite{peraire2008compact}, the LDG SAT and CDG SAT are identical in one space dimension. To see this, consider an arbitrary global vector, $ \bm{g} $, pointing to the right and two elements ordered from left to right, $ \Omega_k $ and $ \Omega_v $, respectively; then $ \beta_{\gamma k} = 1 $, $ \beta_{\gamma v} = 0 $, and $ \beta_{\epsilon k} =0 $. These values of the switches nullify all $ \T^{(5)} $ and $ \T^{(6)} $ SAT coefficients which cast the LDG SAT stable with $ (1/2) \T_{\gamma k}^{(1)} $ in \cref{prop:LDG SAT}, thus the CDG SAT and LDG SAT become identical. 
\end{remark}
\subsection{BO SAT: The Baumann-Oden method}
Unlike the schemes presented so far, the BO method \cite{baumann1999discontinuous} leads to neither a symmetric stiffness matrix nor an adjoint consistent discretization. The numerical fluxes for the BO method are $ \widehat{\fnc{U}}=\avg{\fnc{U}}+\bm{n}\cdot\jump{\fnc{U}_{h}} $, and $ \widehat{\vecfnc{W}}=\avg{\lambda\nabla\fnc{U}} $ \cite{arnold2002unified}, and the residual is given by
\begin{align}\label{eq:BO Residual}
		\fnc{R}\left(\fnc{U}_{h},\fnc{V}\right)&=-\int_{\Omega}\lambda\nabla\fnc{U}_{h}\cdot\nabla\fnc{V}\dd{\Omega}+\int_{\Omega}\fnc{V}\fnc{F}\dd{\Omega}-\int_{\Gamma^{I}}\jump{\fnc{U}_{h}}\cdot\left\{ \lambda\nabla\fnc{V}\right\} -\jump{\fnc{V}}\cdot\avg{\lambda\nabla\fnc{U}_{h}}\dd{\Gamma}-\int_{\Gamma^{D}}(\fnc{U}_{h}-\fnc{U}_{D})\lambda\nabla\fnc{V}\cdot\bm{n}\dd{\Gamma}
		\nonumber
		\\&\quad
		+\int_{\Gamma^{D}}\fnc{V}\left(\lambda\nabla\fnc{U}_{h}\right)\cdot\bm{n}\dd{\Gamma}+\int_{\Gamma^{D}}\fnc{V}\lambda\fnc{S}^{D}(\fnc{U}_{h}-\fnc{U}_{D})\cdot\bm{n}\dd{\Gamma}+\int_{\Gamma^{N}}\fnc{V}\fnc{U}_{N}\dd{\Gamma}. 
\end{align}
Discretization of \cref{eq:BO Residual} and comparison with \cref{eq:Residual 1st form} gives all of the SAT coefficients in \cref{prop:BO SAT} below, except for $ \T_{\gamma}^{(D)} $ which is modified for stability reasons. The coefficient $ 1/\alpha_{\gamma k} $ in $ \T_{\gamma}^{(D)} $ does not appear in the original BO method.
\begin{proposition} \label{prop:BO SAT}
	The BO method is reproduced if the SAT coefficients in \cref{eq:Residual 1st form} are chosen such that 
	\begin{equation*}
		\begin{aligned}
			\T_{\gamma k}^{(1)}&=\T_{\gamma v}^{(1)}=\bm{0},& \T_{\gamma k}^{(2)}&=\T_{\gamma v}^{(2)}=\violet{\frac{1}{2}\B_{\gamma}},&\T_{\gamma k}^{(4)}&=\T_{\gamma v}^{(4)}=\bm{0},\\\T_{\gamma}^{(D)}&=\frac{1}{\alpha_{\gamma k}}\B_{\gamma}\Upsilon_{\gamma\gamma k}\B_{\gamma},&\T_{\gamma k}^{(3)}&=\T_{\gamma v}^{(3)}=\frac{1}{2}\B_{\gamma},&\T_{\gamma\epsilon k}^{(5)}&=\T_{\gamma\epsilon k}^{(6)}=\T_{\gamma\delta v}^{(5)}=\T_{\gamma\delta v}^{(6)}=\bm{0},
		\end{aligned}
	\end{equation*}
	and the discretization arising from using these coefficients is consistent, conservative, and energy stable.
\end{proposition}
\begin{proof}
	It can easily be verified that the SAT coefficients satisfy the conditions for consistency and conservation. The proof for energy stability follows from \cref{thm:Stability Adjoint Inconsistent}.
\end{proof}

\subsection{CNG SAT: The Carpenter-Nordstr{\"o}m-Gottlieb method}
The CNG SAT \cite{carpenter1999stable} was introduced to resolve the multi-domain problem in high-order finite difference methods. Although it was originally presented for advection-diffusion problems, in this work we consider the CNG SAT coefficients that couple the diffusive terms only. The CNG SAT coefficients for multidimensional SBP operators are stated in \cref{prop:CNG SAT} below (see \cite{carpenter1999stable,carpenter2010revisiting,gong2011interface} for analogous coefficients in one-dimensional implementations).
\begin{proposition}\label{prop:CNG SAT}
	The CNG SAT leads to consistent, conservative, and energy stable discretization, and has the coefficients
	\begin{equation*}
		\begin{aligned}
			\T_{\gamma k}^{(1)}&=\T_{\gamma v}^{(1)}=\frac{1}{16}\B_{\gamma}\left[\frac{1}{\alpha_{\gamma k}}\Upsilon_{\gamma\gamma k}+\frac{1}{\alpha_{\gamma v}}\Upsilon_{\gamma\gamma v}\right]\B_{\gamma},&\T_{\gamma k}^{(3)}&=\T_{\gamma v}^{(3)}=\frac{1}{2}\B_{\gamma},&\T_{\gamma}^{(D)}&=\frac{1}{\alpha_{\gamma k}}\B_{\gamma}\Upsilon_{\gamma\gamma k}\B_{\gamma},\\\T_{\gamma\epsilon k}^{(5)}&=\T_{\gamma\epsilon k}^{(6)}=\T_{\gamma\delta v}^{(5)}=\T_{\gamma\delta v}^{(6)}=\bm{0},&\T_{\gamma k}^{(2)}&=\T_{\gamma v}^{(2)}=\bm{0},&\T_{\gamma k}^{(4)}&=\T_{\gamma v}^{(4)}=\bm{0}.
		\end{aligned}
	\end{equation*}	
\end{proposition}
\begin{proof}
	Substituting the SAT coefficients in \cref{eq:Residual 3rd form} and evaluating $ 2R_h(\bm{u}_h, \bm{u}_h) = R_h(\bm{u}_h, \bm{u}_h) +  R_h^T(\bm{u}_h, \bm{u}_h) $, stability follows if $ \A+\A^T \succeq 0$, where $ \A $ is given by \cref{eq:Matrix A for compact SATs} (with $ \T_{\gamma k}^{(2)}=\T_{\gamma v}^{(2)}=\bm{0} $). From \cref{thm:Positive semi-definiteness}, $ \A+\A^T $ is positive semidefinite if 
	\begin{equation}
		\begin{aligned}
			\left[\begin{array}{cc}
			2\T_{\gamma k}^{(1)} & -2\T_{\gamma k}^{(1)}\\
			-2\T_{\gamma k}^{(1)} & 2\T_{\gamma k}^{(1)}
			\end{array}\right]-\left[\begin{array}{cc}
			-\frac{1}{2}\B_{\gamma}\C_{\gamma k} & \frac{1}{2}\B_{\gamma}\C_{\gamma v}\\
			\frac{1}{2}\B_{\gamma}\C_{\gamma k} & -\frac{1}{2}\B_{\gamma}\C_{\gamma v}
			\end{array}\right]\left[\begin{array}{cc}
			2\alpha_{\gamma k}\Lambda_{k}^{*} & 0\\
			0 & 2\alpha_{\gamma v}\Lambda_{v}^{*}
			\end{array}\right]^{-1}\left[\begin{array}{cc}
			-\frac{1}{2}\C_{\gamma k}^{T}\B_{\gamma} & \frac{1}{2}\C_{\gamma k}^{T}\B_{\gamma}\\
			\frac{1}{2}\C_{\gamma v}^{T}\B_{\gamma} & -\frac{1}{2}\C_{\gamma v}^{T}\B_{\gamma}
			\end{array}\right]&\succeq0,
		\end{aligned}
	\end{equation}
	which, after simplification, gives 
	\begin{equation}
		\left[\begin{array}{cc}
		1 & -1\\
		-1 & 1
		\end{array}\right]\otimes\left(2\T_{\gamma k}^{(1)}-\B_{\gamma}\left[\frac{1}{8\alpha_{\gamma k}}\Upsilon_{\gamma\gamma k}+\frac{1}{8\alpha_{\gamma k}}\Upsilon_{\gamma\gamma v}\right]\B_{\gamma}\right)\succeq0.
	\end{equation}
	The stability constraint $ \T_{\gamma k}^{(1)}\succeq\B_{\gamma}[1/(16\alpha_{\gamma k})\Upsilon_{\gamma\gamma k}+1/(16\alpha_{\gamma v})\Upsilon_{\gamma\gamma v}]\B_{\gamma} $ is, therefore, satisfied by the proposed SAT coefficient. The stability constraint on $ \T_{\gamma}^{(D)} $ is the same as the other methods presented earlier. Finally, it follows from \cref{thm:Consistency,thm:Conservation} that the coefficients in \cref{prop:CNG SAT} lead to consistent and conservative discretizations.
\end{proof}

\cref{tab:SATs} summarizes the SAT coefficients corresponding to eight different methods. The coefficients for BR2 and SIPG were first presented in \cite{yan2018interior}, and the analysis therein shows that the methods lead to consistent, conservative, adjoint consistent, and energy stable discretizations. The nonsymmetric interior penalty Galerkin (NIPG) SAT is obtained by modifying $ \T_{\gamma k}^{(1)} $ and $ \T_{\gamma}^{(D)} $ in \cref{prop:BO SAT}, and it leads to consistent, conservative, and energy stable discretizations. \blue{For the implementations of the SATs in \cref{tab:SATs}, we used the facet weight parameter, $ \alpha_{\gamma k} $, provided in \cite{yan2018interior}, 
\begin{equation} \label{eq:alpha_gamma k}
	\alpha_{\gamma k}=\begin{cases}
	\frac{{\cal A}\left(\gamma\right)}{{\cal A}\left(\Gamma_{k}^{I}\right)+2{\cal A}\left(\Gamma_{k}^{D}\right)}, & \text{if }\gamma\in\Gamma^{I},\\
	\frac{2{\cal A}\left(\gamma\right)}{{\cal A}\left(\Gamma_{k}^{I}\right)+2{\cal A}\left(\Gamma_{k}^{D}\right)}, & \text{if }\gamma\in\Gamma^{D},
	\end{cases}
\end{equation}
where $ {\cal A}(\gamma) $ computes the length of facet $ \gamma $ in 2D (and area of $ \gamma $ in 3D).} \violet{Moreover, for the LDG and CDG SATs, we set the arbitrary global vector in \cref{eq:LDG switch with g}  as $ \bm{g}=[\frac{\pi}{2}, e]^T $}.

\begin{table*}[t!]
	\small
	\centering
	\begin{threeparttable}
	\caption{\label{tab:SATs} Interior facet SAT coefficients for diffusion problems. In the cases where only one entry is provided, the two coefficients in the column heading are equal, \eg, for the BR1 SAT, we have $ \T_{\gamma k}^{(2)} = \T_{\gamma v}^{(2)} = -\frac{1}{2}\B_{\gamma} $. All the SATs considered have $ \T_{\gamma k}^{(4)}=\T_{\gamma v}^{(4)}=\bm{0} $.}
	\setlength{\tabcolsep}{.5em}
	\renewcommand\cellgape{\Gape[\jot]}
		\begin{tabular}{l | c| c| c| c }
			\toprule 
			\makecell[l]{SAT} & \makecell[c]{$\T_{\gamma k}^{(1)}$, $\T_{\gamma v}^{(1)}$} & \makecell[c]{$\T_{\gamma k}^{(2)}$, $\T_{\gamma v}^{(2)}$} & \makecell[c]{$\T_{\gamma k}^{(3)}$, $\T_{\gamma v}^{(3)}$}& \makecell[c]{$\T_{\gamma\epsilon k}^{(5)}$, $\T_{\gamma\delta v}^{(6)}$}\\
			\midrule
			\makecell[l]{BR1} & \makecell[c]{$\frac{1}{2}\B_{\gamma}\left[\frac{1}{\alpha_{\gamma k}}\Upsilon_{\gamma\gamma k}+\frac{1}{\alpha_{\gamma v}}\Upsilon_{\gamma\gamma v}\right]\B_{\gamma}$} & \makecell[c]{$-\frac{1}{2}\B_{\gamma}$} & \makecell[c]{$\frac{1}{2}\B_{\gamma}$} & \makecell[c]{$\T_{\gamma\epsilon k}^{(5)}=\frac{1}{16}\B_{\gamma}\Upsilon_{\gamma\epsilon k}\B_{\epsilon}$\vspace{3pt}\\
			$\T_{\gamma\delta v}^{(6)}=\frac{1}{16}\B_{\gamma}\Upsilon_{\gamma\delta v}\B_{\delta}$}\\
			\midrule
			\makecell[l]{BR2} & \makecell[c]{$\frac{1}{4}\B_{\gamma}\left[\frac{1}{\alpha_{\gamma k}}\Upsilon_{\gamma\gamma k}+\frac{1}{\alpha_{\gamma v}}\Upsilon_{\gamma\gamma v}\right]\B_{\gamma}$} & \makecell[c]{$-\frac{1}{2}\B_{\gamma}$} & \makecell[c]{$\frac{1}{2}\B_{\gamma}$} & \makecell[c]{$\bm{0}$}\\
			\midrule
			\makecell[l]{SIPG} & \makecell[c]{$\frac{\left(\lambda_{\max}\right)_{k}\parallel\B_{\gamma}^{\frac{1}{2}}\R_{\gamma k}\H_{k}^{-\frac{1}{2}}\parallel_{2}^{2}}{4\alpha_{\gamma k}}+\frac{\left(\lambda_{\max}\right)_{v}\parallel\B_{\gamma}^{\frac{1}{2}}\R_{\gamma v}\H_{v}^{-\frac{1}{2}}\parallel_{2}^{2}}{4\alpha_{\gamma v}}\B_{\gamma}$} & \makecell[c]{$-\frac{1}{2}\B_{\gamma}$} & \makecell[c]{$\frac{1}{2}\B_{\gamma}$} & \makecell[c]{$\bm{0}$}\\
			\midrule
			\makecell[l]{LDG} & \makecell[c]{$\B_{\gamma}\left[\frac{\beta_{\gamma k}-\beta_{\gamma v}+1}{\alpha_{\gamma k}}\Upsilon_{\gamma\gamma k}+\frac{\beta_{\gamma v}-\beta_{\gamma k}+1}{\alpha_{\gamma v}}\Upsilon_{\gamma\gamma v}\right]\B_{\gamma}$} & \makecell[c]{$\T_{\gamma k}^{(2)}=\frac{\beta_{\gamma v}-\beta_{\gamma k}-1}{2}\B_{\gamma}$\vspace{3pt}\\
			$\T_{\gamma v}^{(2)}=\frac{\beta_{\gamma k}-\beta_{\gamma v}-1}{2}\B_{\gamma}$} & \makecell[c]{$\T_{\gamma k}^{(3)}=\frac{\beta_{\gamma v}-\beta_{\gamma k}+1}{2}\B_{\gamma}$ \vspace{3pt}\\$\T_{\gamma v}^{(3)}=\frac{\beta_{\gamma k}-\beta_{\gamma v}+1}{2}\B_{\gamma}$} & \makecell[c]{$\T_{\gamma\epsilon k}^{(5)}=\frac{\left[1+\beta_{\gamma k}-\beta_{\gamma v}\right]\left[1+\beta_{\epsilon k}-\beta_{\epsilon g}\right]}{16}\B_{\gamma}\Upsilon_{\gamma\epsilon k}\B_{\epsilon}$\vspace{3pt}\\
			$\T_{\gamma\delta v}^{(6)}=\frac{\left[\beta_{\gamma k}-\beta_{\gamma v}-1\right]\left[1+\beta_{\delta v}-\beta_{\delta q}\right]}{16}\B_{\gamma}\Upsilon_{\gamma\delta v}\B_{\delta}$}\\
			\midrule
			\makecell[l]{CDG} & \makecell[c]{$\frac{1}{2}\B_{\gamma}\left[\frac{\beta_{\gamma k}-\beta_{\gamma v}+1}{\alpha_{\gamma k}}\Upsilon_{\gamma\gamma k}+\frac{\beta_{\gamma v}-\beta_{\gamma k}+1}{\alpha_{\gamma v}}\Upsilon_{\gamma\gamma v}\right]\B_{\gamma}$} & \makecell[c]{$\T_{\gamma k}^{(2)}=\frac{\beta_{\gamma v}-\beta_{\gamma k}-1}{2}\B_{\gamma}$ \vspace{3pt}\\
			$\T_{\gamma v}^{(2)}=\frac{\beta_{\gamma k}-\beta_{\gamma v}-1}{2}\B_{\gamma}$} & \makecell[c]{$\T_{\gamma k}^{(3)}=\frac{\beta_{\gamma v}-\beta_{\gamma k}+1}{2}\B_{\gamma}$
			\vspace{3pt}\\ $\T_{\gamma v}^{(3)}=\frac{\beta_{\gamma k}-\beta_{\gamma v}+1}{2}\B_{\gamma}$} &\makecell[c]{$\bm{0}$}\\
			\midrule
			\makecell[l]{BO} & \makecell[c]{$\bm{0}$} & \makecell[c]{$\frac{1}{2}\B_{\gamma}$} & \makecell[c]{$\frac{1}{2}\B_{\gamma}$} & \makecell[c]{$\bm{0}$}\\
			\midrule
			\makecell[l]{NIPG} & \makecell[c]{$\frac{\left(\lambda_{\max}\right)_{k}\parallel\B_{\gamma}^{\frac{1}{2}}\R_{\gamma k}\H_{k}^{-\frac{1}{2}}\parallel_{2}^{2}}{4\alpha_{\gamma k}}+\frac{\left(\lambda_{\max}\right)_{v}\parallel\B_{\gamma}^{\frac{1}{2}}\R_{\gamma v}\H_{v}^{-\frac{1}{2}}\parallel_{2}^{2}}{4\alpha_{\gamma v}}\B_{\gamma}$} & \makecell[c]{$\frac{1}{2}\B_{\gamma}$} & \makecell[c]{$\frac{1}{2}\B_{\gamma}$} & \makecell[c]{$\bm{0}$}\\
			\midrule 
			\makecell[c]{CNG} & \makecell[c]{$\frac{1}{16}\B_{\gamma}\left[\frac{1}{\alpha_{\gamma k}}\Upsilon_{\gamma\gamma k}+\frac{1}{\alpha_{\gamma v}}\Upsilon_{\gamma\gamma v}\right]\B_{\gamma}$} & \makecell[c]{$\bm{0}$} & \makecell[c]{$\frac{1}{2}\B_{\gamma}$} & \makecell[c]{$\bm{0}$}\\
			\bottomrule 
		\end{tabular}
	\begin{tablenotes}
		 \item[Note:] The Dirichlet boundary SAT coefficient is given by $ \T_{\gamma}^{(D)}={1}/{\alpha_{\gamma k}}\B_{\gamma}\Upsilon_{\gamma\gamma k}\B_{\gamma} $ for all the SATs except for the SIPG and NIPG SATs for which $ \T_{\gamma}^{(D)}= {(\lambda_{\max})_{k}}/{\alpha_{\gamma k}}\parallel\B_{\gamma}^{\frac{1}{2}}\R_{\gamma k}\H_{k}^{-\frac{1}{2}}\parallel_{2}^{2}\B_{\gamma} $, where $ (\lambda_{\max})_{k} $ is the largest eigenvalue of $ \Lambda_{k} $.
	\end{tablenotes}
	\end{threeparttable}
\end{table*}

\section{Practical issues} \label{sec:Practial issues}
In addition to the properties presented in \cref{sec:Properties of SBP-SAT discretization}, one needs to consider a few practical issues when deciding which type of SAT to use. In this section, we investigate the relation between the SATs when applied with SBP diagonal-E operators and quantify the sparsity of the system matrix arising from different SBP-SAT discretizations.

\subsection{Equivalence of SATs for diagonal-norm $ \R^{0} $ SBP operators} \label{sec:Equivalence of SATs}
\begin{figure}[t!]
	\centering
	\begin{subfigure}{0.3\textwidth}
		\centering
		\includegraphics[scale=0.25]{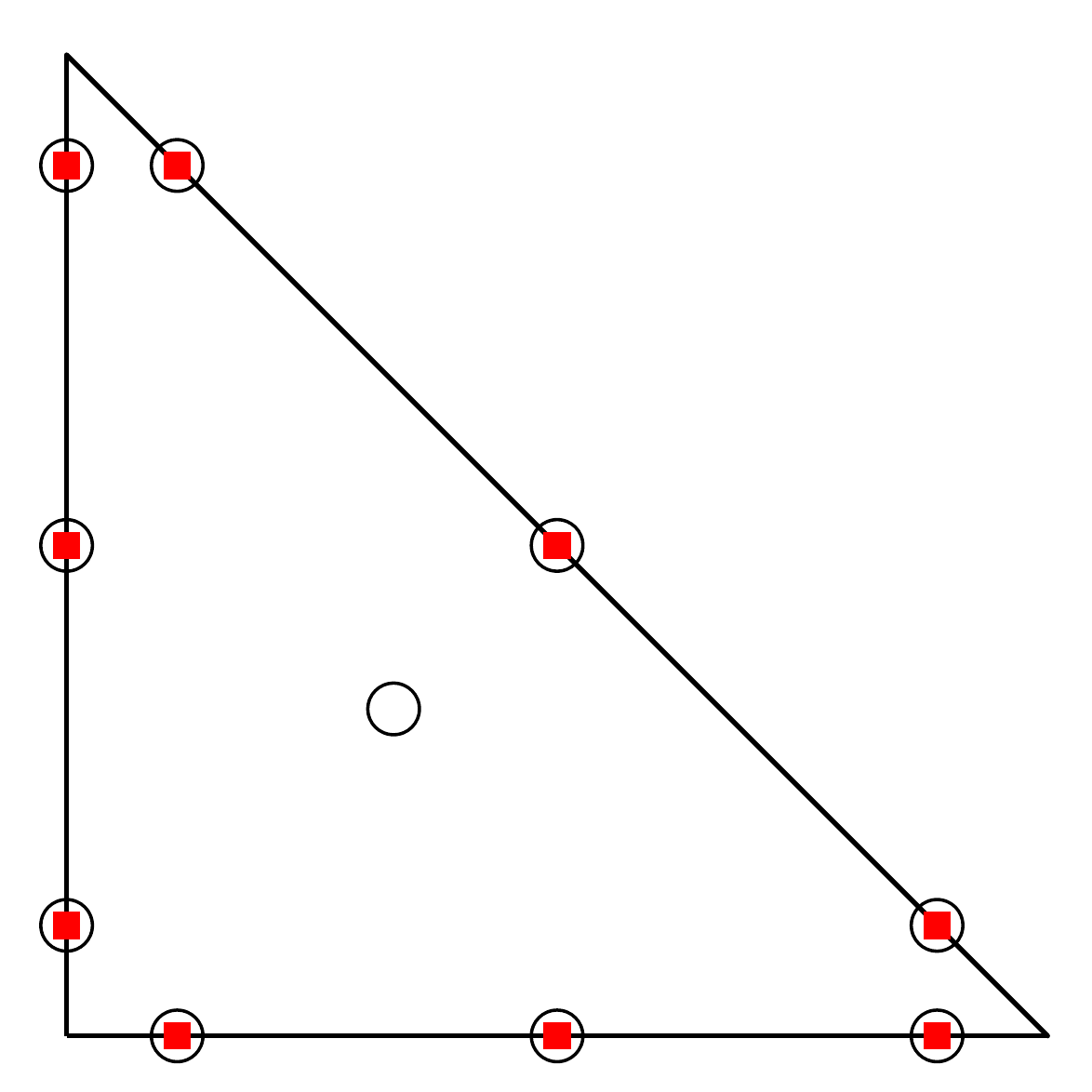}
		\caption{\label{fig:R0 family}$ \R^0 $ (SBP-E), $p=2$}
	\end{subfigure}\hfill
	\begin{subfigure}{0.3\textwidth}
		\centering
		\includegraphics[scale=0.25]{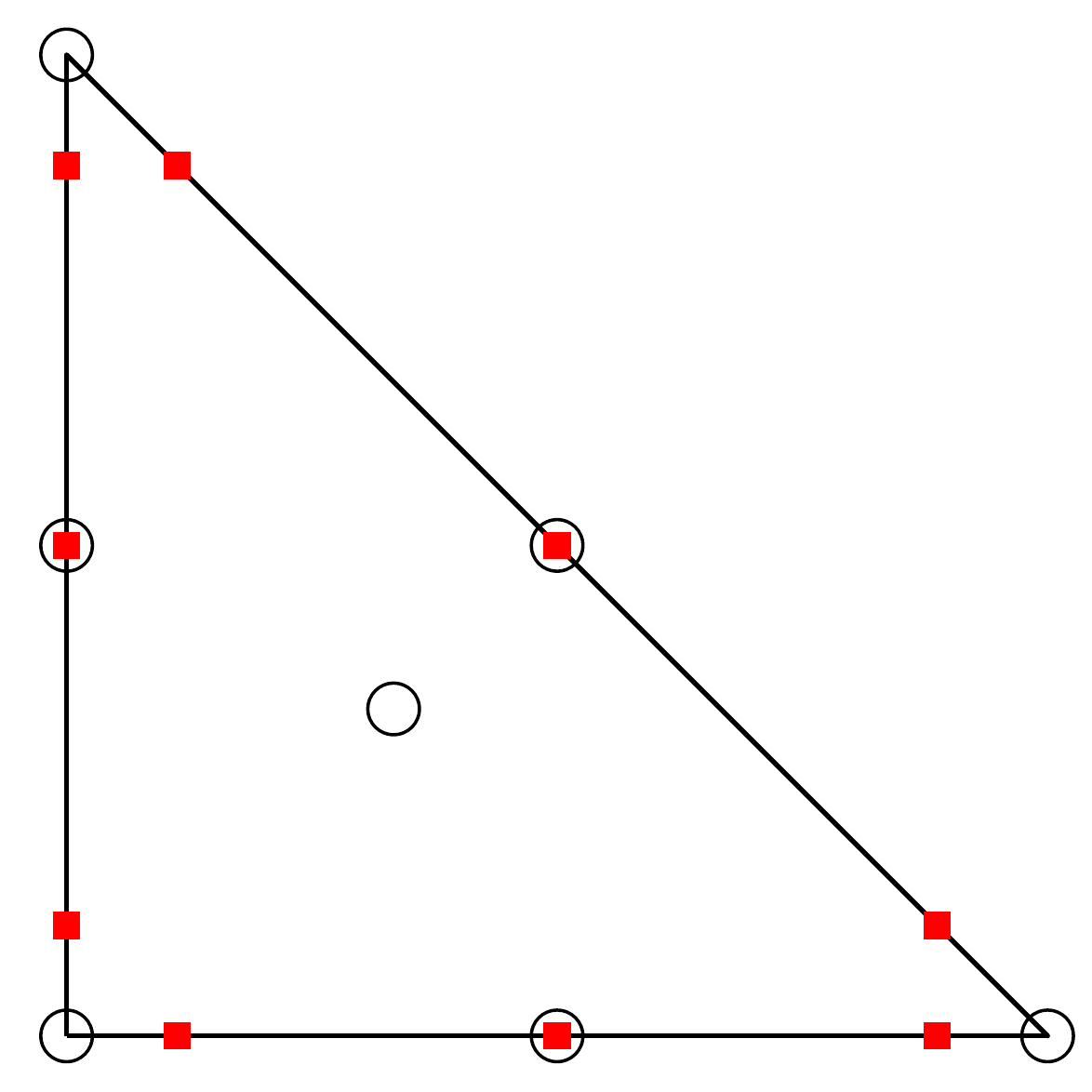}
		\caption{\label{fig:Rd1 family}$ \R^{d-1} $ (SBP-$ \Gamma $), $p=2$}
	\end{subfigure}\hfill
	\begin{subfigure}{0.3\textwidth}
		\centering
		\includegraphics[scale=0.25]{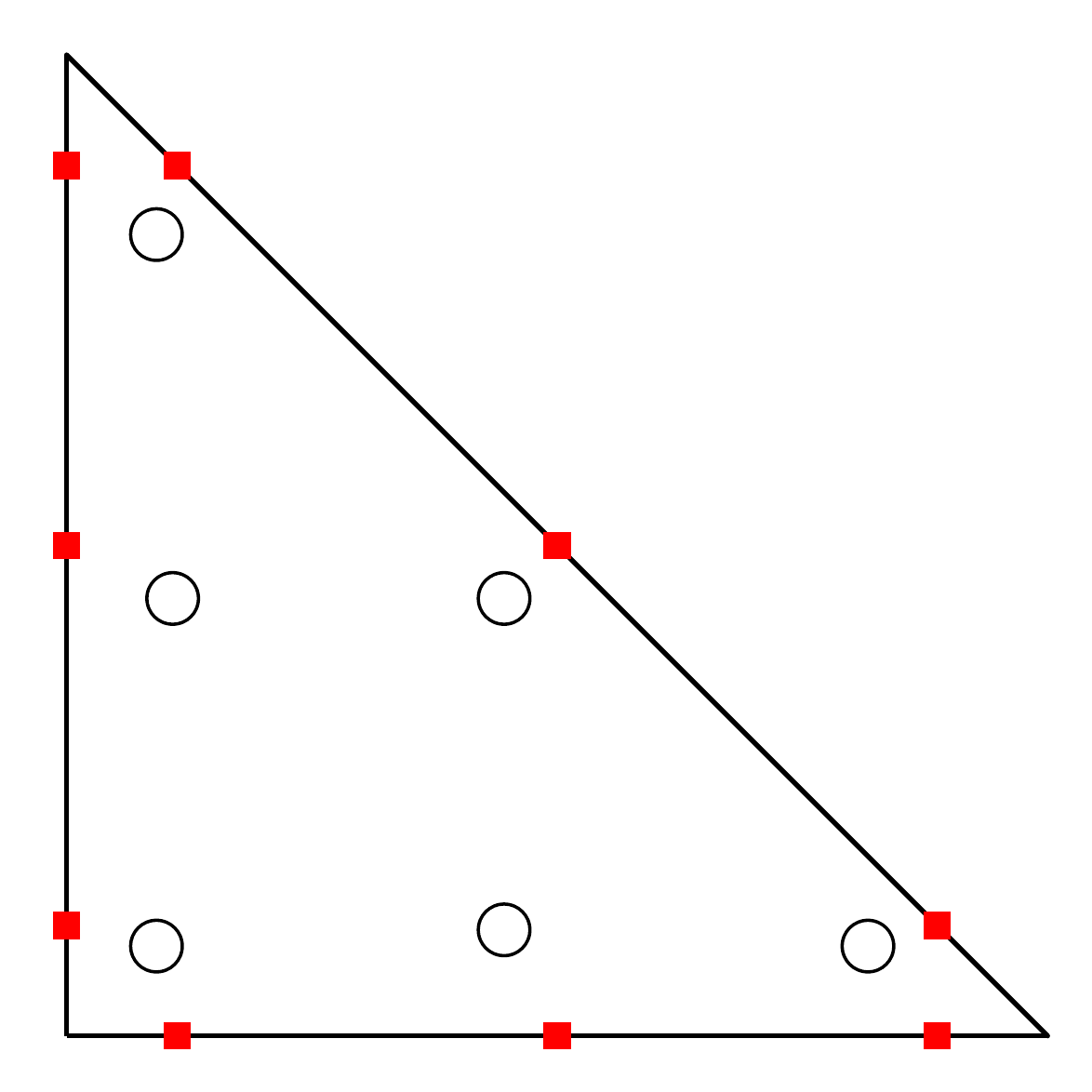}
		\caption{\label{fig:Rd family}$ \R^{d} $ (SBP-$ \Omega $), $p=2$}
	\end{subfigure}
	\caption{\label{fig:SBP families}Examples of degree two SBP operators in the $ \R^{0} $, $ \R^{d-1} $, and $ \R^{d} $ families on the reference triangle. The circles indicate locations of volume nodes, and the squares indicate locations of facet quadrature nodes.}
\end{figure}

Classification of SBP operators based on the dimensions spanned by the extrapolation matrix generalizes the SBP operator families introduced in \cite{fernandez2018simultaneous,chen2017entropy}. The SBP-$ \Omega $ operators \cite{fernandez2018simultaneous}, which fall under the $ \R^{d} $ SBP family, are characterized by having no volume nodes on element facets, \eg, the Legendre-Gauss (LG) operator. In general, however, the $ \R^{d} $ operator family allows volume nodes to be positioned on element facets as long as the $ \R $ matrix spans $ d $ dimensions \cite{marchildon2020optimization}. SBP-$ \Gamma $ \cite{hicken2016multidimensional,fernandez2018simultaneous} operators require $ {p+d-1 \choose d-1} $  volume nodes on each facet that are not collocated with facet quadrature nodes. In contrast, the $ \R^{d-1} $ family allows more volume nodes per facet; hence, it includes operators that cannot be categorized under the SBP-$ \Gamma $ family. SBP operators that have collocated volume and facet nodes on each facet, \eg, Legendre-Gauss-Lobatto (LGL), are classified under the SBP diagonal-E \cite{chen2017entropy} family which is equivalent to the $ \R^{0} $ SBP family. Examples of two-dimensional SBP operators in the $ \R^{0} $, $ \R^{d-1} $, and $ \R^{d} $ families are depicted in \cref{fig:SBP families}.  

For diagonal-norm $ \R^{0} $ SBP operators, the matrix $ \Upsilon_{abk} $, defined in \cref{eq:Upsilon definition}, exhibits an interesting property. In particular, it vanishes if $ a \neq b $, where $ a,b\in\{\gamma,\epsilon_{1},\epsilon_{2},\delta_{1},\delta_{2}\} $. Note that for this operator family, the extrapolation matrix simply picks out the volume nodes that are collocated with facet quadrature nodes, \ie,
\begin{equation}
	\left[\R_{f_{n}k}\right]_{ij}=\begin{cases}
	1 & \text{if}\;i+(n-1)n_f=j,\\
	0 & \text{if}\;i+(n-1)n_f\neq j,
	\end{cases}
\end{equation}
where $ n \in \{1,2,3\}$ is the facet number, $ i\in\{ 1,\dots,n_{f}\}  $, and $ j\in\{ 1,\dots,n_{p}\} $. Therefore, we have 
\begin{equation} \label{eq:equivalence due to Rgk}
	[{\R}_{\gamma k}{\H}_{k}^{-1}\Lambda_{xx,k}{\R}^T_{\epsilon_{1} k}]_{ij} = \sum_{m=1}^{n_p}\sum_{n=1}^{n_p} [{\R}_{\gamma k}]_{in} [{\H}_{k}^{-1}\Lambda_{xx,k}]_{nm} [{\R}^T_{\epsilon_{1} k}]_{mj} = \sum_{m=1}^{n_p} [{\R}_{\gamma k}]_{im} [{\H}_{k}^{-1}\Lambda_{xx,k}]_{mm} [{\R}_{\epsilon_{1} k}]_{jm} = 0,
\end{equation}
where the penultimate equality is a result of the fact that $ {\H}_{k}^{-1}\Lambda_{xx,k} $ is diagonal and $ [{\R}^T_{\epsilon_{1} k}]_{mj}=[{\R}_{\epsilon_{1} k}]_{jm} $. The last equality holds since $ {\R}_{\gamma k} $ and $ {\R}_{\epsilon_{1} k} $ do not contain $ 1 $ in the same column index, as $ \gamma $ and $ \epsilon_{1} $ have different facet numbers. This implies that for the $ \R^{0} $ operator family with diagonal norm matrix, $ \Upsilon_{\gamma \epsilon_{1} k}
=\N_{\gamma k}^{T}\bar{\R}_{\gamma k}\bar{\H}_{k}^{-1}\Lambda_{k}\bar{\R}_{\epsilon_{1} k}^{T}\N_{\epsilon_{1} k} = \bm{0} $, and more generally  
\begin{equation}
	\begin{aligned}
		\Upsilon_{abk}
		=\N_{ak}^{T}\bar{\R}_{a k}\bar{\H}_{k}^{-1}\Lambda_{k}\bar{\R}_{b k}^{T}\N_{bk}=\bm{0}\quad  \text{if}\; a\neq b, \; \text{where}\;  a,b\in\{\gamma,\epsilon_{1},\epsilon_{2},\delta_{1},\delta_{2}\}. 
	\end{aligned}
\end{equation}
Hence, the coefficients $ \T^{(5)} $ and $ \T^{(6)} $ for BR1 SAT in \cref{prop:BR1 SAT} and LDG SAT in \cref{prop:LDG SAT} vanish. We, therefore, have proven the following statement.
\begin{theorem}
	When implemented with diagonal-norm $ \R^{0} $ SBP operators, the BR1, BR2 and SIPG SATs are equivalent in the sense that they can be reproduced by considering $ \sigma_1 \T_{\gamma k}^{(1)} $ and $ \sigma_D \T_{\gamma}^{(D)} $ in \cref{prop:BR1 SAT} for $ \sigma_1,\sigma_D > 0 $. Similarly, the LDG and CDG SATs are equivalent for this family of operators.
\end{theorem}
The equivalence of the SIPG SAT and BR2 SAT is established in \cite{yan2018interior}. Similarly, it is shown in \cite{manzanero2018bassi} that the BR1 and SIPG methods are equivalent when the discretization is restricted to the LGL nodal points, and the LGL quadrature is used to approximate integrals. In the same paper, this property is exploited to find a sharper estimate of the minimum penalty coefficient for stability of the SIPG method. Gassner \etal \cite{gassner2018br1} reported that most drawbacks of the BR1 method are not observed when the Navier-Stokes equations are solved using DG discretization with LGL nodal points and quadrature. Since discretizations with LGL nodal point and quadrature satisfy the SBP property \cite{gassner2013skew,fernandez2014generalized,fernandez2014review} and the operator is in the diagonal-norm $ \R^{0} $ family, \violet{for one-dimensional implementations, the flux used} in \cite{gassner2018br1} can also be regarded as the BR2 SAT implemented with $ \T_{\gamma k}^{(1)}= (1/4)\B_{\gamma}[\Upsilon_{\gamma\gamma k}+\Upsilon_{\gamma\gamma v}]\B_{\gamma} $ and $ \T_{\gamma}^{(D)}=\B_{\gamma}\Upsilon_{\gamma \gamma k}\B_{\gamma} $  (or with the stabilization parameter for the BR2 flux set to $ \eta_0=1 $ in \cite{arnold2002unified}). \violet{For tensor-product implementations of the LGL operators in multiple dimensions, it can be shown, using the structure of the extrapolation matrices as in \cref{eq:equivalence due to Rgk}, that the BR1 SAT is not equivalent to the BR2 SAT. Despite this, when coupled with the BR1 SAT, the LGL operators lead to a smaller stencil width compared to operators that do not have nodes at the boundaries.}

diagonal-norm SBP operators in the $ \R^{0} $ family have also found important application in nonlinear stability analyses. They simplify entropy stability analyses \cite{chen2020review,shadpey2020entropy,crean2018entropy}, and are computationally less expensive than operators in the $ \R^{d} $ family on conforming grids \cite{chan2019efficient,fernandez2019staggered}. However, they exhibit lower solution accuracy and have a larger number of degrees of freedom compared to operators of the same degree in the $ \R^{d} $ SBP family \cite{chen2020review}.

\begin{remark}
	When implemented with the diagonal-norm $ \R^{0} $ SBP operators, the BR1 and LDG SATs are stable with the $ \T_{\gamma k}^{(1)} $ coefficients specified for the BR2 and CDG SATs in \cref{tab:SATs}, respectively. If such a modification is applied, then the BR1 and BR2 SATs as well as the LDG and CDG SATs become identical. We have not implemented this modification for the numerical results presented in \cref{sec:Numerical Results}.
\end{remark}

\subsection{Sparsity and storage requirements}
It is desirable to reduce the number of nonzero entries of the matrix resulting from a spatial discretization of \cref{eq:diffusion problem} to minimize storage requirements and take advantage of efficient sparse matrix algorithms for implicit time-marching methods. More generally, fewer nonzero entries lead to fewer floating point operations, and thus lower computational cost. The sparsity of a matrix is equal to one minus the density of the matrix, which is defined as the ratio of the number of nonzero entries to the total number of entries. 

The linear system of equations resulting from the SBP-SAT discretization on the RHS of \cref{eq:SBP-SAT discretization} is assembled in a global system matrix. This matrix is equivalent to the product of the inverse of the global mass matrix and the global stiffness matrix in the DG framework. An estimate of the number of nonzero entries of the system matrix depends on the type of SBP operator and SAT used. We first note that it has diagonal blocks of size $ n_p^2 $ associated with each element in the domain. Furthermore, for SBP-$ \Omega $ operators the $ \R $ matrix is dense since it spans $ d $ dimensions. Therefore, the number of nonzero entries of the off-diagonal block matrices containing terms such as $ \Rgk^T \T_{\gamma k}^{(1)}\Rgv $ are dense, \ie, they contain $ n_p^2 $ nonzero entries. Assuming simplices are used to tessellate the domain, each element has at most $ d+1 $ immediate neighbors. Thus, we can write an upper bound on the number of nonzero entries of the system matrix arising from the use of SBP-$ \Omega $ operators and any of the compact SATs as
\begin{equation}\label{eq:nnz omega compact}
nnz = n_e \left(n_p^2 + (d+1)n_p^2\right)= (d+2)n_p^2 n_e,
\end{equation}
where $ nnz $ denotes the number of nonzero entries. When SBP-$ \Omega $ operators are implemented with the BR1 SAT, each element is coupled with $ d^2 + 2d + 1 $ elements. Therefore, we have 
\begin{equation}\label{eq:nnz omega BR1}
nnz = n_e\left(n_p^2 + ( d^2 + 2d + 1) n_p^2\right) = ( d^2 + 2d + 2)n_p^2 n_e.
\end{equation}

For the LDG SAT, the number of elements coupled with a target element depends on the switch function. The choice of $ \bm{\beta} $ in \cref{eq:LDG switch} and \cref{eq:LDG switch with g} ensures that there is no element for which all switches point inwards or outwards simultaneously \cite{sherwin20062d}. Using this fact with the expressions for $ \T_{\gamma \epsilon k}^{(5)} $ and $ \T_{\gamma \delta v}^{(6)} $ in \cref{prop:LDG SAT}, it can be shown that the maximum number of elements coupled with a target element by the LDG SAT is $ d^2 + 1 $. Moreover, for every element coupled with $ d^2 + 1 $ neighbors there are $ d $ number of neighbors that will be coupled with less than $ d^2 + 1$ elements when $ d>1 $. Therefore, the number of elements that can have $ d^2 + 1 $ neighbors is limited to $ \lceil{n_e/(d+1)} \rceil$, where $ \lceil \cdot \rceil $ denotes the ceiling operator. Thus, an upper estimate of the number of nonzero entries of the system matrix resulting from the LDG SAT implemented with SBP-$ \Omega $ operator is given by
\begin{equation}\label{eq:nnz omega LDG}
nnz = \left\lceil\frac{n_e}{d+1}\right\rceil(d^2 + 2)n_p^2 + \left\lfloor\frac{n_e d}{d+1}\right\rfloor(d^2 + 1)n_p^2 = \left[\left\lceil\frac{n_e}{d+1}\right\rceil(d^2+2) + \left\lfloor\frac{n_e d}{d+1}\right\rfloor(d^2 + 1)\right]n_p^2,
\end{equation}
where $ \lfloor\cdot\rfloor $ denotes the floor operator. We used affine mapping (or straight-edged elements) to obtain \cref{eq:nnz omega LDG}; otherwise, the LDG SATs may result in more nonzero entries than the estimate in \cref{eq:nnz omega LDG} since the switch function $ \bm{\beta} $ varies along curved facets.

Since the $ \R $ matrix spans $ d-1 $ dimensions for SBP-$ \Gamma $ operators, it has $ n_f $ nonzero columns. Therefore, for implementations with SBP-$ \Gamma $ operators, blocks containing terms such as $ \Dgk^T \T_{\gamma k}^{(2)}\Rgk $ have $ n_f $ nonzero columns. Similarly, blocks containing terms such as $ \Rgk^T \T_{\gamma k}^{(3)}\Dgk $ have $ n_f $ nonzero rows. Thus, the sum $ \Dgk^T \T_{\gamma k}^{(2)}\Rgk +  \Rgk^T \T_{\gamma k}^{(3)}\Dgk$ has $2n_p n_f - n_f^2 $ nonzero entries. Identifying the structure of terms in blocks of the system matrix in a similar manner and using the number of coupled elements, we calculate upper estimates of the number of nonzero entries for different SBP-SAT discretizations of the Poisson problem. The estimates obtained are shown in \cref{tab:nnz}; similar results for DG implementation of the BR2 and CDG fluxes are presented in \cite{peraire2008compact}. In deriving the estimates, we assumed that all elements in the domain are interior; consequently, the number of nonzero entries is overestimated. This assumption, however, implies that the estimates in \cref{tab:nnz} get better with an increasing ratio of the number of interior to number of boundary elements in the domain.

\begin{table*}[!t]
	\small
	\caption{\label{tab:nnz}Estimates of the number of nonzero entries of system matrices resulting from different SBP-SAT discretizations of \cref{eq:Poisson problem}. For LDG and CDG SATs, straight-edged elements are used. If $ d=1 $, the estimates for the SBP-E family apply for operators that have volume nodes on their facets, and the estimates for LDG SAT should be replaced by those presented for CDG SAT. The operators $ \lceil\cdot\rceil $ and $ \lfloor\cdot\rfloor $ denote the ceiling and floor functions, respectively.}
	\centering
	\setlength{\tabcolsep}{1em}
	\renewcommand\cellgape{\Gape[\jot]}
	\begin{tabular}{l|l|l|l}
		\toprule
		\makecell[l]{SAT} & \makecell[c]{SBP-$\Omega$} & \makecell[c]{SBP-$\Gamma$} & \makecell[c]{SBP-E} \\
		\midrule
		\makecell[l]{BR1} & \makecell[l]{$(d^{2}+2d+2)n_{p}^{2}n_{e}$} & \makecell[l]{$\left[n_{p}^{2}+(d+1)(2n_{p}n_{f} - n_{f}^2)+(d^2+d)n_{f}^{2}\right]n_{e}$} & \makecell[l]{$\left[n_{p}^{2}+(d+1)(2n_{p}n_{f}-n_{f}^{2})\right]n_{e}$}\\
		\midrule
		\makecell[l]{BR2, SIPG, \\ BO, NIPG}
		& \makecell[l]{$\left(d+2\right)n_{p}^{2}n_{e}$} & \makecell[l]{$\left[n_{p}^{2}+(d+1)(2n_{p}n_{f}-n_{f}^{2})\right]n_{e}$} & \makecell[l]{$\left[n_{p}^{2}+(d+1)(2n_{p}n_{f}-n_{f}^{2})\right]n_{e}$}\\
		\midrule
		\makecell[l]{CDG, CNG} & \makecell[l]{$(d+2)n_{p}^{2}n_{e}$} & \makecell[l]{$\left[n_{p}^{2}+(d+1)n_{p}n_{f}\right]n_{e}$} & \makecell[l]{$\left[n_{p}^{2}+(d+1)n_{p}n_{f}\right]n_{e}$}\\
		\midrule
		\makecell[l]{LDG} & \makecell[l]{$\left\lceil\frac{n_e}{d+1}\right\rceil(d^2+2)n_p^2$ \vspace{2pt}\\$+ \left\lfloor\frac{n_e d}{d+1}\right\rfloor(d^2 + 1)n_p^2$} &
		\makecell[l]{$\left[n_{p}^{2}+(d+1)n_{p}n_{f}\right] n_e$ \vspace{2pt}\\$+\left[(d^2+1)-\left\lceil\frac{n_e}{d+1}\right\rceil(d+1) - \left\lfloor\frac{n_e d}{d+1}\right\rfloor(d+2)\right]n_f^2$} &
		\makecell[l]{$\left[n_{p}^{2}+(d+1)n_{p}n_{f}\right]n_{e}$}\\
		\bottomrule 
	\end{tabular}
\end{table*}

From \cref{tab:nnz} it can be deduced that the BR1 SAT yields the largest number of nonzero entries for a given type of SBP operator. In contrast, the CDG and CNG SATs give the smallest number of nonzero entries. While it is fairly easy to rank the SATs based on the number of nonzero entries they produce for a given type of operator, such a comparison involving different types of SBP operator is not straightforward due to varying number of volume nodes, $ n_p $. 

\section{Numerical results} \label{sec:Numerical Results}
To verify the theoretical analyses presented in the previous sections, we consider the two-dimensional Poisson problem
\begin{equation}\label{eq:Poisson problem}
	\begin{aligned}
		-\nabla\cdot\left(\lambda\nabla\fnc{U}\right)&=\fnc{F}&&\text{in}\;\Omega=\left[0,20\right]\times\left[-5,5\right],\\\bm{n}\cdot\left(\lambda\nabla\fnc{U}\right)&=\fnc{U}_{N}&&\text{on}\;\Gamma^{N}=\left\{ \left(x,y\right):y\in\left[-5,5\right],x=20\right\} ,\\\fnc{U}&=\fnc{U}_{D}&&\text{on}\;\Gamma^{D}=\Gamma\backslash\Gamma^{N}, 
	\end{aligned}
\end{equation}
where $ \lambda= \bigl[ \begin{smallmatrix} 4x+1 & y \\
y & y^{2} + 1 \end{smallmatrix} \bigr]$, and the source term and boundary conditions are determined via the method of manufactured solution, \ie, we choose the exact solution to be
\begin{equation}\label{eq:Exact solution}
	\fnc{U} = \sin(\frac{\pi}{8} x) \sin(\frac{\pi}{8} y),
\end{equation}
and evaluate $ \fnc{F} $, $ \fnc{U}_N $, $ \fnc{U}_D $ from \cref{eq:Poisson problem}. Similarly, we specify the exact adjoint solution as 
\begin{equation}\label{eq:Exact adjoint}
	\psi = x+y
\end{equation}
and evaluate the source term and boundary conditions associated with the adjoint problem from 
\begin{equation} \label{eq:Poisson adjoint problem}
	\begin{aligned}
		-\nabla\cdot\left(\lambda\nabla\fnc{\psi}\right)&=\fnc{G}&&\text{in}\;\Omega=\left[0,20\right]\times\left[-5,5\right],\\\bm{n}\cdot\left(\lambda\nabla\fnc{\psi}\right)&=\fnc{\psi}_{N}&&\text{on}\;\Gamma^{N}=\left\{ \left(x,y\right):y\in\left[-5,5\right],x=20\right\} ,\\\fnc{\psi}&=\fnc{\psi}_{D}&&\text{on}\;\Gamma^{D}=\Gamma\backslash\Gamma^{N}.
	\end{aligned}
\end{equation}
Finally, a linear functional of the form\footnote{Note that $ \psi_{D} $ and $ \psi_{N} $ are evaluated from $ \psi $, but usually there is no need to know the adjoint solution. Thus, the functional simply contains $ \fnc{U} $ as an unknown, and the values of $ \psi_{D} $ and $ \psi_{N} $ are given as coefficients or functions in the expression for the functional.} \cref{eq:Functional}, \ie,
\begin{equation*}
\fnc{I} (\fnc{U})= \int_{\Omega}\fnc{G}\fnc{U}\dd{\Omega} 
- \int_{\Gamma^D}\psi_{D} (\lambda\nabla\fnc{U})\cdot\bm{n}\dd{\Gamma}
+ \int_{\Gamma^N}\psi_{N}\fnc{U}\dd{\Gamma},
\end{equation*}
is considered. Since we know the primal solution, the adjoint, the boundary conditions, and the source terms, the linear functional can be evaluated exactly, and its value, accurate to fifteen significant figures, is \violet{$ \fnc{I}(\fnc{U}) = -27.0912595377575 $}. 

\begin{figure}[!t]
	\centering
	\begin{subfigure}{0.5\textwidth}
		\centering
		\includegraphics[scale=0.52]{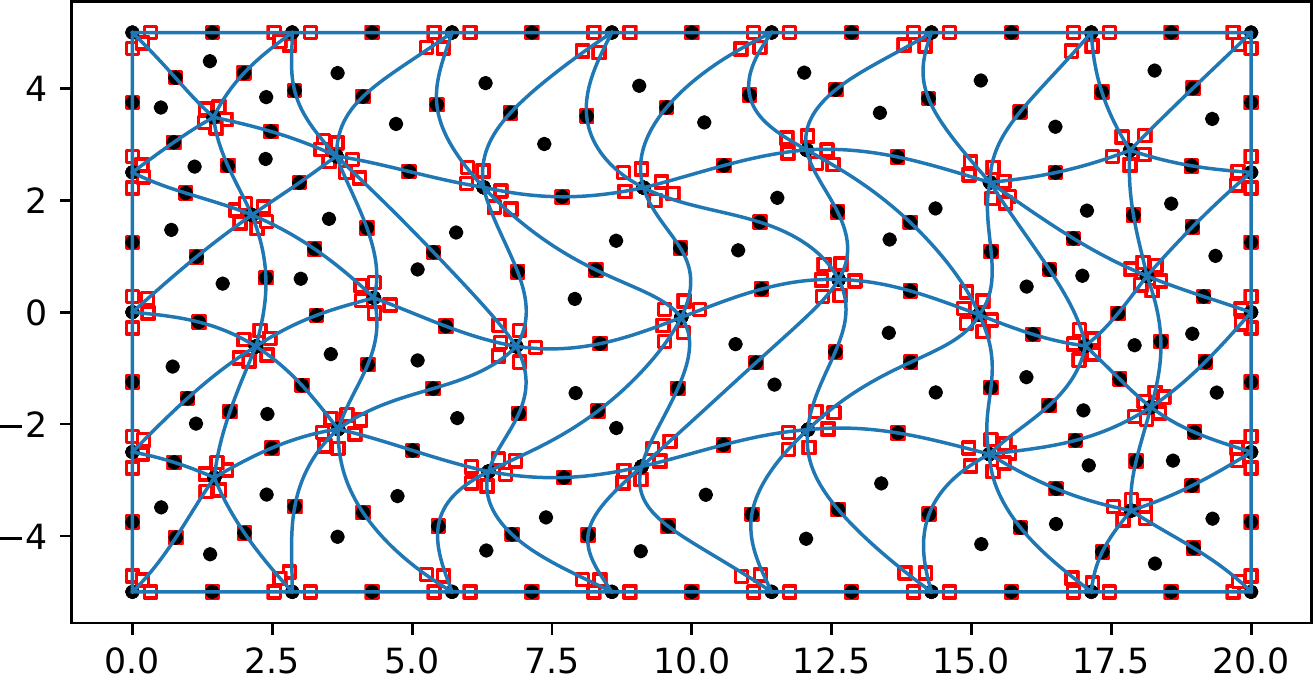}
		\caption{\label{fig:mesh Rd1 family} SBP-$ \Gamma $ operator, $p=2$, $p_{\rm map}=2$}
	\end{subfigure}\hfill
	\begin{subfigure}{0.5\textwidth}
		\centering
		\includegraphics[scale=0.52]{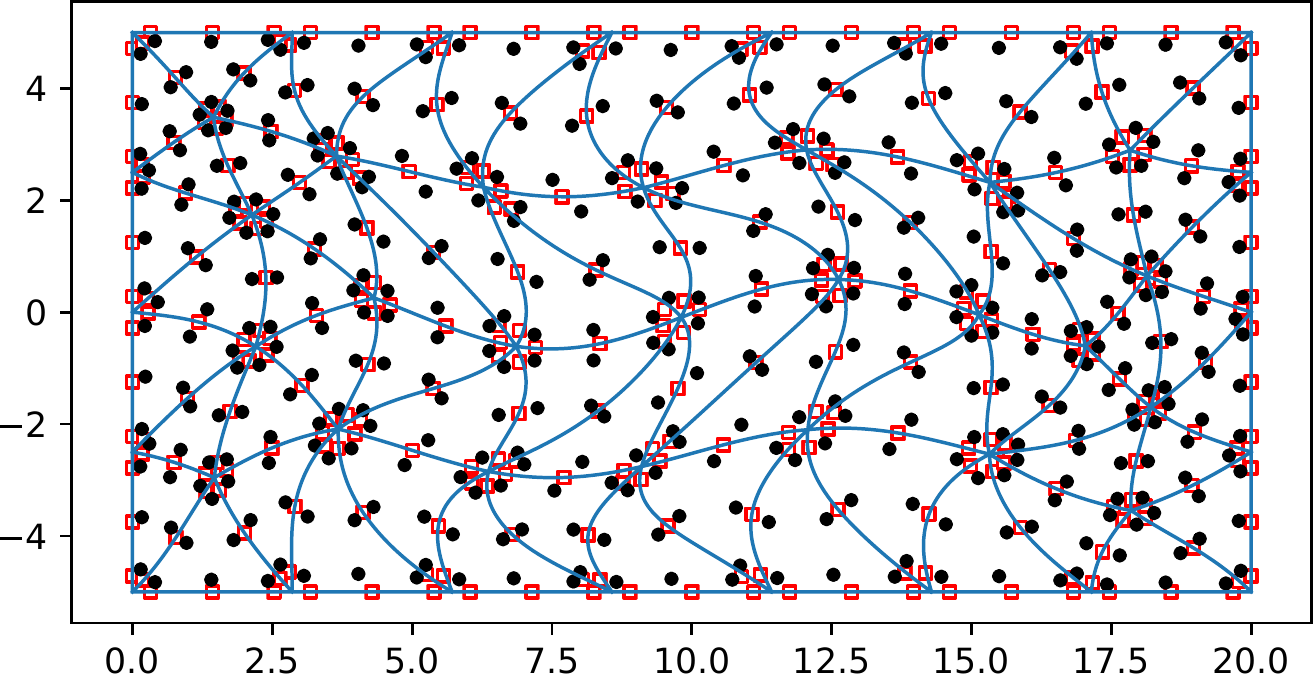}
		\caption{\label{fig:mesh Rd family} SBP-$ \Omega $ operator, $p=2$, $p_{\rm map}=2$}
	\end{subfigure}
	\caption{\label{fig:mesh SBP families} Physical domain tessellated with 68 curved triangular elements. The circles and squares indicate, respectively, the locations of volume nodes and facet quadrature points obtained with a degree two curvilinear mapping, and the lines define the facets of each element due to the perturbation given by \cref{eq:Mesh purturbation}.}
\end{figure}

The physical domain is tessellated with triangular elements, and the $ \alpha $-optimized Lagrange interpolation nodes on the reference element are mapped through an affine mapping to the physical elements. Then the triangular elements are curved by perturbing the coordinates of the $ \alpha$-optimized Lagrange interpolation nodes, $ \bm{x}_{L,k} $ and $ \bm{y}_{L,k} $, using the functions \cite{chan2019efficient} 
\begin{equation}\label{eq:Mesh purturbation}
	\begin{aligned}
		\tilde{\bm{x}}_{k}&=\bm{x}_{L,k}+\frac{5}{4}\cos\left(\frac{\pi}{20}\bm{x}_{L,k}-\frac{\pi}{2}\right)\cos\left(\frac{3\pi}{10}\bm{y}_{L,k}\right),& \quad \tilde{\bm{y}}_{k}&=\bm{y}_{L,k}+\frac{5}{8}\sin\left(\frac{\pi}{5}\tilde{\bm{x}}_{k}-2\pi\right)\cos\left(\frac{\pi}{10}\bm{y}_{L,k}\right).
	\end{aligned}
\end{equation}
The mesh Jacobian remain positive for each element under the curvilinear transformation. Examples of curvilinear grids with degree two SBP-$ \Gamma $ and SBP-$ \Omega $ operators are shown in \cref{fig:mesh SBP families}. A mapping degree of two is used for all numerical results presented. \blue{In all cases, the numerical solutions are obtained by solving the discrete equations using a direct method; specifically, the ``spsolve" function from the SciPy sparse linear algebra library in Python is used.}

\subsection{Accuracy} \label{sec:Accuracy}
The errors in the primal and adjoint solutions are computed, respectively, as 
\begin{equation} \label{eq:Error}
	\begin{aligned}
		\norm{\bm{u}_h - \bm{u}}_{\H} =\sqrt{ \sum_{\Omega_k \in \fnc{T}_h}\left(\bm{u}_{h,k} - \bm{u}_k\right)\H_k \left(\bm{u}_{h,k} - \bm{u}_k\right)},\quad && 
		\norm{\bm{\psi}_h - \bm{\psi}}_{\H} = \sqrt{\sum_{\Omega_k \in \fnc{T}_h}\left(\bm{\psi}_{h,k} - \bm{\psi}_k\right)\H_k \left(\bm{\psi}_{h,k} - \bm{\psi}_k\right)},
	\end{aligned}
\end{equation}
and the functional error is calculated as $ \abs{I_h(\bm{u}_h) - \fnc{I}(\fnc{U})} $. To study the accuracy and convergence properties of the primal solution, adjoint solution, and functional under mesh refinement, we consider four successively refined grids with 68, 272, 1088, 4352 elements. The nominal element size is calculated as $ h \equiv 20/\sqrt{n_e} $.

\begin{figure}[!t]
	\centering
	\begin{subfigure}{0.33\textwidth}
		\centering
		\includegraphics[scale=0.2]{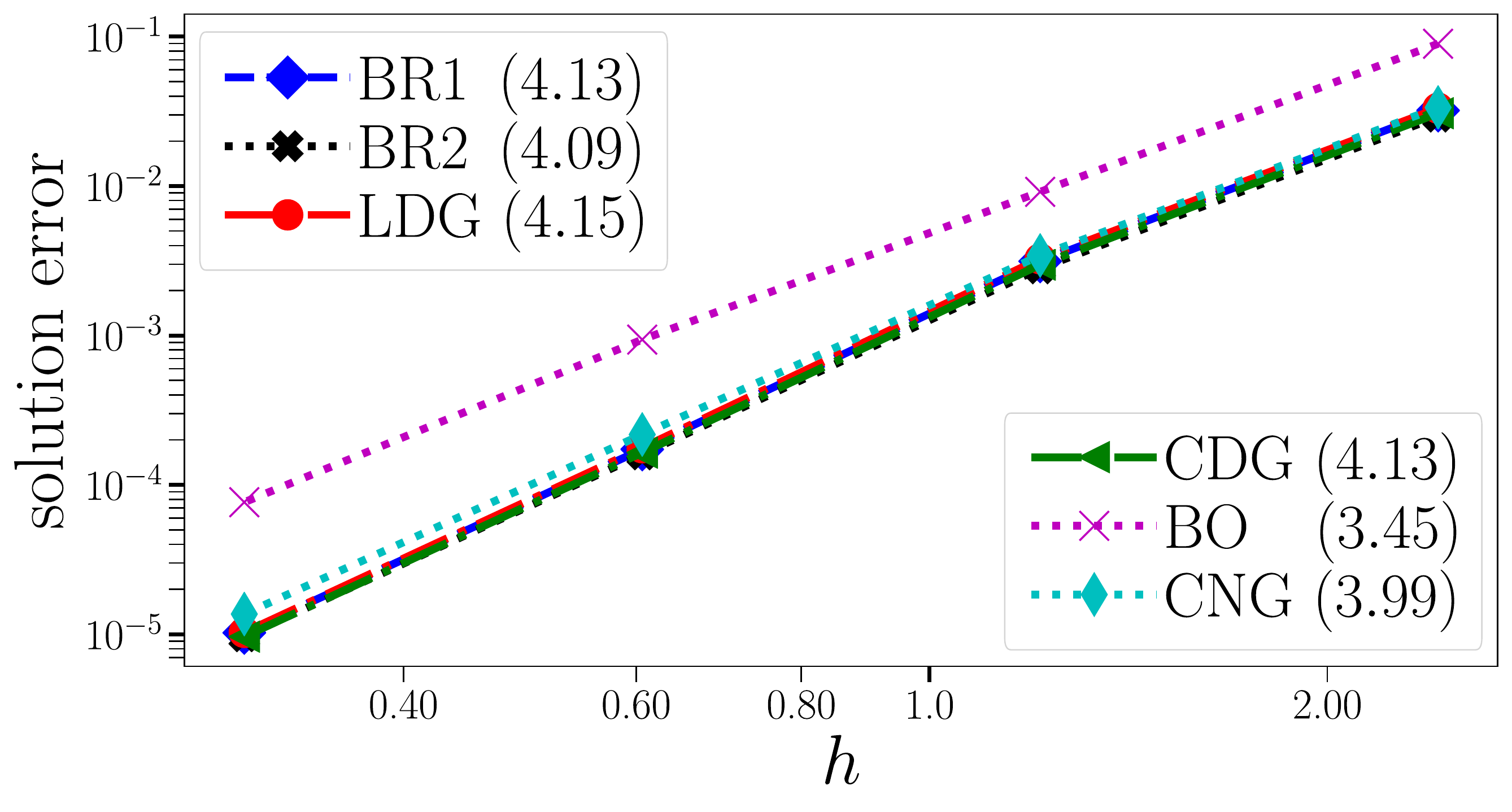}
		\caption{\label{fig:Omega soln p=3} SBP-$ \Omega $ operator, $p=3$}
	\end{subfigure}\hfill
	\begin{subfigure}{0.33\textwidth}
		\centering
		\includegraphics[scale=0.2]{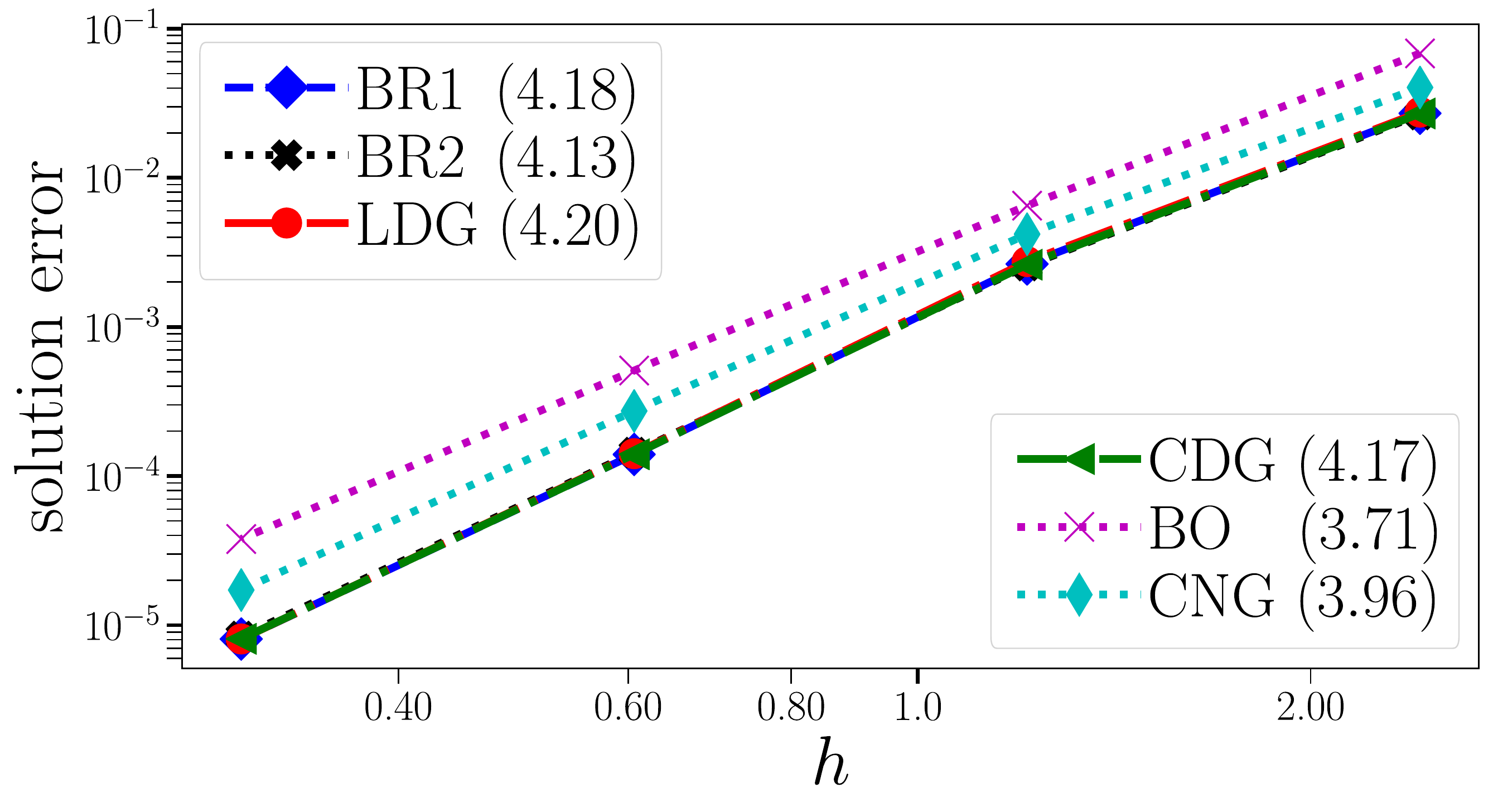}
		\caption{\label{fig:Gamma soln p=3} SBP-$ \Gamma $ operator, $p=3$}
	\end{subfigure} 
	\begin{subfigure}{0.33\textwidth}
		\centering
		\includegraphics[scale=0.2]{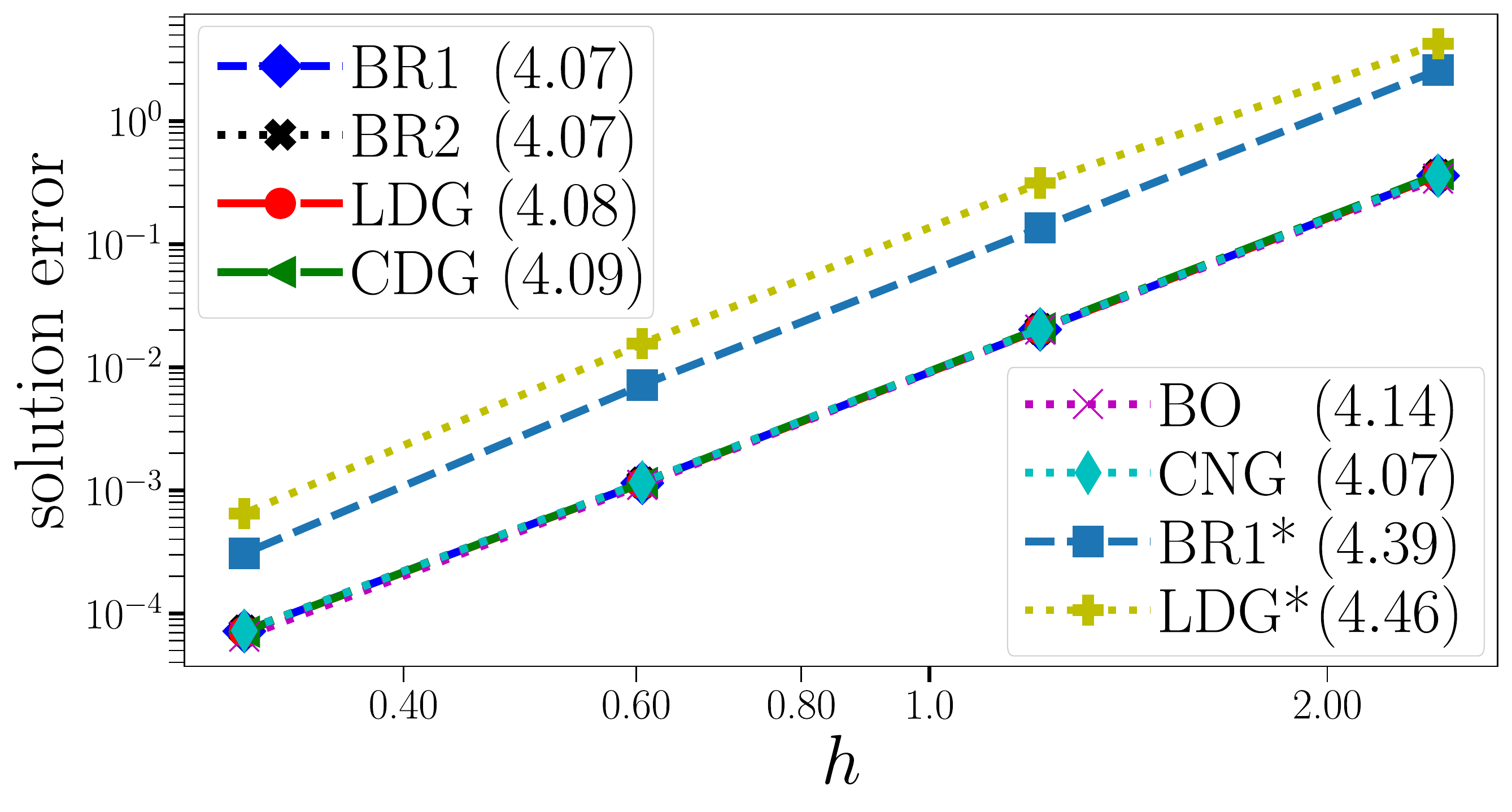}
		\caption{\label{fig:diage soln p=3} SBP-E operator, $p=3$}
	\end{subfigure}
	\\
	\begin{subfigure}{0.33\textwidth}
		\centering
		\includegraphics[scale=0.2]{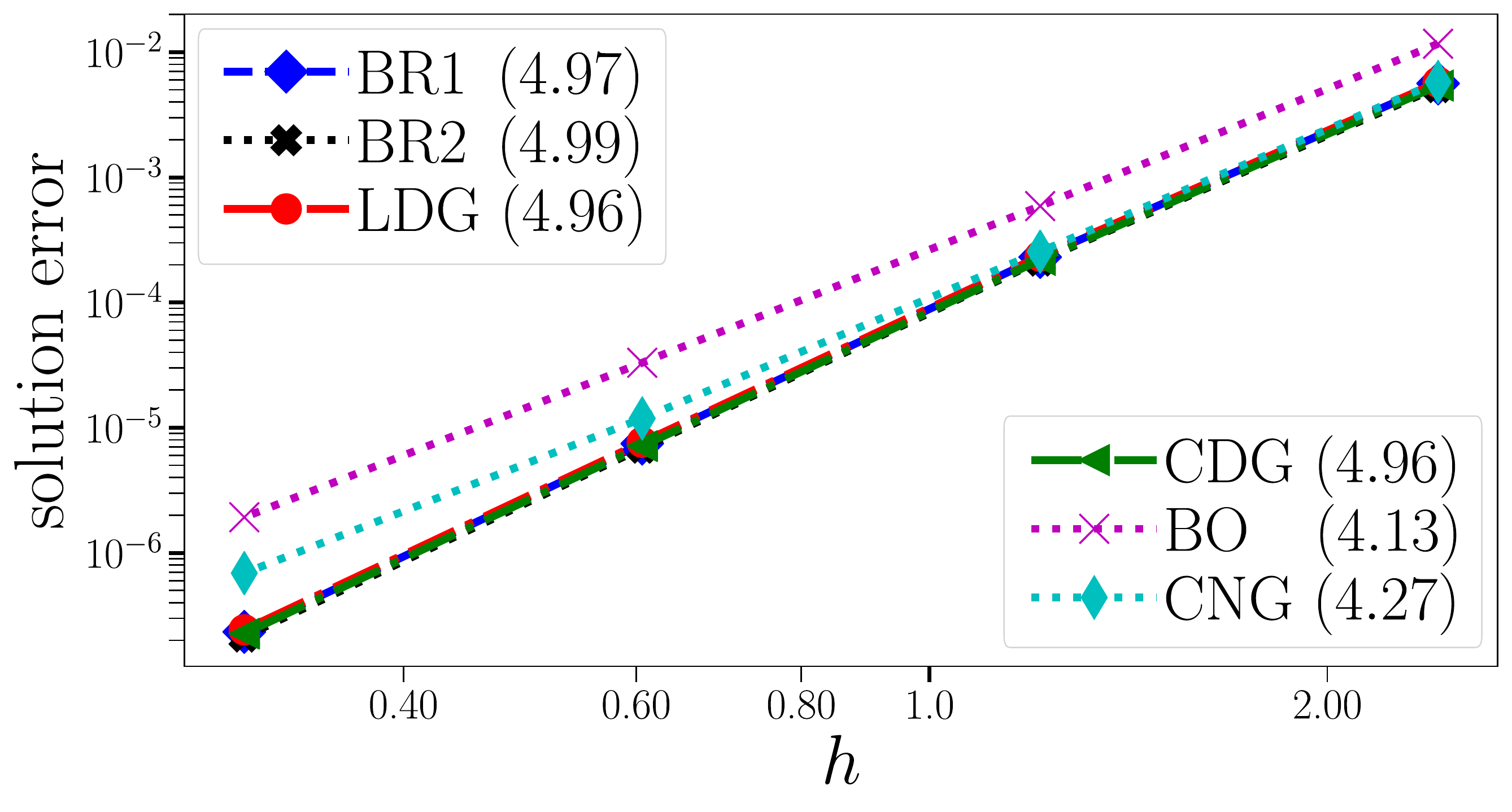}
		\caption{\label{fig:Omega soln p=4} SBP-$ \Omega $ operator, $p=4$}
	\end{subfigure}\hfill
	\begin{subfigure}{0.33\textwidth}
		\centering
		\includegraphics[scale=0.2]{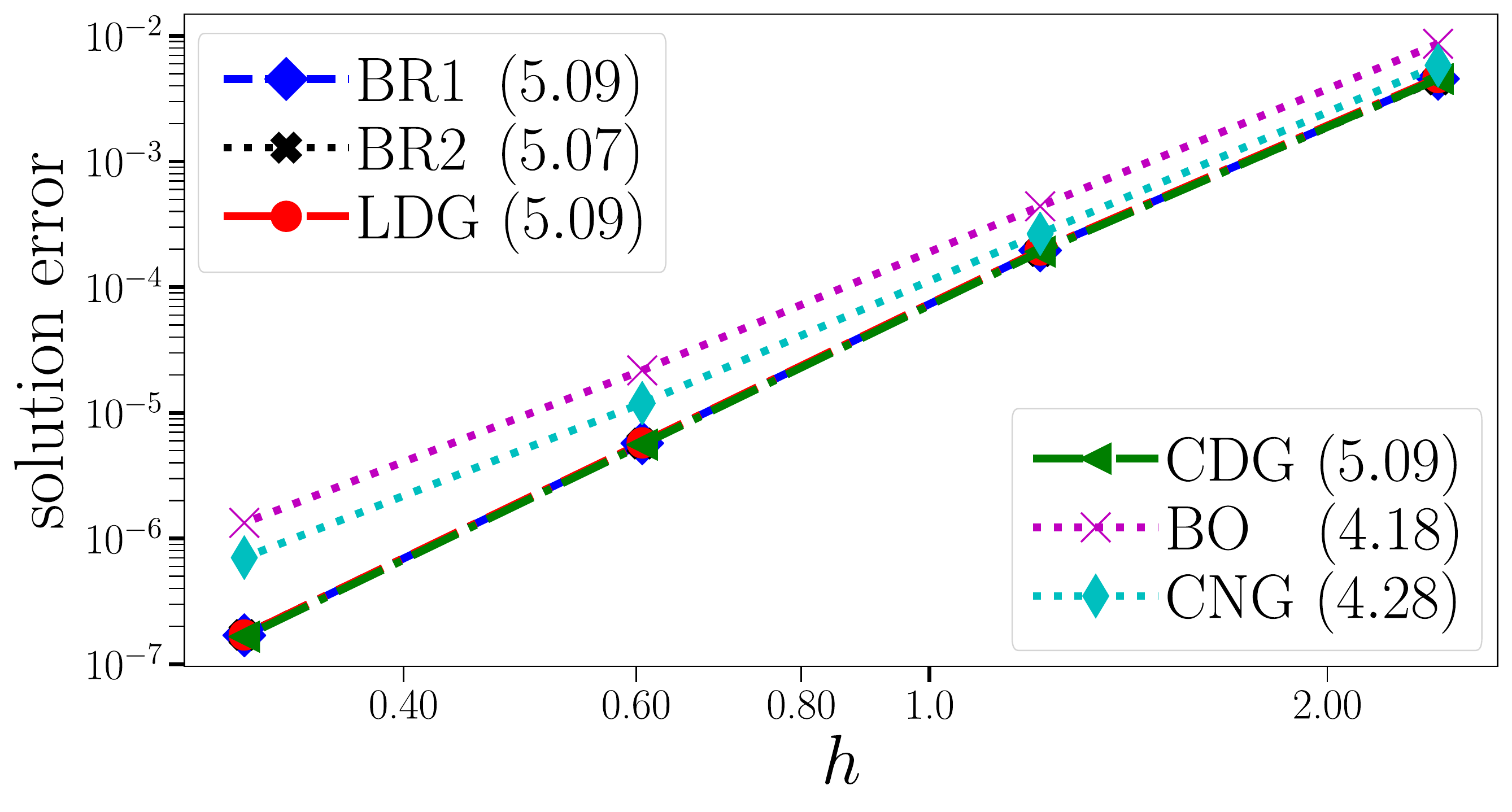}
		\caption{\label{fig:Gamma soln p=4} SBP-$ \Gamma $ operator, $p=4$}
	\end{subfigure} 
	\begin{subfigure}{0.33\textwidth}
		\centering
		\includegraphics[scale=0.2]{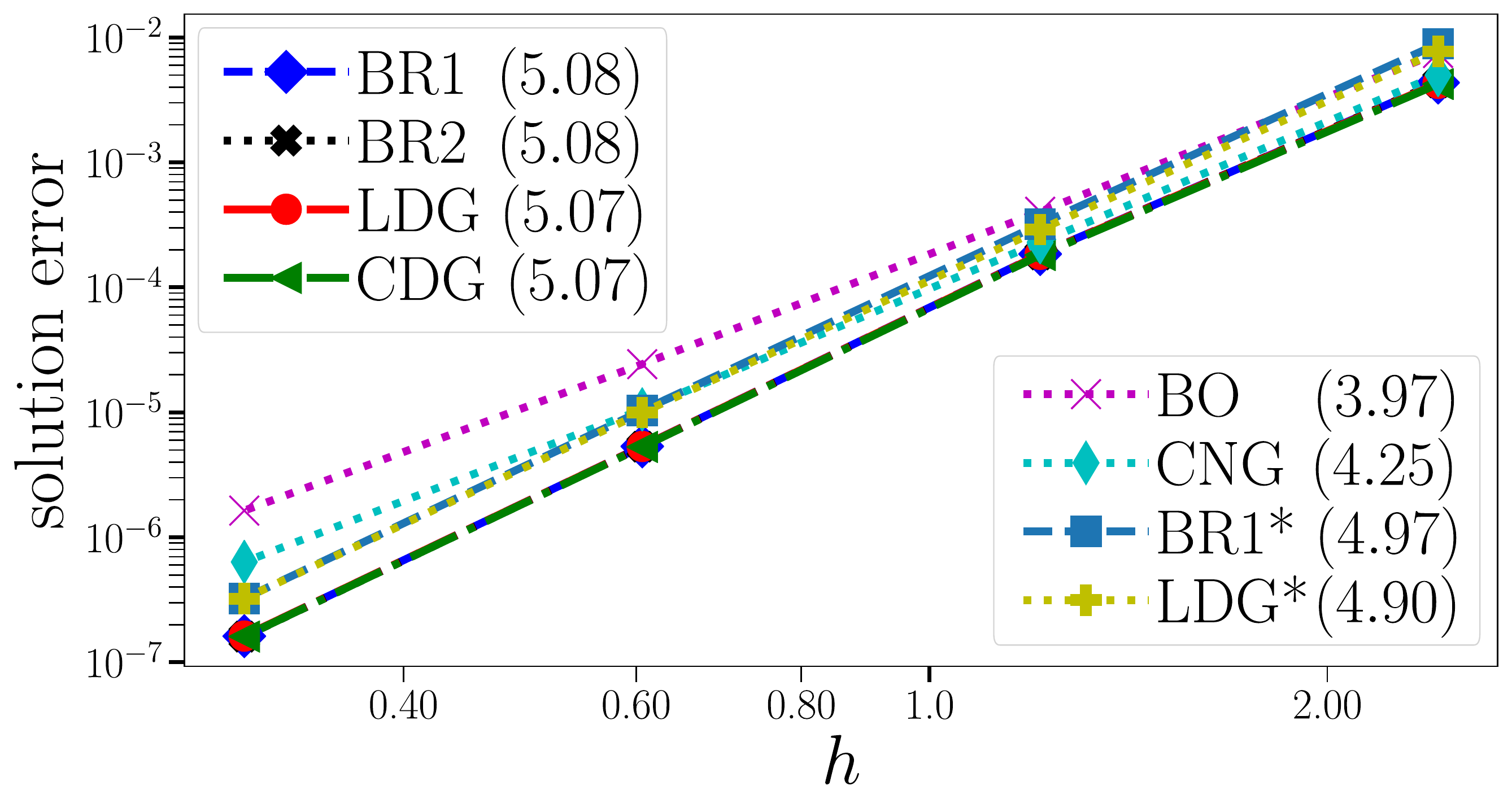}
		\caption{\label{fig:diage soln p=4} SBP-E operator, $p=4$}
	\end{subfigure}
	\caption{\label{fig:Solution error} Solution error under grid refinement. Solution convergence rates (shown in parenthesis) are calculated by fitting a line through the last three error values on the refined meshes. \violet{The BR1* and LDG* SATs represent the unmodified BR1 and LDG SATs, which are implemented with the SBP diagonal-E operators only}.}
\end{figure}

Figure \ref{fig:Solution error} shows the solution errors and convergence rates under mesh refinement for six types of SAT implemented with three types of SBP operators. Schemes with the BR1, BR2, LDG, and CDG SATs display solution convergence rates of $ p+1 $ and achieve very similar solution error values. In contrast, schemes with the BO and CNG SATs exhibit an even-odd convergence phenomenon; schemes with odd degree SBP operators converge at rates of $ p+1 $ while those with even degree operators converge at reduced rates of $ p $. The even-odd convergence property of the BO method is well-known, \eg, see \cite{shu2001different,kirby2005selecting,carpenter2010revisiting}. Furthermore, schemes with the BO SAT exhibit the largest solution error values in almost all cases considered (except for the case with degree three SBP diagonal-E operator). 

Numerical experiments in the literature with odd degree, one-dimensional operators show that the BR1 flux results in suboptimal solution convergence rate of $ p $ \cite{kirby2005selecting,shu2001different,bassi2005discontinuous,hesthaven2007nodal}. However, as can be seen from \cref{fig:Solution error}, this characteristic is not observed when the BR1 SAT is implemented with SBP operators on unstructured triangular meshes. \violet{For the BR1 and LDG SATs, if $ \T_{\gamma k}^{(1)} $ and $ \T_{\gamma}^{(D)} $ are not modified and the extended boundary SATs are not included, then discretizations with the SBP-$ \Omega $ and SBP-$ \Gamma $ operators produce system matrices that have eigenvalues with positive real parts. For the unmodified\footnote{\violet{The unmodified BR1 and LDG SATs are denoted by BR1* and LDG*, respectively, in all figures and tables. If used without a qualifier, the names BR1 and LDG refer to the modified versions of the BR1 and LDG SATs.}} BR1 and LDG SATs, which include extended boundary SATs, positive eigenvalues are not produced with all types of the SBP operators. Despite being stable, however, functional superconvergence is not observed for the unmodified BR1 and LDG SATs except when used with the SBP diagonal-E operators. As noted in \cref{sec:Equivalence of SATs}, when used with SBP diagonal-E operators, the BR1 and LDG SATs (both modified and unmodified) have compact stencil width, and they are adjoint consistent for problems  with non-homogeneous Dirichlet boundary conditions. When the unmodified LDG SAT is implemented with $ \mu =0$, suboptimal solution convergence rates are observed for some of the cases; hence, we implemented the unmodified LDG SAT with $ \T_{\gamma k}^{(D)} =\frac{3}{2}\B_{\gamma}\Upsilon_{\gamma\gamma k}\B_{\gamma} $, which corresponds to a nonzero value of $ \mu $ at Dirichlet boundary facets. It can be seen from \cref{fig:Solution error} that the unmodified BR1 and LDG SATs lead to solution convergence rates of $ p+1 $.}

\begin{figure}[!t]
	\centering
	\begin{subfigure}{0.5\textwidth}
		\centering
		\includegraphics[scale=0.25]{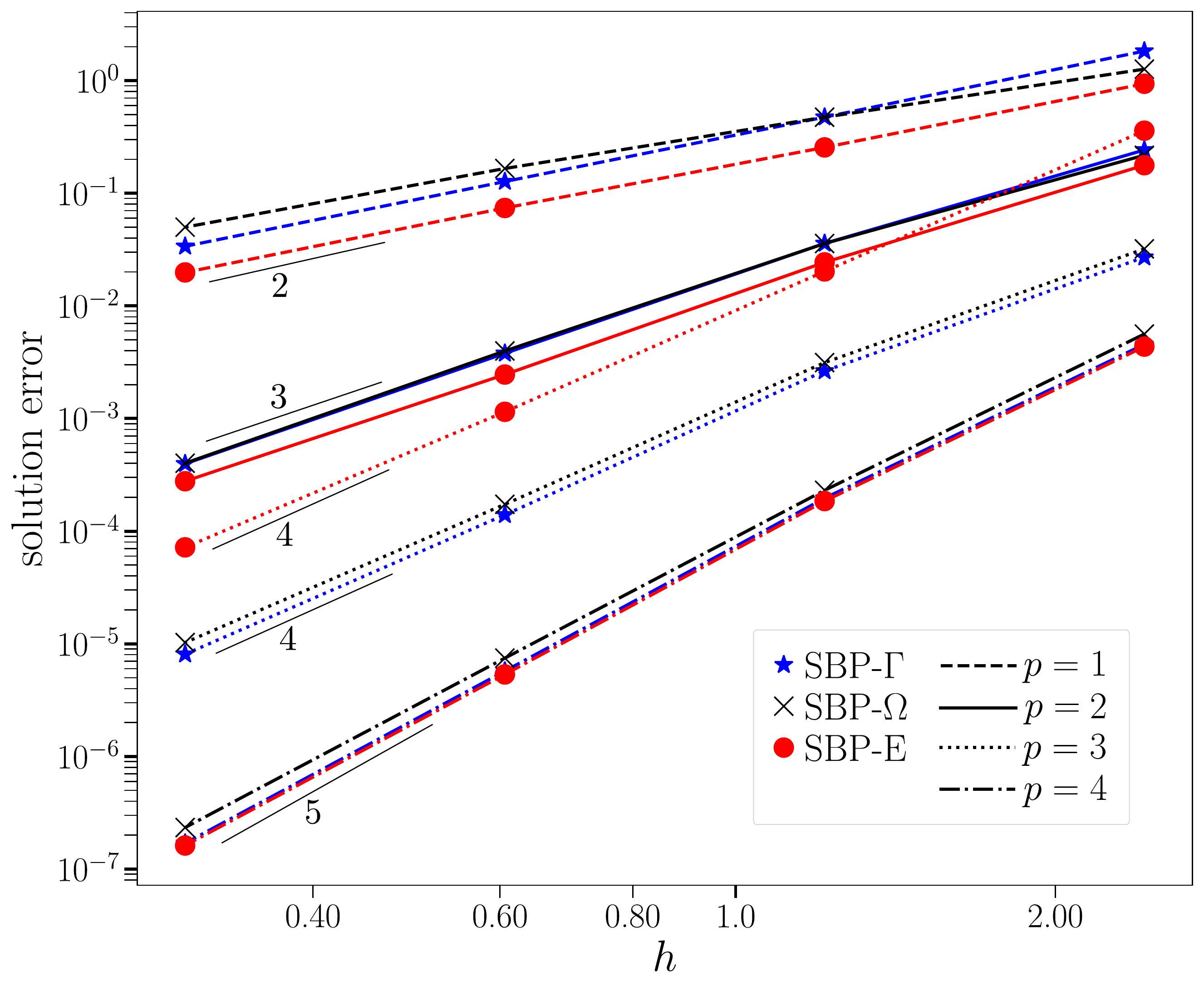}
		\caption{\label{fig:Soln error all BR1 SAT} Solution error with BR1 SAT}
	\end{subfigure}\hfill
	\begin{subfigure}{0.5\textwidth}
		\centering
		\includegraphics[scale=0.25]{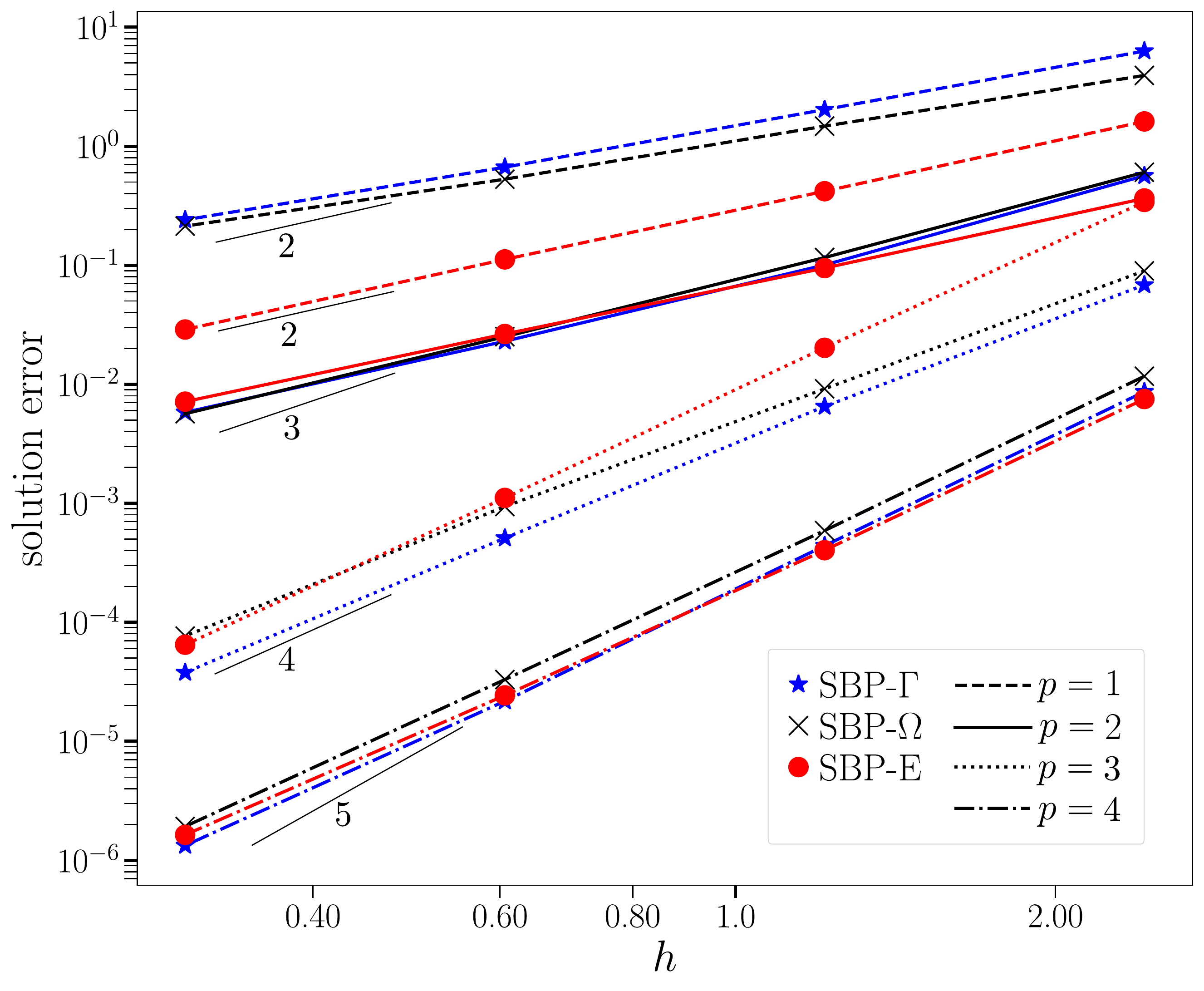}
		\caption{\label{fig:Soln error all BO SAT} Solution error with BO SAT}
	\end{subfigure} 
	\caption{\label{fig:Soln error all} Variation of solution error under grid refinement with respect to three types of SBP operators. Slopes corresponding to $ p+1 $ convergence rates are shown by short, thin lines.}
\end{figure}

Figure \ref{fig:Soln error all} shows the errors produced by the three types of SBP operator when implemented with the BR1 and BO SATs. In general, solution error is not very sensitive to the type of SBP operator used except in a few cases, \eg, the cases where the degree three SBP diagonal-E operator is implemented with SATs other than the BO SAT. Except for the BO SAT, all of the other SATs show very similar solution error convergence behavior as that of the BR1 SAT.

\begin{figure}[!t]
	\centering
	\begin{subfigure}{0.33\textwidth}
		\centering
		\includegraphics[scale=0.2]{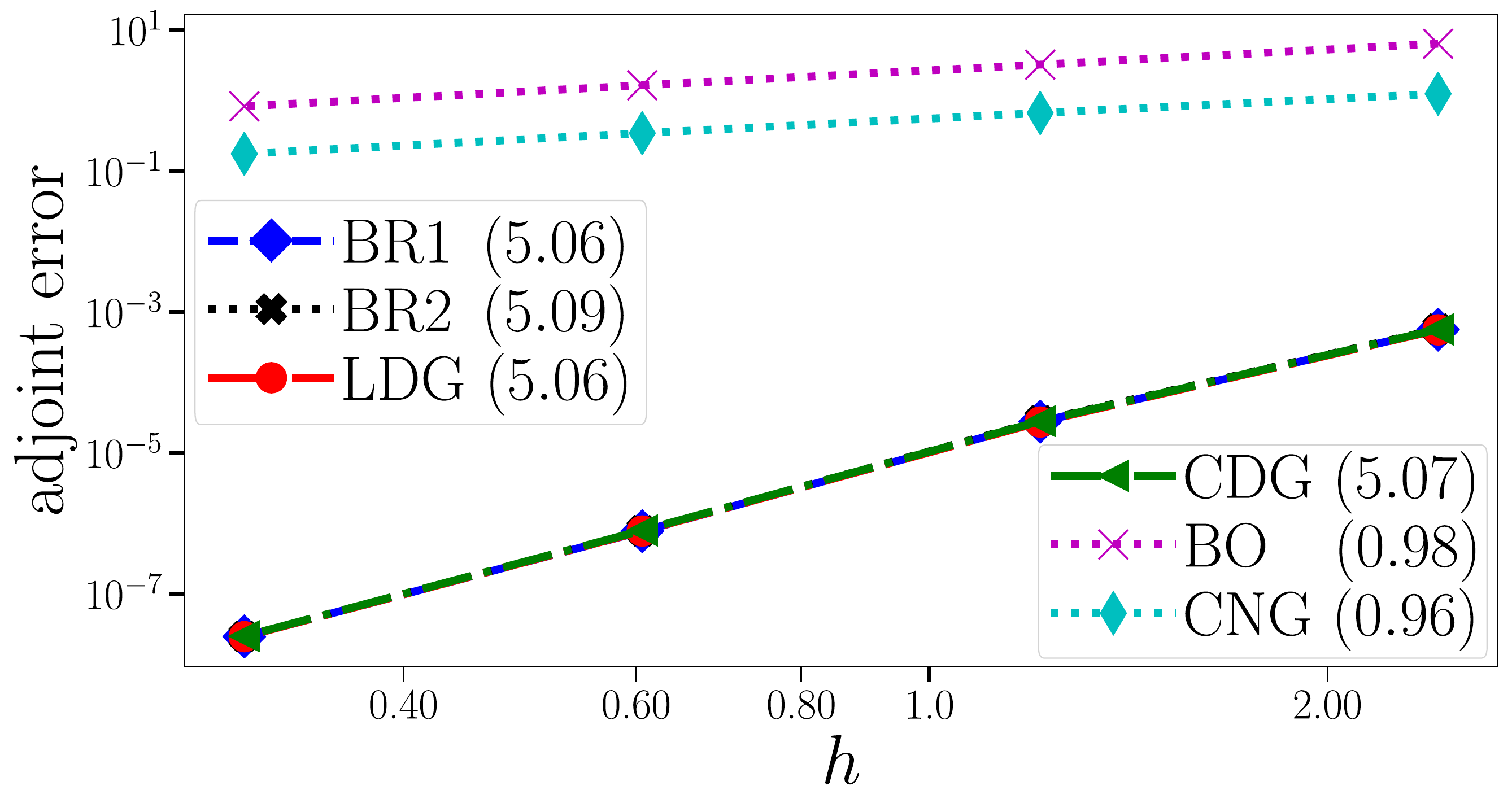}
		\caption{\label{fig:Omega adj p=3}SBP-$ \Omega $, $p=3$}
	\end{subfigure}\hfill
	\begin{subfigure}{0.33\textwidth}
		\centering
		\includegraphics[scale=0.2]{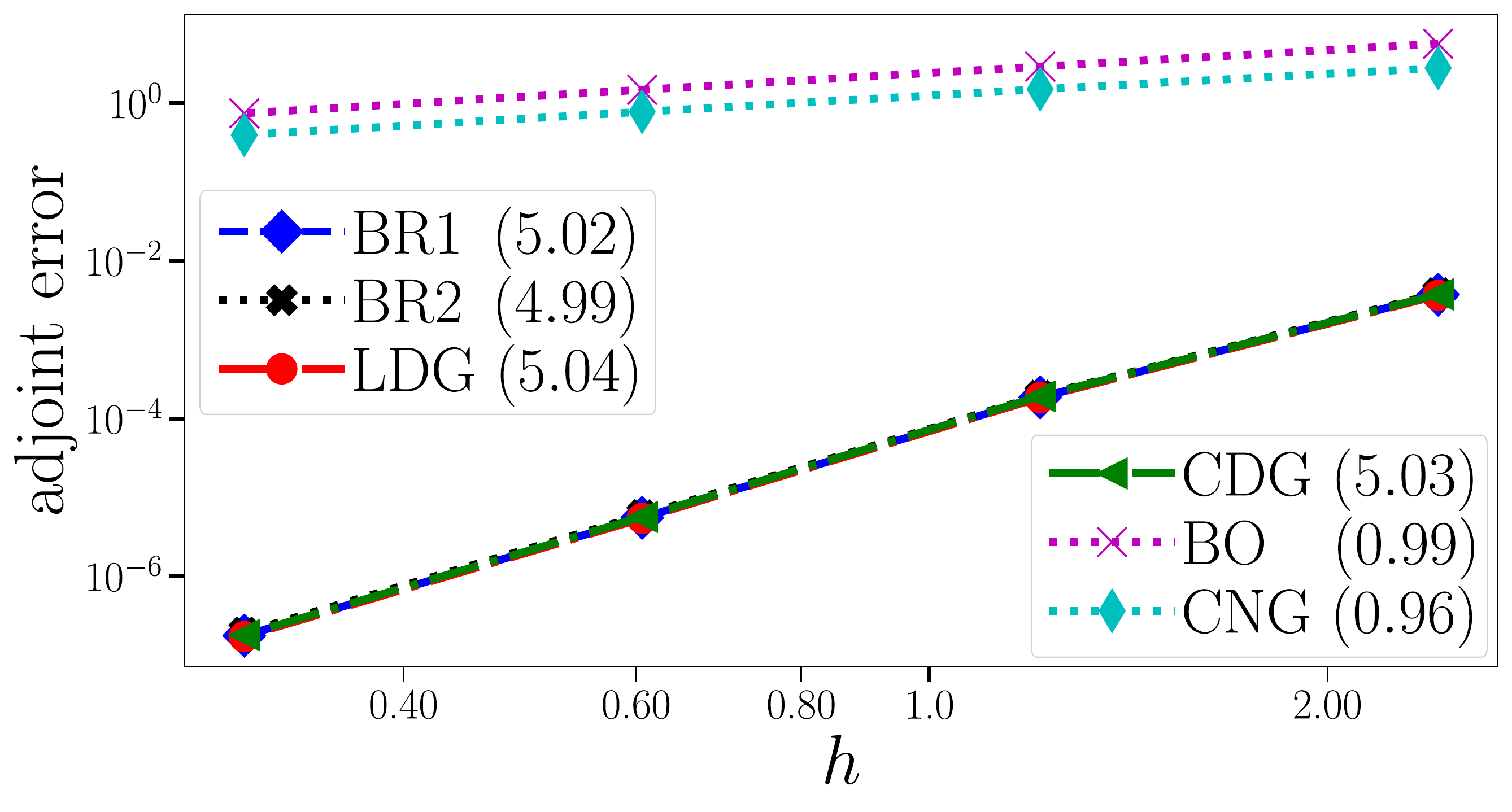}
		\caption{\label{fig:Gamma adj p=3} SBP-$ \Gamma $, $p=3$}
	\end{subfigure} 
	\begin{subfigure}{0.33\textwidth}
		\centering
		\includegraphics[scale=0.2]{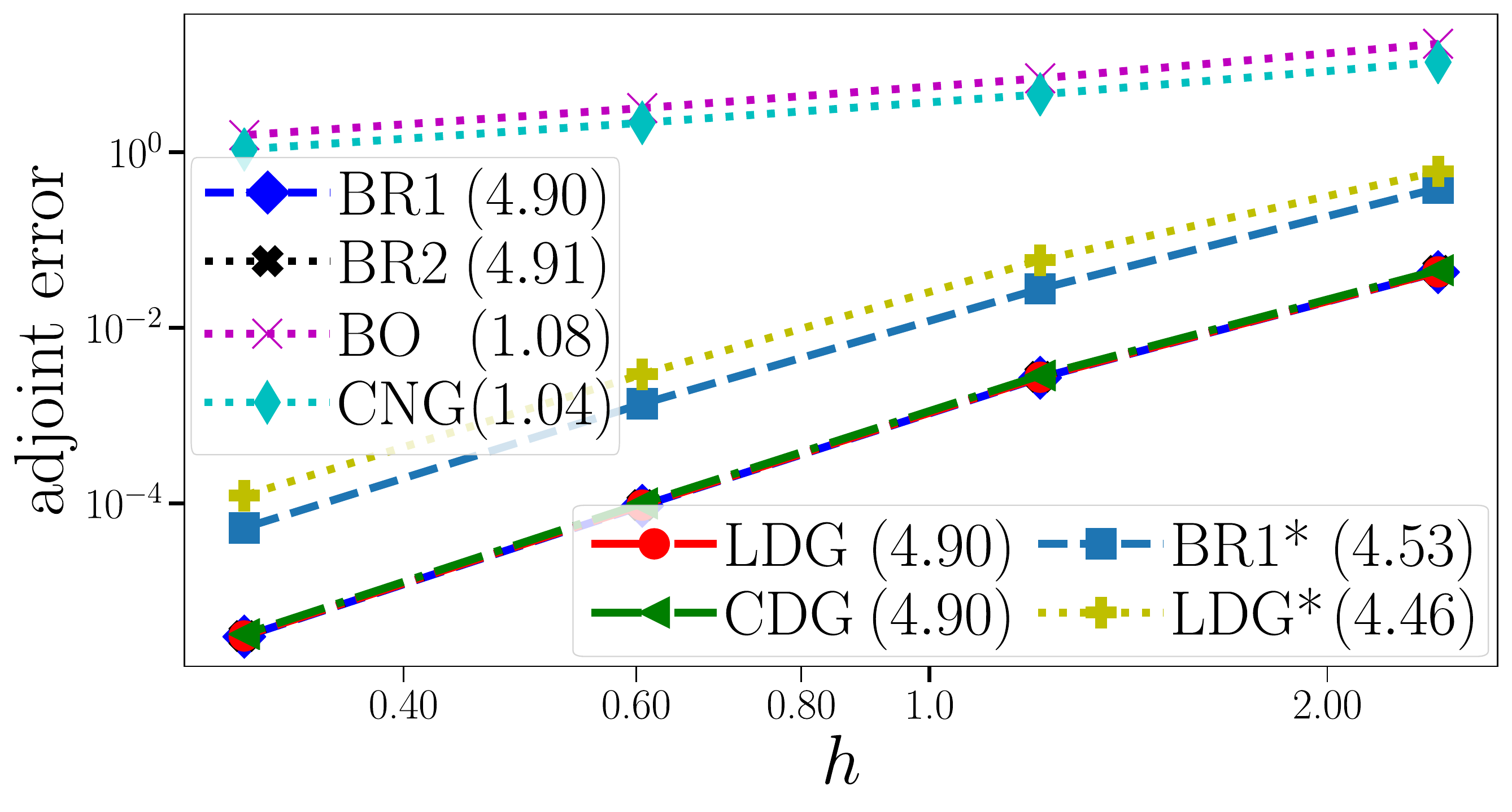}
		\caption{\label{fig:diage adj p=3} SBP-E, $p=3$}
	\end{subfigure}
	\\
	\begin{subfigure}{0.33\textwidth}
		\centering
		\includegraphics[scale=0.2]{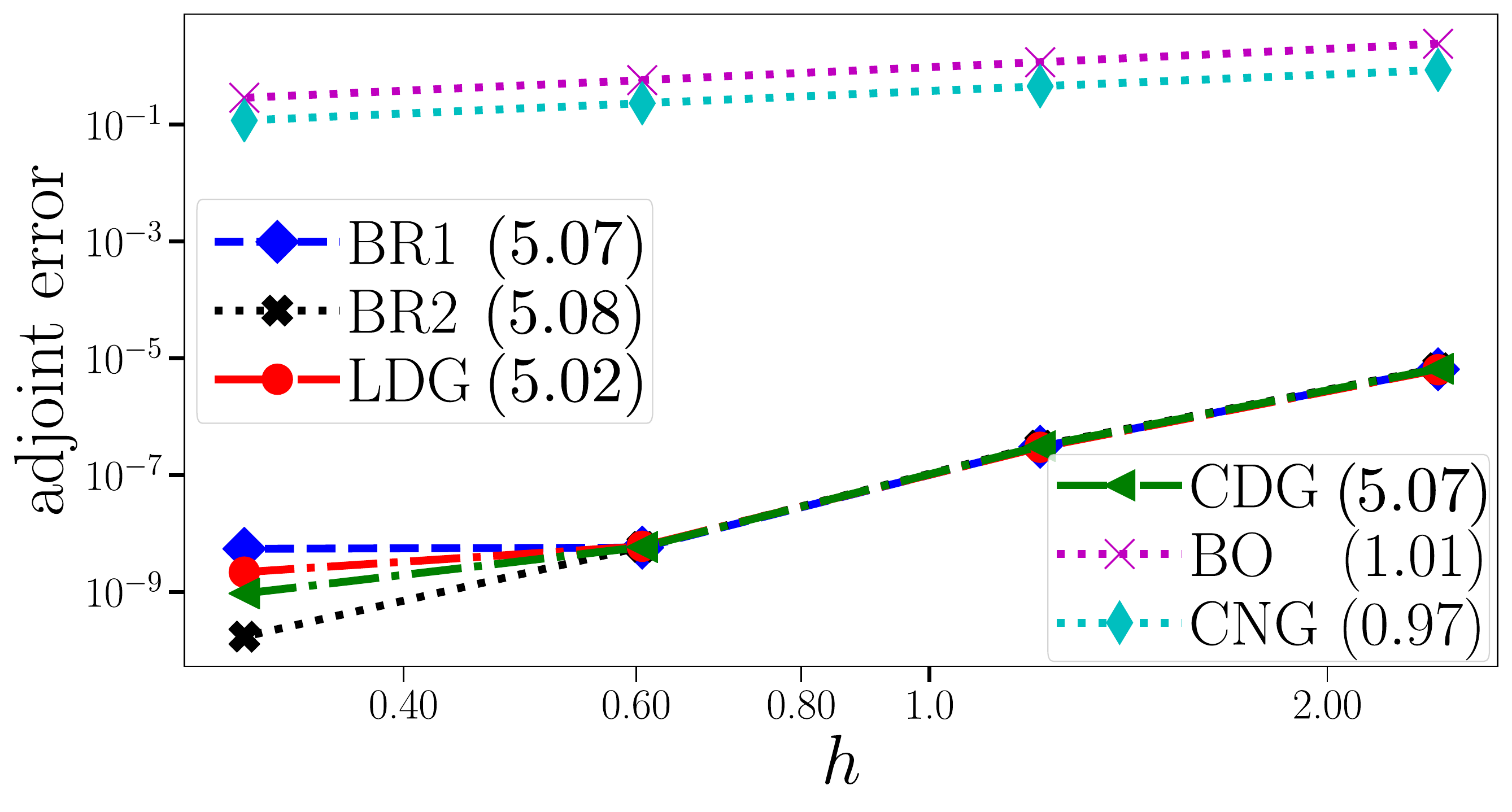}
		\caption{\label{fig:Omega adj p=4} SBP-$ \Omega $, $p=4$}
	\end{subfigure}\hfill
	\begin{subfigure}{0.33\textwidth}
		\centering
		\includegraphics[scale=0.2]{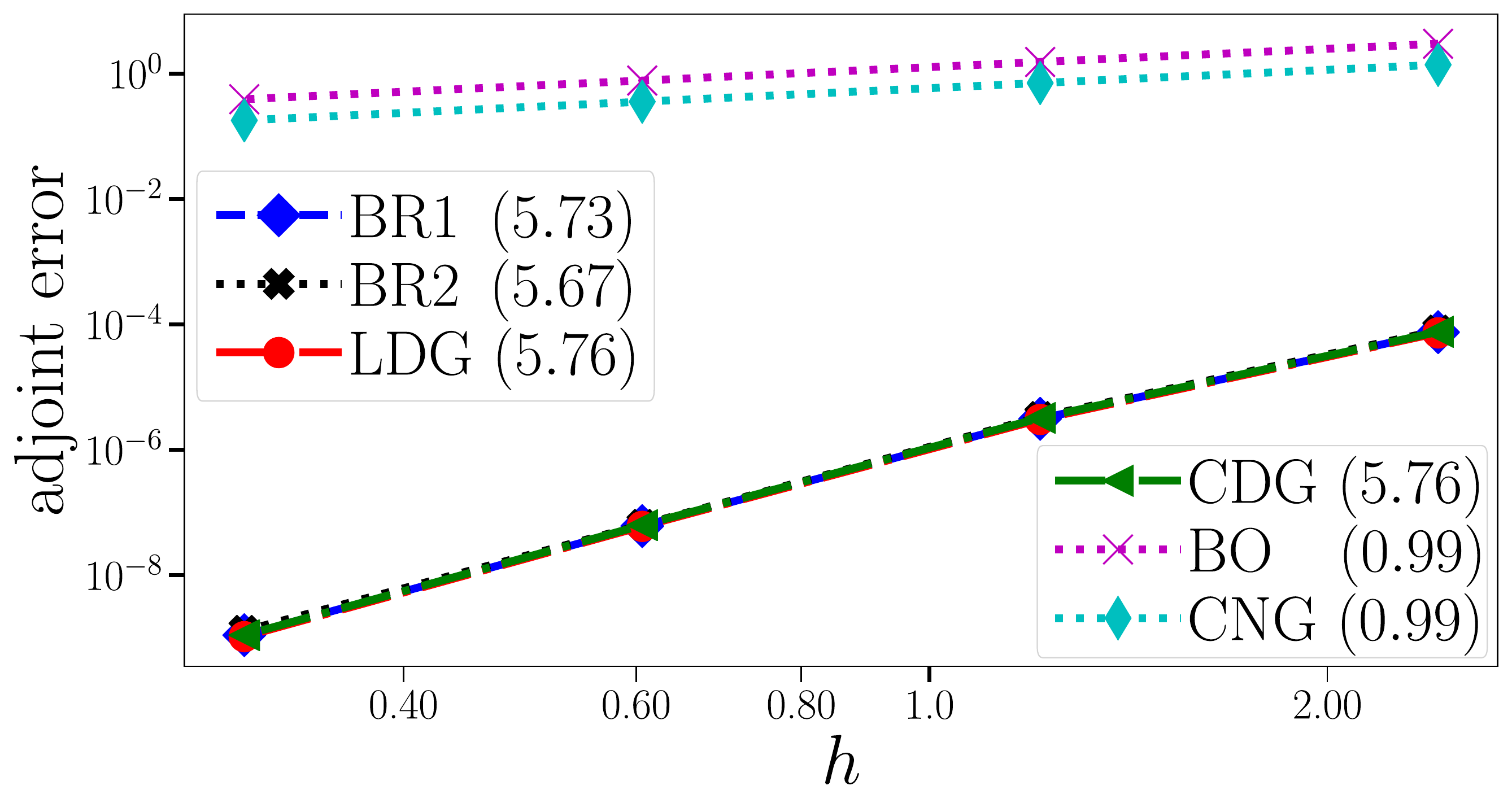}
		\caption{\label{fig:Gamma adj p=4} SBP-$ \Gamma $, $p=4$}
	\end{subfigure} 
	\begin{subfigure}{0.33\textwidth}
		\centering
		\includegraphics[scale=0.2]{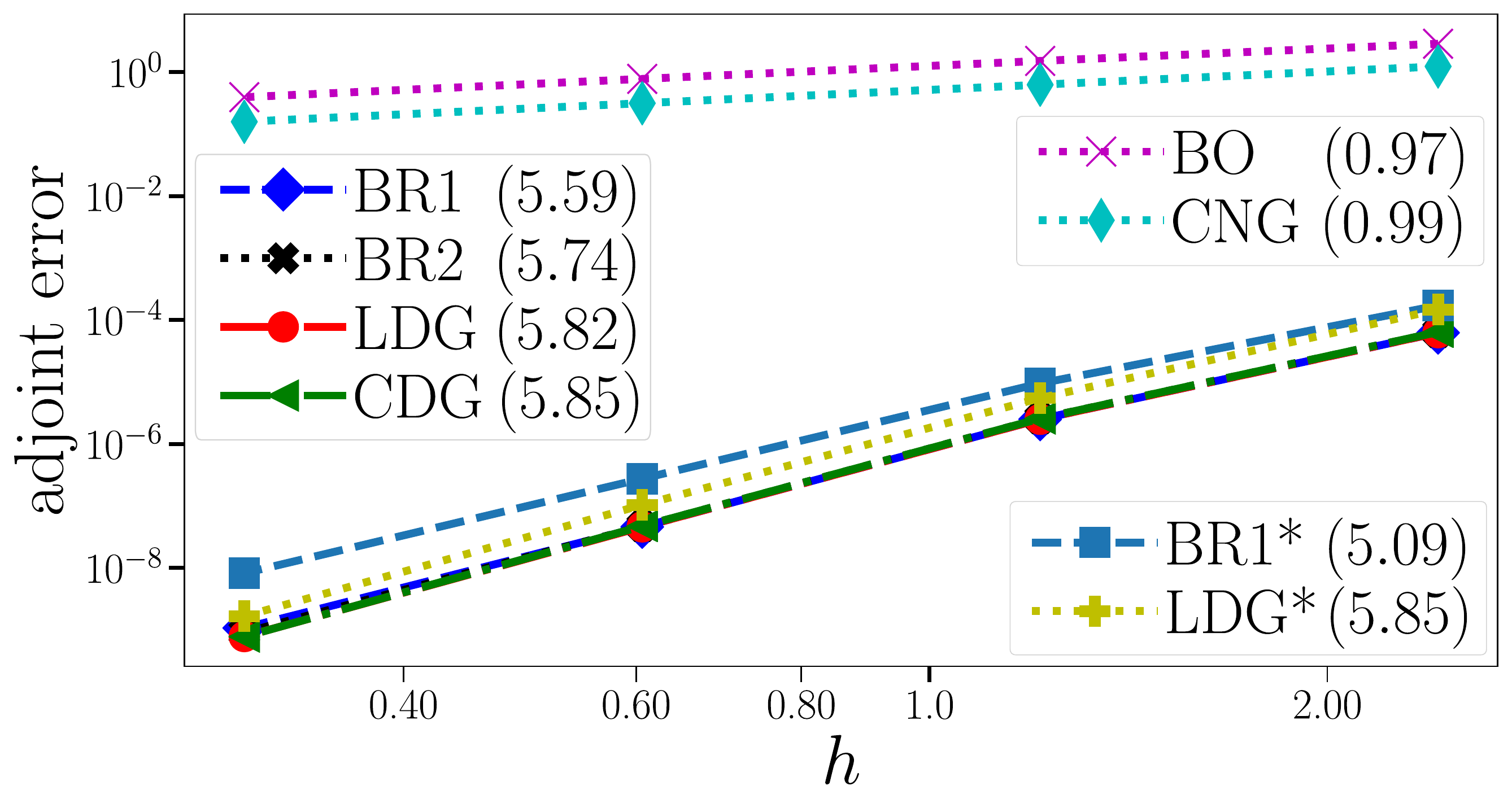}
		\caption{\label{fig:diage adj p=4} SBP-E, $p=4$}
	\end{subfigure}
	\caption{\label{fig:Adjoint error} Adjoint error under grid refinement. Adjoint convergence rates (shown in parenthesis) are calculated by fitting a line through the last three error values on the refined meshes except for the adjoint consistent SATs with SBP-$ \Omega $ operator of degree $ p=4 $ for which the first 3 grids are used. \violet{The BR1* and LDG* SATs represent the unmodified BR1 and LDG SATs}.}
\end{figure}

The errors and convergence rates of the adjoint solution under mesh refinement are presented in \cref{fig:Adjoint error}. All of the adjoint consistent SATs lead to schemes that converge to the exact adjoint at a rate of $p+1$ \violet{or larger}. In contrast, schemes with the BO and CNG SATs have error values of $ \fnc{O}(1) $. Similar properties as with the primal solution are observed regarding the sensitivity of the adjoint error values to the type of SBP operator used.

\begin{figure}[!t]
	\centering
\ignore{
	\begin{subfigure}{0.33\textwidth}
		\centering
		\includegraphics[scale=0.2]{errs_func_VarOper_omega_p1.pdf}
		\caption{\label{fig:Omega func p=1} SBP-$ \Omega $, $p=1$}
	\end{subfigure}\hfill
	\begin{subfigure}{0.33\textwidth}
		\centering
		\includegraphics[scale=0.2]{errs_func_VarOper_gamma_p1.pdf}
		\caption{\label{fig:Gamma func p=1} SBP-$ \Gamma $, $p=1$}
	\end{subfigure} 
	\begin{subfigure}{0.33\textwidth}
		\centering
		\includegraphics[scale=0.2]{errs_func_VarOper_diage_p1.pdf}
		\caption{\label{fig:diage func p=1} SBP-E, $p=1$}
	\end{subfigure}
	\\
	\begin{subfigure}{0.33\textwidth}
		\centering
		\includegraphics[scale=0.2]{errs_func_VarOper_omega_p2.pdf}
		\caption{\label{fig:Omega func p=2} SBP-$ \Omega $, $p=2$}
	\end{subfigure}\hfill
	\begin{subfigure}{0.33\textwidth}
		\centering
		\includegraphics[scale=0.2]{errs_func_VarOper_gamma_p2.pdf}
		\caption{\label{fig:Gamma func p=2} SBP-$ \Gamma $, $p=2$}
	\end{subfigure} 
	\begin{subfigure}{0.33\textwidth}
		\centering
		\includegraphics[scale=0.2]{errs_func_VarOper_diage_p2.pdf}
		\caption{\label{fig:diage func p=2} SBP-E, $p=2$}
	\end{subfigure}
	\\
}
	\begin{subfigure}{0.33\textwidth}
		\centering
		\includegraphics[scale=0.2]{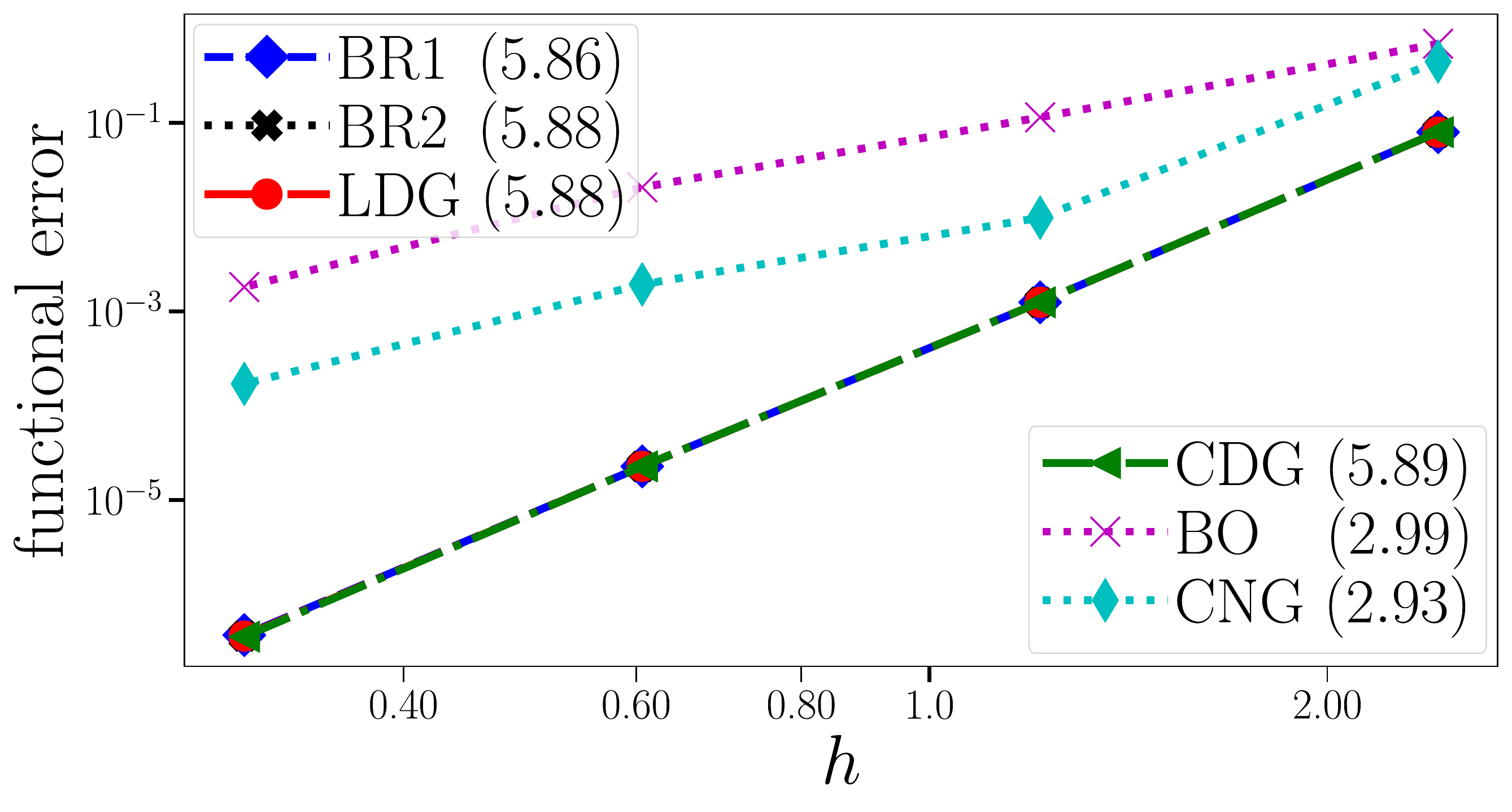}
		\caption{\label{fig:Omega func p=3} SBP-$ \Omega $, $p=3$}
	\end{subfigure}\hfill
	\begin{subfigure}{0.33\textwidth}
		\centering
		\includegraphics[scale=0.2]{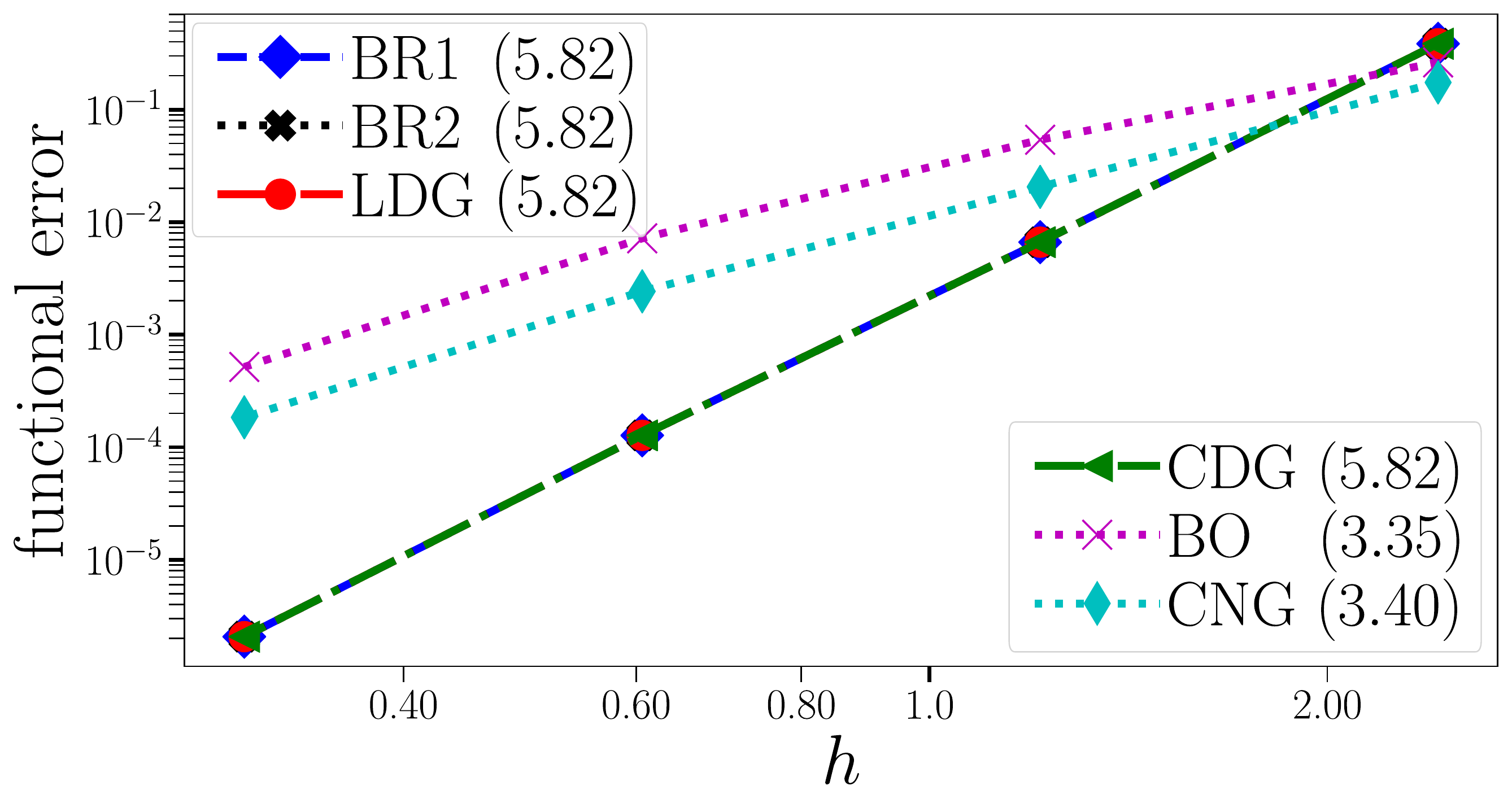}
		\caption{\label{fig:Gamma func p=3} SBP-$ \Gamma $, $p=3$}
	\end{subfigure} 
	\begin{subfigure}{0.33\textwidth}
		\centering
		\includegraphics[scale=0.2]{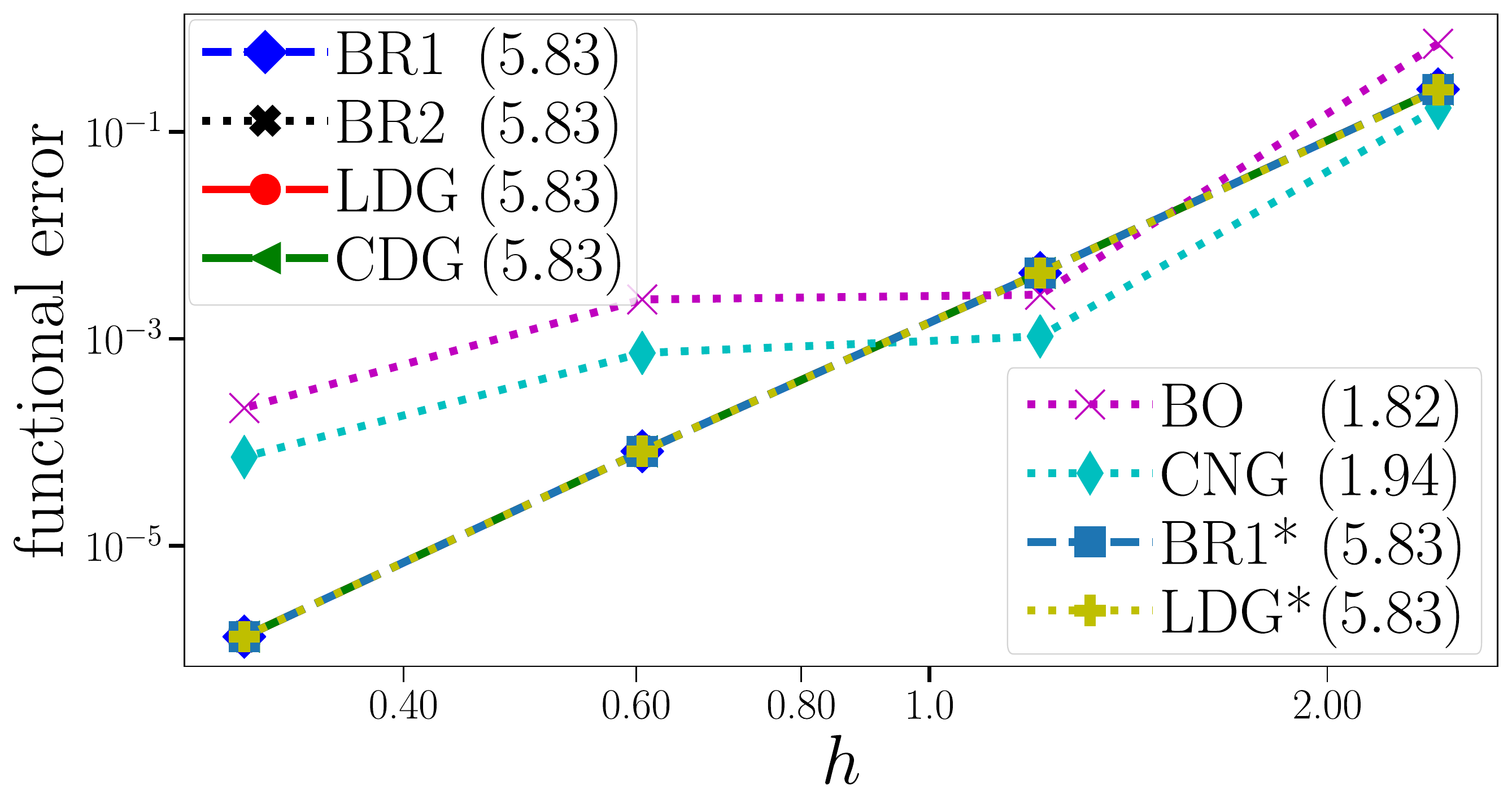}
		\caption{\label{fig:diage func p=3} SBP-E, $p=3$}
	\end{subfigure}
	\\
	\begin{subfigure}{0.33\textwidth}
		\centering
		\includegraphics[scale=0.2]{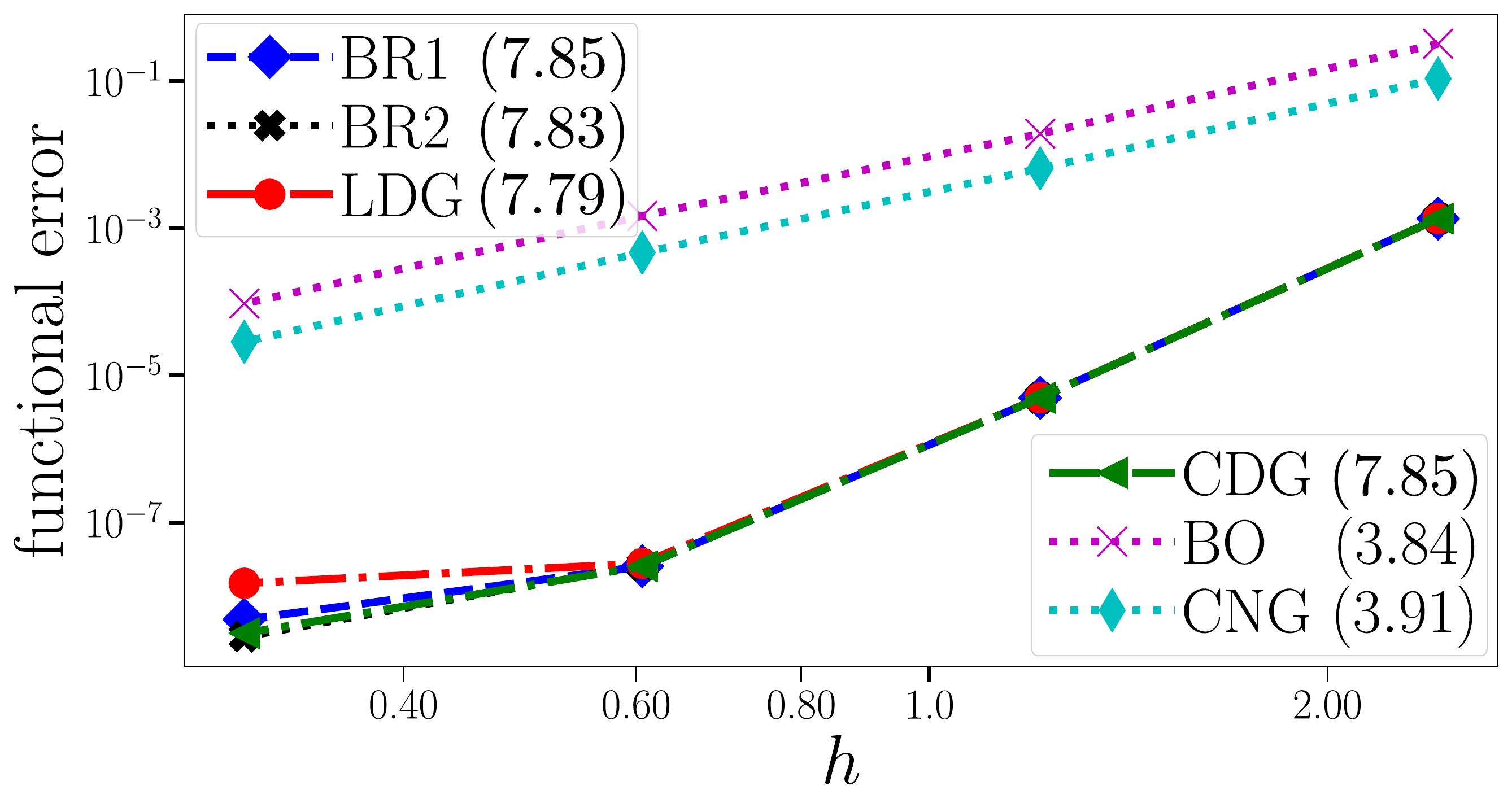}
		\caption{\label{fig:Omega func p=4} SBP-$ \Omega $, $p=4$}
	\end{subfigure}\hfill
	\begin{subfigure}{0.33\textwidth}
		\centering
		\includegraphics[scale=0.2]{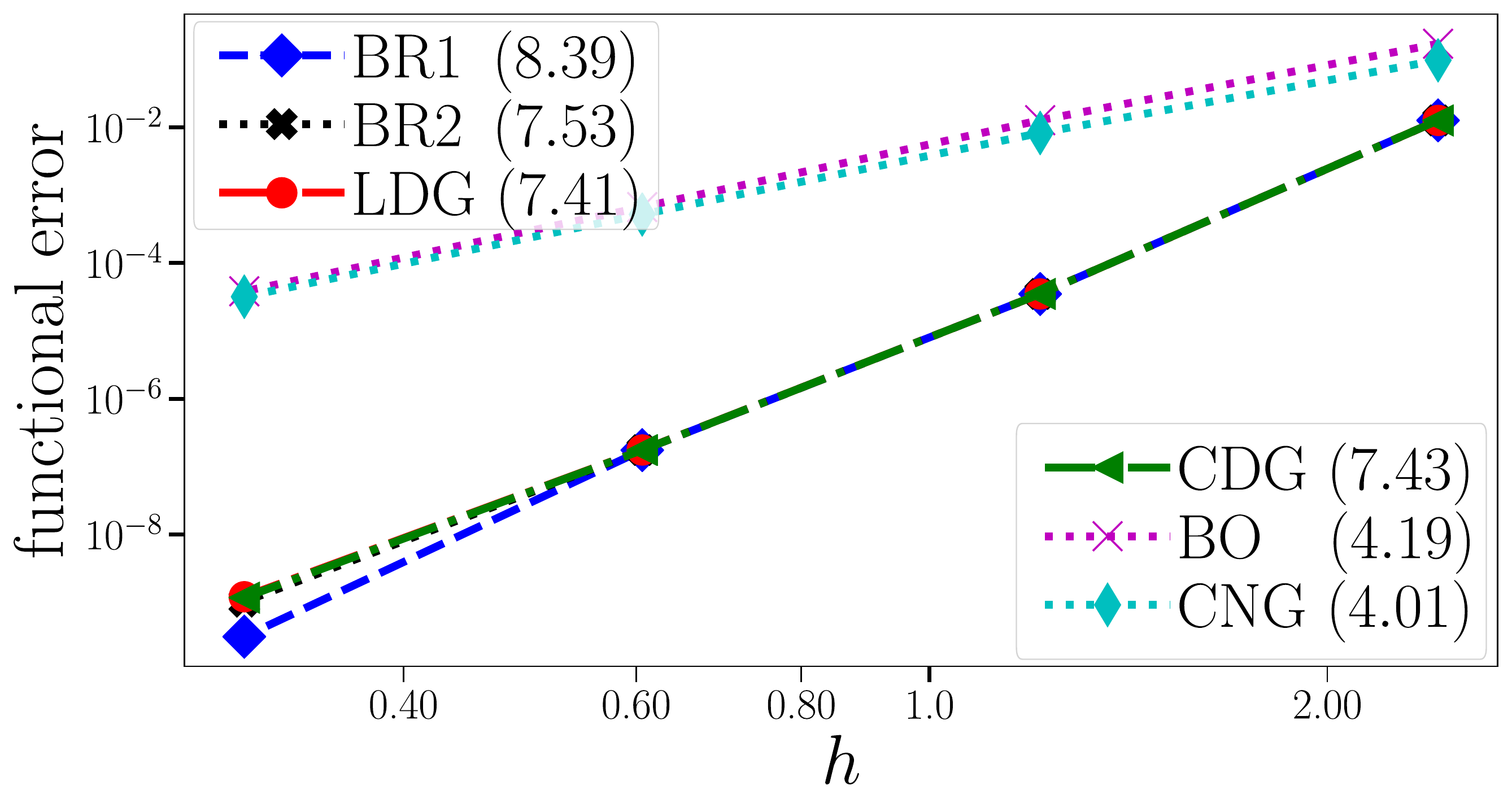}
		\caption{\label{fig:Gamma func p=4} SBP-$ \Gamma $, $p=4$}
	\end{subfigure} 
	\begin{subfigure}{0.33\textwidth}
		\centering
		\includegraphics[scale=0.2]{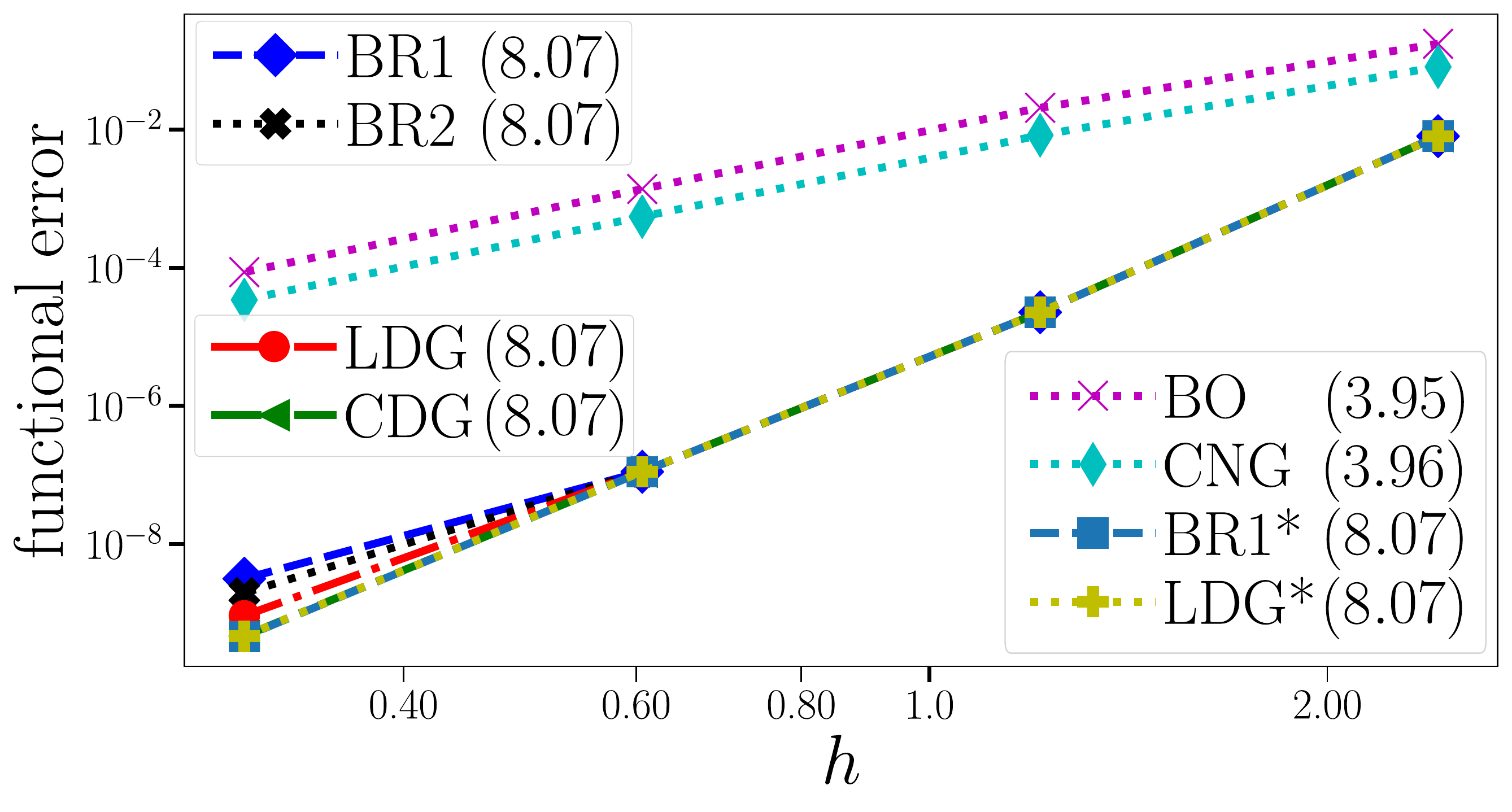}
		\caption{\label{fig:diage func p=4} SBP-E, $p=4$}
	\end{subfigure}
	\caption{\label{fig:Functional error} Functional error under grid refinement. Functional convergence rates (shown in parenthesis) are calculated by fitting a line through the last three error values on the refined meshes except for the adjoint consistent SATs with SBP-$ \Omega $ and SBP-E operators of degree $ p=4 $ for which the first 3 grids are used. \violet{The BR1* and LDG* SATs represent the unmodified BR1 and LDG SATs}.}
\end{figure}

Functional errors and convergence rates are displayed in \cref{fig:Functional error}. As expected, functional superconvergence rates of $ 2p $ are observed for schemes with primal and adjoint consistent SATs. The adjoint inconsistent SATs, BO and CNG, do not display functional superconvergence rates of $ 2p $. While the adjoint consistent schemes achieve comparable functional error values, the CNG SAT outperforms the BO SAT in this regard in most cases. 

\subsection{Eigenspectra} \label{sec:Eigenspectra}
\begin{figure}[!t]
	\begin{subfigure}{0.33\textwidth}
		\includegraphics[scale=0.16, right]{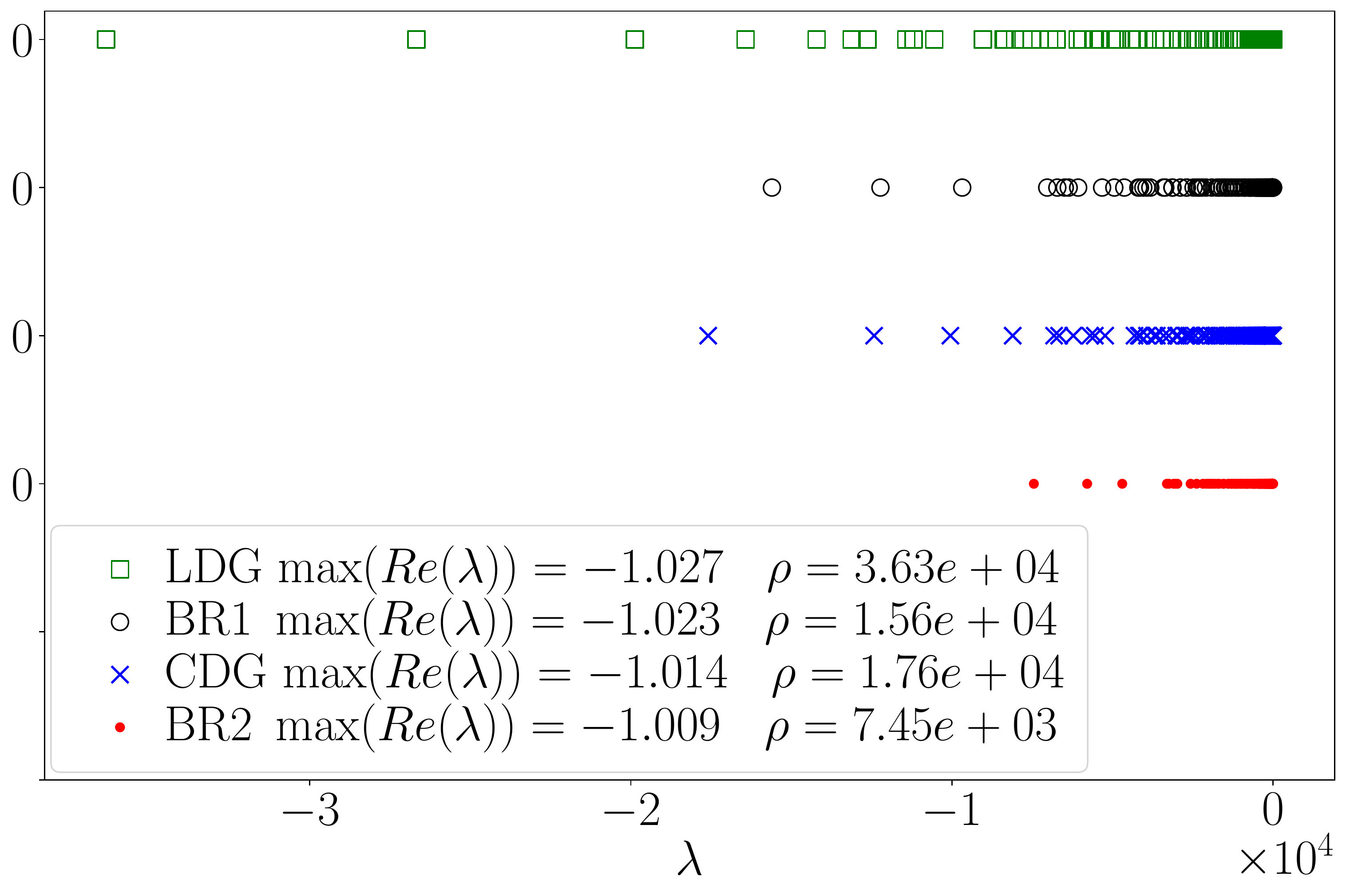}
		\caption{\label{fig:Omega symmetric spectra p=3} SBP-$ \Omega $, $p=3$}
	\end{subfigure}\hfill
	\begin{subfigure}{0.33\textwidth}
		\includegraphics[scale=0.16,  right]{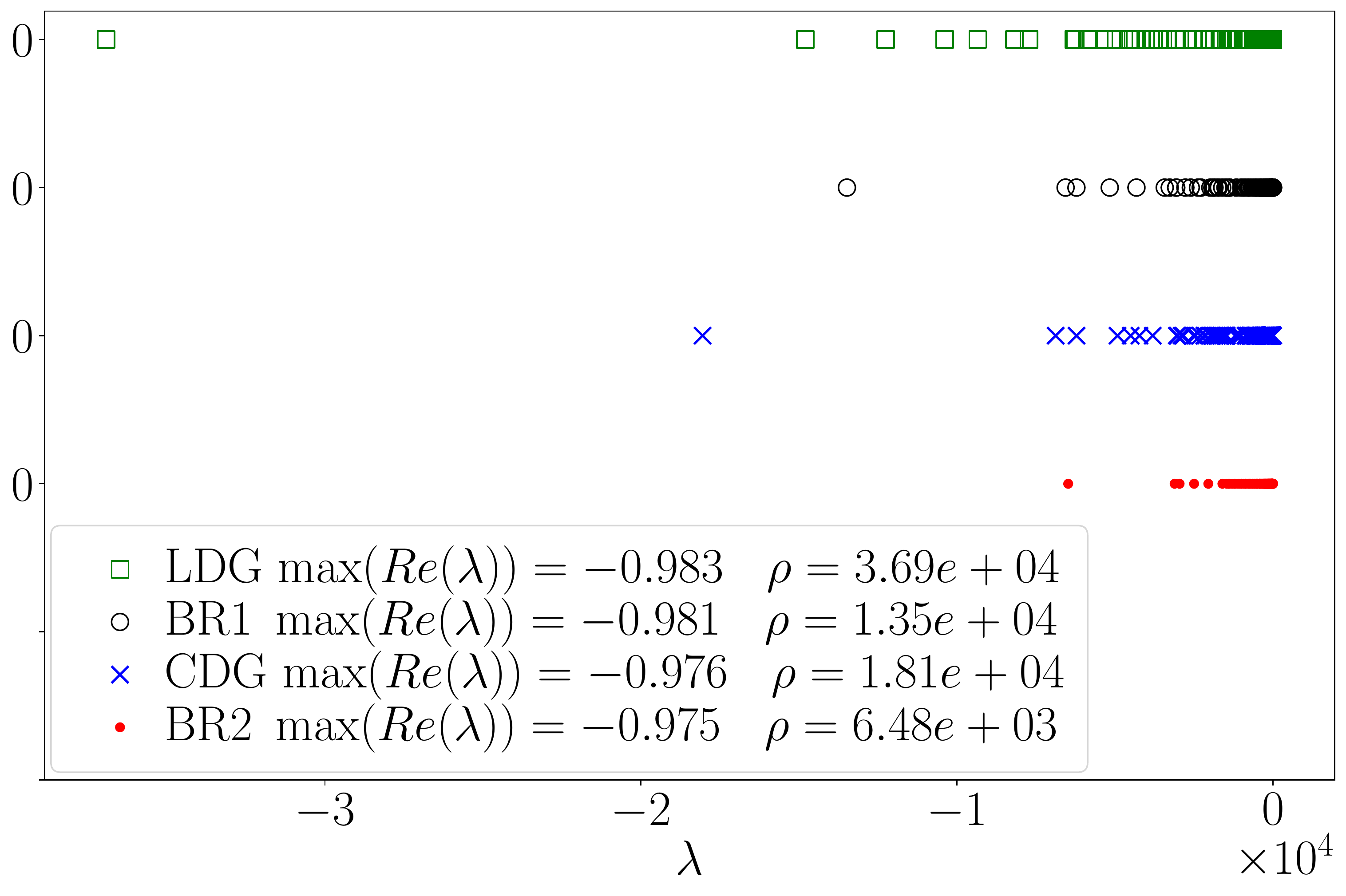}
		\caption{\label{fig:Gamma symmetric spectra p=3} SBP-$ \Gamma $, $p=3$}
	\end{subfigure} \hfill
	\begin{subfigure}{0.33\textwidth}
		\includegraphics[scale=0.16,right]{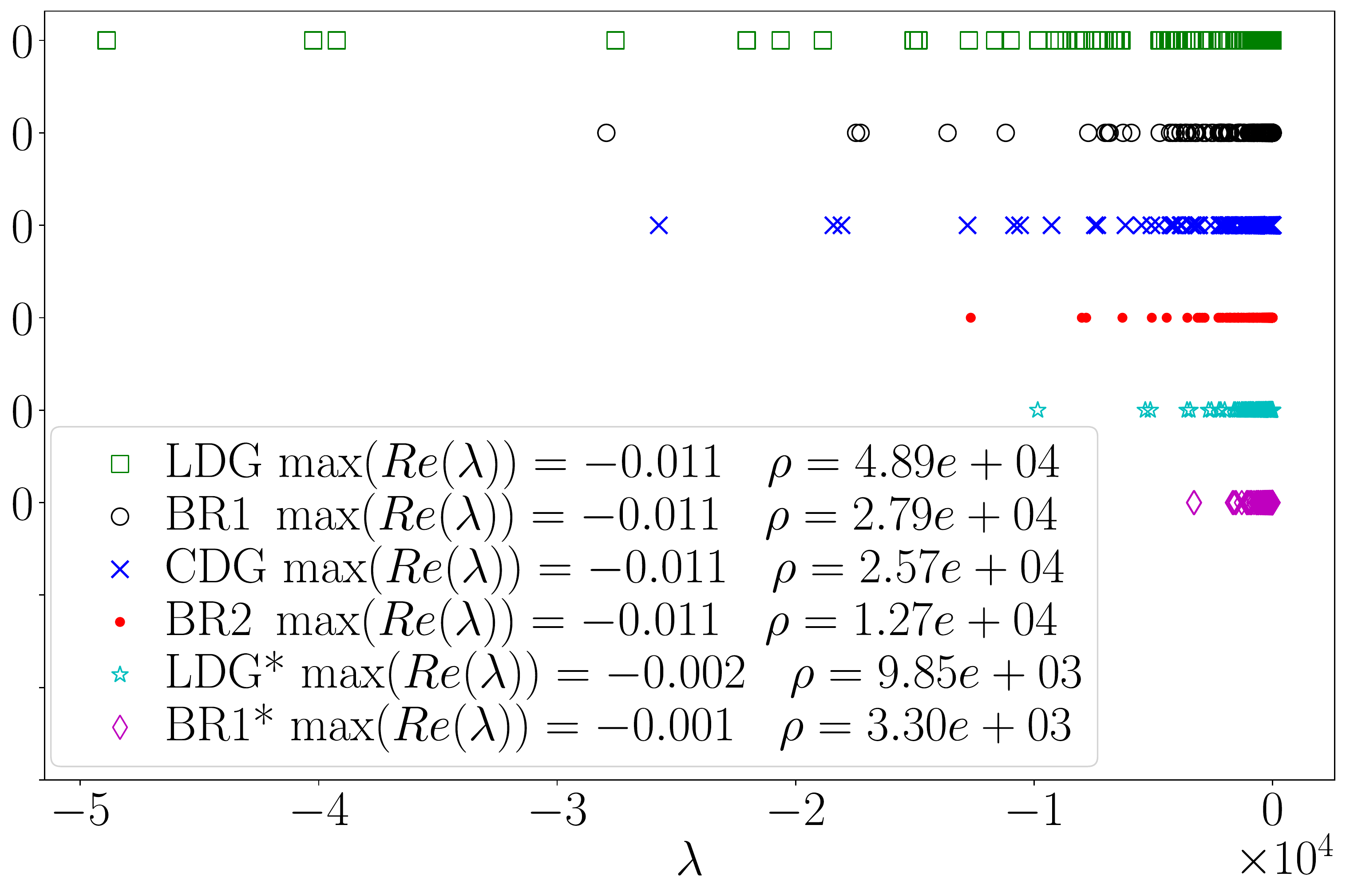}
		\caption{\label{fig:diage symmetric spectra p=3} SBP-E, $p=3$}
	\end{subfigure}
	\\
	\begin{subfigure}{0.33\textwidth}
		\includegraphics[scale=0.16,right]{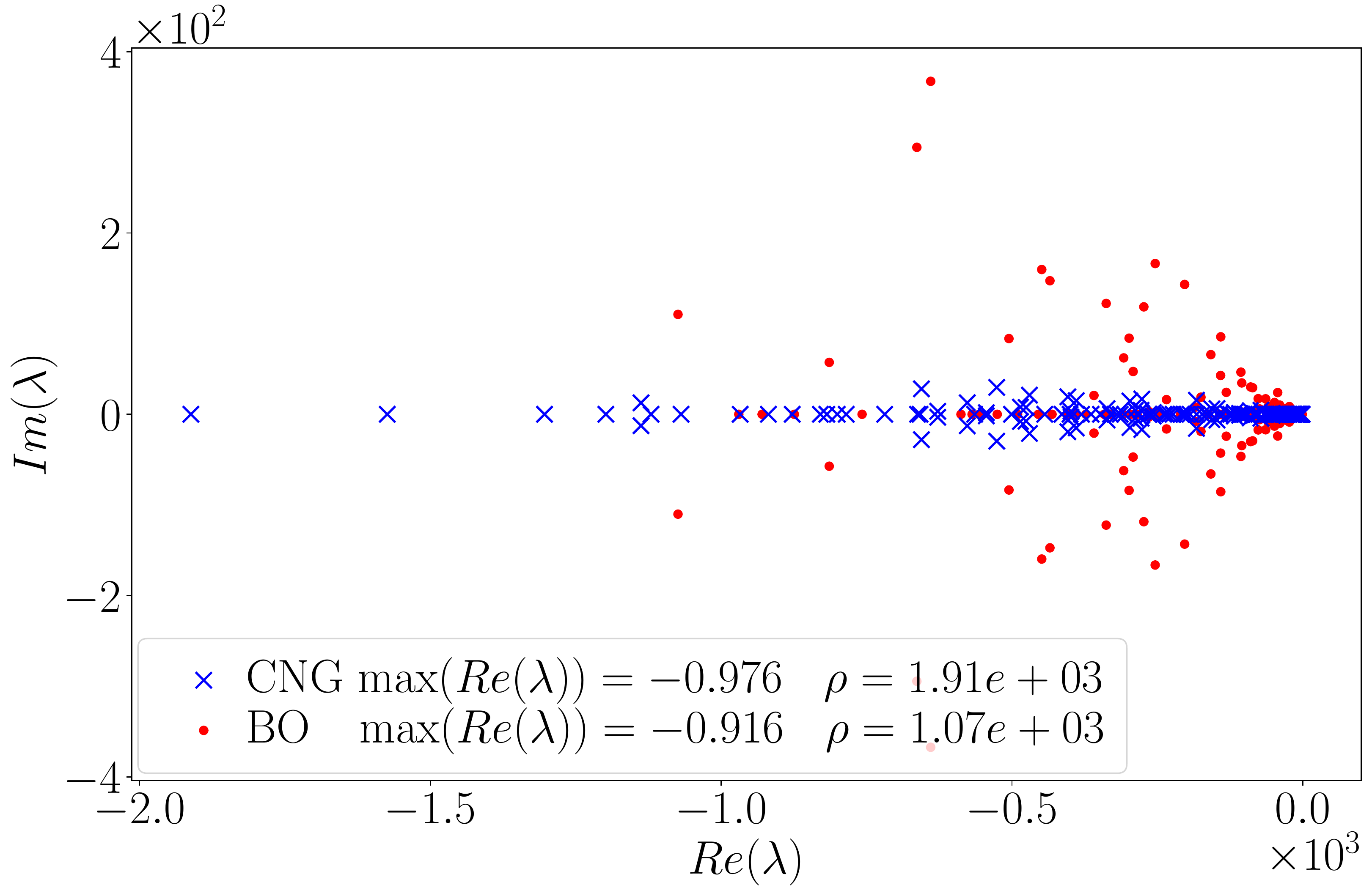}
		\caption{\label{fig:Omega asymmetric spectra p=3} SBP-$ \Omega $, $p=3$}
	\end{subfigure}\hfill
	\begin{subfigure}{0.33\textwidth}
		\includegraphics[scale=0.16,right]{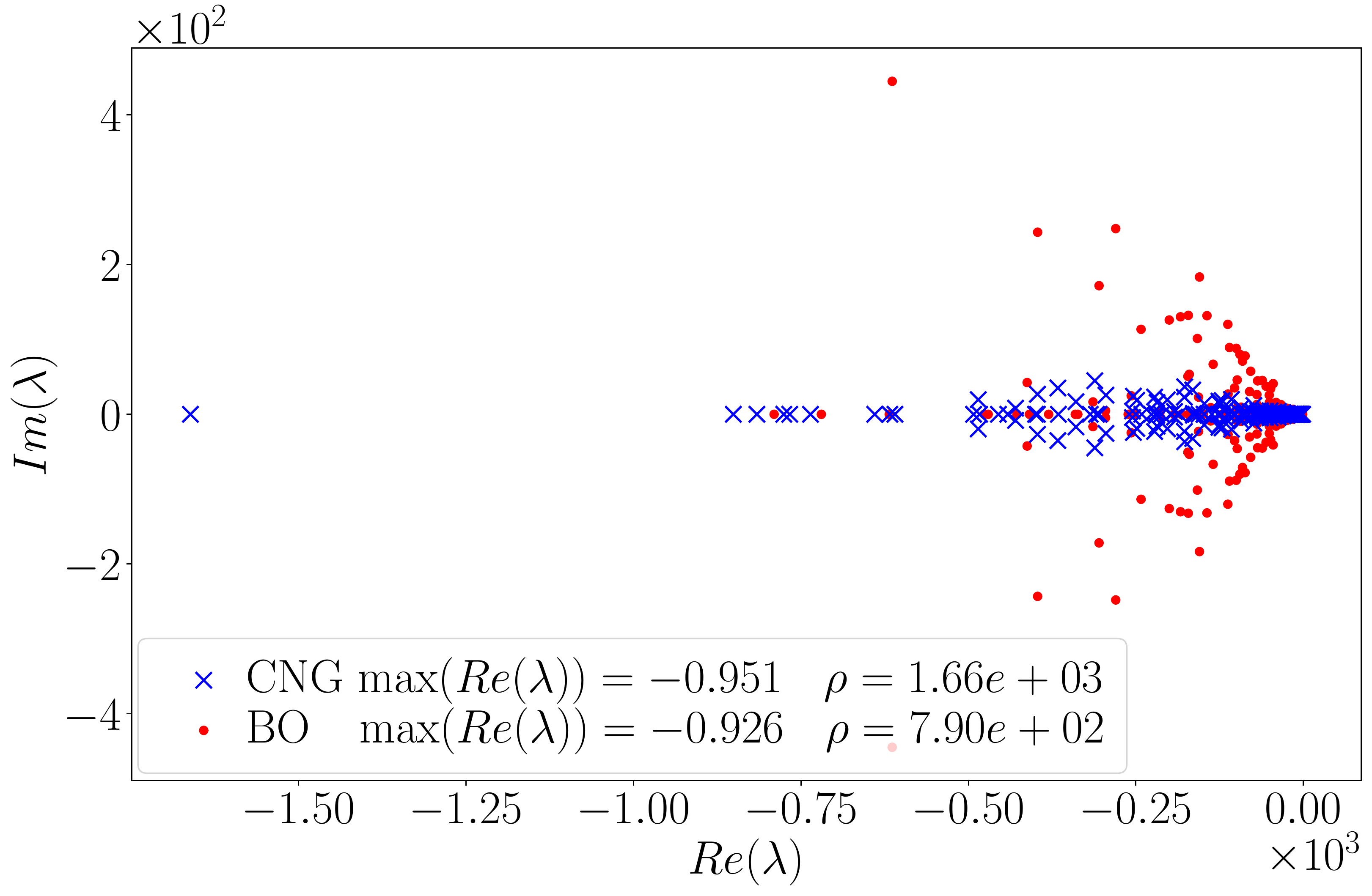}
		\caption{\label{fig:Gamma asymmetric spectra p=3} SBP-$ \Gamma $, $p=3$}
	\end{subfigure} \hfill
	\begin{subfigure}{0.33\textwidth}
		\includegraphics[scale=0.16,right]{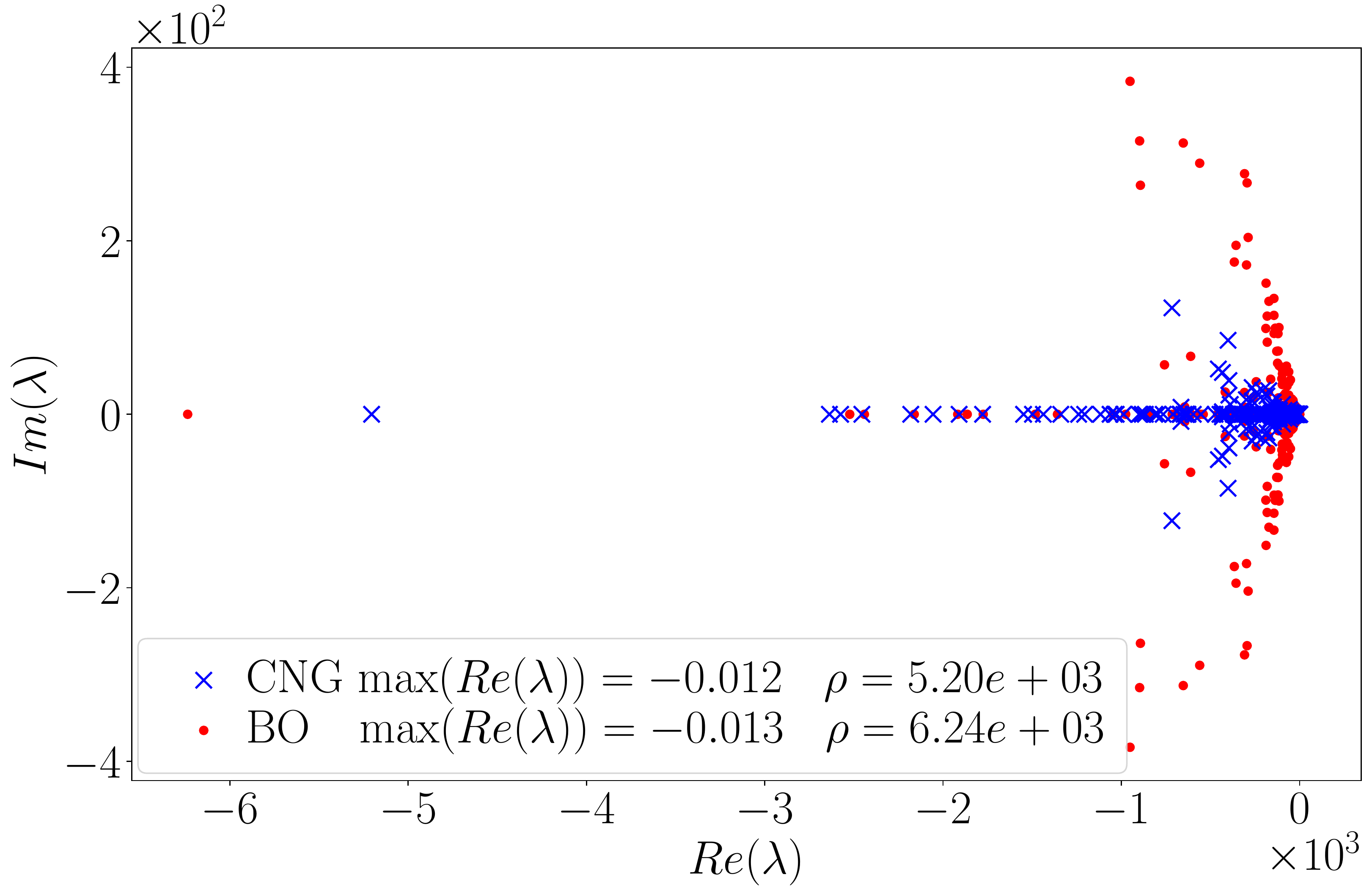}
		\caption{\label{fig:diage asymmetric spectra p=3} SBP-E, $p=3$}
	\end{subfigure}
	\\
	\begin{subfigure}{0.33\textwidth}
		\includegraphics[scale=0.16,right]{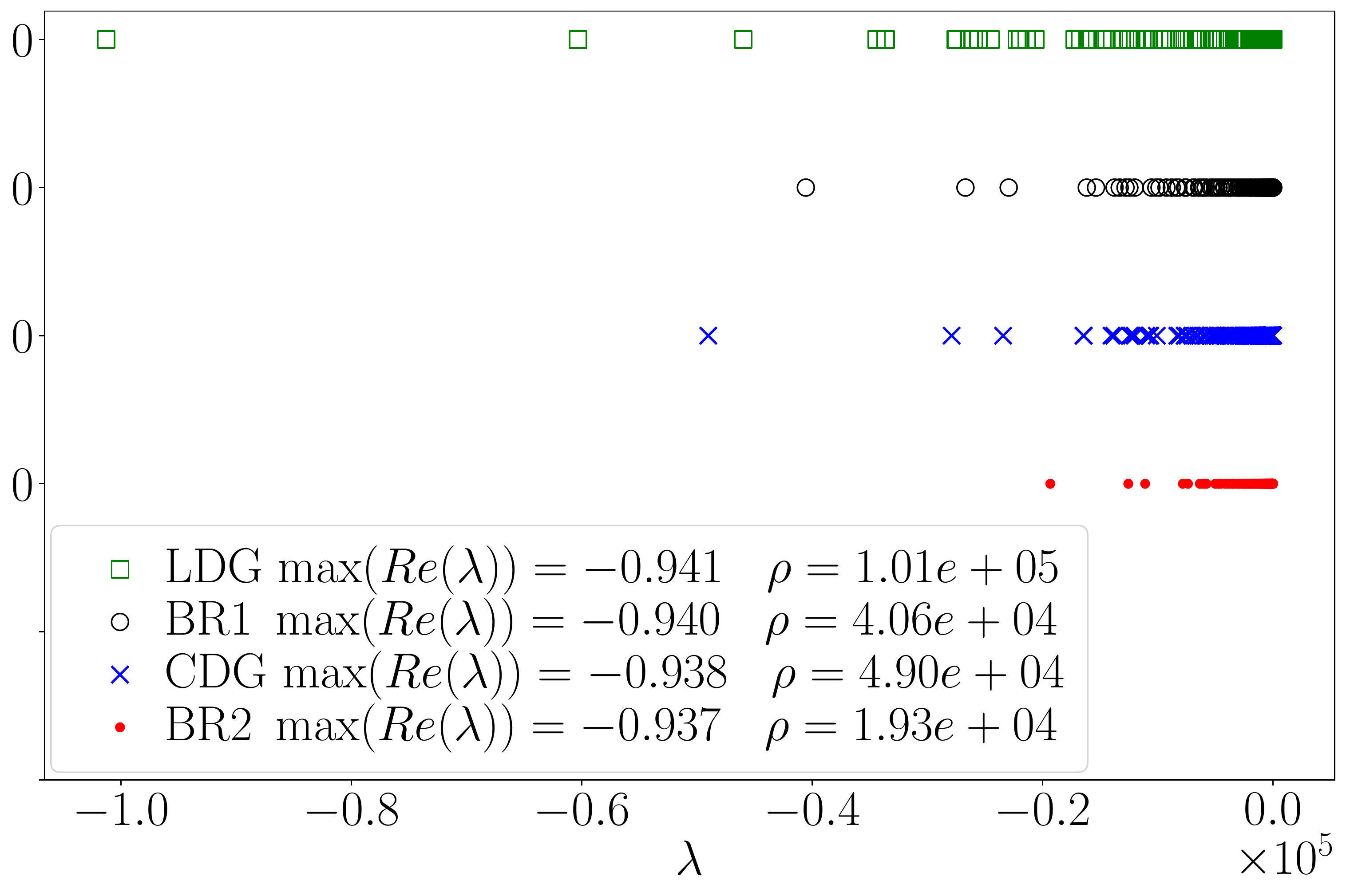}
		\caption{\label{fig:Omega symmetric spectra p=4} SBP-$ \Omega $, $p=4$}
	\end{subfigure}\hfill
	\begin{subfigure}{0.33\textwidth}
		\includegraphics[scale=0.16,right]{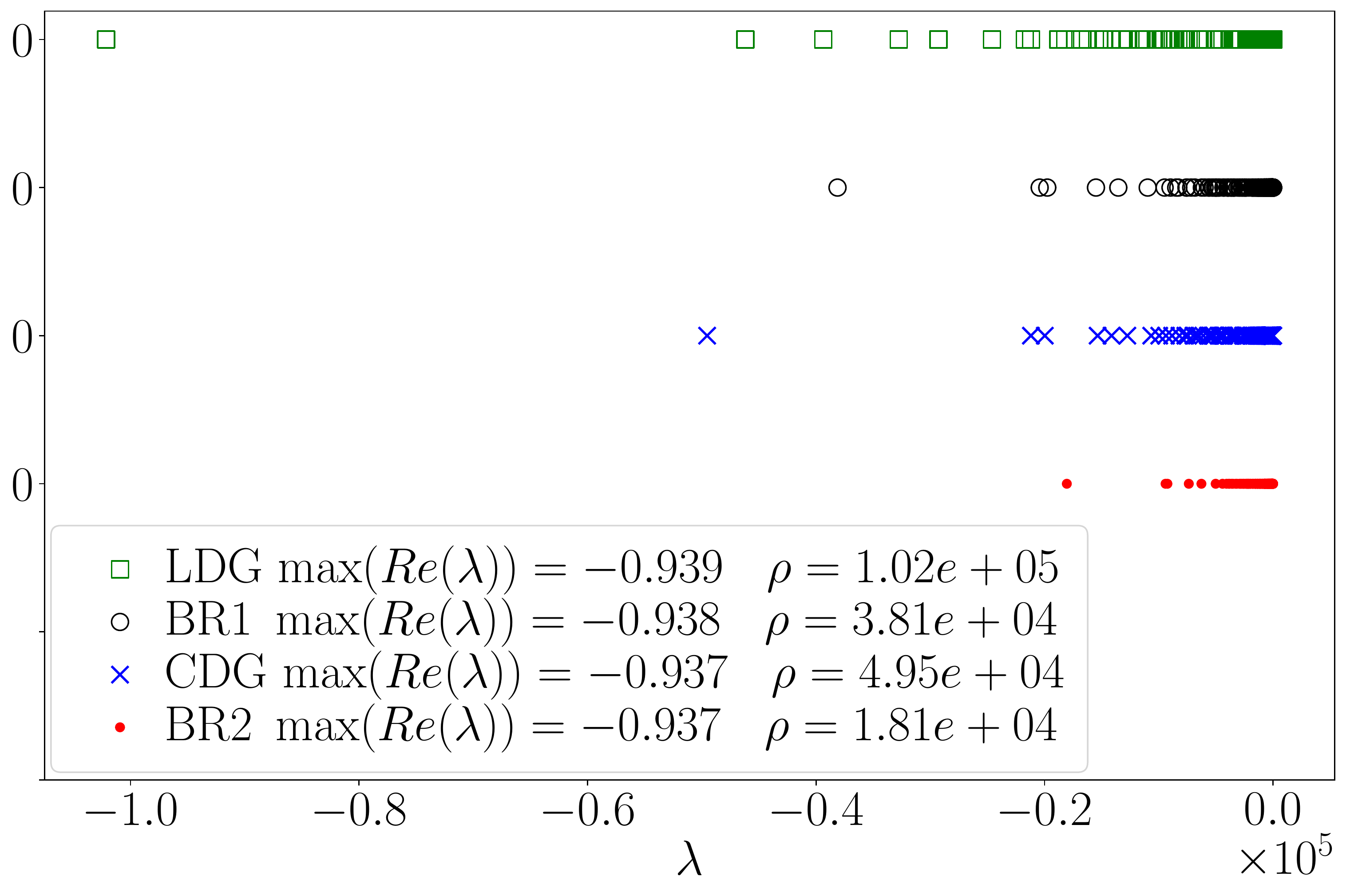}
		\caption{\label{fig:Gamma symmetric spectra p=4} SBP-$ \Gamma $, $p=4$}
	\end{subfigure} \hfill
	\begin{subfigure}{0.33\textwidth}
		\includegraphics[scale=0.16,right]{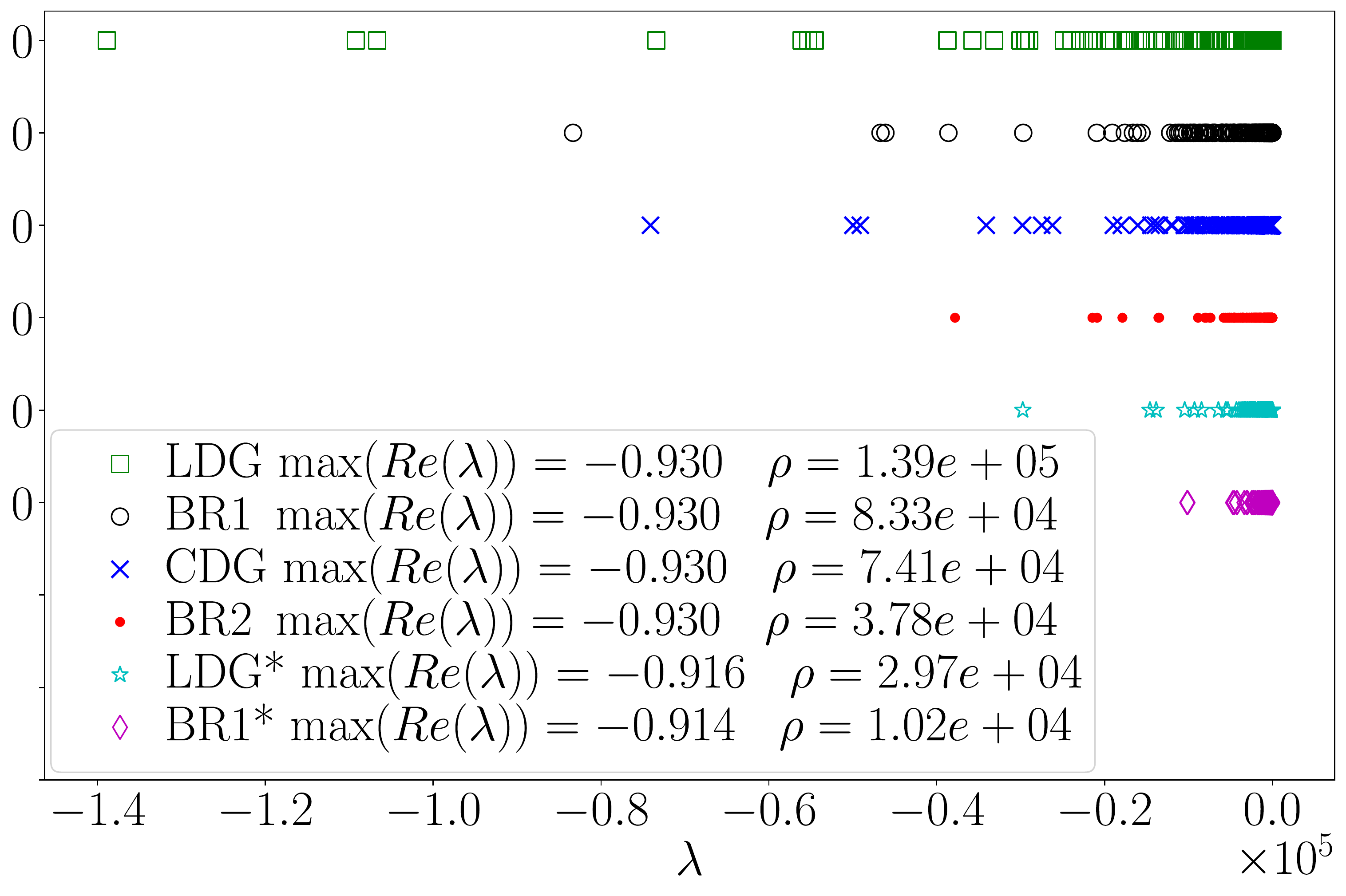}
		\caption{\label{fig:diage symmetric spectra p=4} SBP-E, $p=4$}
	\end{subfigure}
	\\
	\begin{subfigure}{0.33\textwidth}
		\includegraphics[scale=0.16,right]{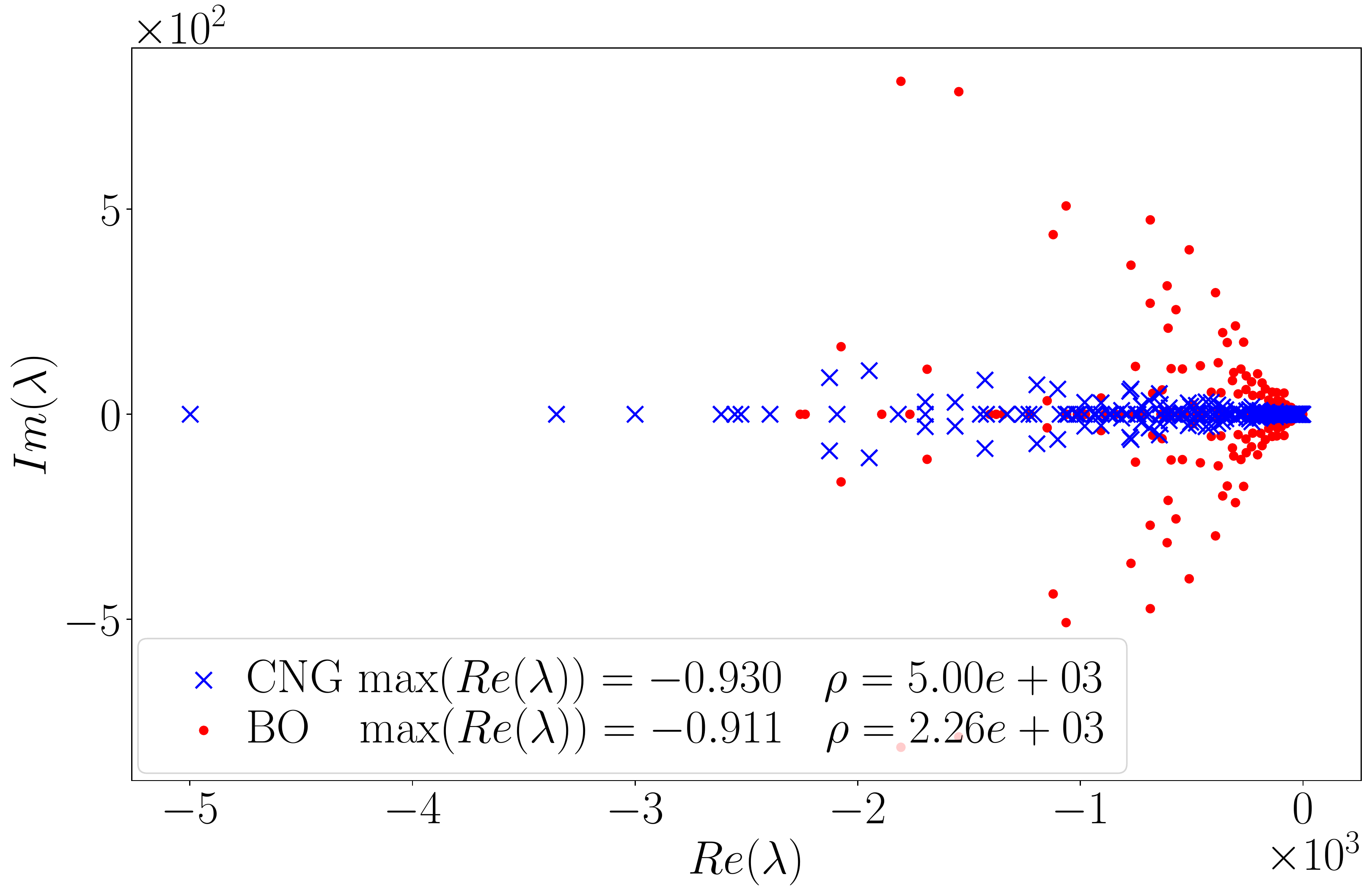}
		\caption{\label{fig:Omega asymmetric spectra p=4} SBP-$ \Omega $, $p=4$}
	\end{subfigure}\hfill
	\begin{subfigure}{0.33\textwidth}
		\includegraphics[scale=0.16,right]{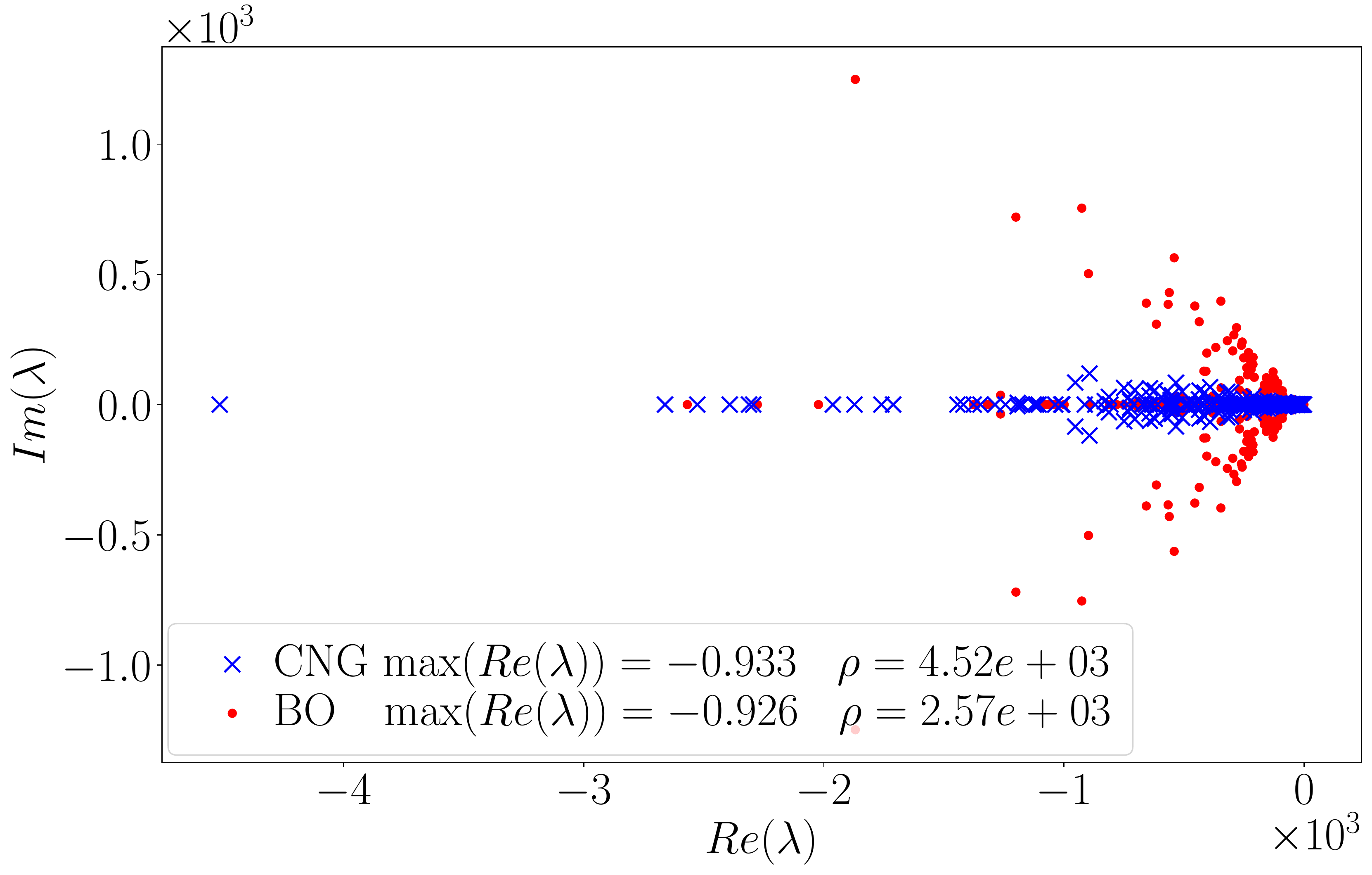}
		\caption{\label{fig:Gamma asymmetric spectra p=4} SBP-$ \Gamma $, $p=4$}
	\end{subfigure} \hfill
	\begin{subfigure}{0.33\textwidth}
		\includegraphics[scale=0.16,right]{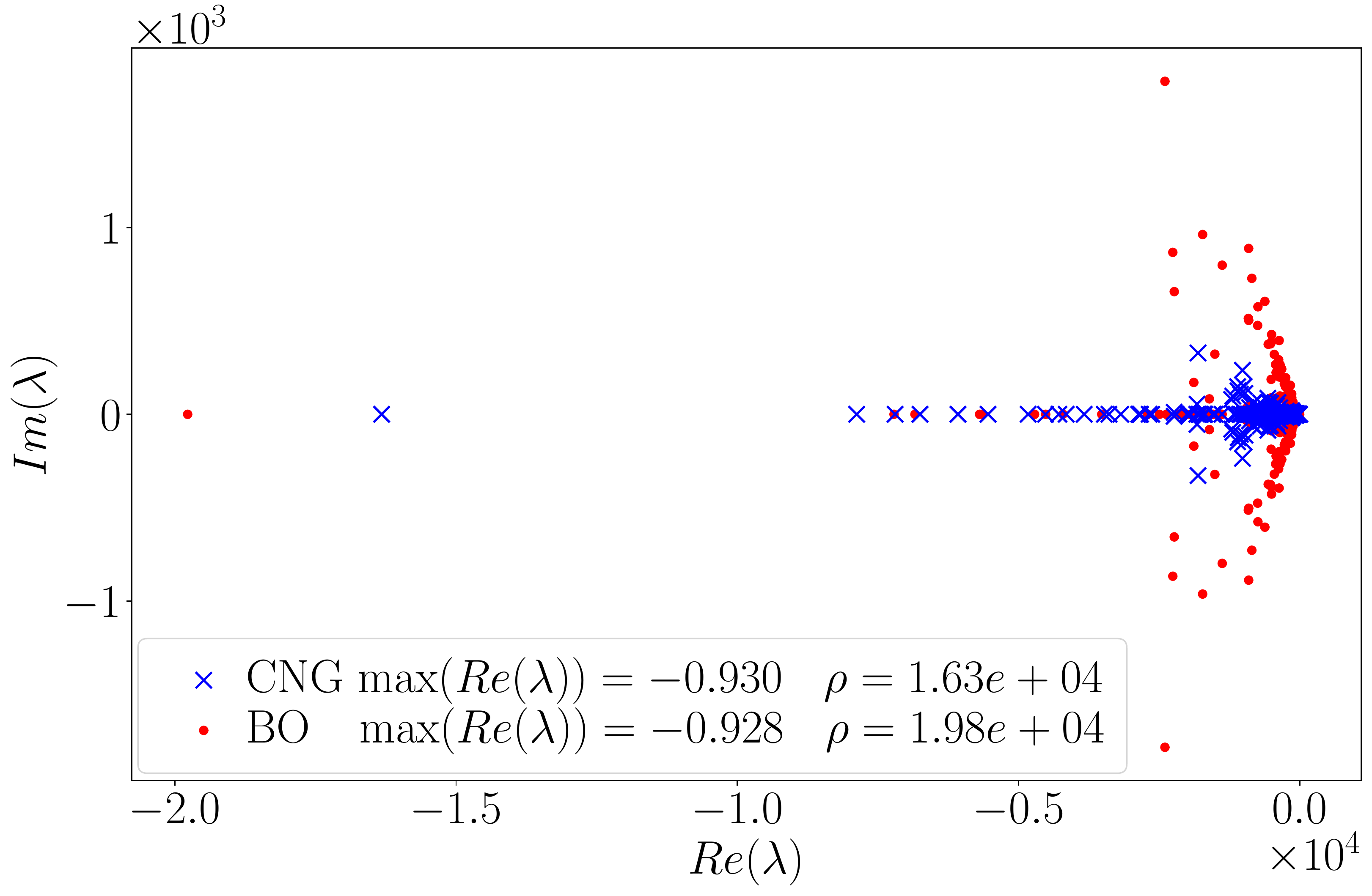}
		\caption{\label{fig:diage asymmetric spectra p=4} SBP-E, $p=4$}
	\end{subfigure}
	\caption{\label{fig:Eigenspectra} Eigenspectra of the system matrix resulting from SBP-SAT discretization of \cref{eq:Poisson problem} with $ n_e = 14 $ elements. \violet{The BR1* and LDG* SATs represent the unmodified BR1 and LDG SATs}.} 
\end{figure}

The maximum time step that can be used with explicit time-marching schemes depends on the spectral radius of the system matrix. Figure \ref{fig:Eigenspectra} shows the eigenspectra of the system matrices arising from the SBP-SAT discretizations of \cref{eq:Poisson problem}. While the BO and CNG SATs produce eigenvalues with imaginary parts, all of the adjoint consistent SATs have eigenvalues on the negative real axis. The BO SAT leads to the smallest spectral radius, $ \rho $, except when used with SBP diagonal-E operators. SBP diagonal-E operators achieve their smallest spectral radii when used with the \violet{unmodified BR1 SAT}. The \violet{modified} LDG SAT produces the largest spectral radius regardless of the type of SBP operator it is used with. In fact, the spectral radius obtained with the LDG SAT is about four times larger than the spectral radius obtained with the BR2 SAT. In comparison, the BR1 and CDG SATs yield spectral radii about twice as large as that of the BR2 SAT. The spectral radii of the BR1, LDG and CDG SATs can be reduced by approximately a factor of $ {1}/{\sigma_{1}} $ if $ \T_{\gamma k}^{(1)} $ is multiplied by $ 0<\sigma_{1} < 1 $, but this would compromise the stability of the discretizations. \violet{The unmodified BR1 and LDG SATs have smaller $ \T_{\gamma k}^{(1)} $ coefficients compared to the rest of the adjoint consistent SATs, and they produce smaller spectral radii, as can be seen from \cref{fig:diage symmetric spectra p=3,fig:diage symmetric spectra p=4}.}

The variation of the spectral radius with respect to the SBP operators can also be inferred from \cref{fig:Eigenspectra}. The SBP-$ \Omega $ and SBP-$ \Gamma $ operators show comparable spectral radii in all cases. In contrast, the SBP diagonal-E operator produces larger spectral radii than the SBP-$ \Omega $ and SBP-$ \Gamma $ operators. It also exhibits the largest ratio of the magnitudes of the smallest to the largest eigenvalues for the $ p=3 $ case. It can be seen from \cref{fig:Eigenspectra} that the eigenvalue with the smallest magnitude for the case with the $ p=3 $ SBP diagonal-E operator has a magnitude approximately two orders of magnitude smaller than those produced with the SBP-$ \Omega $ and SBP-$ \Gamma $ operators. This is also reflected in the condition number of the system matrix presented in \cref{tab:Condition number}.

\subsection{Conditioning}
The condition number of a system matrix affects the solution accuracy and the convergence rate of iterative solvers for implicit methods. \cref{tab:Condition number} shows the condition numbers of the system matrices resulting from various SBP-SAT discretizations of \cref{eq:Poisson problem}. The LDG SAT produces the largest condition numbers, and the BR2 SAT yields the smallest condition numbers among the adjoint consistent SATs. It can also be inferred from \cref{tab:Condition number} that compared to the LDG SAT, the BO and CNG SATs yield approximately an order of magnitude smaller condition numbers. They also give significantly smaller condition numbers compared to the BR1, BR2, and CDG SATs. \violet{In contrast, the unmodified LDG and BR1 SATs yield smaller condition numbers than the rest of the adjoint consistent SATs when used with all but the degree three SBP diagonal-E operators.} A comparison of the condition numbers in \cref{tab:Condition number} by the type of SBP operator reveals that the SBP diagonal-E operators lead to larger condition numbers than the SBP-$ \Gamma $ and SBP-$ \Omega $ operators. As noted in \cref{sec:Accuracy}, the solution and adjoint errors are considerably larger for the case with $ p=3 $ SBP diagonal-E operator compared to the solutions with the same degree SBP-$ \Gamma $ and SBP-$ \Omega $ operators which yield system matrices with significantly smaller condition numbers as can be seen from \cref{tab:Condition number}. 

\begin{figure}[!t]
	\centering
	\begin{subfigure}{0.33\textwidth}
		\centering
		\includegraphics[scale=0.2,right]{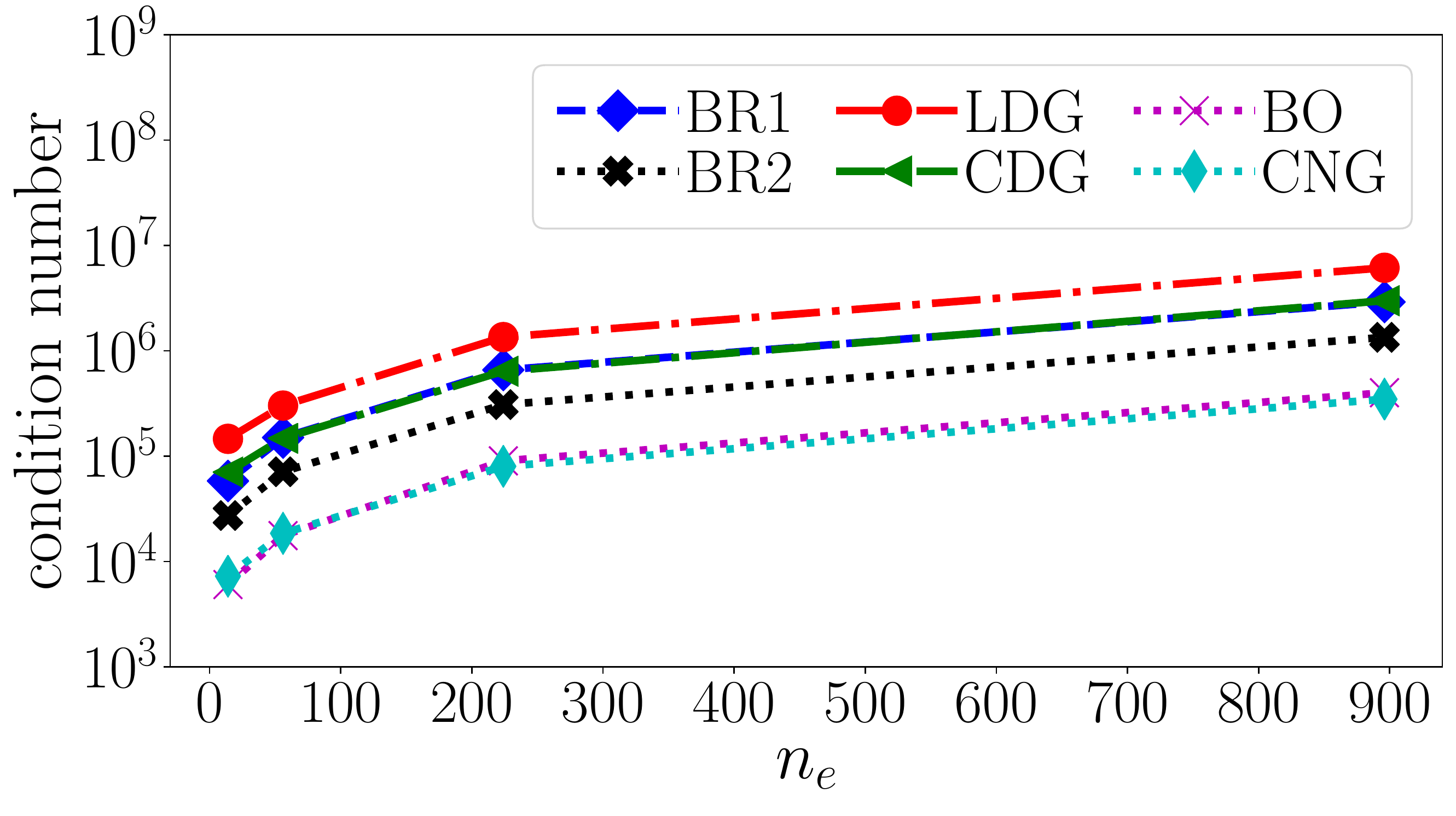}
		\caption{\label{fig:cond num omega p=4} SBP-$ \Omega $, $p=4$}
	\end{subfigure}
	\begin{subfigure}{0.33\textwidth}
		\centering
		\includegraphics[scale=0.2,right]{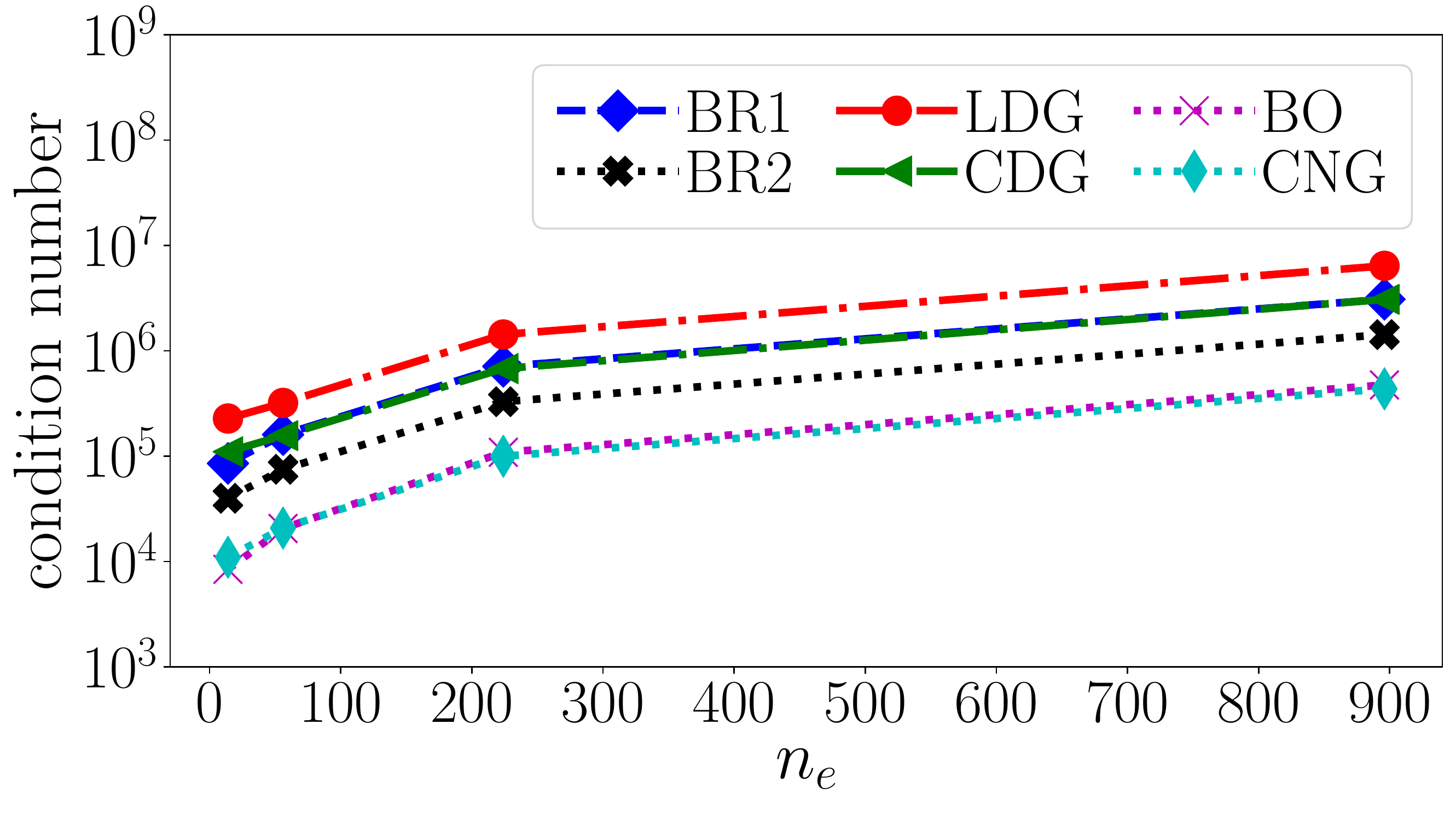}
		\caption{\label{fig:cond num gamma p=4} SBP-$ \Gamma $, $p=4$}
	\end{subfigure} 
	\begin{subfigure}{0.33\textwidth}
		\centering
		\includegraphics[scale=0.2,right]{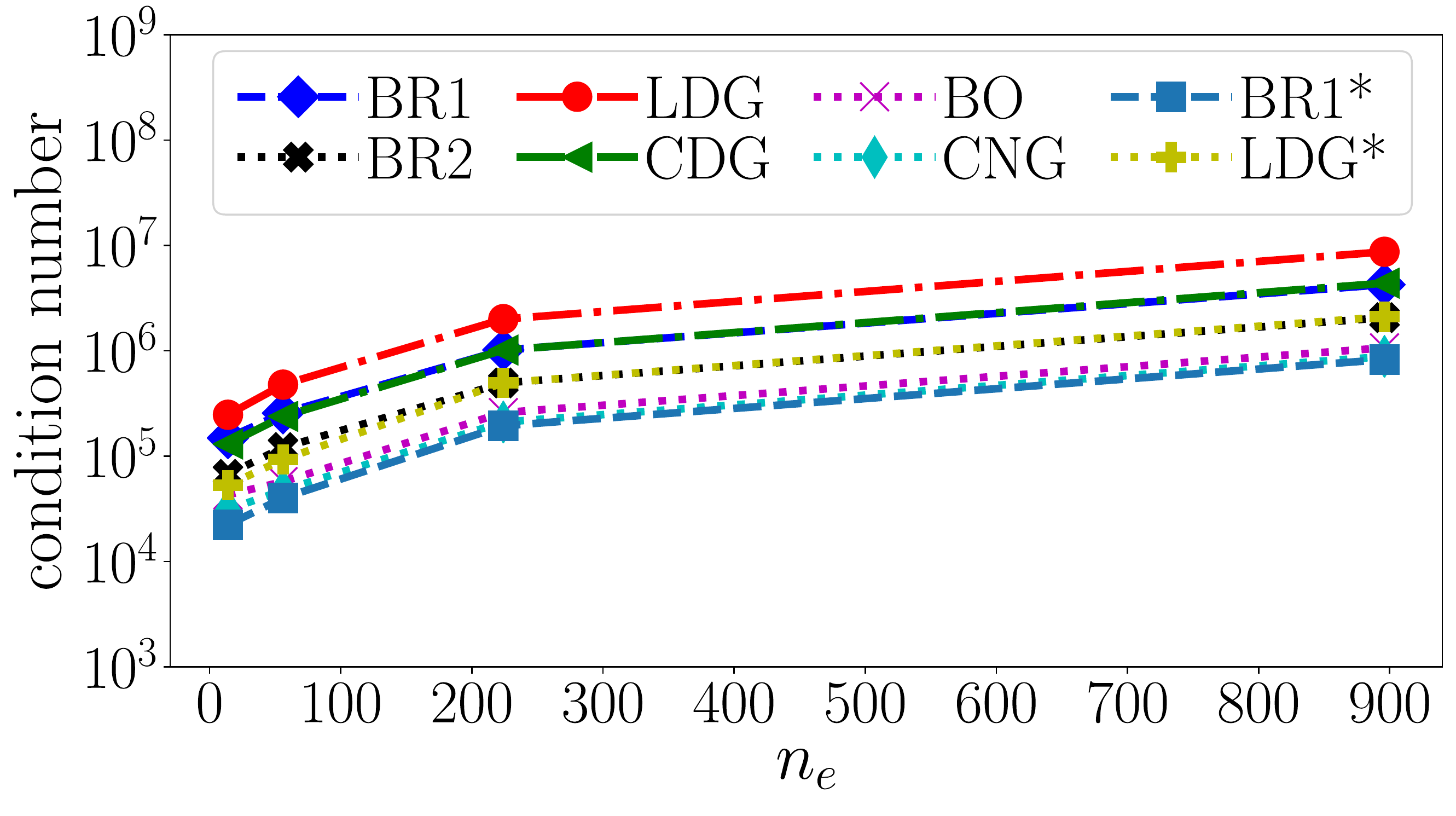}
		\caption{\label{fig:cond num diage p=4} SBP-E, $p=4$}
	\end{subfigure}
	\caption{\label{fig:Condition number scaling} Growth of condition number with respect to mesh refinement. \violet{The BR1* and LDG* SATs represent the unmodified BR1 and LDG SATs}.}
\end{figure}

\begin{figure}[!t]
	\centering
	\begin{subfigure}{0.33\textwidth}
		\centering
		\includegraphics[scale=0.2,right]{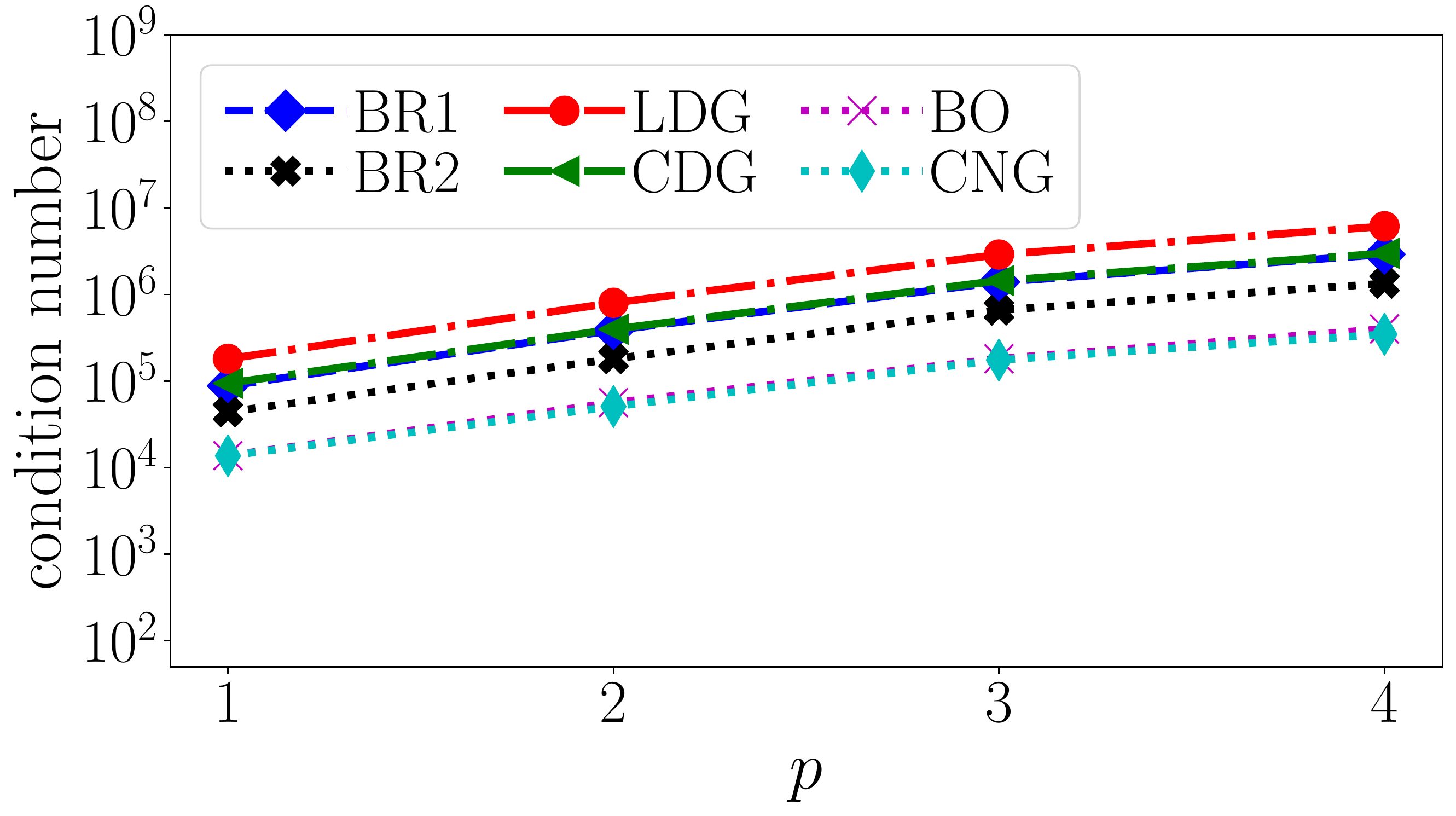}
		\caption{\label{fig:cond num omega} SBP-$ \Omega $, $n_e=896$}
	\end{subfigure}
	\begin{subfigure}{0.33\textwidth}
		\centering
		\includegraphics[scale=0.2,right]{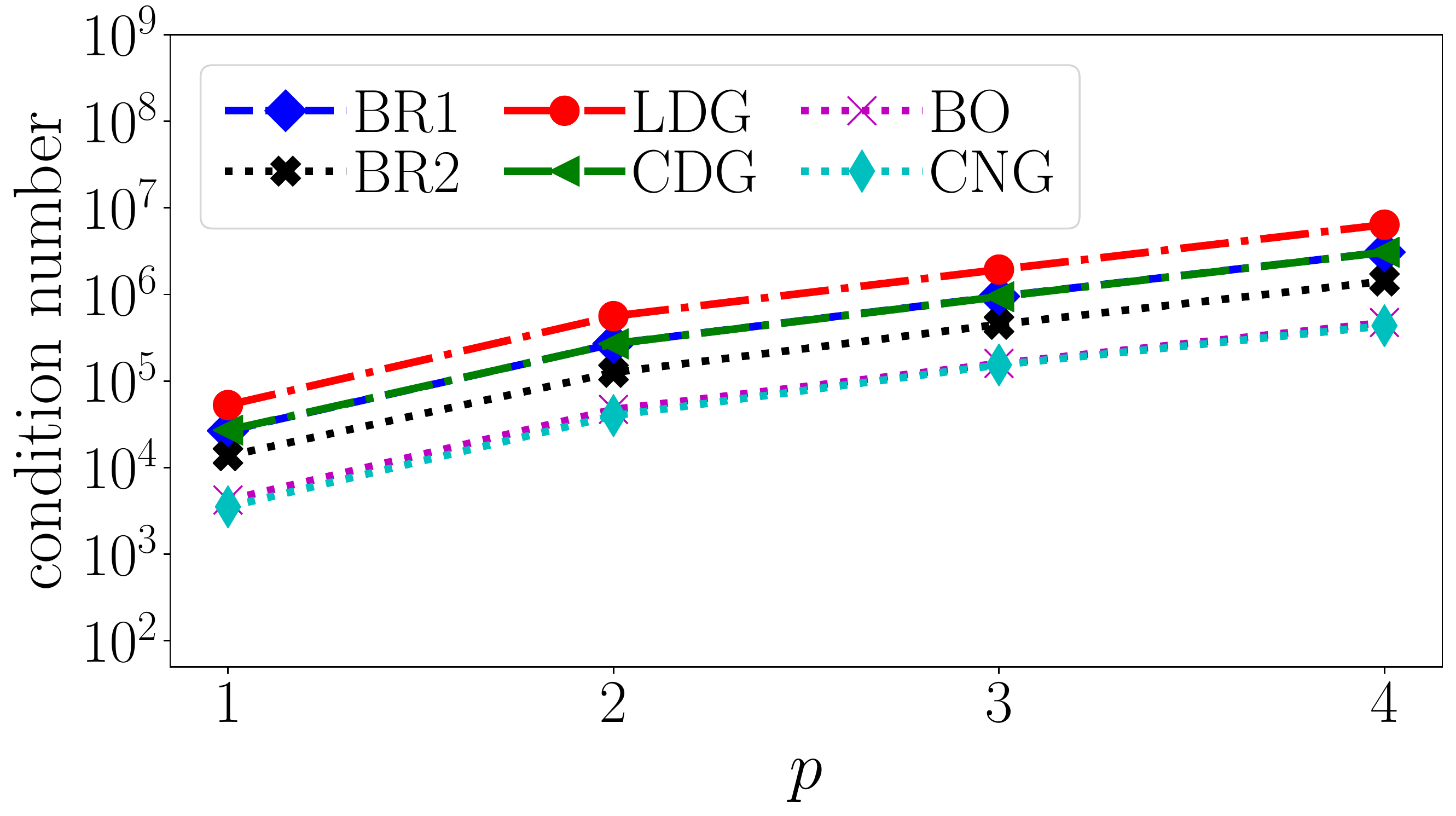}
		\caption{\label{fig:cond num gamma} SBP-$ \Gamma $, $n_e=896$}
	\end{subfigure} 
	\begin{subfigure}{0.33\textwidth}
		\centering
		\includegraphics[scale=0.2,right]{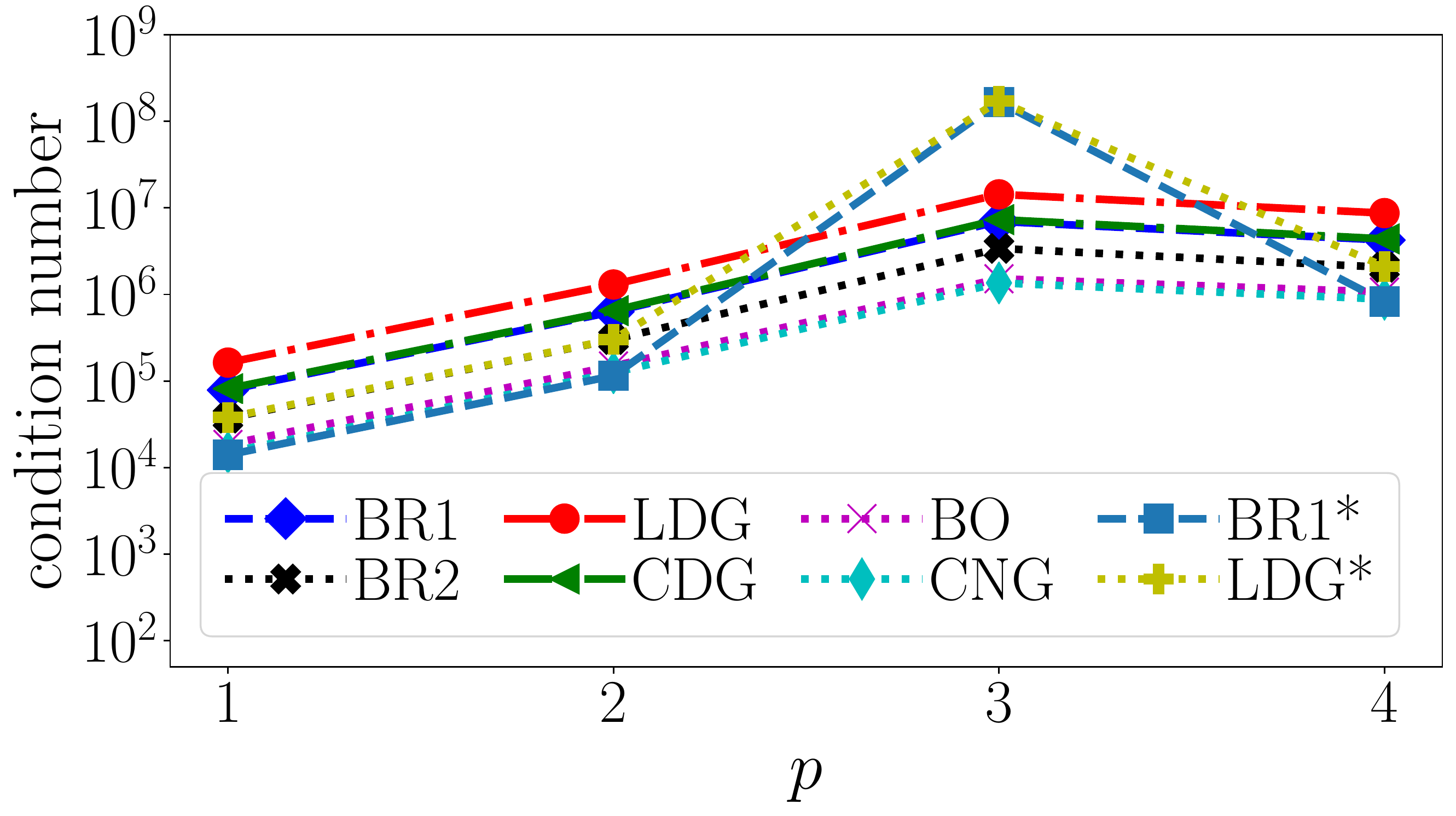}
		\caption{\label{fig:cond num diage} SBP-E, $n_e=896$}
	\end{subfigure}
	\caption{\label{fig:Condition number scaling with p} Growth of condition number with degree of SBP operators. \violet{The BR1* and LDG* SATs represent the unmodified BR1 and LDG SATs}.}
\end{figure}

The growth of the condition number with mesh refinement for degree four SBP operators is depicted in \cref{fig:Condition number scaling}. The figure shows that the scaling factors between the condition numbers resulting from the use of the different types of SAT remain roughly the same under mesh refinement. This holds for the lower degree SBP operators as well. \blue{Similarly, \cref{fig:Condition number scaling with p} shows that, for the SATs considered, the condition number scales at approximately the same rate as the degree of the operators increases. For SBP-$ \Omega $ and SBP-$ \Gamma $ operators, the increase in condition number with the polynomial degree of the operators is roughly linear; a similar observation was made for DG operators in \cite{kirby2005selecting}. From \cref{fig:cond num diage}, we see that the condition number for the degree three SBP diagonal-E operator is larger than that of the degree four operator, and this trend is observed in a more pronounced manner for smaller mesh sizes.}

\begin{table*} [!t]
	\small
	\caption{\label{tab:Condition number} Condition number of the system matrix arising from discretization of \cref{eq:Poisson problem} using $ n_e=14 $ elements. \violet{The BR1* and LDG* SATs represent the unmodified BR1 and LDG SATs}.}
	\centering
	\renewcommand*{\arraystretch}{1.1}
\begin{tabular}{ccccccccccc}
	\toprule
		$p$& 	Operator     &      BR1 &     BR1* &      BR2 &      LDG &     LDG* &      CDG &       BO &      CNG \\
	\midrule
		   &    SBP-$\Omega$ & 5.09e+02 &     --   & 2.55e+02 & 1.01e+03 &     --   & 4.96e+02 & 1.32e+02 & 7.60e+01 \\
	     1 &    SBP-$\Gamma$ & 5.05e+02 &     --   & 3.02e+02 & 1.29e+03 &     --   & 6.88e+02 & 1.01e+02 & 1.14e+02 \\  
	       &    SBP-E        & 9.13e+02 & 3.65e+02 & 4.12e+02 & 2.01e+03 & 6.88e+02 & 9.62e+02 & 2.36e+02 & 1.57e+02 \\
	\midrule
	       &    SBP-$\Omega$ & 3.88e+03 &     --   & 1.90e+03 & 8.59e+03 &     --   & 4.30e+03 & 3.91e+02 & 5.33e+02 \\
	     2 &    SBP-$\Gamma$ & 6.30e+03 &     --   & 3.00e+03 & 1.70e+04 &     --   & 8.36e+03 & 5.83e+02 & 8.28e+02 \\
	       &    SBP-E        & 1.06e+04 & 2.09e+03 & 4.82e+03 & 2.49e+04 & 3.84e+03 & 1.16e+04 & 2.86e+03 & 2.08e+03 \\
	\midrule
	       &    SBP-$\Omega$ & 1.85e+04 &     --   & 8.86e+03 & 4.30e+04 &     --   & 2.11e+04 & 1.98e+03 & 2.42e+03 \\
	     3 &    SBP-$\Gamma$ & 2.66e+04 &     --   & 1.26e+04 & 7.22e+04 &     --   & 3.54e+04 & 2.76e+03 & 3.58e+03 \\
	       &    SBP-E        & 2.65e+06 & 3.63e+06 & 1.22e+06 & 6.69e+06 & 4.76e+06 & 3.15e+06 & 6.60e+05 & 5.26e+05 \\
	\midrule
	       &    SBP-$\Omega$ & 5.81e+04 &     --   & 2.73e+04 & 1.46e+05 &     --   & 7.05e+04 & 6.04e+03 & 7.26e+03 \\
	     4 &    SBP-$\Gamma$ & 8.54e+04 &     --   & 3.98e+04 & 2.27e+05 &     --   & 1.10e+05 & 8.43e+03 & 1.09e+04 \\
	       &    SBP-E        & 1.49e+05 & 2.23e+04 & 6.71e+04 & 3.86e+05 & 5.30e+04 & 1.77e+05 & 4.29e+04 & 3.10e+04 \\
	\bottomrule
\end{tabular}
\end{table*}

\subsection{Verification of sparsity and storage requirement estimates} \label{sec:Varification of sparsity}
We verify estimates of the number of nonzero entries presented in \cref{tab:nnz} for system matrices resulting from different SBP-SAT discretizations of \cref{eq:Poisson problem}. The accuracy of the estimates is measured by the percent error with respect to the actual number of nonzero entries obtained numerically. We also compute the relative densities of the system matrices using the density due to the BR1 SAT as a reference for normalization. The results for degree four SBP operators are shown in \cref{tab:nnz numerical}. \violet{The largest errors in the estimated number of nonzero entries of the system matrices resulting from discretizations with SBP-$ \Omega $, SBP-$ \Gamma $, and SBP diagonal-E operators are $ 8.34\% $, $ 2.10\% $, and $ 0.74\% $, respectively}. For fewer elements the errors increase (\eg, $ \approx 20\% $ with 68 elements) because the ratio of the number of interior elements to boundary elements decreases.

\begin{table*}[!t]
	\small
	\caption{\label{tab:nnz numerical} Number of nonzero entries of system matrices resulting from SBP-SAT discretization of \cref{eq:Poisson problem} with $ 4352 $ degree $ p=4 $ curved SBP elements, percent error of estimated number of nonzero entries compared to actual number of nonzero entries, and relative densities (rel. density) of system matrices with respect to nonzero entries obtained with BR1 SAT.}
	\centering
	\renewcommand*{\arraystretch}{1.1}
	\begin{tabular}{l c c c c c c c c c}
		\toprule
		\multirow{2}{*}{SAT} & \multicolumn{3}{c}{SBP-$ \Omega $}  &  \multicolumn{3}{c}{SBP-$ \Gamma $} & \multicolumn{3}{c}{SBP-E} \\ \cmidrule(lr){2-4} \cmidrule(lr){5-7} \cmidrule(lr){8-10}
		    &   $nnz$   & \%error& rel. density & $nnz$ & \%error& rel. density & $nnz$ & \%error & rel. density \\ \cmidrule(lr){1-1} \cmidrule(lr){2-4} \cmidrule(lr){5-7} \cmidrule(lr){8-10}
		BR1 &  $ 9,594,000 $ &   $+2.06$ & $ 1.0000 $ &  $ 4,041,648 $ &  $+1.11$ & $ 1.0000 $ &  $ 4,617,968 $ & $+0.74$ & $ 1.0000 $ \\
		BR2 &  $ 3,877,200 $ &   $+1.02$ & $ 0.4041 $ &  $ 3,406,448 $ &  $+0.80$ & $ 0.8428 $ &  $ 4,617,968 $ & $+0.74$ & $ 1.0000 $ \\
		LDG &  $ 4,820,400 $ &   $+8.34$ & $ 0.5024 $ &  $ 2,674,048 $ &  $+2.10$ & $ 0.6616 $ &  $ 3,523,168 $ & $+0.55$ & $ 0.7629 $ \\
		CDG &  $ 3,877,200 $ &   $+1.02$ & $ 0.4041 $ &  $ 2,569,248 $ &  $+0.62$ & $ 0.6357 $ &  $ 3,523,168 $ & $+0.55$ & $ 0.7629 $ \\
		BO  &  $ 3,877,200 $ &   $+1.02$ & $ 0.4041 $ &  $ 3,406,448 $ &  $+0.80$ & $ 0.8428 $ &  $ 4,617,968 $ & $+0.74$ & $ 1.0000 $ \\
		CNG &  $ 3,877,200 $ &   $+1.02$ & $ 0.4041 $ &  $ 2,569,248 $ &  $+0.62$ & $ 0.6357 $ &  $ 3,523,168 $ & $+0.55$ & $ 0.7629 $ \\
		\bottomrule
	\end{tabular}
\end{table*}

The relative densities in \cref{tab:nnz numerical} show that the storage requirements for discretizations with SBP-$ \Omega $ operators coupled with the BR1 SAT can be reduced by up to $ \approx 60\% $ if instead compact SATs are used. For SBP-$ \Gamma $ operators density reductions up to $ \approx 35\% $ are observed if CNG, CDG, or LDG SATs are used. The same set of SATs yield $ \approx 23\% $ reduction in density when used with SBP diagonal-E operators. Compared to BR2 and BO SATs, we observe $\approx 20\% $ reduction in density when LDG, CDG and CNG SATs are used with the SBP-$ \Gamma $ operator. 

\begin{figure}[!t]
	\centering
		\begin{subfigure}{0.23\textwidth}
			\centering
			\includegraphics[scale=0.20]{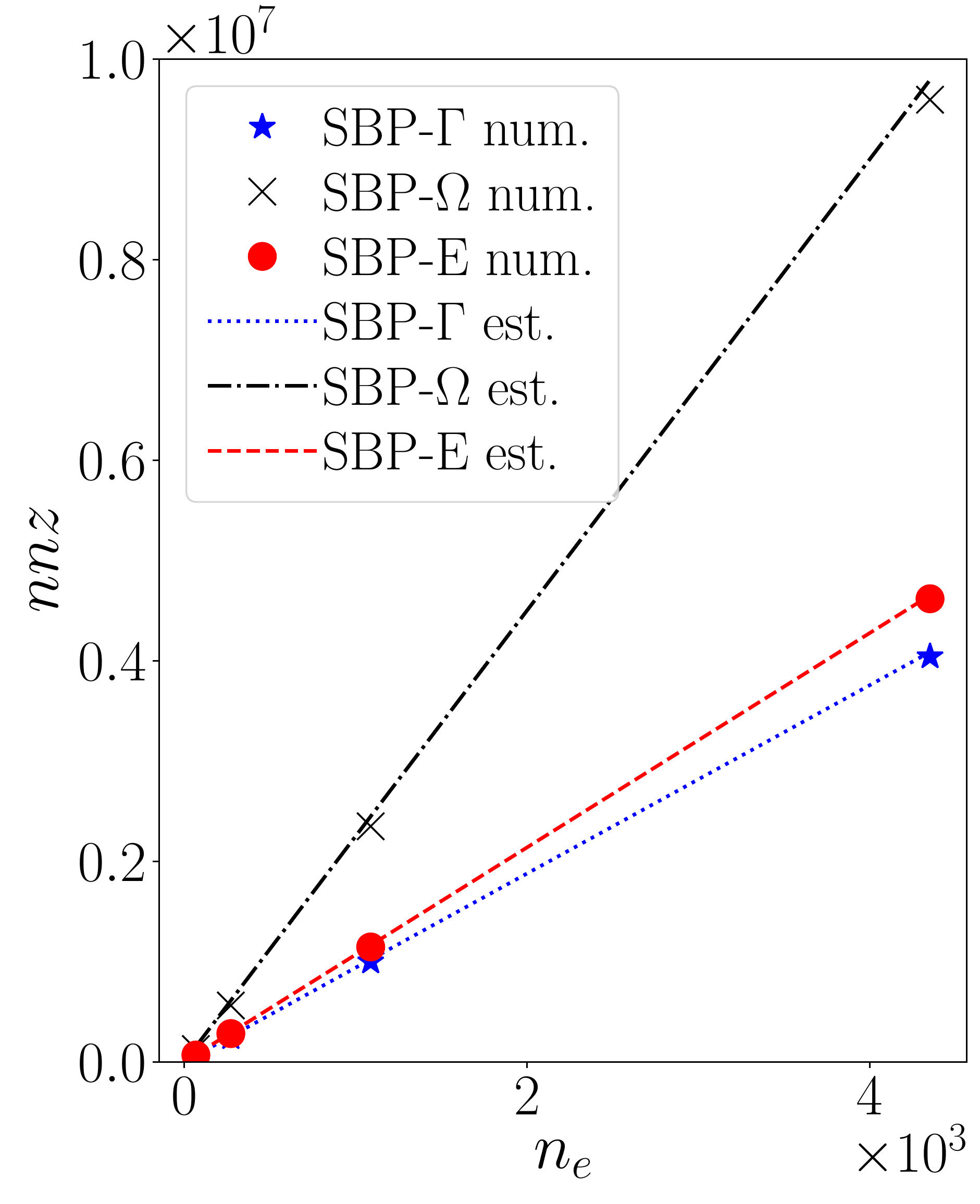}
			\caption{\label{fig:nnz_BR1_p4} BR1 SAT, $p=4$}
		\end{subfigure}\quad
		\begin{subfigure}{0.23\textwidth}
			\centering
			\includegraphics[scale=0.20]{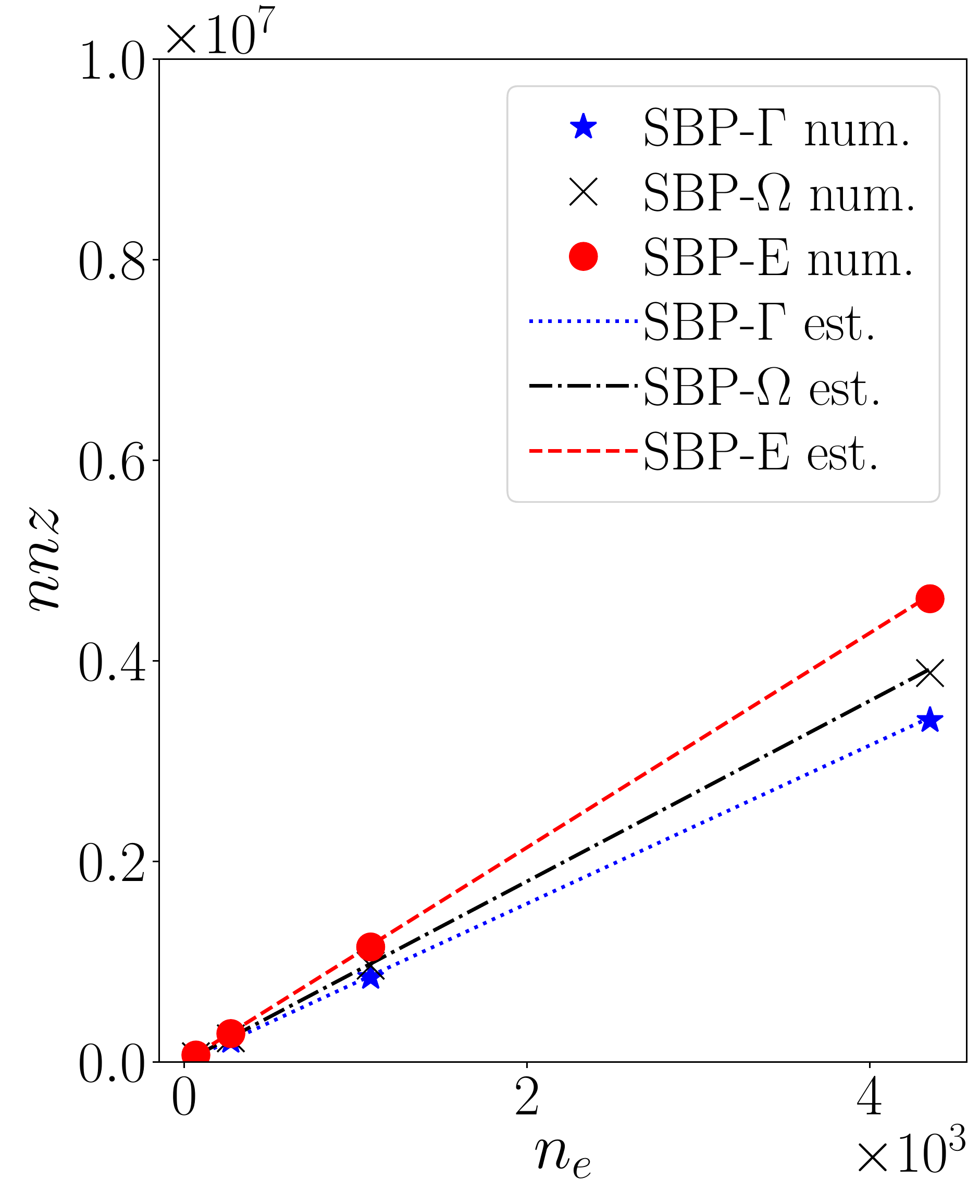}
			\caption{\label{fig:nnz_BR2_p4} BR2 SAT, $p=4$}
		\end{subfigure}\quad
		\begin{subfigure}{0.23\textwidth}
			\centering
			\includegraphics[scale=0.20]{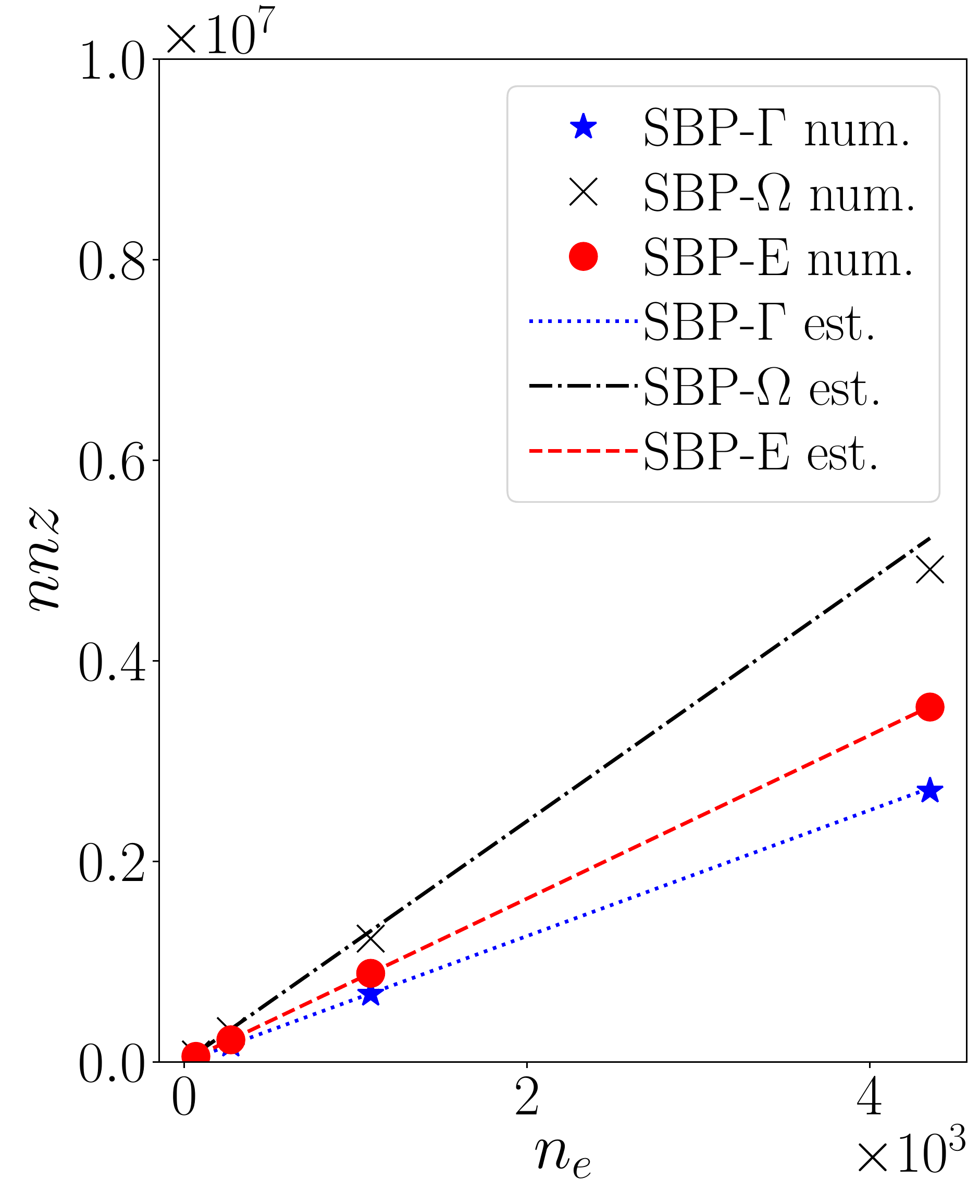}
			\caption{\label{fig:nnz_LDG_p4} LDG SAT, $p=4$}
		\end{subfigure}\quad
		\begin{subfigure}{0.23\textwidth}
			\centering
			\includegraphics[scale=0.20]{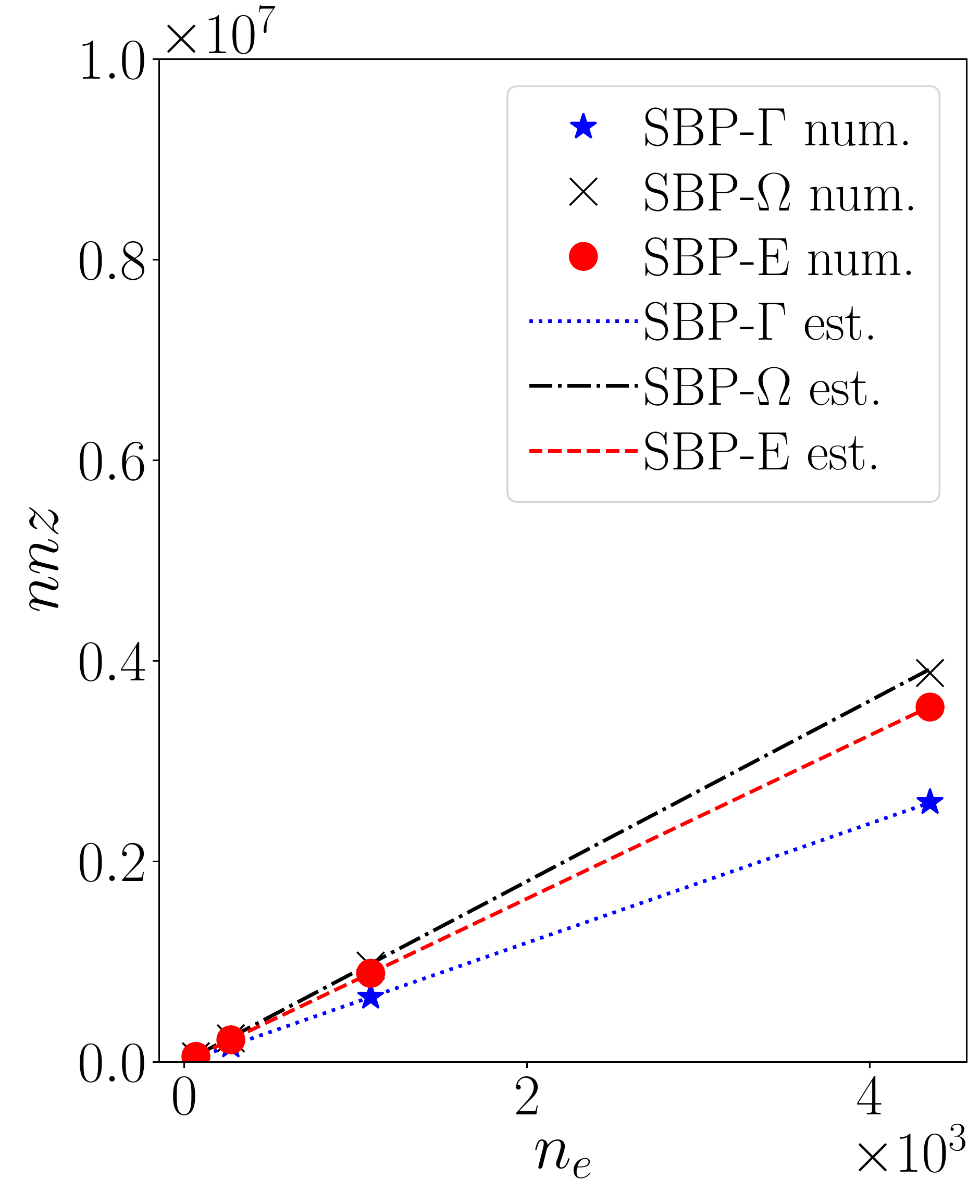}
			\caption{\label{fig:nnz_CDG_p4} CDG SAT, $p=4$}
		\end{subfigure}
	\caption{\label{fig:nnz} Comparison of number of nonzero entries when different types of SBP operators are implemented with similar SAT. Values obtained from numerical experiment are denoted by ``num." and those estimated are denoted  by ``est." in the legends.}
\end{figure}

Figure \ref{fig:nnz} shows the variation of the number of nonzero entries due to implementation of different types of degree four SBP operator with a single type of SAT. The SBP-$ \Gamma $ operator produces the fewest nonzero entries regardless of the choice of SAT. This trend is observed for lower degree operators as well, but for implementations with the BR2 SAT, the SBP-$ \Gamma $ and SBP-$ \Omega $ operators produce very similar numbers of nonzero entries. For SBP-$ \Omega $ and SBP diagonal-E operators, no conclusive statement can be made regarding which operator produces a smaller number of nonzero entries when implemented with the same SAT. Combining the observations from \cref{fig:nnz} and \cref{tab:nnz numerical}, we can conclude that the minimum number of nonzero entries (highest sparsity and lowest storage requirement) is obtained when SBP-$ \Gamma $ operators are used with the CNG SAT. While the CDG SAT also produces the same number of nonzero entries, it requires storing values of the switch functions for each facet in the discretization.

\section{Conclusions}\label{sec:conclusions}

Using a general framework, we have analyzed the numerical properties of discretizations of diffusion problems \blue{with diagonal-norm multidimensional SBP operators and various types of SAT}. The framework enables implementation of SATs without writing diffusion problems as a first-order systems of equations. This offers flexibility to switch from one type of SAT to another with a simple parameter selection of the SAT coefficients. The main theoretical results can be summarized as follows.
\begin{itemize}
	\item Conditions required for consistency, conservation, adjoint consistency, and energy stability of SBP-SAT discretizations of diffusion problems with multidimensional SBP operators are established.
	\item A functional error of order $ h^{2p} $ is attained when primal and adjoint consistent SATs are used with degree $ p $ multidimensional SBP operators in curvilinear coordinates. 
	\item Several types of SAT that correspond to known DG fluxes, including those leading to extended stencils, are identified. Instability issues observed with some these SATs are addressed by modifying the SAT coefficients. 
	\item It is shown that the BR1, BR2, and SIPG SATs are equivalent when implemented with diagonal-norm $ \R^0 $ (SBP diagonal-E) operators, which include the frequently used LGL operator in one space dimension. For the same family of operators, the LDG and CDG SATs are shown to be equivalent.
	\item Upper bounds on the number of nonzero entries in the system matrices arising from $ d $-dimensional SBP-SAT discretizations are derived. 
\end{itemize}

Numerical experiments with the two-dimensional Poisson problem were conducted to study the accuracy, eigenspectra, conditioning, and sparsity of various SBP-SAT discretizations. The adjoint consistent SATs display primal and adjoint solution convergence rates of $ p+1 $. In contrast, the adjoint inconsistent SATs, BO and CNG, show solution convergence rates of $ p+1 $ and $ p $ for odd and even degree operators, respectively. Functional superconvergence rates of $ 2p $ are attained with the adjoint consistent SATs, while the adjoint inconsistent SATs converge at lower rates. The reduction in the functional error values is more notable than the reduction in the solution error values when adjoint consistent SATs are used instead of adjoint inconsistent SATs. We summarize the rest of our observations as follows.
\begin{itemize}
	\item \violet{When used with SBP-$ \Omega $ and SBP-$ \Gamma $ operators}, the BR1 and LDG SATs couple second neighbor elements; hence, they are less amenable for code parallelization than the other types of SAT.
	\item \violet{When used with the SBP-$ \Omega $ and SBP-$ \Gamma $ operators}, the BR2 SAT leads to a system matrix with the smallest spectral radius compared to the rest of the adjoint consistent SATs. In contrast, the LDG SAT leads to a system matrix with the largest spectral radius. 
	\item \violet{When used with the SBP diagonal-E operators, the unmodified BR1 and LDG SATs are compact, adjoint consistent, and energy stable. Except for the $ p=3 $ operator, they lead to smaller condition numbers compared to the other adjoint consistent SATs. Furthermore, the unmodified BR1 SAT leads to system matrices with the smallest spectral radius while the unmodified LDG SAT produces the sparsest system matrices.}
	\item The BR2 SAT yields about half as large a spectral radius and condition number as the CDG SAT, but the CDG SAT produces system matrices with up to $ 25\% $ fewer nonzero entries when implemented with SBP-$ \Gamma $ and SBP diagonal-E operators.
	\item Compared to the adjoint consistent SATs \violet{other than the unmodified BR1 and LDG SATs}, the BO and CNG SATs lead to system matrices with significantly smaller spectral radii. This is also reflected in the conditioning of their system matrices, which have $ 1.5$ to $20 $ times smaller condition numbers.
	\item The CNG and CDG SATs produce system matrices with $ 20\%$ to $60\% $ fewer number of nonzero entries compared to the other types of SAT.  
	\item If functional superconvergence is not a priority, the CNG SAT offers interesting properties such as a reduced condition number (about half of that of the BR2 SAT), a larger time step, and a sparse system matrix. However, the scheme suffers from larger solution error and even-odd solution convergence behavior. 
\end{itemize}

We acknowledge that the choice of SATs to solve diffusion problems is not straightforward due to the competing numerical properties, which can be problem dependent, but our observations indicate that \violet{when used with SBP-$ \Omega $ and SBP-$ \Gamma $ operators}, the BR2 and CDG SATs offer superior numerical properties in most cases, and the CNG SAT is a better alternative in some cases. \violet{For the SBP diagonal-E operators, the unmodified BR1 and LDG SATs show significantly better numerical properties compared to the rest of the SATs.}  \blue{It is possible that other types of SAT with better numerical properties fall under the general framework, and this may be studied in the future.} 

\appendix 
\section{Construction of SBP operators on curved elements} \label{sec:Curvilinear Transformation}  High-order methods require accurate enough representation of curved geometries to achieve optimal solution convergence rates \cite{hesthaven2007nodal,bassi1997higheuler}. One approach to generate curved elements is to reposition facet nodes of linear meshes generated on physical elements such that they coincide with facet quadrature points on curved physical boundaries, and propagate the curvature to volume nodes \cite{hesthaven2007nodal}. In this work, however, we assume that a curvilinear mesh is available or an analytical relation is known such that coordinates of $ \alpha $-optimized Lagrange interpolation nodes discussed in \cite{hesthaven2007nodal} are accessible on each curved physical element. We then apply polynomial interpolation to find the SBP nodal locations and grid metrics in the physical space.

Crean \etal \cite{crean2018entropy} showed that SBP operators on curved physical elements preserve design order accuracy, freestream flow, and the SBP property if the curvilinear mapping satisfies \cref{assu: mapping}.
The geometric mapping from a point in the reference element, $ (\xi,\eta) \in \hat{\Omega}$, to a point in the physical element, $ (x,y) \in \Omega_k$, is defined by 
\begin{align} \label{eq:mapping}
(x,y) = \fnc{M}_k(\xi,\eta) \equiv \sum_{j=1}^{n_s^*} c_j \hat{\phi}_j(\xi,\eta), 
\end{align} 
where $ c_j $ is the coordinate of $ j $-th Lagrange interpolation node on the the physical element, $ \hat{\phi}_j \in \polyref{p_{\rm map}}$ is the $ j $-th Lagrange polynomial basis associated with the $ j $-th node on the reference element, and $ n_s^* = {p_{\rm map}+d \choose d} $ is the cardinality of the polynomial basis for the mapping. At the Lagrange nodes of the reference element, $ \hat{\phi}_{j} $ satisfies
\begin{equation} \label{eq:Lagrange basis}
\hat{\phi}_j(\xi_i, \eta_i) = \sum_{\ell=1}^{n_s^*} k_\ell^{(j)}\hat{\varphi}_\ell (\xi_i, \eta_i)=\delta_{ij}\qquad \text{for}\;j=1,\dots,n_s^*,
\end{equation}
where $ k_{\ell}^{(j)} \in \IR{}$, $ \hat{\varphi} $ is another basis function, and $ \delta_{ij} $ is the Kronecker delta operator. The basis $ \hat{\varphi} $ is chosen to be the orthonormalized canonical basis given in \cite{hesthaven2007nodal} as, 
\begin{equation} \label{eq:orthonormal basis}
\hat{\varphi}_m(\xi, \eta)=\sqrt{2}\fnc{P}_i(a)\fnc{P}_j^{(2i+1, 0)}(b)(1-b)^i,
\end{equation}
where  $ \fnc{P}_n^{(\alpha,\beta)}$ is the $ n $-th order Jacobi polynomial, $ m= j+(p_{\rm map}+1)i+1 - i/2(i-1) $, $ ij\ge 0 $, $ i+j\le p_{\rm map} $, $ a=2({1+\xi})/({1-\eta}) - 1 $, and $ b=\eta $.
Writing \cref{eq:Lagrange basis} in matrix form we have $ \hat{\V}_L \K = \I_{n_s} $ which yields $ \K = \hat{\V}_L^{-1} $, where the coefficient matrix, $ \K \in \IRtwo{n_s^*}{n_s^*}$, contains the coefficient $ k_\ell^{(j)} $ in the $ j$-th column and $ \ell$-th row, and $ \hat{\V}_L $ is the Vandermonde matrix constructed using the orthonormal basis in \cref{eq:orthonormal basis} and the $ \alpha $-optimized Lagrange nodes, $ \hat{S}_L = \{ \xi_i,\eta_i \}_{i=1}^{n_s^*}$, presented in \cite{hesthaven2007nodal}. The $ \alpha $-optimized Lagrange nodes minimize the Lebesgue constant and ensure the Vandermonde matrix is well-behaved \cite{hesthaven2007nodal}. Using the matrix forms of \cref{eq:mapping} and \cref{eq:Lagrange basis}, the coordinates of the SBP volume nodes, $ \bm{x}_k $, $ \bm{y}_k $, and facet node, $ \bm{x}_\gamma $, $ \bm{y}_\gamma $, in the physical element can be calculated as
\begin{equation} \label{eq:sbp nodes on physical element}
\begin{aligned}
\bm{x}_k = \hat{\V}_{\Omega} \hat{\V}_L^{-1}\tilde{\bm{x}}_{k},  && \bm{y}_k = \hat{\V}_{\Omega} \hat{\V}_L^{-1}\tilde{\bm{y}}_{k}, &&
\bm{x}_{\gamma} = \hat{\V}_{\gamma} \hat{\V}_L^{-1}\tilde{\bm{x}}_{k}, &&
\bm{y}_{\gamma} = \hat{\V}_{\gamma}\hat{\V}_L^{-1}\tilde{\bm{y}}_{k},
\end{aligned}
\end{equation}
where $ \tilde{\bm{x}}_k,\; \tilde{\bm{y}}_k \in \IR{n_s^*} $ are vectors of the $ x $ and $ y $ coordinates of the Lagrange interpolation nodes in $ \Omega_k $. Using the derivatives of \cref{eq:mapping} with respect to the reference coordinates,
\begin{equation} \label{eq:derivative of polynomial in basis}
\begin{aligned}
\pdv{\fnc{M}}{\xi} &= \sum_{j=1}^{n_{s}^*}\sum_{i=1}^{n_{s}^*}c_{j}k_i^{(j)}\pdv{\hat{\varphi}_{i}}{\xi},
&& 
\pdv{\fnc{M}}{\eta} = \sum_{j=1}^{n_{s}^*}\sum_{i=1}^{n_{s}^*}c_{j}k_i^{(j)}\pdv{\hat{\varphi}_{i}}{\eta},
\end{aligned}
\end{equation}
we compute the exact grid metrics by forming the derivatives of the Vandermonde matrix on $ \hat{S} $, \ie,
\begin{equation} \label{eq:grid metrics volume}	
\begin{aligned}
\bm{x}_{\xi,k} = \hat{\V}_{\xi,\Omega} \hat{\V}_L^{-1}\tilde{\bm{x}}_{k}, && 
\bm{y}_{\xi,k} = \hat{\V}_{\xi,\Omega} \hat{\V}_L^{-1}\tilde{\bm{y}}_{k}, &&
\bm{x}_{\eta,k} = \hat{\V}_{\eta,\Omega} \hat{\V}_L^{-1}\tilde{\bm{x}}_{k}, &&
\bm{y}_{\eta,k} = \hat{\V}_{\eta,\Omega} \hat{\V}_L^{-1}\tilde{\bm{y}}_{k},
\end{aligned}
\end{equation}
where the subscripts $ \xi $ and $ \eta $ denote partial derivatives with respect to $ \xi $ and $ \eta $, \eg, $ \bm{x}_{\xi,k} $ is the restriction of $ {\partial x}/{\partial \xi} $ on to the nodes $ S_k $. Similarly, the facet grid metrics are computed as
\begin{equation} \label{eq:grid metrics facet}	
\begin{aligned}
\bm{x}_{\xi,\gamma k} = \hat{\V}_{\xi,\gamma} \hat{\V}_L^{-1}\tilde{\bm{x}}_{k}, &&
\bm{y}_{\xi,\gamma k} = \hat{\V}_{\xi,\gamma} \hat{\V}_L^{-1}\tilde{\bm{y}}_{k}, &&
\bm{x}_{\eta,\gamma k} = \hat{\V}_{\eta,\gamma} \hat{\V}_L^{-1}\tilde{\bm{x}}_{k}, && 
\bm{y}_{\eta,\gamma k} = \hat{\V}_{\eta,\gamma} \hat{\V}_L^{-1}\tilde{\bm{y}}_{k}.
\end{aligned}
\end{equation}

The mapping Jacobian matrices for the volume, $ \fnc{J}_k: \hat{\Omega} \rightarrow \IRtwo{d}{d} $, and facets, $ \fnc{J}_{f_\gamma}: \hat{l}_\gamma \rightarrow \IRtwo{d}{(d-1)} $, are given, respectively, by
\begin{equation}
\fnc{J}_{k} = \left[\begin{array}{cc}
\pdv{x}{\xi} & \pdv{x}{\eta}\\
\pdv{y}{\xi} & \pdv{y}{\eta}
\end{array}\right], \quad 
\text{and} \quad
\fnc{J}_{f_\gamma} = \left[\begin{array}{c}
\pdv{x}{s}\\
\pdv{y}{s}
\end{array}\right] , 
\end{equation} 
where $ s=s(\xi,\eta) $ is the parametric equation of the line, $ \hat{l}_\gamma $, connecting the end points of facet $ \gamma $ on the reference element. The outward pointing unit normal vectors on facet $ \gamma $ of element $ \Omega_k $ are given by 
\begin{equation}
\bm{n}_{\gamma k} = \frac{\myabs{\fnc{J}_k}}{\myabs{\fnc{J}_{f_\gamma}}} \fnc{J}_k^{-T}\hat{\bm{n}}_\gamma
= \frac{1}{\myabs{\fnc{J}_{f_\gamma}}} \left[\begin{array}{cc}
\pdv{y}{\eta} & -\pdv{y}{\xi}\\
-\pdv{x}{\eta} & \pdv{x}{\xi}
\end{array}\right]\left[\begin{array}{c}
{\hat{n}}_{\gamma \xi}\\
{\hat{n}}_{\gamma \eta}
\end{array}\right],
\end{equation}
where $ \myabs{\fnc{J}_k} $ is the determinant of the Jacobian, and $ \myabs{\fnc{J}_{f_\gamma}} = \sqrt{[(\fnc{J}_{f_\gamma})_{1}]^2 + [(\fnc{J}_{f_\gamma})_{2}]^2}$. We evaluate $ \myabs{\fnc{J}_k} $  and $ \myabs{\fnc{J}_{f_\gamma}} $ at the volume and facet nodes of $ \Omega_k $ as
\begin{equation}
\begin{aligned}
\J_k &= \mydiag\qty(\myabs{\bm{x}_{\xi,k} \circ \bm{y}_{\eta,k} - \bm{x}_{\eta,k} \circ \bm{y}_{\xi,k} }),
&&
\J_{f_1} = \mydiag\qty(\frac{1}{\sqrt{2}}\sqrt{(\bm{x}_{\eta,\gamma k}-\bm{x}_{\xi,\gamma k})^2 + (\bm{y}_{\eta,\gamma k}-\bm{y}_{\xi,\gamma k})^2}),
\\
\J_{f_2} &= \mydiag\qty(\sqrt{(-\bm{x}_{\eta,\gamma k})^2 + (-\bm{y}_{\eta,\gamma k})^2}),
&&
\J_{f_3} = \mydiag\qty(\sqrt{(\bm{x}_{\xi,\gamma k})^2 + (\bm{y}_{\xi,\gamma k})^2}),
\end{aligned}
\end{equation} 
respectively, where $ \circ $ denotes the Hadamard (element-wise) product of vectors, and the operator $ \mydiag(\cdot) $ takes in a vector and creates a diagonal matrix with the vector placed in the main diagonal. 

The SBP operators on the physical element are constructed following \cite{shadpey2020entropy,crean2018entropy}. The norm matrices on the physical element read
\begin{equation} \label{eq:norm matrices}
\begin{aligned}
\H_k &= \J_k \hat{\H}, 
& &
\B_\gamma = \J_{f_\gamma} \hat{\B}_\gamma.
\end{aligned}
\end{equation}
The normals at facet $ \gamma $ are stored in the diagonal matrices
\begin{equation} \label{eq:normal matrices}
\begin{aligned}
\Nxgk &= \J_{f_\gamma}^{-1}[\mydiag(\bm{y}_{\eta,\gamma k}) \hat{n}_{\gamma \xi} - \mydiag(\bm{y}_{\xi,\gamma k}) \hat{n}_{\gamma \eta}],
&&
\Nygk = \J_{f_\gamma}^{-1}[-\mydiag(\bm{x}_{\eta,\gamma k}) \hat{n}_{\gamma \xi} + \mydiag(\bm{x}_{\xi,\gamma k}) \hat{n}_{\gamma \eta}].
\end{aligned}
\end{equation}
The surface integral matrices in the $ x $ and $ y $ directions are given by
\begin{equation} \label{eq:Ex and Ey matrices}
\begin{aligned}
\Exk &= 	\sumfk \Rgk^T\B_\gamma \Nxgk \Rgk,
&&
\Eyk = 	\sumfk \Rgk^T\B_\gamma \Nygk \Rgk,
\end{aligned}
\end{equation}
and the skew-symmetric matrices are constructed as
\begin{equation} \label{eq:Sxk and Syk}
\begin{aligned}
\Sxk &= \frac{1}{2}\qty(\mydiag(\bm{y}_{\eta,\gamma k})\Qxi - \Qxi^T\mydiag(\bm{y}_{\eta,\gamma k})) + \frac{1}{2} \qty(-\mydiag(\bm{y}_{\xi,\gamma k})\Qeta + \Qeta^T\mydiag(\bm{y}_{\xi,\gamma k})),
\\
\Syk &= \frac{1}{2}\qty(-\mydiag(\bm{x}_{\eta,\gamma k})\Qxi + \Qxi^T\mydiag(\bm{x}_{\eta,\gamma k})) + \frac{1}{2} \qty(\mydiag(\bm{x}_{\xi,\gamma k})\Qeta - \Qeta^T\mydiag(\bm{x}_{\xi,\gamma k})).
\end{aligned}
\end{equation}
Finally, the derivative operators are computed as
\begin{equation} \label{eq:Dx and Dy}
\begin{aligned}
\Dxk &= \H_k^{-1}\left(\Sxk + \frac{1}{2}\Exk\right),
&&
\Dyk = \H_k^{-1}\left(\Syk + \frac{1}{2}\Eyk\right).
\end{aligned}
\end{equation}

\section{Summary of notation}
\blue{The analysis presented in this work is notation heavy; hence, we tabulate some of the important notation in \cref{tab:notation} for quick referencing.} 

\begin{table*}[!t]
	\small
	\caption{\label{tab:notation}Summary of important notation}
	\centering
	\setlength{\tabcolsep}{1em}
	\renewcommand{\arraystretch}{1.25}
	\begin{tabular}{l l l}
		\toprule
		\makecell[l]{Notation} & \makecell[l]{Equation} & \makecell[l]{Description}
		\\ \midrule 
		\makecell[l]{$ \D_{xk} $, $ \D_{yk} $} & \makecell[l]{\cref{eq:Dx and Dy}} & \makecell[l]{First derivative operators in the $ x $ and $ y $ directions}
		\\
		\makecell[l]{$ \D_k^{(2)} $} & \makecell[l]{\cref{eq:D2 1st form}} & \makecell[l]{Second derivative operator approximating $ \nabla \cdot (\lambda \nabla) $ on element $ \Omega_k $}
		\\
		\makecell[l]{$ \H_{k} $, $ \B_{\gamma} $} & \makecell[l]{\cref{eq:norm matrices}} & \makecell[l] {Diagonal norm matrices of element $ \Omega_k $ and facet $ \gamma $, respectively }
		\\
		\makecell[l]{$\E_{x,k}$, $ \E_{y,k} $} & \makecell[l]{\cref{eq:Ex and Ey matrices}} & \makecell[l] {Surface integral matrix in the $ x $ and $ y $ directions}
		\\
		\makecell[l]{$ \R_{\gamma k} $} & \makecell[l]{\cref{eq:extrapolation matrix}} &\makecell[l]{Extrapolation matrix from volume nodes in element $ \Omega_k $ to facet nodes on facet $ \gamma $}
		\\
		\makecell[l]{$ \D_{\gamma k} $} & \makecell[l]{\cref{eq:D_gamma k}} & \makecell[l]{Normal derivative operator approximating $ \bm{n}\cdot(\lambda\nabla) $ on facet $ \gamma $ of element $ \Omega_k $ }
		\\
		\makecell[l]{$ \N_{x \gamma k} $, $ \N_{y \gamma k} $} & \makecell[l]{\cref{eq:normal matrices}} & \makecell[l]{Diagonal matrices with the $ x $ and $ y $ components of the normal vector on face $ \gamma $ of element $ \Omega_k $}
		\\
		\makecell[l]{$ \M_k $} & \makecell[l]{\cref{eq:M_k matrix}} & \makecell[l]{A positive semidefinite matrix used for the approximation $\bm{v}_k^T\M_k\bm{u}_k \approx \int_{\Omega} \nabla\fnc{V}\cdot (\lambda\nabla \fnc{U}) \dd \Omega$}
		\\
		\makecell[l]{$ \Lambda_k $} & \makecell[l]{\cref{eq:Lambda}} & \makecell[l]{A block matrix containing the diffusivity coefficients in all combinations of directions}
		\\
		\makecell[l]{$ \T_{ak}^{(i)} $, $ \T_{abk}^{(j)} $} & \makecell[l]{\cref{eq:Interface SATs,eq:Boundary SATs}} & \makecell[l]{SAT coefficient matrices for facets $ a,b\in\{\gamma,\epsilon,\delta\} $, $ i=\{1,2,3,4,D\}$, and $ j=\{5,6\} $}
		\\
		\makecell[l]{$ \Upsilon_{abk} $} & \makecell[l]{\cref{eq:Upsilon definition}} & \makecell[l]{Component of the SAT coefficient matrices defined for facet $ a,b\in\{\gamma,\epsilon,\delta\} $}
		\\
		\makecell[l]{$ \alpha_{\gamma k} $} & \makecell[l]{\cref{eq:alpha_gamma k}} & \makecell[l]{A facet weight parameter satisfying $ \sum_{\gamma \subset \Gamma_k} \alpha_{\gamma k} = 1$ }
		\\
		\makecell[l]{$ \beta_{\gamma k} $} & \makecell[l]{\cref{eq:LDG switch with g}} & \makecell[l]{Switch function defined at facet $ \gamma $ of element $ \Omega_k $}
		\\
		\bottomrule 
	\end{tabular}
\end{table*}

\section*{Declaration of competing interest}
The authors declare that they have no known competing financial interests or personal relationships that could have appeared to influence the work reported in this paper.

\section*{Acknowledgments}
The authors would like to thank Professor Masayuki Yano for his insights on DG fluxes for elliptic problems, David Craig Penner for his helpful feedback on the functional superconvergence proofs, and the anonymous referees for their valuable comments. All figures are produced using Matplotlib \cite{hunter2007matplotlib}.
\addcontentsline{toc}{section}{Acknowledgments}

\bibliographystyle{model1-num-names}
\bibliography{references}
\addcontentsline{toc}{section}{\refname}

\end{document}